\documentclass[11pt,tS^*oside]{amsart}
\usepackage{amsmath}
\usepackage{amssymb}
\usepackage{amsfonts}
\usepackage{amscd}
\usepackage{color}
\usepackage[small,nohug,heads=vee]{diagrams}
\diagramstyle[labelstyle=\scriptstyle]

\newtheorem{thm}{Theorem}[section]

\newtheorem{lemma}[thm]{Lemma}

\newtheorem{cor}[thm]{Corollary}

\newtheorem{claim}[thm]{Claim}

\theoremstyle{definition}
\newtheorem{defn}[thm]{Definition}

\newtheorem{example}[thm]{Example}

\newtheorem{notation}[thm]{Notation}

\newtheorem{rem}[thm]{Remark}

\newtheorem{conj}[thm]{Conjecture}

\newcommand\op{\operatorname}
\newcommand\m{\medskip}

\newcommand\n{\noindent}

\newcommand{\AB}{\op{ab}}
\newcommand\ab{\op{ab}}
\newcommand\abs{\op{abs}}
\newcommand\ample{\op{ample}}

\newcommand{\an}{{\op{an}}}

\newcommand\ad{\op{ad}}
\newcommand\as{\op{as}}

\newcommand\Aut{\op{Aut}}

\newcommand\chara{\op{char}}

\newcommand\codim{\op{codim}}

\newcommand\Cone{\op{Cone}}
\newcommand\comm{\op{comm}}

\newcommand\Conv{\op{Conv}}

\newcommand\Cut{\op{Cut}}

\newcommand\cutlog{\op{cutlog}}

\newcommand\DD{\op{DD}}
\newcommand\diag{\op{diag}}

\newcommand\e{{\mathbf e}}

\newcommand\End{\op{End}}
\newcommand\ex{\op{ex}}
\newcommand\ext{\op{ext}}
\newcommand\Ext{\op{Ext}}

\newcommand\Fan{\op{Fan}}

\newcommand\FC{\op{FC}}

\newcommand\fppf{\op{fppf}}

\newcommand\Gal{\op{Gal}}

\newcommand\Hom{\op{Hom}}

\newcommand\id{\op{id}}
\newcommand\Image{\op{Im}}

\newcommand\Impart{\op{Im}}

\newcommand\inv{\op{inv}}

\newcommand\init{\op{in}}

\newcommand{\LCM}{\op{LCM}}

\newcommand\ma{\op{ma}}

\newcommand\Map{\op{Map}}

\newcommand\new{\op{new}}

\newcommand\NeFC{\text{N\'eFC}}

\newcommand\Orb{\op{Orb}}

\newcommand\ord{\op{ord}}

\newcommand\Pic{\op{Pic}}
\newcommand\pr{\op{pr}}

\newcommand\Proj{\op{Proj}}

\newcommand\red{\op{red}}

\newcommand\Res{\op{Res}}

\newcommand\rank{\op{rank}}
\newcommand{\rig}{{\op{rig}}}

\newcommand\Sch{\op{Sch}}

\newcommand\Semi{\op{Semi}}
\newcommand\Sing{\op{Sing}}
\newcommand\Sk{\op{Sk}}

\newcommand\su{\op{su}}

\newcommand\specialization{\op{sp}}
\newcommand\Spec{\op{Spec}}
\newcommand\Spf{\op{Spf}\,}

\newcommand\spl{\op{spl}}

\newcommand\Star{\op{Star}}

\newcommand\Sup{\op{Sup}}

\newcommand\symm{\op{symm}}

\newcommand\univ{\op{univ}}
\newcommand\unr{\op{unr}}

\newcommand\Vor{\op{Vor}}
\newcommand\vor{\op{vor}}
\newcommand\wt{\op{wt}}

\newcommand\Zar{\op{Zar}}

\newcommand{\Scheme}{\cS ch}

\newcommand{\tc}{{\tilde c}}

\newcommand{\tildeh}{{\tilde h}}

\newcommand{\tx}{{\tilde x}}

\newcommand{\tz}{{\tilde z}}
\newcommand{\tzeta}{{\tilde\zeta}}

\newcommand{\tpi}{{\tilde\pi}}
\newcommand{\tpsi}{{\tilde\psi}}

\newcommand{\tmu}{{\tilde\mu}}

\newcommand{\tdelta}{{\tilde\delta}}

\newcommand{\tA}{{\tilde A}}

\newcommand{\tC}{{\tilde C}}
\newcommand{\tD}{{\tilde D}}

\newcommand{\tG}{{\tilde G}}

\newcommand{\tP}{{\tilde P}}
\newcommand{\tQ}{{\tilde Q}}
\newcommand{\tS}{{\tilde S}}
\newcommand{\tT}{{\tilde T}}
\newcommand{\tU}{{\tilde U}}

\newcommand{\wtW}{{\widetilde W}}
\newcommand{\tX}{{\tilde X}}

\newcommand{\tZ}{{\tilde Z}}

\newcommand{\tbL}{{\tilde\bL}}
\newcommand{\tcL}{{\tilde\cL}}

\newcommand{\tcN}{{\tilde\cN}}

\newcommand{\underbarf}{\underline{f}}

\newcommand{\underbarW}{\underline{W}}

\newcommand{\underbarZ}{\underline{Z}}

\newcommand{\barf}{{\bar{f}}}

\newcommand{\baru}{{\bar{u}}}

\newcommand{\barx}{{\bar x}}

\newcommand{\barZ}{\bar Z}

\newcommand{\bartau}{\bar\tau}

\newcommand{\bA}{{\mathbf A}}

\newcommand{\bC}{{\mathbf C}}
\newcommand{\bD}{{\mathbf D}}
\newcommand{\be}{{\mathbf e}}

\newcommand{\bF}{{\mathbf F}}
\newcommand{\bG}{{\mathbf G}}

\newcommand{\bL}{{\mathbf L}}
\newcommand{\bM}{{\mathbf M}}
\newcommand{\bN}{{\mathbf N}}
\newcommand{\bP}{{\mathbf P}}
\newcommand{\bQ}{{\mathbf Q}}
\newcommand{\bR}{{\mathbf R}}

\newcommand{\bW}{{\mathbf W}}
\newcommand{\bZ}{{\mathbf Z}}
\newcommand{\cA}{{\mathcal A}}
\newcommand{\cB}{{\mathcal B}}

\newcommand{\cD}{{\mathcal D}}

\newcommand{\cF}{{\mathcal F}}
\newcommand{\cG}{{\mathcal G}}
\newcommand{\cH}{{\mathcal H}}

\newcommand{\cI}{{\mathcal I}}

\newcommand{\cK}{{\mathcal K}}
\newcommand{\cL}{{\mathcal L}}
\newcommand{\cM}{{\mathcal M}}
\newcommand{\cN}{{\mathcal N}}

\newcommand{\cO}{{\mathcal O}}
\newcommand{\cP}{{\mathcal P}}

\newcommand{\cR}{{\mathcal R}}
\newcommand{\cS}{{\mathcal S}}

\newcommand{\cW}{{\mathcal W}}

\newcommand{\fm}{{\mathfrak m}}

\newcommand{\fp}{{\mathfrak p}}

\ifx\undefined\pf
\newcommand{\pf}{\proof}
\fi
\ifx\undefined\cal
\newcommand{\cal}{\mathcal}
\fi
\ifx\undefined\Bbb
\newcommand{\Bbb}{\mathbf}
\fi
\ifx\undefined\bold
\newcommand{\bold}{\mathbf}
\fi
\ifx\undefined\frak
\newcommand{\frak}{\mathfrak}
\fi
\begin{document}

\title[Relative compactification]
{Relative compactification of semiabelian N\'eron models, II}
\author[Iku Nakamura]{I.~Nakamura}
\address{Department of Mathematics, 
Hokkaido University, Sapporo, 060, Japan}
\email{nakamura@math.sci.hokudai.ac.jp}
\thanks{The author was partially 
supported by the Grants-in-aid 
 for Scientific Research (C) (No. 17K05188, 22K03261), JSPS. 
This work was partially supported also by the Research Institute 
for Mathematical Sciences, an International Joint Usage/Research Center located in Kyoto University.}
\thanks{2000 {\it Mathematics Subject Classification}.
Primary 14K05; Secondary 14J10,  14K99.}
\thanks{{\it Key words and phrases}. 
Abelian varieties, N\'eron model, relative compactification.}
\begin{abstract}Let $R$ be a complete discrete valuation ring, 
$k(\eta)$ its fraction field, $S=\Spec R$, 
$(G_{\eta},\cL_{\eta})$ a polarized  
abelian variety  over $k(\eta)$ 
with $\cL_{\eta}$ symmetric ample cubical    
 and $\cG$ the N\'eron model of $G_{\eta}$ over $S$. 
Suppose that $\cG$ is semiabelian over $S$. 
Then there exists a {\it unique} relative compactification 
$(P,\cN)$ of $\cG$ such that 
($\alpha$) $P$ is Cohen-Macaulay with $\codim_P(P\setminus\cG)=2$ and  
($\beta$) $\cN$ is ample invertible 
with $\cN_{|\cG}$ cubical and $\cN_{\eta}=\cL^{\otimes n}_{\eta}$ 
for some positive integer $n$. The totally degenerate case 
has been studied in \cite{MN24}. We discuss here 
first the partially degenerate case and then the case 
where $R$ is a Dedekind domain.
\end{abstract}
\maketitle
\setcounter{tocdepth}{1}
\setcounter{section}{0}
\vskip -1cm
\tableofcontents

\section{Introduction}
\label{sec:introduction}
This is a continuation of \cite{MN24}. 
Let $R$ be a complete discrete valuation ring (abbr. CDVR),  
$I$ the maximal ideal of $R$, $k(\eta)$ (resp. $k(0)$)
the fraction (resp. residue) field of $R$, $S:=\Spec R$,  
$\eta$ (resp. $0$) the generic (resp. closed) point  of $S$ and 
$\Omega:=\overline{k(\eta)}$ an algebraic closure of $k(\eta)$. 
Let $(G_{\eta},\cL_{\eta})$ be a polarized 
abelian variety over $k(\eta)$ with $\cL_{\eta}$ 
symmetric ample cubical, 
$\cG$ the N\'eron model 
of $G_{\eta}$ 
and $G:=\cG^0$ the identity component of $\cG$. 
A triple $(P,i,\cN)$ is called a 
{\it relative compactification of $\cG$} 
(abbr. {\it compactification of $\cG$})
(extending $(G_{\eta},\cL_{\eta})$) if 
\begin{enumerate}
\item[(rc1)]$P$ is an irreducible proper flat $S$-scheme;  
\item[(rc2)]
$i:\cG\hookrightarrow P$ is an open immersion with  
$P_{\eta}=i(\cG_{\eta})=i(G_{\eta})$; 
\item[(rc3)]$\cN$ is an ample invertible $\cO_P$-module 
with $i^*\cN_{\eta}\simeq\cL^{\otimes l}_{\eta}$ for some $l\geq 1$.
\end{enumerate}

 The purpose of this article is to prove the following theorems. 
\begin{thm}\label{thm:main thm}
If $\cG$ is semiabelian  over $S$,  
then 
there exists a relative compactification 
$(P,i,\cN)$ of $\cG$ 
 such that 
\begin{enumerate} 
\item[(a)] $P$ is Cohen-Macaulay;
\item[(b)] $i(\cG)=P\setminus\Sing(P_0)$ with   
$\codim_{P}\Sing(P_0)\geq 2$ where $\Sing(P_0)$ denotes 
 the singular locus of $P_0$;
\item[(c)] $i^*\cN$ is cubical with $i_*i^*\cN=\cN$;
\item[(d)] 
$\cG$ acts on $P$ so that $i$ is $\cG$-equivariant.
\end{enumerate} 
\end{thm}

\begin{thm}\label{thm:P=Proj A(cG cN)/intro}
If $(P',i',\cN')$ is another relative 
compactification of $\cG$ subject to (a)-(c), 
then $P\simeq P'$. 
\end{thm}

Theorem~\ref{thm:main thm} (resp. Theorem~\ref{thm:P=Proj A(cG cN)/intro}) 
is proved in \S\S~\ref{sec:str of Pl in pd case}-\ref{sec:proof of main thm} 
in the partially degenerate case 
(resp. \S~\ref{sec:uniqueness of P} in any case). 
See \cite{MN24} for Theorem~\ref{thm:main thm} in the totally degenerate case. 
We call $(P,i,\cN)$ in Theorem~\ref{thm:main thm} 
(as well as Theorem~\ref{thm:main thm Dedekind}) 
a {\it cubical compactification of $\cG$.}

The article is organized as follows: In \S\S~\ref{sec:partially degenerate case}-\ref{sec:deg data pd case}, we recall basic notions and theorems 
mainly from \cite{FC90}. 
In \S~\ref{sec:central extensions}, we study central extensions 
of a finite group $S$-sceme by $\bG_{m,S}$.  
In \S~\ref{sec:Deg data of Neron model}, we construct 
degeneration data of the N\'eron model $(\cG,i^*\cN)$ 
after a suitable base change.  In \S~\ref{sec:twisted Mumford families},  
we construct a compactification $P_{l}$ of $\cG$ 
by using the degeneration data of $(\cG,i^*\cN)$.   
In \S~\ref{sec:str of Pl in pd case}, 
we study the structure of the compactification $P_{l}$. 
In \S~\ref{sec:proof of main thm},
we prove Theorem~\ref{thm:main thm} 
in the partially degenerate case. 
In \S\S~\ref{sec:totally deg case revisited}-\ref{sec:cohomology of Pl}, we 
compute the cohomology groups of $P$. 
In \S~\ref{sec:uniqueness of P}, 
we prove Theorem~\ref{thm:P=Proj A(cG cN)} 
($=$ Theorem~\ref{thm:P=Proj A(cG cN)/intro} in any case). 
In \S~\ref{sec:rel compactification Dedekind}, we prove 
Theorem~\ref{thm:main thm Dedekind}\ 
($=$Theorems~\ref{thm:main thm}/\ref{thm:P=Proj A(cG cN)/intro} over a Dedekind domain) 
using \S~\ref{sec:cohomology of Pl}. 
In \S~\ref{sec:complex case}, 
we explain the complex case. 
In Appendix~\ref{sec:appendix algebraizability of hatA}
we prove algebraizability of the abelian part  
of the formal completion of $G$. 
In Appendix~\ref{sec:appendix deg data},  
we extend the degeneration data 
of $(G,\cL)$ from $(X,Y)$ to $(X,X)$. 
Throughout this article, we use the notation in \cite{MN24} 
unless otherwise mentioned. See \cite[\S\S~2.3/4.1/10.3/Notation~5.7]{MN24}. 
See also Notation~\ref{notation:Notation Rinit Sinit keta=Kmin(xi) for zeta}/\ref{notation:theta expansion in pd case}.

{\small{\it Acknowledgement.}\quad We are very grateful to Professor K.~Mitsui for his collaboration in \cite{MN24}.}

\section{The partially degenerate case}
\label{sec:partially degenerate case}

Let $R$, $I$, $k(\eta)$, $k(0)$ and $S$ be the same as in 
\S~\ref{sec:introduction}, 

Throughout this article, 
 let $(G,\cL)$ be a semiabelian $S$-scheme with $\cL$ 
{\it symmetric ample cubical} invertible sheaf  
and $G_0:=G\times_S0$.    
By definition, $G_0$ is an extension of an abelian variety $A_0$ 
by a $k(0)$-torus $T_0$ and $(G_{\eta},\cL_{\eta})$ 
is a polarized abelian variety 
over $k(\eta)$. In this section  we will
discuss the partially degenerate case
where the abelian part $A_0$ of $G_0$ is nontrivial. 
\S\S~\ref{subsec:Raynaud extensions}-\ref{subsec:Fourier series general}   
are based on \cite[II, 5]{FC90}. For an affine $S$-scheme $W$ 
of finite type, we denote by $T_{X,W}$ (resp. $\bG_{m,W}$) 
 a split $W$-torus $\Spec \Gamma(\cO_W)[w^x;x\in X]$ 
for some lattice $X$ (resp. $\Spec \Gamma(\cO_W)[w^x;x\in \bZ]$). 
Let $T$ be any $W$-torus. Then 
the pullback $T_U$ is a split $U$-torus 
for some finite \'etale cover $U$ of $W$ 
\cite[X, \S~1]{SGA3}. 
We write $\bG_{m,\bZ}$ by $\bG_m$ if no confusion is possible.

\subsection{Raynaud extensions}
\label{subsec:Raynaud extensions}
Let $S_n:=\Spec R/I^{n+1}$, 
$G_n:=G\times_SS_n$ and $\cL_n:=\cL\times_SS_n$.
Hence associated to $(G,\cL)$ are the group $S^{\wedge}$-scheme 
$G^{\wedge}=\projlim G_n$ and an invertible sheaf 
$\cL^{\wedge}=\projlim \cL\otimes_R R/I^{n+1}$. 
Since $T_0$ is 
a $k(0)$-torus, $T_0$ lifts to a unique 
multiplicative closed subgroup scheme $T_n$ of $G_n$ 
flat over $S_n$ for every $n$ 
by \cite[I, 2.2]{FC90}/\cite[IX, 3.6]{SGA3}, 
and hence 
$\hat{T}=\projlim T_n$ is a formal $S^{\wedge}$-torus. 
Since $T_n$ is a closed subgroup scheme of $G_n$ flat over $S_n$, 
we have a quotient group scheme $A_n:=G_n/T_n$ flat over $S_n$ 
and an exact sequence 
$0\to T_n\overset{\iota_n}{\to} G_n\overset{\pi_n}{\to} A_n\to 0$ 
by \cite[I,V$_{\new}$,10.1.2/10.1.3]{SGA3}. 
Since $A_0$ is an abelian $k(0)$-scheme by assumption, 
$A_n$ is a group $S_n$-scheme flat smooth 
and proper over $S_n$. 
 Hence $A_n$ is an abelian $S_n$-scheme.  
Let $\hat{A}:=\projlim A_n$. Then  
$\hat{A}$ is algebraizable:
\footnote{See Lemma~\ref{lemma:A is alg for nonsplit T0}. 
See also \cite[II, p.~34, lines 2-3]{FC90}. }
\begin{equation}\label{eq:A is alg for nonsplit T0}
\text{there exists an abelian $S$-scheme $A$ 
such that $A^{\wedge}\simeq\hat{A}$.}
\end{equation}  The formal $S^{\wedge}$-torus $\hat{T}$ 
is also algebraizable by 
$T:=\Hom(\underline{X}(T),\bG_{m,S})$, where 
$\underline{X}(T):=\underline{X}(\hat{T})$ is the unique lifting of 
$\underline{X}(T_0)$ as a locally constant \'etale sheaf on $S^{\wedge}$. 
Namely $T$ is an $S$-torus such that $T^{\wedge}\simeq \hat{T}$.

The group $S^{\wedge}$-scheme 
$G^{\wedge}$ fits into 
an exact sequence
\begin{equation}\label{eq:exact seq of Gwedge}
1\to T^{\wedge}\overset{\hat\iota}{\to} G^{\wedge}\, 
\overset{\hat\pi}{\to}\, A^{\wedge}\to 0.
\end{equation} 

Let $A^t:=\Pic^0_{A/S}$ be the dual abelian $S$-scheme. 
By \cite[VIII, 1.6]{SGA7}, there exist isomorphisms: 
$$\op{Ext}(A^{\wedge},T^{\wedge})\simeq 
\Hom(\underline{X}(T^{\wedge}),A^{t,\wedge})\simeq 
\Hom(\underline{X}(T),A^t)\simeq \op{Ext}(A,T).$$  
Hence the sequence (\ref{eq:exact seq of Gwedge}) 
is algebraizable by \cite[p.~34]{FC90}, that is,  
there exist a semiabelian $S$-scheme 
$\tG$ with $\tG^{\wedge}\simeq G^{\wedge}$ 
and a closed immersion $\iota_T:T\to\tG$ 
fitting into an exact sequence 
\begin{equation}\label{eq:exact seq of tG}
1\to T\overset{\iota_T}{\to} \tG\, 
\overset{\pi}{\to}\, A\to 0 
\end{equation}such that $\iota_T^{\wedge}=\hat\iota$ 
and $\pi^{\wedge}=\hat\pi$. Namely (\ref{eq:exact seq of Gwedge}) is 
the $I$-adic completion of (\ref{eq:exact seq of tG}).
We call $\tG$ (or Eq.~(\ref{eq:exact seq of tG})) 
the {\it Raynaud extension of $G$.}  
The Raynaud extension Eq.~(\ref{eq:exact seq of tG}) 
is functorial in $G$. See \cite[p.~35]{FC90}.

Let $\lambda(\cL_{\eta}):G_{\eta}\to G_{\eta}^t$ 
be the polarization morphism of $\cL_{\eta}$, 
$K(\cL_{\eta}):=\ker(\lambda(\cL_{\eta}))$ 
and $K(\cL)$ the closure of $K(\cL_{\eta})$ in $G$. 
Since $K(\cL)$ is flat quasi-finite over $S$, we have by \cite[IV, 4.C]{An73} 
the quotient $G/K(\cL)$ denoted by $G^t$, which is a semiabelian $S$-scheme 
whose closed fiber $G^t_0$ is an extension of an abelian $k(0)$-variety $A_0^t$ by a $k(0)$-torus $T_0^t$. Moreover $G^t$ is the connected N\'eron model 
of $G^t_{\eta}$. Hence $\lambda(\cL_{\eta})$ 
extends to a morphism $\lambda(\cL):G\to G^t$.  
In the same manner as 
Eqs.~(\ref{eq:exact seq of Gwedge})/(\ref{eq:exact seq of tG}), 
we have another exact sequence 
\begin{equation}\label{eq:exact seq of t(Gt)wedge}
1\to T^{t,\wedge}\overset{\hat{\iota}^{t}}{\to} 
G^{t,\wedge}\to A^{t,\wedge}\to 0, 
\end{equation}which is uniquely algebraized into  
\begin{equation}\label{eq:exact seq of t(Gt)}
1\to T^t\overset{\iota_{T^t}}{\to} \tG^t\, 
\overset{\pi^t}{\to}\, A^t\to 0. 
\end{equation} See Eq.~(\ref{diagram:commu diag of exact seqs}).  We call $\tG^t$ the {\it Raynaud extension of $G^t$.}

The quasi-finite morphism $\lambda(\cL):G\to G^t$ 
 induces an $S_n$-morphism 
$G_n\to G_n^t$ and hence 
an $S_n$-morphism $T_n\to A_n^t$, which  
 turns out to be trivial because 
an abelian scheme $A_n^t$ over $S_n$ has no rational curves in it. 
Hence $\lambda(\cL)$ induces morphisms $\lambda_{A_n}:A_n\to A_n^t$ and 
$\lambda_{T_n}:T_n\to T_n^t$. Thus we eventually obtain  
a commutative diagram of exact sequences:
\begin{equation}\label{diagram:formal commu diag of exact seqs}
\begin{diagram}
0&\rTo&T^{\wedge}&\rTo^{\iota_T^{\wedge}}&G^{\wedge}&\rTo^{\pi^{\wedge}}&A^{\wedge}&\rTo&0\\
&&\dTo^{\lambda_{T^{\wedge}}}&&\dTo^{\lambda_{G^{\wedge}}}&&\dTo^{\lambda_{A^{\wedge}}}\\
0&\rTo&T^{t,\wedge}&\rTo^{\iota_{T^t}^{\wedge}}&G^{t,\wedge}&\rTo^{\pi^{t,\wedge}}&A^{t,\wedge}&\rTo&0.
\end{diagram}
\end{equation}

Since Eqs.~(\ref{eq:exact seq of Gwedge})/(\ref{eq:exact seq of t(Gt)wedge}) 
are uniquely algebraizable, 
we have a commutative diagram of exact sequences which algebraizes   
Eq.~(\ref{diagram:formal commu diag of exact seqs}):
\begin{equation}\label{diagram:commu diag of exact seqs}
\begin{diagram}
0&\rTo&T&\rTo^{\iota_T}&\tG&\rTo_{\pi}&A&\rTo&0\\
&&\dTo^{\lambda_T}&&\dTo^{\lambda_{\tG}}&&\dTo^{\lambda_A}\\
0&\rTo&T^t&\rTo^{\iota_{T^t}}&\tG^t&\rTo^{\pi^t}&A^t&\rTo&0.
\end{diagram}
\end{equation}
Now we have two extensions of $A$ by $T^t$ arising 
from Eq.~(\ref{diagram:commu diag of exact seqs}):
\begin{enumerate}
\item[]the pushout $\tG':=\tG\times_S T^t/(\iota_T,\lambda_T)(T)$ 
of Eq.~(\ref{eq:exact seq of tG}) by $\lambda_T:T\to T^t$;
\item[]the pullback $\tG'':=\tG^t\times_{A^t}A$ of Eq.~(\ref{eq:exact seq of t(Gt)}) 
by $\lambda_A:A^\to A^t$.
\end{enumerate}where the quotient $\tG'$ is 
an $S$-scheme by \cite[IV, 4.C]{An73} because $(\iota_T,\lambda_T)(T)$ is 
a closed $S$-flat subgroup scheme of $\tG\times_S T^t$. 
We have a natural morphism $\Lambda:\tG'\to\tG''$ sending 
$(g,v)\mapsto(\lambda_{\tG}(g),\pi(g))$ for functorial points $g\in \tG$ and 
$v\in T^t$.
Hence we have a commutative diagram of exact sequences:
\begin{equation}\label{diagram:commu diag of ext of A by Tt}
\begin{diagram}
0&\rTo&T^t&\rTo^{e_{\tG}\times\id_T}&\tG'&\rTo_{\pi\circ\pr_1}&A&\rTo&0\\
&&\dTo^{\|}&&\dTo^{\Lambda}&&\dTo^{\|}\\
0&\rTo&T^t&\rTo^{\iota_{T^t}\times e_A}&\tG''&\rTo^{\pi^t\circ\pr_1}&A&\rTo&0
\end{diagram}
\end{equation}where $e_{\tG}$ (resp. $e_A$) is the identity element of 
$\tG$ (resp. $A$). Since $\lambda(\cL)$ is quasi-finite surjective, 
so are $\lambda_{G^{\wedge}}$, $\lambda_T$ and $\lambda_A$. 
Hence $\Lambda$ is quasi-finite surjective and hence bijective by 
Eq.~(\ref{diagram:commu diag of ext of A by Tt}). 
It follows that $\Lambda$ is an isomorphism: 
\begin{equation}\label{eq:isom of two extensions}
\Lambda:\tG\times_S T^t/(\iota_T,\lambda_T)(T)=\tG'
\simeq\tG''=\tG^t\times_{A^t}A.
\end{equation}

\subsection{Raynaud extensions in the split case}
\label{subsec:Raynaud extensions split case}
By taking a finite \'etale Galois cover 
$f:S^*\to S$ of $S$ if necessary, we may assume 
$T_0\times_{S_0}{0^*}$ is a split $k(0^*)$-torus, 
where $0^*$ is the closed point of $S^*$. 
In what follows by choosing $S^*$ as $S$, 
we assume that 
\begin{equation}\label{eq:T0 split}
\text{$T_0$ is a split $k(0)$-torus.}
\end{equation}
Then every $T_n$ is a split $S_n$-torus 
and hence $T^{\wedge}$ is a split $S^{\wedge}$-torus. 
Therefore there exists a lattice $X$ such that 
$T^{\wedge}=\Spf R[X]^{\wedge}$ such that $T^{\wedge}$ 
is algebraizable by a split $S$-torus $T:=\Spec R[X]$. 

By \cite[2.13~(2)]{MN24}
/\cite[II, 1.1 (ii)]{MB85}
there exists a unique cubical structure 
of $\cL$, which induces a cubical structure 
of $\cL^{\wedge}$.  Since $T^{\wedge}$ 
is a split torus, 
$\cL^{\wedge}_{|T^{\wedge}}$ is trivial, so that  
$\cL^{\wedge}$ descends by \cite[I, 7.2.2]{MB85} 
to a unique ample cubical
invertible sheaf $\cM^{\wedge}$ 
on $A^{\wedge}$, that is, 
$\cL^{\wedge}\simeq\pi^{\wedge,*}(\cM^{\wedge})$. 
Hence there exists an ample invertible sheaf $\cM$ on 
the abelian $S$-scheme $A$ whose $I$-adic completion is  
$\cM^{\wedge}$. 
Moreover $\underline{X}(T)=X_S$ and 
$\Ext(A,T)\simeq\Hom(X_S,A^t)$. Hence attached to 
the sequence (\ref{eq:exact seq of tG}), 
we have a classifying morphism ($=$ extension class) 
$c\in\Hom(X_S,A^t)$ which induces   
 a classifying morphism $c^{\wedge}$ 
of the sequence (\ref{eq:exact seq of Gwedge}) 
via the $I$-adic completion $\tG^{\wedge}\simeq G^{\wedge}$. \par 
Let $x\in X$. Since $X=\Hom_{\text{gr.sch.}}(T,\bG_{m,S})$, 
we have the pushout (pushforward)  torsor 
$\bF^{\times}_x:=\bG_{m,S}\times_S\tG/(x,\iota_T)(T)$ 
of the diagram (\ref{eq:exact seq of tG}) by $x$, that is,  
the quotient of $\bG_{m,S}\times_S\tG$ by 
the $S$-flat closed subgroup scheme $(x,\iota_T)(T)\ (\simeq T)$.  
Hence we have a commutative diagram of exact sequences:
\begin{equation}\label{eq: pushforward of exact seq of tG}
\begin{diagram}
1&\rTo&T         &\rTo^{\iota_T}&\tG &\rTo^{\pi}&A&\rTo& 0\\
 &    &\dTo^{x}&              &\dTo&          &\dTo^{\|}&&\\
1&\rTo&\bG_{m,S} &\rTo^{i_x}&\bF^{\times}_x&\rTo&A&\rTo& 0
\end{diagram}
\end{equation} Let $\bF_x$ be the line bundle on $A$ associated 
with $\bF_x^{\times}$. Since the classifying morphism 
of the $\bG_{m,S}$-torsor $\bF_x^{\times}$ is $c(x)$, 
we have $\bF_x=c(x)\in A^t(S)$ as a line bundle on $A$. 
By the diagram 
(\ref{eq: pushforward of exact seq of tG}), we have an isomorphism 
of $A$-schemes:  
\begin{equation}\label{eq:isom of tG}
\tG\simeq \bigoplus_{x\in X}\bF^{\times}_x. 
\end{equation}

Since $T_0$ is a split $k(0)$-torus, so is $T_0^t$ and hence $T^t$ 
(resp. $T^{t,\wedge}$) is a split $S$-torus 
(resp. a split $S^{\wedge}$-torus). Therefore there exists a lattice $Y$ such that $T^t=\Spec R[Y]$ and $T^{t,\wedge}=\Spf R[Y]^{I\op{-adic}}$. 
Hence $\underline{X}(T^t)=Y_S$ and 
$\Ext(A^t,T^t)\simeq\Hom(Y_S,A)$. Thus we have 
a classifying morphism ($=$ extension class) 
$c^t\in\Hom(Y_S,A)$ of Eq.~(\ref{eq:exact seq of t(Gt)}), 
whose $I$-adic completion is a classifying morphism $c^{t,\wedge}$ 
of Eq.~(\ref{eq:exact seq of t(Gt)wedge}).\par

There is a canonical morphism $\lambda(\cL):G\to G^t$ 
extending $\lambda(\cL_{\eta})$, which eventually induces 
morphisms $\lambda_A=\lambda_A(\cM):A\to A^t$ and 
$\lambda_T:T\to T^t$, each being finite flat surjective. 
Then 
$\lambda_T$ induces an injective homomorphism $\phi:Y\to X$ 
with finite cokernel. Now we recall the isomorphisms:
\begin{equation}\label{eq:c_ext and c_hom}
\begin{aligned}
c_{\ext}\in\Ext(A, T)&\simeq\Hom(X_S,{\cal Ext}^1(A,\bG_m))\simeq\Hom(X_S,A^t)\ni c,\\
c^t_{\ext}\in\Ext(A^t, T^t)&\simeq\Hom(Y_S,{\cal Ext}^1(A^t,\bG_m))\simeq\Hom(Y_S, A)\ni c^t
\end{aligned}
\end{equation}where $c_{\ext}$ (resp, $c^t_{\ext})$ is the class corresponding 
to $c$ (resp $c^t$) under the isomorphism. 
Note that ${\cal Ext}^1(A,\bG_m)$ (resp. ${\cal Ext}^1(A^t,\bG_m)$) is an fppf sheaf over $S$ representable by the $S$-scheme $A^t$ (resp. $A$). 
\footnote{See \S~\ref{subsec:fppf presheaf  hW is sheaf} for fppf sheaves.} The group 
$\Ext(A, T)$ is the set of isomorphism classes of commutative group $S$-schemes which are extensions of $A$ by $T$. The extension $c_{\ext}\in\Ext(A, T)$ induces, as a pushout (\ref{eq: pushforward of exact seq of tG}), 
a commutative group $S$-scheme 
$c(x)=\bF_x^{\times}\in{\cal Ext}^1(A,\bG_m)$ $(\forall x\in X(T))$ 
which is an extension 
of $A$ by $\bG_{m,S}$. Since $\bF_x^{\times}$ determines 
a unique line bundle $\bF_x\in A^t(S)$ on $A$, 
we may view $c(x)\in A^t(S)$. This is 
what Eq.~(\ref{eq:c_ext and c_hom}) means. 
By Eq.~(\ref{eq:isom of two extensions}), we have 
$\Hom(A,\lambda_T)(c_{\ext})=\Hom(\lambda_A,T^t)(c^t_{\ext})\in\Ext(A,T^t)$, 
which is translated into 
\begin{equation}\label{eq:a_phi=lambda_ct}
c\circ\phi=\lambda_A\circ c^t\in\Hom(Y_S,A^t).
\end{equation}

\subsection{Fourier series}
\label{subsec:Fourier series general}  
The ample cubical sheaf 
$\cL$ on $G$ induces a sheaf 
$\cL^{\wedge}$ on  $G^{\wedge}$, which  
descends to an ample cubical sheaf $\cM^{\wedge}$ 
on $A^{\wedge}$. Then   
$(A^{\wedge},\cM^{\wedge})$ is algebraized 
into a polarized abelian $S$-scheme $(A,\cM)$ with 
$\cM$ cubical.  Hence 
$\cL^{\wedge}\simeq(\pi^{\wedge})^*\cM^{\wedge}=(\pi^*\cM)^{\wedge}$ 
via the isomorphism $G^{\wedge}\simeq\tG^{\wedge}$. 
Since $T\subset\tG$, $T$ acts on 
$\tG$ keeping $\pi$ invariant. 
 Let $x\in X$ and 
\begin{equation}\label{eq:cO^x}
\cO^x:=\{h\in\cO_{\tG}; \tT_t^*h=x(t)h\ \text{for any functorial point 
$t$ of $T$}\},
\end{equation}where $\tT_t$ is translation of $\tG$ by $t$.
 Since $\cO^{x}$ is $\cO_A$-invertible by 
$\cO_A=\cO_{\tG}^{T\op{-inv}.}$, 
there exists an invertible 
$\cO_A$-module $\cO_{x}$ such that $\pi^{-1}\cO_{x}=\cO^{x}$. 
Since $T$ acts on $\bF_x$ by $t\cdot f=x(t)^{-1}f$ for any functorial point $t$ of $T$, we have $O_{-x}\simeq\cO_{A}(\bF_x)$\ $(\forall x\in X)$, so that 
$\pi_{*}(O_{\tG})\simeq\bigoplus_{x\in X}\cO_{x}$ 
by Eq.~(\ref{eq:isom of tG}). 
By the sign change of $c$, we have 
\begin{equation}\label{eq:-Fx=c(x)=Ox}
c(x)=\cO_{x}=\cO_A(-\bF_x)\in A^t(S), 
\end{equation} If we choose a local trivialization of 
$\tG$ as a split $U$-torus 
over an affine open subscheme $U$ of $A$ such that 
$\tG_U=\Spec \Gamma(\cO_U)[w^x_U; x\in X]$ 
$(w^x_Uw^y_U=w^{x+y}_U, \forall x,y\in X)$, then 
$\cO_{x,U}\simeq \cO^x_U$ is a subsheaf 
of $\cO_{\tG_U}$ generated by $w^{x}_U$.  
By the sign change of $c^t$, we have  
$\lambda_A\circ c^t=c\circ\phi$.

Via the isomorphism $G^{\wedge}\simeq\tG^{\wedge}$, we have 
$(\pi^{\wedge})_*\cO_{G^{\wedge}}\simeq(\pi_*\cO_{\tG})^{I\text{-adic}}
\simeq(\bigoplus_{x\in X}\cO_{x})^{I\text{-adic}}
$ as $\cO_{A^{\wedge}}$-modules.  
Let $\cM_{x}:=\cM\otimes_{\cO_A}\cO_{x}$\ $(x\in X)$. 
Then 
\begin{equation}\label{eq:isom of Lwedge and sum of Mxwedge}
(\pi^{\wedge})_*\cL^{\wedge}\simeq (\pi_*\pi^*\cM)^{I\text{-adic}}
=(\cM\otimes_{\cO_A}\pi_*\cO_{\tG})^{I\text{-adic}}
\simeq(\bigoplus_{x\in X}\cM_{x})^{I\text{-adic}},
\end{equation}which induces the composite homomorphism 
\begin{equation}\label{eq:composite homom}
\begin{aligned}
\Gamma(G,\cL)&\hookrightarrow\Gamma(G^{\wedge},\cL^{\wedge})=
\Gamma(A^{\wedge},(\pi^{\wedge})_*\cM^{\wedge})\\
&=(\bigoplus_{x\in X}\Gamma(A,\cM_{x}))^{I\text{-adic}}
\subset \prod_{x\in X}\Gamma(A^{\wedge},\cM^{\wedge}_{x}),
\end{aligned}
\end{equation}
whose $x$-th projection induces an $R$-homomorphism 
$$\sigma_x:
\Gamma(G,\cL)\hookrightarrow\Gamma(G^{\wedge},\cL^{\wedge})
\to\Gamma(A^{\wedge},\cM^{\wedge}_{x})
\overset{\simeq}{\to}\Gamma(A,\cM_{x}).$$
Every $\theta\in\Gamma(G,\cL)$
 is uniquely expressed as  
$\theta=\sum_{x\in X}\sigma_{x}(\theta)$ for some   
$\sigma_{x}(\theta)\in\Gamma(A,\cM_{x})$. 
\begin{defn}\label{defn:weight of cOx}
We define the $X$-weight of $\sigma_x$ by $\wt(\sigma_x)=x$. 
Hence $\wt(\cO_x^{\times})=\wt(\cM_x^{\times})
=\wt((\cM^{\otimes m})^{\times}_x)=x$, 
while $\wt(0)=0$. 
\end{defn}

Let 
$T_{\alpha}:A\to A$ be the translation of $A$ by $\alpha\in A(S)$ 
and $y\in Y$. By Eq.~(\ref{eq:a_phi=lambda_ct}), $T_{c^t(y)}^*\cM\simeq 
\cM\otimes_{\cO_A}\cO_{\phi(y)}$, while 
$T_{c^t(y)}^*\cO_{x}\simeq\cO_{x}$ by $\cO_{x}\in A^t(S)$. 
Hence as rigidified sheaves, 
\begin{equation}\label{eq:tSy linearized isom}
\begin{aligned}
T^*_{c^t(y)}(\cM)&=\cM_{\phi(y)}\otimes_{\cO_A} \cM(c^t(y)),\\ 
T^*_{c^t(y)}(\cO_{x})&=\cO_{x}
\otimes_{\cO_A} \cO_{x}(c^t(y)),\\  
T^*_{c^t(y)}(\cM_{x})&= 
\cM_{x+\phi(y)}\otimes_{\cO_A} \cM_{x}(c^t(y)).
\end{aligned}
\end{equation}

Since $\cM_{\eta}$ is a cubical 
$\bG_{m,\eta}$-torsor on $A$, 
so is the pullback $(c^t)^*\cM_{\eta}$ on $Y$. 
A trivial cubic $\bG_{m,\eta}$-torsor 
on $Y$ (resp. a trivial biextension of $Y\times X$ by 
$\bG_{m,\eta}$) is defined to be    
$1_Y=Y\times\bG_{m,\eta}$ 
(resp. $1_{Y\times X}=Y\times X\times\bG_{m,\eta}$). 
The Poincar\`e sheaf $\cP$ of $A$
is by definition the universal invertible sheaf on $A^t\times A$ 
for the functor $\Pic^0_A$ rigidified  along 
$(e_{A^t}\times_S A)\cup(A^t\times_S e_A)$ where $e_A$ and $e_{A^t}$ are 
the identity elements of $A$ and $A^t$ respectively. 
It is an invertible sheaf on 
$A^t\times_S A$ such that 
$\cP_{a\times A}\in\Pic^0(A)=A^t$ corresponds to $a\in A^t$. 
The $\bG_m$-torsor associated with $\cP$ is 
also a biextension of $A^t\times A$ by $\bG_{m,S}$, so that 
 $(c^t\times c)^*\cP_{\eta}^{\otimes(-1)}$ is a biextension of
$Y\times X$ by $\bG_{m,\eta}$. 

In what follows, for simplicity, we often identify 
an invertible sheaf with the $\bG_m$-torsor associated with it.

\section{Degeneration data} 
\label{sec:deg data pd case}

\S\S~\ref{subsec:deg data semiabelian}-\ref{subsec:split obj in DDample} of this section is based on \cite[II-III]{FC90}. See also \cite{R70}. 
\subsection{Degeneration data of a semiabelian scheme}
\label{subsec:deg data semiabelian} 
\begin{thm} \cite[II, 5.1, p.~44]{FC90}
\label{thm:degeneration data2} 
The $R$-module homomorphisms 
$\sigma_{x}:\Gamma(G,\cL)\to \Gamma (A,\cM_x)$ $(x\in X)$ 
have the following properties: 
\begin{enumerate}
\item\label{item:nonzero sigma} $\sigma_{x}\neq 0$\quad $(\forall x\in X)$;
\item\label{item:sigma}there exist a cubic 
trivialization {\small$\psi : 1_Y \to 
(c^t)^*\cM_{\eta}^{\otimes(-1)}$} and a trivialization of biextension 
{\small$\tau :1_{Y\times X}\to (c^t\times c)^*\cP_{\eta}^{\otimes(-1)}$}  
such that 
{\small\begin{equation}\label{eq:sigma psi tau}
\sigma_{x+\phi(y)}(\theta)=
\psi(y)\tau(y,x)T^*_{c^t(y)}(\sigma_{x}(\theta))\  
(\forall \theta\in \Gamma(G,\cL), \forall x\in X, \forall y\in Y);
\end{equation}} 
in particular, $\psi(0)=\tau(0,x)=\tau(y,0)=1$;  
\item\label{item:tau} $\tau(y,x)$ is bilinear in $y$ and $x$; 
\item\label{item:psi/tau} $\tau(y,\phi(z))=\tau(z,\phi(y))
=\psi(y+z)\psi(y)^{-1}\psi(z)^{-1}$
 $(\forall y,z\in Y)$;
\item\label{item:positive tau} $\tau(y,\phi(y))\in I$ $(\forall y\in Y\setminus\{0\})$;
\item\label{item:positive psi} for every $n\geq 0$, 
$\psi(y)\in I^n$ for all but finitely many $y\in Y$;  
\item\label{item:Gamma Geta L general}  
$\Gamma(G_{\eta},\cL_{\eta})$ is identified 
via Eq.~(\ref{eq:composite homom}) with 
the $k(\eta)$-vector subspace consisting of Fourier series 
$\theta=\sum_{x\in X}\sigma_x(\theta)$ satisfying 
\begin{equation}\label{eq:Gamma(geta,cLeta) psi tau}
\sigma_{x+\phi(y)}(\theta)=
\psi(y)\tau(y,x)T^*_{c^t(y)}(\sigma_{x}(\theta))\ \ 
(\forall x\in X, \forall y\in Y)
\end{equation} 
with $\sigma_{x}(\theta)\in \Gamma(\cA,\cM_x)\otimes_Rk(\eta)$.
\end{enumerate}
\end{thm}

\begin{defn}
\label{defn:sigma psi tau mgeq1}
For $m\geq 1$, we consider a semiabelian $S$-schme 
$(G,\cL^{\otimes m})$. By Theorem~\ref{thm:degeneration data2}, we have 
\begin{gather*}
\sigma^{(m)}:\Gamma(G,\cL^{\otimes m})\to \Gamma(A,\cM^{\otimes m}\otimes _{\cO_A}\cO_x),\ \phi^{(m)}=m\phi:Y\to X,\\
\psi^{(m)}:1_Y\to (c^t)^*\cM_{\eta}^{\otimes(-m)},\ 
\tau^{(m)}:1_{Y\times X}\to (c^t)^*\cP_{\eta}^{\otimes(-1)}
\end{gather*} 
 such that $(\sigma,\psi,\phi,\tau)=(\sigma^{(1)},\psi^{(1)},\phi^{(1)},\tau^{(1)})$ and  
{\small\begin{equation}\label{eq:sigmam psim taum}
\sigma^{(m)}_{x+m\phi(y)}(\theta)=
\psi^{(m)}(y)\tau^{(m)}(y,x)T^*_{c^t(y)}(\sigma^{(m)}_{x}(\theta))
\end{equation}} 
$(\forall \theta\in \Gamma(G,\cL^{\otimes m})\otimes_Rk(\eta), \forall x\in X, \forall y\in Y)$. \end{defn}

We will give $\sigma^{(m)}$, $\psi^{(m)}$ and $\tau^{(m)}$ explicitly 
in Corollary~\ref{cor:psim,taum}. 

\subsection{Split objects in $\DD_{\ample}$}
\label{subsec:split obj in DDample}

\begin{defn}\label{defn:split obj in DDample}
A {\it split semiabelian $S$-scheme} is 
a pair $(G,\cL)$ of a semiabelian $S$-scheme $G$ 
and a symmetric ample cubical 
invertible sheaf $\cL$ on $G$ 
such that 
$G_0$ is an extension of an abelian variety $A_0$ 
by a {\it split torus} $T_0$. 
 Let $G^t_{\eta}$ be the dual abelian variety of $G_{\eta}$,
$\lambda(\cL_{\eta}):G_{\eta}\to G^t_{\eta}$ the polarization (homo)morphism, 
$K(\cL_{\eta}):=\ker(\lambda(\cL_{\eta}))$, $K(\cL):=$ the closure of 
$K(\cL_{\eta})$ in $G$ and $G^t:=G/K(\cL)$. By \cite[p.~35, line~5]{FC90}/\cite[IV, 4.C]{An73}, $G^t$ is a semiabelian $S$-scheme and we have a surjective morphism $\lambda_G:G\to G^t$. By \cite[7.4/3]{BLR90}, $G^t$ is a connected N\'eron model of $G^t_{\eta}$.  
\end{defn}

To avoid repeating long definitions, let $\DD_{\ample}$ be   
the category of data consisting of \cite[III, \S~2, (1)-(7), p.~57]{FC90} 
satisfying a certain condition 
($=$ Definition~\ref{defn:split obj zeta in DDample}~(7)).
Here we only recall those split objects in $\DD_{\ample}$.  
\begin{defn}\label{defn:split obj zeta in DDample}
A {\it split object $\zeta$ over $S$ or $R$} 
in $\DD_{\ample}$ is 
a set of data
\begin{equation}\label{eq:split obj zeta in DDample}
\zeta=(\tG,A,T,X,Y,c,c^t,\iota,\lambda,\phi,\tau,\tcL,\psi,\cM)
\end{equation}
such that 
\begin{enumerate}
\item\label{item:tG c A At} 
$A$ is an abelian $S$-scheme, $T$ is a split $S$-torus and $\tG$ is
an extension of $A$ by $T$ with $\pi:\tG\to A$ the natural projection; 
the extension $\tG$ determines 
and is determined by a homomorphism $c:X\to A^t$ 
via the {\it negative} of push-out where $X=\Hom(T,\bG_{m,S})$ is 
the character group of $T$ with $T=\Spec R[X]$ 
and $A^t$ is the dual abelian scheme of $A$;
\item\label{item:phi inj} 
$\phi:Y\to X$ is a monomorphism of free abelian groups of equal  
rank;
\item\label{item:ct} $c^t:Y\to A$ is a homomorphism which determines and is determined by 
an extension $\tG^t$ of $A^t$ by $T^t$ via push-out 
where $T^t:=\Spec R[Y]$; 
\item\label{item:tcL and cM} $\cM$ is a cubical ample invertible sheaf on $A$ with $\tcL=\pi^*\cM$ and $\lambda\circ c^t=c\circ\phi$ for 
its polarization morphism 
$\lambda:=\lambda(\cM):A\to A^t$; 
\item\label{item:iota lifts ct} 
$\iota:Y\to \tG_{\eta}$ is a homomorphism lifting $c^t_{\eta}$;\ $\iota$ determines and is determined by 
a trivialization of biextension 
$\tau:1_{Y\times X}\simeq (c^t\times c)^*\cP_{A,\eta}^{\otimes (-1)}$;
\item\label{item:symmetric triv} $Y$ acts on $\tcL_{\eta}$ over $\iota$ in the way compatible with $\phi$;\footnote{See 
Eq.~(\ref{eq:action Sy on sigma_x}).}\ the action 
of $Y$ determines and is determined by 
a cubic trivialization 
$\psi:1_Y\simeq \iota^*\tcL_{\eta}^{\otimes (-1)}$ 
compatible with $\tau\circ(\id\times\phi):
1_{Y\times Y}\simeq (c^t\times c\circ\phi)^*\cP_{A,\eta}^{\otimes(-1)}$;
\item\label{item:positivity} $\tau(y,\phi(y))\in I$ 
$(\forall y\in Y\setminus\{0\})$\  
(positivity condition).  
\end{enumerate} 
\end{defn}

\begin{defn}\label{defn:symm}
A cubic trivialization $\psi:1_Y\simeq\iota^*\tcL_{\eta}$ is 
called {\it symmetric} if  $\psi(-y)=\psi(y)$ $(\forall y\in Y)$. 
If $\psi$ is symmetric, then $\psi(y)^2=\tau(y,\phi(y))$ $(\forall y\in Y)$. 
We say $\zeta$ is {\it symmetric} if 
$\cM$ is symmetric and $\psi$ is symmetric. 
Note that $\cM$ is symmetric iff $\tcL$ is symmetric because 
$\tcL=\pi^*\cM$ and $\tcL^{\wedge}\simeq
(\bigoplus_{x\in X}\cM_x)^{I\op{-adic}}$ 
with $[-\id_{A^{\wedge}}]^*\cM_x^{\wedge}=\cM_{-x}^{\wedge}$.  If $(\tG,\tcL)$ is 
the Raynaud extension of a semiabelian scheme $(G,\cL)$ 
with $\cL$ is symmetric, then $\cL^{\wedge}\simeq\tcL^{\wedge}$, 
$\tcL$ and $\cM$ are symmetric. 
\end{defn}

\begin{defn}\label{defn:FC ample datum} 
Let $(G,\cL)$ be 
a split semiabelian $S$-scheme. Then we obtain a split object $\zeta$ over $S$ 
 in $\DD_{\ample}$ by \S\S~\ref{subsec:Raynaud extensions split case}-\ref{subsec:Fourier series general}, where $\cM$ can be chosen symmetric 
by choosing $\cL\otimes (-\id_{\cG})^*\cL$ for $\cL$
and $\cM\otimes (-\id_A)^*\cM$ for $\cM$ if necessary. 
 Similarly by Theorem~\ref{thm:degeneration data2}, we define 
$\zeta_{m}$ by 
\begin{equation}\label{eq:split obj zetam in DDample}
\zeta_{m}:=(\tG,A,T,X,Y,c,c^t,\iota,\lambda^{(m)},
\phi^{(m)},\tau^{(m)},\tcL^{\otimes m},
\psi^{(m)},\cM^{\otimes m})
\end{equation}
where $\zeta_{1}=\zeta$.
It is also a split object  over $S$ in $\DD_{\ample}$, 
which encodes the data of $(G,\cL^{\otimes m})$. 
We denote $\zeta_m$ by $\FC(G,\cL^{\otimes m})$ and call it 
the {\it FC datum of $(G,\cL^{\otimes m})$}.  
In what follows, we assume that {\it $\cM$ and $\zeta_m$ are symmetric.} 
If $(G,\cL)$ is totally degenerate, then $\zeta_{m}$ is reduced to 
$(X,Y,\phi^{(m)},\tau^{(m)},\psi^{(m)})$ with $A=0$. 
See \cite[3.7/3.12]{MN24}. 
\end{defn}

\subsection{More about split objects}
\label{subsec:more about split objects}

Let $\zeta$ be the split object 
Eq.~(\ref{eq:split obj zeta in DDample}). 
Now we revisit  
Definition~\ref{defn:split obj zeta in DDample}~(\ref{item:iota lifts ct})-(\ref{item:symmetric triv}). First we note that for $y\in Y$,
\begin{align*}(c^t(y)\times c)^*\cP_{\eta}^{\otimes(-1)}
=\prod_{x\in X}\cO_{x}(c^t(y))_{\eta}^{\otimes(-1)}=y\times X\times 
\bG_{m,\eta}.
\end{align*}  Hence a trivialization of biextension 
\begin{equation}\label{eq:tau trivialization of biextension}
\tau:1_{Y\times X}:=(Y\times X)\times\bG_{m,\eta}\overset{\simeq}{\to} (c^t\times c)^*\cP_{\eta}^{\otimes(-1)}
\end{equation}
is given by $\tau(y,x):
=\tau(y\times x\times 1)\in\cO_{x}(c^t(y))_{\eta}^{\otimes(-1)}$ 
$(x\in X, y\in Y)$. 

Since $\pi_*(\cO_{\tG})=\bigoplus_{x\in X}\cO_x$, we have 
$\tG_{\eta}=\Spec(\bigoplus_{x\in X}\cO_{x,\eta})$, so that 
the homomorphism $\iota:Y\to \tG_{\eta}$ 
lifting $c^t_{\eta}$ 
is  given by $\coprod_{y\in Y}\iota(y)$ 
where $\iota(y)^*\in(\bigoplus_{x\in X}\cO_x(c^t(y))_{\eta}^{\otimes (-1)}$:
\begin{equation}\label{eq:iotay}\iota(y)^*=
\bigoplus_{x\in X}\tau(y,x)
\in(c^t(y)\times c)^*\cP_{\eta}^{\otimes(-1)}. 
\end{equation} Since $A$ is proper (and hence separated) over $S$, 
$c^t_{\eta}$ extends to $c^t:Y_S\to A$, which 
is uniquely determined by $c^t_{\eta}$ and hence by $\iota$.
Let $\psi:1_Y:=Y\times\bG_{m,\eta}\simeq (c^t)^*\cM_{\eta}^{\otimes(-1)}$ be 
the cubical trivialization given by $\zeta$. 
Since $\iota^*\tcL_{\eta}=\iota^*\pi^*\cM_{\eta}
=\coprod_{y\in Y}\cM(c^t(y))_{\eta}$,
we have $\psi(y)\in \cM(c^t(y))_{\eta}^{\otimes(-1)}
=\tcL(\iota(y))_{\eta}^{\otimes(-1)}$  $(y\in Y)$, so that 
$\psi:1_Y\simeq \iota^*\tcL_{\eta}^{\otimes(-1)}$. 

By Theorem~\ref{thm:degeneration data2} and 
Eq.~(\ref{eq:tSy linearized isom}), 
we define an isomorphism of $\cO_A$-modules
$\tS_{y}:T_{c^t(y)}^*\left(\bigoplus_{x\in X}\cM_{x,\eta}\right)\to
\bigoplus_{x\in X}\cM_{x,\eta}$ in terms of 
$\psi$ and $\tau$ by 
\begin{equation}\label{eq:action Sy on sigma_x}
\tS_{y}T_{c^t(y)}^*\left(\bigoplus_{x\in X}\sigma_{x}\right)
=\bigoplus_{x\in X}\sigma'_{x},\ 
\sigma'_{x+\phi(y)}=\psi(y)\tau(y,x)T_{c^t(y)}^*\sigma_{x} 
\end{equation}
where $\sigma_{x},\sigma'_{x}\in\cM_{x,\eta}$.   
Thus $\tS_{y}$ is compatible 
with $\phi$, $\psi$ and $\tau$ over $k(\eta)$ 
in the sense of Eq.~(\ref{eq:action Sy on sigma_x}). 
This is essentially the same as  
Definition~\ref{defn:split obj zeta in DDample}~(\ref{item:symmetric triv}). 
The constant $\psi(y)\tau(y,x)$ in Eq.~(\ref{eq:action Sy on sigma_x}) 
is derived from Theorem~\ref{thm:degeneration data2}
~(\ref{item:Gamma Geta L general}) 
because $\sigma_{x}(s)\in\Gamma(A,\cM_{x})$ and 
$\sigma_{x+\phi(y)}(s)=\psi(y)\tau(y,x)T_{c^t(y)}^*\sigma_{x}(s)$ 
if $s\in\Gamma(G,\cL)$.

Next we define an isomorphism of $\cO_A$-modules 
$\tS_{y}:T_{c^t(y)}^*\pi_*\cO_{\tG_{\eta}}\to\pi_*\cO_{\tG_{\eta}}$ 
in a similar fashion by 
\begin{equation}\label{eq:transf tSy on fx}
\tS_{y}T_{c^t(y)}^*(f_x)=\tau(y,x)T_{c^t(y)}^*f_x\ 
(f_x\in\cO_x), 
\end{equation}which makes the second equality 
of Eq.~(\ref{eq:tSy linearized isom}) precise.  This is  
derived from Eq.~(\ref{eq:action Sy on sigma_x}) as follows.
If $s\in\Gamma(G,\cL)$, then Eq.~(\ref{eq:action Sy on sigma_x}) 
implies  
\begin{gather}
\tS_{y}\tT_{c^t(y)}^*\left(\sum_{x\in X}\sigma_{x}(s)\right)
=\sum_{x\in X}\sigma_{x+\phi(y)}(s),\label{eq:action Sy on sigma_x(s)}\\ 
\sigma_{x+\phi(y)}(s)=\psi(y)\tau(y,x)T_{c^t(y)}^*(\sigma_{x}(s)). 
\label{eq:relation of sigma_x(s)}
\end{gather}
Assume $\sigma_z(s)\neq 0$ $(z=0,x)$. 
By Eq.~(\ref{eq:relation of sigma_x(s)})  
$\sigma_{z+\phi(y)}(s)\neq 0$ $(z=0,x)$. 
Let $V(\sigma_0(s)\sigma_{\phi(y)}(s))$ be the zero locus of 
the section $\sigma_0(s)\sigma_{\phi(y)}(s)$. Then  
on the open subset 
$U:=A\setminus V(\sigma_0(s)\sigma_{\phi(y)}(s))$ of $A$,  we have 
\begin{equation}\label{eq:transf on sigma_x sigma0}
\begin{aligned}
\tS_y T_{c^t(y)}^*\left(\frac{\sigma_x(s)}{\sigma_0(s)}\right)&:=
\frac{\tS_y T_{c^t(y)}^*\sigma_x(s)}{\tS_y T_{c^t(y)}^*\sigma_0(s)}
=\frac{\sigma_{x+\phi(y)}(s)}{\sigma_{\phi(y)}(s)}\\
&=\tau(y,x)T_{c^t(y)}^*\left(\frac{\sigma_x(s)}{\sigma_0(s)}\right). 
\end{aligned}
\end{equation}  Since $\cM_0$ is ample, 
$\Gamma(A,\cM_x)\neq 0$ $(\forall x\in X)$. 
Since $\cO_x$ is torsion free,   
Eq.~(\ref{eq:transf on sigma_x sigma0}) implies 
Eq.~(\ref{eq:transf tSy on fx}) everywhere on $A$. 
Hence the definition  of 
$\tS_{y}^*$ on $\bigoplus_{x\in X}\cM_{x}$  
by Eq.~(\ref{eq:action Sy on sigma_x}) 
 is consistent with 
that of $\tS_y$  
on $\bigoplus_{x\in X}\cO_x$ 
by Eq.~(\ref{eq:transf tSy on fx}).

\subsection{Split objects as degeneration data}
\label{subsec:split objects as deg data}
\begin{thm}
\label{thm:G' isom G'' iff FC same}
Let $(G',\cL')$ and $(G'',\cL'')$ be 
split semiabelian $S$-schemes. 
If $\zeta:=\FC(G',\cL')=\FC(G'',\cL'')$, then 
$(G',\cL')\simeq (G'',\cL'')$.
\end{thm}
\begin{proof}Let $(G^{(1)},\cL^{(1)}):=(G',\cL')$\ 
and $(G^{(2)},\cL^{(2)}):=(G'',\cL'')$. 
Let $\zeta:=(\tG,A,T,X,Y,c,c^t,\iota,\phi,\tau,\cL,\psi,\cM)$. 
Let $\pi:\tG\to A$ be the natural projection. Then we have an isomorphism
$\pi_*\cO_{\tG}\simeq\bigoplus_{x\in X}\cO_x$.
Since $(G^{(\nu)})^{\wedge}\simeq(\tG)^{\wedge}$, 
we have formal morphisms 
$\varpi^{(\nu)}:G^{(\nu),\wedge}\to A^{\wedge}$ $(\nu=1,2)$, 
so that 
\begin{equation}
\label{eq:isom Gformal G'formal TXformal}
\Spf (\varpi^{(\nu)})_*\cO_{G^{(\nu),\wedge}}=G^{(\nu),\wedge}
\simeq \tG^{\wedge}
=\Spf(\pi_*\cO_{\tG})^{\wedge}.
\end{equation}   Since $\zeta=\FC(G',\cL')=\FC(G'',\cL'')$, $c$ and $\cM$ 
are shared by $\tG'$ and $\tG''$, so that there exists an isomorphism 
$\hat\Phi:(G')^{\wedge}\simeq (G'')^{\wedge}$ such that 
$\hat\Phi^*(\varpi'')^*\cM^{\wedge}
=(\varpi')^*\cM^{\wedge}$ and 
$\hat\Phi^*(\varpi'')^{-1}\cO_x^{\wedge}
=(\varpi')^{-1}\cO_x^{\wedge}$ $(\forall x\in X)$.

Since $\lambda^{(m)}=\lambda_A(\cL^{\otimes m})=m\lambda$, 
we have $\phi^{(m)}(y)=m\phi(y)$ $(\forall m\geq 1, \forall y\in Y)$. 
By \cite[p.~50, line~6]{FC90}
 $\psi^{(2)}(y)=\psi(y)^2$ 
and $\tau^{(2)}(y,x)=\tau(y,x)$ 
because:
$$\psi(y)^2\tau(y,x)\tau(y,x')=\psi^{(2)}(y)
\tau^{(2)}(y,x+x')\ \ (\forall x,x'\in X, \forall y\in Y).
$$

Let $m=2^{l}$ $(l\geq 2)$ and $V^{(\nu)}_m:=\Gamma(G^{(\nu)}_{\eta},(\cL^{(\nu)}_{\eta})^{\otimes m})$. Then $\psi^{(m)}(y)=\psi(y)^{m}$ 
and $\tau^{(m)}(y,x)=\tau(y,x)$. 
By Theorem~\ref{thm:degeneration data2}~
(\ref{item:Gamma Geta L general}), 
$V^{(\nu)}_m$
 is identified with  
{\small 
\begin{equation*}
\label{eq:Geta Letam explicit}
\begin{aligned}
{}&\left\{
\theta^{(\nu)}=\sum_{x\in X}\sigma^{(\nu,m)}_x(\theta^{(\nu)});
\begin{matrix}
\sigma^{(\nu,m)}_{x+my}(\theta^{(\nu)})
=\psi(y)^m\tau(y,x)T^*_{c^t(y)}(\sigma^{(\nu,m)}_x(\theta^{(\nu)})) 
\\  
\sigma^{(\nu,m)}_x(\theta^{(\nu)})\in\Gamma(A^{\wedge},
(\cM^{\otimes m}\otimes_{\cO_A}\cO_x)^{\wedge})\otimes_Rk(\eta)\\
 (\forall x\in X, \forall y\in Y)
\end{matrix}\right\}.
\end{aligned}
\end{equation*}}

Hence  $\hat\Phi$ induces an isomorphism
$\Psi_m:V''_m\to V'_m$ with  
$\sigma^{(1,m)}_x(\Psi_m(\theta''))=\sigma^{(2,m)}_x(\theta'')$\ 
$(\forall \theta''\in V''_m, \forall x\in X)$.  
By the same proof as \cite[3.14]{MN24} from here, there is an isomorphism  
$\Phi:(G'_{\eta},(\cL'_{\eta})^{\otimes m})\simeq (G''_{\eta},(\cL''_{\eta})^{\otimes m})$ because $(\cL^{(\nu)})^{\otimes m}$ is very ample by $m\geq 4$. 
Since $G^{(\nu)}$ is the connected N\'eron model of $G^{(\nu)}_{\eta}$, 
$G'\simeq G''$. Since $\cL^{(\nu)}$ is a unique extension of $\cL^{(\nu)}_{\eta}$ to $G^{(\nu)}$ with $\Phi^*(\cL'_{\eta})=\cL''_{\eta}$, we have $(G',\cL')\simeq (G'',\cL'')$. 
\end{proof}

\begin{lemma}\label{lemma:H0(G,Lm) as pd Fourier}
Let $(G,\cL)$ be a split 
semiabelian $S$-scheme and  $\zeta:=\FC(G,\cL)$ 
as in Eq.~(\ref{eq:split obj zeta in DDample}). 
Then $\Gamma(G_{\eta},\cL_{\eta}^{\otimes m})$ is identified with 
the following $k(\eta)$-vector space for any $m\geq 1$:
{\small 
\begin{equation}
\label{eq:Geta Letam explicit 2nd}
\begin{aligned}
\left\{
\theta=\sum_{x\in X}\sigma^{(m)}_x(\theta);
\begin{matrix}\sigma^{(m)}_{x+my}(\theta)
=\psi(y)^m\tau(y,x)T^*_{c^t(y)}(\sigma^{(m)}_x(\theta)) 
\\  
\sigma^{(m)}_x(\theta)\in\Gamma(A^{\wedge},
(\cM^{\otimes m}\otimes_{\cO_A}\cO_x)^{\wedge})\otimes_Rk(\eta)\\
(\forall x\in X, \forall y\in Y)
\end{matrix}\right\}.
\end{aligned}
\end{equation}}
\end{lemma}
\begin{proof}This is a partially degenerate counterpart 
(abbr. pd-counterpart) 
of \cite[4.21]{MN24}. 
We choose a flat projective $S$-scheme $P$ with $P\supset G$ and 
$P_{\eta}=G_{\eta}$ constructed 
by \cite[\S~3]{AN99}/\cite[\S~4]{Nakamura99}. 
Recall that for proving \cite[4.21]{MN24}, we needed 
\cite[4.16/4.19-4.21]{MN24} and \cite[3.9]{Nakamura99}, 
all of whose pd-counterparts   
including \cite[4.10]{Nakamura99}
are proved in parallel to the totally degenerate case if  
we choose $\zeta=\FC(G,\cL)$ for $\xi$ of \cite[4.16]{MN24}. 
To prove Lemma, it suffices to check the pd--counterpart 
of \cite[4.20]{MN24} by applying \cite[3.9/4.10~(3)]{Nakamura99}.
Let $t$ be a uniformizer of $R$. 
With the notation of \cite[4.20]{MN24}, 
we set $\eta_m=\zeta_{x-m\alpha,\alpha}(\xi^*_{\alpha})^m$, 
$\gamma_m(x)=t^{v_t(\eta_m(x))}$ and 
$c_x(\theta)=\sigma_x(\theta)/\gamma_m(x)$, where 
we regard $c_x(\theta)$ (resp. $\gamma_m(x)$) 
as an element of $X$-weight $x$ (resp. zero). 
Let $\bar\psi(y):=t^{-v_t(\psi(y))}\psi(y)$ and  
$\bartau(y,x):=t^{-v_t(\tau(y,x))}\tau(y,x)$. As in 
\cite[4.20]{MN24}, 
\begin{align*} 
&c^{(m)}_{x+my}(\theta)=\bar\psi(y)^m
\bartau(y,x)T^*_{c^t(y)}(c^{(m)}_x(\theta))
\\
&\Leftrightarrow 
c^{(m)}_{x+my}(\theta)\gamma_m(x+my)\vartheta^m
=T^*_{c^t(y)}(c^{(m)}_x(\theta)\gamma_m(x)\vartheta^m)\\
&\Leftrightarrow \sigma_{x+my}(\theta)=\psi(y)^m\tau(y,x)T^*_{c^t(y)}(\sigma_x(\theta)).
\end{align*}  
By \cite[4.19~(3)]{Nakamura99} and Nakayama's lemma, 
we obtain the equality $\overset{*}{=}$ below, 
so that we see via Eq.~(\ref{eq:composite homom})
\begin{align*}
\Gamma&(G{_\eta},\cL_{\eta}^{\otimes m})
=\Gamma(P,\cL^{\otimes m})\otimes_Rk(\eta)=\Gamma(P^{\wedge},(\cL^{\wedge})^{\otimes m})\otimes_Rk(\eta)\\
&\overset{*}{=}\left\{
\theta=\sum_{x\in X}c^{(m)}_x(\theta)\eta_m(x); 
\begin{matrix} c^{(m)}_{x+my}(\theta)
=\bar\psi(y)^m\bartau(y,x)T^*_{c^t(y)}(c^{(m)}_x(\theta)) 
\\  
c^{(m)}_x(\theta)\in\Gamma(A^{\wedge},
(\cM^{\otimes m}\otimes_{\cO_A}\cO_x)^{\wedge})\otimes_Rk(\eta)\\
(\forall x\in X, \forall y\in Y)
\end{matrix}\right\}\\
&=\left\{
\theta=\sum_{x\in X}\sigma^{(m)}_x(\theta);
\begin{matrix}\sigma^{(m)}_{x+my}(\theta)
=\psi(y)^m\tau(y,x)T^*_{c^t(y)}(\sigma^{(m)}_x(\theta)) 
\\  
\sigma^{(m)}_x(\theta)\in\Gamma(A^{\wedge},
(\cM^{\otimes m}\otimes_{\cO_A}\cO_x)^{\wedge})\otimes_Rk(\eta)\\
 (\forall x\in X, \forall y\in Y)
\end{matrix}\right\}
\end{align*}
Hence we have Lemma.\end{proof}

\begin{cor}\label{cor:psim,taum}
$\psi^{(m)}=\psi^m$, $\tau^{(m)}=\tau$ and $\phi^{(m)}=m\phi$ 
$(\forall m\geq 1)$.
\end{cor}
\begin{proof}This follows from Theorem~\ref{thm:degeneration data2}/Definition~\ref{defn:sigma psi tau mgeq1}/Lemma~\ref{lemma:H0(G,Lm) as pd Fourier}. 
\end{proof}

\section{Central extensions}
\label{sec:central extensions}
\subsection{Extensions of a constant finite group}
\label{subsec:central ext}
Let $\Lambda$ be a commutative group and 
$H$ be a finite commutative group acting 
on $\Lambda$. We denote the action of $H$ on $\Lambda$ 
by $\tau:(x,a)\mapsto x(a)$ $(x\in H,a\in\Lambda)$.
\begin{defn}\label{defn:extension of H by Lambda}
Let $E$ be a (not necessarily commutative) group. 
Then $E$ is called an {\it extension of $H$ by $\Lambda$} if 
\begin{enumerate}
\item[($\alpha$)]\ $\Lambda$ is a normal subgroup of $E$ such that 
$E/\Lambda\simeq H$; 
\item[($\beta$)]\ if $\zeta(x)$ be a lifting of $x\in H$ 
to $E$ with $\zeta(0)=1_E$ the identity element of $E$, 
then $\ad(w)=\tau$, that is, $\ad(w)(x,a):=\zeta(x)\cdot a\cdot\zeta(x)^{-1}=x(a)=\tau(x,a)$ $(x\in H,a\in\Lambda)$, where $a\cdot b$ is 
the product of $a,b\in E$. 
\end{enumerate} If $E$ is an extension of $H$ by $\Lambda$, then 
$E\simeq H\times\Lambda$ as sets. If moreover $\Lambda$ is contained in the center of $E$, then $E$ is called a {\it central extension of $H$ by $\Lambda$.} 
If $H$ acts trivially on $\Lambda$, then 
$E$ is a central extension of $H$ by $\Lambda$ in view of 
($\beta$).
\end{defn}

\begin{defn}\label{defn:Ext(H,Lambda)}
Two extensions $E$ and $E'$ of $H$ by $\Lambda$ are isomorphic (as extensions) if there is an isomorphism $f:E\to E'$ as groups such that 
$f_{|\Lambda}=\id_{\Lambda}$ and the homomorphism $\barf:H\overset{\simeq}{\to} E/\Lambda\overset{\simeq}{\to} H$ induced from $f$ is $\id_H$.  
We define $\Ext(H,\Lambda)$ to be the set of all isomorphism classes of extensions of $H$ by $\Lambda$. 
\end{defn}
  
\begin{thm}\label{thm:Ext = H^2}Let $H$ be a finite commutative group and 
$\Lambda$ a commutative group with $H$-action on it. Then  
$\Ext(H,\Lambda)\simeq H^2(H,\Lambda)\ \ \mathrm{(}bijective\mathrm{)}.$ 
\end{thm}
\begin{proof}
Here we recall how to define a map  
$\sigma:\Ext(H,\Lambda)\to H^2(H,\Lambda)$. 
See \cite[1.2.4]{NSW99} for a complete proof of Theorem~\ref{thm:Ext = H^2} including the case where $H$ is a possibly non-commutative group. 
In what follows we denote the product of $H$ additively.

Let $E$ be an extension of $H$ by $\Lambda$. 
Then associated with $E$,  
we define a $2$-cocycle $\sigma(E)=(\phi(x,y):x,y\in H)\in H^2(H,\Lambda)$ 
with $\phi(x,y)\in\Lambda$
as follows. We choose a lifting $\zeta(x)$ to $E$ of each 
$x\in H$ with $\zeta(0)=1_E$. 
Then there exists $\phi(x,y)\in\Lambda$  such that 
$\zeta(x)\cdot \zeta(y)=\phi(x,y)\cdot \zeta(x+y)$ $(\forall x,y\in H)$.  
Since the group law of $E$ is associative, we have 
$(\zeta(x)\cdot \zeta(y))\cdot \zeta(z)
=\zeta(x)\cdot (\zeta(y)\cdot \zeta(z))$ 
$(\forall x,y,z\in H)$. Hence we have by 
Definition~\ref{defn:extension of H by Lambda}~($\beta$)  
\begin{equation}\label{eq:two cocycles}
\begin{aligned}
\phi(x,y)\cdot\phi(x+y,z)&=x(\phi(y, z))\cdot\phi(x, y+z),\\
\phi(x,0)=\phi(0,&x)=1\ \ (\forall x,y,z\in H).
\end{aligned}
\end{equation}Hence $\sigma(E)$ is a $2$-cocycle in the sense that 
$\phi$ satisfies Eq.=(\ref{eq:two cocycles}). 
Let $\sigma(E):=(\phi(x,y); x,y\in H)$. 
If $E$ splits, that is, $E\simeq H\times\Lambda$ as groups, then there exists 
a $1$-cochain $\psi:=
(\psi(x);x\in H)$ with $\psi(x)\in\Lambda$ such that 
$\tzeta(x):=\psi(x)\zeta(x)$ 
trivializes $E$ in the sense that 
$\tzeta(x)\cdot\tzeta(y)=\tzeta(x+y)$ $(\forall x,y\in H)$, 
whence $E\simeq H\times\Lambda$ as groups. 
In this case, $\phi$ is a 
$2$-coboundary:
\begin{equation}\label{eq:2 coboundary}
\phi(x,y)=\psi(x+y)\cdot\psi(x)^{-1}\cdot x(\psi(y))^{-1}.
\end{equation} 

The second group cohomology $H^2(H,\Lambda)$is defined by 
$$H^2(H,\Lambda)=Z^2(H,\Lambda)/B^2(H,\Lambda),
$$
where 
\begin{align*}
Z^2(H,\Lambda)&:=\left\{(\phi(x,y);x,y\in H); 
\begin{matrix}
\phi(x,y)\in\Lambda\ (\forall x,y\in H)\\
 \text{Eq.~(\ref{eq:two cocycles}) holds}\end{matrix}
\right\},\\
B^2(H,\Lambda)&:=\left\{
\begin{pmatrix}\psi(x+y)\cdot\psi(x)^{-1}\cdot(x(\psi(y))^{-1}\\ 
x,y\in H
\end{pmatrix}; \psi(x)\in\Lambda\ (\forall x\in H) 
\right\}.
\end{align*}For $\phi\in Z^2(H,\Lambda)$, we denote 
the cohomology class of $\phi$ in $H^2(H,\Lambda)$ by $[\phi]$.  

Thus we define a map $\sigma:\Ext(H,\Lambda)\to H^2(H,\Lambda)$ 
by $$\sigma(E)=(\phi(x,y);x,y\in H)\mod B^2(H,\Lambda),$$ 
which is a bijection by \cite[1.2.4]{NSW99}. 
\end{proof}

Now we consider the case where $H=\bZ/n\bZ$.  

\begin{thm}
\label{thm:H2(H,Lam)}
Let $H=\bZ/n\bZ$ be a cyclic group of order $n$ with $h$ generator, 
and $\Lambda$ a commutative group with $H$-action on it. 
Let $\Lambda^H$ be the subgroup of $\Lambda$ 
of all $H$-invariants of $\Lambda$, 
and $N_H(\Lambda):=\{\prod_{x\in H}x(\lambda); \lambda\in\Lambda\}$. Then 
$\Omega:H^2(H,\Lambda)\simeq\Lambda^H/N_H(\Lambda)$, where 
$\Omega$ and its inverse $\Psi$ are given as follows:
\begin{gather*}
\Omega(\phi)=\sum_{\nu=0}^{n-1}\phi(\nu h,h),\ 
\Psi(x)(ih,jh)=([\frac{i+j}{n}]-[\frac{i}{n}]-[\frac{j}{n}])x\\ 
(i,j\in[0,n-1], x\in\Lambda^H).
\end{gather*} 
\end{thm}

See \cite[pp.~105-108]{Saito97}.  
See also \cite[Chap.~I, 1.9.11]{NSW99}. 

\begin{cor}\label{cor:H^2(GmS)}
Let $R$ be a CDVR and $S=\Spec R$. 
Let $H=\bZ/n\bZ$ and $\Lambda=\bG_m(S)$. 
If $H$ acts on $\bG_m(S)$  trivially, then  
\begin{enumerate}
\item\label{item:isom GmS/GmS_n} $H^2(H, \bG_m(S))=\bG_m(S)/\bG_m(S)^n$ 
where $\bG_m(S)^n:=\{a^n;a\in\bG_m(S)\}$;
\item\label{item:pullback trivialized} let 
$E\in\Ext(H_S,\bG_{m,S})$, $\phi_E$ a 2-cocycle of $E$, 
$a:=\Omega(\phi_E)\in\bG_m(S)$, $R_L:=R[x]/(x^n-a)$ and $W=\Spec R_L$;
then the pullback $E_W$ of $E$ to $W$ splits, that is, 
$E_W\simeq \bG_{m,W}\oplus H_W$ as group $W$-schemes.
\end{enumerate} 
\end{cor}
\begin{proof}(\ref{item:isom GmS/GmS_n}) follows 
from Theorem~\ref{thm:H2(H,Lam)}. 
We shall prove (\ref{item:pullback trivialized}).  
 Let $\zeta(i)$ be a lifting of $i\in H$ to $E$ with $\zeta(0)=1$ and 
$\phi_E(i,j):=\zeta(i)\zeta(j)/\zeta(i+j)\in\bG_m(S)$. Then   
$\phi_E:=(\phi_E(i,j);i,j\in H)$ is a $2$-cocycle of $E$. 
Since $\Phi(\phi_E)=a$, we may assume $\phi_E(i,j)=\Psi(a)(i,j)$ 
by Theorem~\ref{thm:H2(H,Lam)}. 
Let $\alpha\in R_L$ such that $\alpha^{-n}=a$. Then we have 
\begin{align*}
\phi_E(i,j)&=\Psi(a)(i,j)=a^{[\frac{i+j}{n}]-[\frac{i}{n}]-[\frac{j}{n}]}\\
&=\alpha^{n(\frac{i+j}{n}-[\frac{i+j}{n}])}\alpha^{-n(\frac{i}{n}-[\frac{i}{n}])}\alpha^{-n(\frac{j}{n}-[\frac{j}{n}])}.
\end{align*}  Let $\tzeta(i)=\zeta(i)\alpha^{n(\frac{i}{n}-[\frac{i}{n}])}$. 
Then $\tzeta(i)\tzeta(j)=\tzeta(i+j)$, whence $E_W$ splits. 
\end{proof}

\subsection{Extensions of a finite group $S$-scheme by $\bG_{m,S}$}
\label{subsec:central ext of finite gr scheme by Gm}
Let $R$ be a CDVR and $S=\Spec R$. 
Let $H$ be a finite $S$-scheme, and $E$ a $\bG_{m}$-torsor over $H$. 
That is, $\pi:E\to H$ is a smooth $S$-morphism such that 
there exists a Zariski cover $(U_i;i\in I)$ of $H$ 
with $E_{U_i}\simeq U_i\times_S\bG_{m,S}$.
Since $H$ is finite flat over $S$, this implies 
$E\simeq H\times_S\bG_{m,S}$ as $S$-schemes. 
Hence we have a section $\zeta:H\to E$ as $S$-schemes 
with $\zeta(0)=1_E$ the identity element of $E$, 
which we call a lifting of $H$ to $E$. Since $H$ is affine, so is $E$ because 
$E\simeq H\times_S\bG_{m,S}$ as $S$-schemes. 
Let $D:=\Gamma(\cO_H)$. Then $D$ is the Hopf algebra of $H$.

\begin{defn}\label{defn:extension group scheme of H by GmS} 
Let $H$ be a finite flat commutative group $S$-scheme 
acting on $\bG_{m,S}$, and 
$\tau:H\times_S\bG_{m,S}\to\bG_{m,S}$ the action of $H$ on $\Lambda$.  
A $\bG_m$-torsor $\pi:E\to H$  
is called an {\it extension of $H$ by $\bG_{m,S}$} if 
\begin{enumerate}
\item[(a)]\ $E$ is a group $S$-scheme;
\item[(b)]\ $\bG_{m,S}$ is a closed normal group $S$-subscheme of $E$ 
(flat over $S$)  such that 
the quotient $Q:=E/\bG_{m,S}$ is the group $S$-scheme $H$;\footnote{
To be more precise, the fppf quotient sheaf $Q:=E/\bG_{m,S}$ 
is representable by a group $S$-scheme by \cite[IV, 4.C]{An73}. 
See \S~\ref{subsec:fppf presheaf  hW is sheaf} for fppf shaves. }
\item[(c)] if $w:H\to  E$ is a lifting of $H$ to $E$ 
with $\zeta(0)=1_E$, then $\ad(\zeta)=\tau$, where $\ad(\zeta):H\times_S\bG_{m,S}\to \bG_{m,S}$ is the morphism sending 
$(x,a)\mapsto \zeta(x)\cdot a\cdot \zeta(x)^{-1}$.
\end{enumerate}      
If $\bG_{m,S}$ is contained in the center of $E$, then $E$ is called 
a {\it central extension of $H$ by $\bG_{m,S}$}.  Note that 
if $H$ acts trivially on $\bG_{m,S}$, then   
$E$ is a central extension of $H$ by $\bG_{m,S}$ in view of (c).
\end{defn}

\begin{defn}\label{defn:Ext HS by GmS}
Two extensions $E$ and $E'$ of $H$ by $\bG_{m,S}$ are isomorphic (as extensions) if there is an isomorphism $f:E\to E'$ as group $S$-schemes such that 
$f_{|\bG_{m,S}}=\id_{\bG_{m,S}}$ and the homomorphism $\barf:H\overset{\simeq}{\to} E/\Lambda\overset{\simeq}{\to} H$ induced from $f$ is $\id_H$.  
Let $\Ext(H,\bG_{m,S})$ be the set of isomorphism classes of extensions 
of $H$ by $\bG_{m,S}$.
\end{defn}

Let $H^m:=H\times_S\cdots\times_SH$ ($m$-times), and 
let $C^n(H,\bG_{m,S})$ be the set of all $S$-morphisms 
$\phi:H^n\to \bG_{m,S}$. Hence $C^n(H,\bG_{m,S})=\bigotimes_{i=1}^nD^{\times}$

\begin{defn}\label{defn:2 cocycle coboundary}Let $\phi\in(D\otimes_RD)^{\times}$.We call $\phi\in C^2(H,\bG_{m,S})$ a {\it 2-cocycle} 
if $\phi$ satisfies the condition:for any $S$-scheme $T$, 
\begin{equation}\label{eq:two cocycles x,y,z in H(T)}
\begin{aligned}
\phi(x,y)\cdot\phi(x+y,z)&=x(\phi(y, z))\cdot\phi(x, y+z),\\
\phi(x,0)=\phi(0,&x)=1\ \ (\forall x,y,z\in H(T)).
\end{aligned}
\end{equation}
 Next we call $\phi$ a {\it 2-coboundary} iff 
there exists $\psi\in C^1(H,\bG_{m,S})$ such that 
$\phi=d^1\psi$, that is,  for any $S$-scheme $T$, 
\begin{equation}\label{eq:coboundary d1psi}
\phi(x,y)=d^1\psi(x,y)=\psi(x+y)\psi(x)^{-1}(x(\psi(y)))^{-1}\  
(\forall x,y\in H(T)). 
\end{equation}

We define $Z^2(H,\bG_{m,S})$ (resp. $B^2(H,\bG_{m,S})$) 
to be the subset of $C^2(H,\bG_{m,S})$ of all 2-cocycles 
(resp. 2-coboundaries). We define 
$$H^2(H,\bG_{m,S})=Z^2(H,\bG_{m,S})/B^2(H,\bG_{m,S}).
$$

For $\phi\in Z^2(H,\bG_{m,S})$, we denote the class 
of $\phi$ in $H^2(H,\bG_{m,S})$ by $[\phi]$. 
Since $\Hom_S(H,\bG_{m,S})=\Hom_R(R[w,w^{-1}],D)=D^{\times}$, we have 
$\phi\in C^2(H,\bG_{m,S})=(D\otimes_RD)^{\times}$ and 
$\psi\in C^1(H,\bG_{m,S})=D^{\times}$. 
\end{defn}

\begin{lemma}\label{lemma:equi cond ass and comm}
Let $H$ be a finite flat commutative group $S$-scheme and 
$E$ is a $\bG_{m,S}$-torsor on $H$, 
$\zeta:H\to E$ a lifting of $H$ to $E$ and 
$\phi\in C^2(H,\bG_{m,S})=(D\otimes_RD)^{\times}$.  
Let $T$ any $S$-scheme, $x,y\in H(T)$, $a, b\in\bG_{m,S}(T)$ and  
$\phi(x,y):=(x,y)^*\phi\in\Gamma(\cO_T)^{\times}$. 
We define a multiplication of $E$ by
\begin{equation}\label{eq:mult of E by phi}
(\zeta(x)\cdot a)\cdot (\zeta(y)\cdot b)=\zeta(x+y)\cdot 
\phi(x,y)\cdot a\cdot x(b).
\end{equation}Then 
\begin{enumerate}\item\label{item:E group scheme} the following are equivalent:
\begin{enumerate}
\item\label{item:E associative} 
multiplication Eq.~(\ref{eq:mult of E by phi}) of $E$ is associative;
\item\label{item:phi 2 cocycle}$\phi\in Z^2(H,\bG_{m,S})$;
\item\label{item:2cocyle any T} Eq.~(\ref{eq:two cocycles x,y,z in H(T)}) 
is true for any $S$-scheme $T$;
\item\label{item:2cocyle univ}Eq.~(\ref{eq:two cocycles x,y,z in H(T)}) 
 is true for $T=H$ and $x=y=z=\id_H$;
\end{enumerate}
\item\label{item:E group scheme splits}the following are equivalent:
\begin{enumerate}
\item\label{item:E splits as gr scheme}$E$ splits, that is, $E\simeq H\times\bG_{m,S}$ as group $S$-schemes;
\item\label{item:phi 2 coboundary}$\phi\in B^2(H,\bG_{m,S})$;
\item\label{item:phi coboundary}there exists $\psi\in C^1(H,\bG_{m,S})$ such that Eq.~(\ref{eq:coboundary d1psi}) is true for any $S$-scheme $T$;
\item\label{item:phi coboundary univ}
there exists $\psi\in  C^1(H,\bG_{m,S})$ such that Eq.~(\ref{eq:coboundary d1psi}) is true for $T=H$ and $x=y=\id_H$.
\end{enumerate}
\end{enumerate}
\end{lemma}
\begin{proof}Equivalence of 
(\ref{item:E associative})-(\ref{item:2cocyle any T}) is clear. 
Let $T$ be any $S$-scheme 
and $x,y,z\in H(T)$. The pullback of 
Eq.~(\ref{eq:two cocycles x,y,z in H(T)}) for (\ref{item:2cocyle univ}) 
by the morphism $(x,y,z):T\to H\times_SH\times_SH$ 
 is Eq.~(\ref{eq:two cocycles x,y,z in H(T)}) for (\ref{item:2cocyle any T}). 
Hence (\ref{item:2cocyle univ} implies (\ref{item:2cocyle any T}). 
Since (\ref{item:2cocyle any T}) clearly implies (\ref{item:2cocyle univ}, this proves (\ref{item:E group scheme}). 
(\ref{item:E group scheme splits}) is proved similarly. 
\end{proof}

\begin{defn}\label{defn:univ H2(H Gm)}
Let $D=\Gamma(\cO_H)$ and $n=\rank_RD$. Then we define 
\begin{align*}
Z^2_{\univ}(H,\bG_{m,S})&=\left\{\left(\phi(x,y);x,y\in\bZ\id_H\right);
\begin{matrix}\phi(x,y)\in D^{\times}\\
\text{Eq.~(\ref{eq:two cocycles x,y,z in H(T)})  
holds}
\end{matrix}
\right\},\\
B^2_{\univ}(H,\bG_{m,S})&=\left\{\begin{matrix}
\left(\psi(x+y)\psi(x)^{-1}(x(\psi(y))^{-1}; x,y\in\bZ\id_H\right)\\
 \psi(x)\in D^{\times}
\end{matrix}
\right\},\\
H^2_{\univ}(H,\bG_{m,S})&=Z^2_{\univ}(H,\bG_{m,S})/B^2_{\univ}(H,\bG_{m,S}),
\end{align*}where we note $n\id_H=0$ in $H(H)$ by Deligne's theorem \cite[p.~4]{OT70}.
\end{defn}

The following is proved in the same manner as Theorem~\ref{thm:Ext = H^2}:
\begin{thm}\label{thm:Ext = H^2 group schemes}
Let $H$ be a finite flat commutative group $S$-scheme 
acting on $\bG_{m,S}$. Then  
we have a bijection by the map $E\mapsto\sigma(E)$,
\begin{equation*}
\Ext(H,\bG_{m,S})\simeq H^2(H,\bG_{m,S})\simeq H^2_{\univ}(H,\bG_{m,S})\ \  
\mathrm{(}bijective\mathrm{)}.
\end{equation*}
\end{thm}

\begin{thm}
\label{thm:central extension E splits when}
Let $K=\bZ/n\bZ$, and let $H$ be a finite commutative group $S$-scheme 
of rank $n$ acting on $\bG_{m,S}$ with $D=\Gamma(\cO_H)$. 
Let $\id_H\in H(H)$ be the identity of $H$ 
and assume $n=\ord(\id_E)$. We define an action of $K$ on 
$D=\Hom(H,\bG_{m,S})$ by 
$K\times D\ni (x,f)\mapsto (x\id_H)^*f:=f\circ (x\id_H)$, and 
$N_K(f)=\prod_{x\in K}(x\id_H)^*f$. 
Moreover let $\phi\in Z^2(H,\bG_{m,S})$, $E=E_{\phi}\in\Ext(H,\bG_{m,S})$, and 
$\Omega(\phi):=\prod_{x=0}^{n-1}\phi(x\id_E,\id_E)\in D^{\times}$.
Then 
\begin{enumerate}
\item\label{item:isom GmD/GmDn} 
$\Omega$ induces an isomorphism 
$H^2_{\univ}(H,\bG_{m,S})\simeq (D^{\times})^H/N_HD^{\times}$;
\item\label{item:E splits} $E$ splits iff 
$\Omega(\phi)\in N_HD^{\times}$; 
\item\label{item:EW splits} if $a:=\Omega(\phi)\in D^{\times}$, 
$R^*=R[x]/(x^n-a)$ and $W=\Spec R^*$,  
then the pullback $E_W$ of $E$ splits.
\end{enumerate}
\end{thm}
\begin{proof}
By Corollary~\ref{cor:H^2(GmS)} and  
Definition~\ref{defn:univ H2(H Gm)},
$$H^2_{\univ}(H,\bG_{m,S})\simeq H^2(\bZ/n\bZ,D^{\times})=H^2(K,D^{\times}),$$
where $K$ naturally acts on $D^{\times}$ as follows:  
since $x\in K$ gives $x\id_H\in H(H)$, via the isomorphism 
$D^{\times}=\Hom(H,\bG_{m,S})$, 
this induces a natural action 
$\tau_K:K\times D^{\times}\mapsto D^{\times}$ by 
$\tau_K(x,f)=(x\id_H)^*f$. By this  
we define $N_K(g):=\prod_{x\in K}(x\id_H)^*g$ for $g\in D$.   
Hence 
$H^2(K,D^{\times})\simeq (D^{\times})^{K}/N_KD^{\times}$. 
This proves (\ref{item:isom GmD/GmDn}). 
(\ref{item:E splits})-(\ref{item:EW splits}) 
follows from Corollary~\ref{cor:H^2(GmS)} 
and Theorem~\ref{thm:Ext = H^2 group schemes}. 
\end{proof}

\subsection{Formulation in terms of Hopf algebras}
\begin{lemma}\label{lemma:def of gr scheme E}
Let $H$ be a finite flat commutative group $S$-scheme acting on $\bG_{m,S}$, 
$\tau:H\times\bG_{m,S}\to\bG_{m,S}$ the action of $H$ on $\bG_{m,S}$, 
 $E=H\times\Lambda$ as $S$-schemes, and 
$\phi\in Z^2(H,\bG_{m,S})\subset (D\otimes_R D)^{\times}$. 
Then we define a multiplication of $E$ by
\begin{equation}\label{eq:group law of E}
(x,a)\cdot (y,b)=(x+y, \phi(x,y)\cdot a\cdot x(b))
\end{equation}where $T$ is an $S$-scheme,  
$x,y\in H(T)$ and $a,b\in\bG_{m,S}(T)$ are 
functorial points, and $0(b)=b$. 
Then $E\in\Ext(H,\Lambda)$ with $\sigma(E)=[\phi]$, where 
$\sigma(E)$ is the class of 
the 2-cocycle associated with $E$. We denote $E$ by $E_{\phi}$. 
\end{lemma}
\begin{proof}
Associativity of the group law of $E$
follows from  Eq.~(\ref{eq:two cocycles}). Here $(0,1)$ is 
the identity of $E$ because $x(1)=1$ $(\forall x\in H(T))$, while 
$(x,\lambda)^{-1}=(-x,x(a)a^{-1})$. Indeed, 
let $b:=x(a)a^{-1}$. We compute 
$(-x,b)\cdot(x,a)=(0,x(a)\cdot(-x)(a)\cdot a^{-1})
=(0,\tau(x,a)\cdot\tau(-x,a)\cdot a^{-1})
=(0,\tau(0,a)\cdot a^{-1})=(0,a\cdot a^{-1})=(0,1)$. 
Hence $(-x,b)=(x,a)^{-1}$, and $E$ is a group $S$-scheme.  
This shows Definition~\ref{defn:extension group scheme of H by GmS}~(a).  
The group $\bG_{m,S}$ is identified with 
a normal subgroup $\{(0,a);a\in\bG_{m,S})\}$ of $E$,  
and then $E/\Lambda\simeq H$, which shows Definition~\ref{defn:extension group scheme of H by GmS}~(b). 
A lifting $\zeta$ of $H$ to $E$ is given by $\zeta(x)=(x,1)$, and 
$\zeta(x)(0,a)\zeta(x)^{-1}=(x,1)\cdot(0,a)\cdot(-x,1)
=(0,x(a))$, which 
shows Definition~\ref{defn:extension group scheme of H by GmS}~(c). Hence $E\in\Ext(H,\Lambda)$ with $\sigma(E)=[\phi]$. 
\end{proof}

Now we describe the above group $S$-scheme 
$E:=E_{\phi}$ via Hopf algebras. 
\begin{lemma}\label{lemma:Hopf al formulation}
Let $D=\Gamma(\cO_H)$, $A:=\Gamma(\cO_E)=D[w, w^{-1}]$, 
$\Delta_H$ 
(resp. $\iota_H$, $\epsilon_H$) 
comultiplication (resp. coinverse, counit) of $H$ and 
$\iota_{\bG_{m,S}}$ coinverse of $\bG_{m,S}$. Then 
comultiplication $\Delta$ (resp. coinverse $\iota$, counit $\epsilon$) 
of $E:=E_{\phi}$ is given via Lemma~\ref{lemma:def of gr scheme E} as follows:
\begin{gather*}
\Delta(w)=\phi\cdot w\otimes\tau^*(w),\ 
\iota(w)=\left(\tau^*\cdot(1\otimes\iota_{\bG_{m,S}})\right)(w),\ 
\epsilon(w)=1,\\ 
\Delta(\lambda)=\Delta_H(\lambda),\  
\iota(\lambda)=\iota_H(\lambda),\ 
\epsilon(\lambda)=\epsilon_H(\lambda)\ (\lambda\in D). 
\end{gather*}
\end{lemma}
\begin{proof} 
This is clear from Lemma~\ref{lemma:def of gr scheme E}. 
\end{proof}

Let $H^m:=H\times_S\cdots\times_SH$ ($m$-times), 
and let $\pr_{ij}:H^3\to H^2$ be 
the projection to the  $(i,j)$-component. 
Let $w_i:=w$, $\tau_i:=\tau$ and 
$\phi_{ij}:=\pr_{ij}^*\phi$.
By Lemma~\ref{lemma:Hopf al formulation}, we compute 
\begin{align*}
(\Delta\otimes 1)\Delta(w)&=(\Delta\otimes 1)
(\phi_{23}w_2\otimes\tau_2^*w_3)\\
&=(\Delta_H\otimes 1_D)(\phi_{23})(\Delta(w_2)\otimes(\Delta_H\otimes 1_{R[w^{\pm}]})\tau_2^*w_3)\\
&=\left((\Delta_H\otimes 1_D)\phi_{23})\phi_{12}\right)
\left(w_1\otimes\tau_1^*w_2\otimes\tau_1^*\tau_2^*w_3\right),\\
(1\otimes\Delta)\Delta(w)&=(1\otimes\Delta)(\phi_{12}w_1\otimes\tau^*_1w_2)\\
&=(1_D\otimes\Delta_H)(\phi_{12})(w_1\otimes(\tau_1^*\otimes\tau_1^*)(\Delta(w_2)))\\
&=(1_D\otimes\Delta_H)(\phi_{12})(w_1\otimes(\tau_1^*\otimes\tau_1^*)(\phi_{23}w_2\otimes\tau_2^*w_3))\\
&=\left((1_D\otimes\Delta_H)(\phi_{12})\tau_1^*\phi_{23}\right)
\left(w_1\otimes\tau_1^*w_2\otimes\tau_1^*\tau_{2}^*w_3\right). 
\end{align*} Note 
$(\Delta_H\otimes 1_{R[w^{\pm}]})\tau_2^*=\tau_1^*\tau_2^*$.
It follows 
\begin{equation}\label{eq:phi 2 cocycle Hopf alg}
((\Delta_H\otimes 1_D)\phi_{23})\phi_{12}
=((1_D\otimes\Delta_H)\phi_{12})\tau_1^*\phi_{23},
\end{equation}
which coincides with  Eq.~(\ref{eq:two cocycles x,y,z in H(T)}). 
Similarly Eq.~(\ref{eq:coboundary d1psi}) is translated into 
\begin{equation}\label{eq:phi 2 coboundary Hopf alg}
\phi=d^1\psi=(\Delta_H\psi)\psi_1^{-1}\tau_1^*\psi_2^{-1} 
\end{equation}where $\psi_i=p_i^*\psi$ and 
$p_i:H^2\to H$ is the $i$-th projection.

\subsection{Central extensions of a finite group $S$-scheme}
\label{subsec:applications to Z/nZ or mu_n}
In what follows we consider the case where   
a finite flat commutative group $S$-scheme $H$ 
acts trivialy on $\bG_{m,S}$.  
Hence as is mentioned in 
Definition~\ref{defn:extension group scheme of H by GmS}, 
any $E\in\Ext(H,\bG_{m,S})$ is 
a central extension of $H$ by $\bG_{m,S}$, 
that is, $\bG_{m,S}$ is {\it central}. 
Each $E$ is of the form $E_{\phi}$ for some $\phi\in H^2(H,\bG_{m,S})$.

First we define: 
\begin{defn}\label{defn:radical normal base change}
Let $R$ be a CDVR, $k(\eta)$ its fraction field and $S=\Spec R$. 
Let $R^*$ be a CDVR, $K^*$ the fraction field of $R^*$ and $W:=\Spec R^*$. 
If $K^*$ is the splitting field in the algebraic closure $\Omega$ 
of $k(\eta)$ of a polynomial $\prod_{i=1}^q(x^{n_i}-a_i)$ 
for some $(n_i,a_i)\in\bN\times R^{\times}$ $(i\in[0,q-1])$, 
then we say that $K^*$ is a finite radical normal extension of $k(\eta)$ 
and $\pi:W\to S$ is a {\it finite radical normal cover of} $S$, 
or $W$ is {\it finite radical normal over} $S$. See \S~\ref{subsc:Gal(W/S)}.
\end{defn}

\begin{lemma}
\label{lemma:splitting E_in_Cent_Z/nZ_Gm}
Let $H=(\bZ/n\bZ)_S$ and $E\in\Ext(H,\bG_{m,S})$.
Then there exists a finite radical normal cover $W$ of $S$ such that  
the pullback $E_{W}$ splits. 
\end{lemma}
\begin{proof}Let $D=\Gamma(\cO_H)$. By Theorem~\ref{thm:central extension E splits when}, there exists $\phi\in Z^2(H,\bG_{m,S})$ such that $E=E_{\phi}$. Let  
$a:=\Omega(\phi)\in D^{\times}$.  
By definition,  
 $H$ is a disjoint union of $n$ copies of $S$, say, 
$H=\coprod_{i\in[0,n-1]}S_i$ as $S$-schemes 
where $S_i$ is the $i$-th copy of $S$ and 
$S_0$ the identity component of $H$.  
Hence $D=\Gamma(\cO_H)=\bigoplus_{i=0}^{n-1}R\delta_i$  
where $\delta_i$ is the function on $H$ such that 
$\delta_i$ on $S_j$ equals $\delta_{i,j}$. 
It follows that 
\begin{equation}\label{eq:str of D_times}
\bG_m(D)=D^{\times}=\bigoplus_{i=0}^{n-1}R^{\times}\delta_i,\ 
\delta_i\delta_j=\delta_{i,j}\delta_i\ (\forall i,j\in[0,n-1]).
\end{equation}

Hence $a\in D^{\times}$ is written uniquely as 
$a=\sum_{i=0}^{n-1}a_i\delta_i$ with $a_i\in R^{\times}$. 
Let $K$ be the splitting field in $\Omega$ of the polynomial:
\begin{equation}\label{eq:splitting equations}
\prod_{i=0}^{n-1}(x^n-a_i).
\end{equation}

Let $R_K$ be the integral closure of $R$ in $K$ and $W=\Spec R_K$. 
Since $a_i\in R^{\times}$, there exists $(\alpha_i;i\in[1,n])$ such that 
$\alpha_i\in R_K^{\times}$ and $\alpha_i^n=a_i$. 
Let $\alpha:=\sum_{i=0}^{n-1}\alpha_i\delta_i\in 
\Gamma(\cO_{H_W})^{\times}=(D\otimes_RR_K)^{\times}$. Then $\alpha^n=a$, 
so that the pullback $E_W$ of $E$ splits by 
Theorem~\ref{thm:central extension E splits when}~(\ref{item:E splits}).
\end{proof}

\begin{rem}\label{rem:two theories of central ext are distinct}
Let $H=\bZ/n\bZ$. Then 
$H^2(H,\bG_m(S))=\bG_m(S)/\bG_m(S)^n$ 
can be different from $H^2_{\univ}(H_S,\bG_{m,S})=\bG_m(D)/\bG_m(D)^n$ 
by Eq.~(\ref{eq:str of D_times}).  
\end{rem}

Next we consider the case $H=\mu_{n,S}$. We prove:
\begin{lemma}\label{lemma:splitting E_in_Cent_mu_n_Gm}   
Let $H=\mu_{n,S}$ and $E\in\Ext(H,\bG_{m,S})$.
Then there exists a finite radical normal cover $W$ of $S$ such that  
the pullback $E_{W}$ splits. 
\end{lemma}
\begin{proof}  
Our proof of Lemma~\ref{lemma:splitting E_in_Cent_mu_n_Gm} 
consists of 4 steps:

\n{\it Step 1.}\quad By applying  
Theorem~\ref{thm:central extension E splits when} 
to the case $H=\mu_{n,S}$ and $K=\bZ/n\bZ$, we have 
$H^2(H,\bG_{m,S})=D^{\times}/(D^{\times})^n$. 
If $\chara k(0)=0$ or 
$n$ is prime to $\chara k(0)$, then 
Lemma~\ref{lemma:splitting E_in_Cent_mu_n_Gm} 
 follows from 
Lemma~\ref{lemma:splitting E_in_Cent_Z/nZ_Gm} because 
$\mu_{n,S}\simeq(\bZ/n\bZ)_S$ 
by taking a finite radical normal cover of $S$ for $S$
if necessary. 

\n{\it Step 2.}\quad Next consider the case where $p:=\chara k(0)>0$. 
Let $q$ be the maximal power of $p$ which divides $n$, and $r:=n/q$. 
There are two cases:
\begin{center} (i)\ $\chara R=0$,\quad (ii)\ $\chara R=p$. 
\end{center}

\n{\it Step 3.}\quad Now we consider Case (i).  
Let $D_{m,R}:=\Gamma(\cO_{\mu_{m,S}})$ for $m\in\bN$.
Then $D_{n,R}=D_{r,R}\otimes_RD_{q,R}$. We first compute $D_{r,R}$. 
Since $r$ is prime to $\chara k(0)$ and hence $\chara R$, 
we may assume that $R$ contains a primitive $n$-th root of unity 
by taking a finite radical normal cover $S'$ of $S$ for $S$ 
if necessary.  By Eq.~(\ref{eq:str of D_times}),  
$\mu_{r,S}\simeq(\bZ/r\bZ)_S=\coprod_{i\in\bZ/r\bZ}S_i$ as $S$-schemes 
where $S_i$ is the $i$-th copy of $S$, and hence 
\begin{equation}\label{eq:D_rS_times}
D^{\times}_{r,R}=\bigoplus_{i=0}^{r-1}
R^{\times}\delta_i,\  \delta_i\delta_j=\delta_{i,j}\delta_i\ 
(i,j\in[0,r-1]).
\end{equation}

Next we shall describe $D_{q,S}$ explicitly. 
Let $\zeta_q$ be a primitive $q$-th root of unity of $1$,  
$R':=R[\zeta_q]$ and $S':=\Spec R'$. Since the finite cover $S'\to S$ 
is radical normal, we assume $S'=S$ by taking $S'$ for $S$.
The $S$-scheme $\mu_{q,S}=\Spec R[x]/(x^q-1)$ has $q$-sections 
$s_k:S\to \mu_{q,S}$ defined by $s_k^*x=\zeta_q^k$ $(k\in[0,q-1])$. Let 
$\pi:\tmu_{q,S}\to\mu_{q,S}$ be the normalization of the $S$-scheme 
$\mu_{q,S}$, $S_k$ the $k$-th irreducible component of $\tmu_{q,S}$ 
which is a copy of $S$ with $s_k(S)=\pi(S_k)$, and 
$P:=\Spec k(0)[x]/(x-1)$ the unique closed point of $\mu_{q,S}$. 
Then $\tmu_{q,S}\simeq\coprod_{k=1}^qS_k$ and 
we have an exact sequence of $\cO_{\mu_{q,S}}$-modules:
\begin{equation}\label{eq:resolution of D}
0\to D_{q,S}\overset{\pi^*}{\to} \Gamma(\cO_{\tmu_{q,S}})
=\bigoplus_{k=0}^{q-1}R\sigma_k
\overset{\oplus_{k=1}^{q-1}\lambda_k}{\longrightarrow}\bigoplus_{k\in[2,q]}k(0)
\to 0,
\end{equation}
where $\sigma_k\in\Gamma(\cO_{\tmu_{q,S}})$ with 
${\sigma_k}_{|S_j}=\delta_{k,j}$, $\sigma_k\sigma_l=\delta_{k,l}\sigma_k$\ 
$(k,l\in[0,q-1])$ and 
$\lambda_l(\bigoplus_{k=0}^{q-1}a_k\delta_k)=a_l(P)-a_{l-1}(P)$\ 
$(l\in[1,q-1])$. 
It follows that 
\begin{equation}\label{eq:unit of D_nS}
D_{n,S}^{\times}=(D_{r,S}\otimes_RD_{q,S})^{\times}
=\bigoplus_{i=0}^{r-1}D_{q,S}^{\times}\delta_i.  
\end{equation}

Let $E:=E_{\phi}$ for some $\phi\in Z^2(H,\bG_{m,S})$ and 
$\Omega(\phi))=a\in D_{n,S}^{\times}$. 
Then $a$ is uniquely written as 
$a=\sum_{i=0}^{q-1}a_i\delta_i$ with $a_i\in D_{q,S}^{\times}$. 
By Eq.~(\ref{eq:resolution of D}), each $a_i$ is uniquely written as 
$a_i=\sum_{k=0}^{q-1}a_{i,k}\sigma_k$ for some $a_{i,k}\in R^{\times}$ 
such that $a_{i,1}(P)=a_{i,l}(P)$ $(\forall l\in[1,q-1])$.

By (the proof of) Lemma~\ref{lemma:splitting E_in_Cent_Z/nZ_Gm} 
there exists a finite radical normal cover $W$ of $S$ 
with $R^*:=\Gamma(\cO_W)$ 
and $\alpha_{i,k}\in (R^*)^{\times}$ $(k\in[0,q-1])$ 
such that $\alpha_{i,k}^{n}=a_{i,k}$. 
Since $\alpha_{i,1}(P^*)^n=\alpha_{i,k}(P^*)^n$, whence 
$\alpha_{i,1}(P^*)^r=\alpha_{i,k}(P^*)^r$ where $P^*$ 
is the unique closed point of $W$ lying over $P$.  
Since we have assumed $\mu_r\subset R$, we may assume 
$\alpha_{i,1}(P^*)=\alpha_{i,k}(P^*)$ $(\forall k\in[1,q-1])$ 
by multiplying $\alpha_{i,k}$ by some $r$-th root of unity. 
Hence $\alpha_i:=\sum_{i=0}^{q-1}\alpha_{i,k}\sigma_k\in 
(D_{q,S}\otimes_RR^*)^{\times}=D_{q,W}^{\times}$ $(\forall i\in[0,r-1])$. 
Let $\alpha=\sum_{i=0}^{r-1}\alpha_i\delta_i\in D_{n,W}^{\times}$. 
Then $a=\sum_{i=0}^{r-1}\alpha_i^n\delta_i=\alpha^n$. 
Hence $E_{W}$ splits by 
Theorem~\ref{thm:central extension E splits when}~(\ref{item:E splits}). 
This proves Lemma in Case (i).

\n{\it Step 4.}\quad In Case (ii), 
by taking a finite radical normal cover of $S$ if necessary, we have 
$S=\Spec R$, $R=k(0)[[t]]$ and
\begin{gather*}D_{q,R}=k(0)[[t]][x]/(x^q-1)
=k(0)[[t]][\epsilon]/(\epsilon^q)\ (\epsilon:=x-1),\\ 
D_{q,R}^{\times}=(k(0)[[t]])^{\times}
\bigoplus\bigoplus_{k=0}^{q-1}k[[t]]\epsilon^k,\\
D_{n,R}^{\times}=(D_{r,R}\otimes_RD_{q,S})^{\times}
=\bigoplus_{i=0}^{r-1}D_{q,R}^{\times}\delta_i.
\end{gather*}

Let $Z'=(\mu_{n,S})_{\red}\simeq\mu_{r,S}$ and $Z:=(Z')^0$: 
the identity component of $Z'$. Then $Z\simeq S$, 
$\Gamma(\cO_Z)=R$ and $\Gamma(\cO_{Z'})=D_{r,R}$.
Let $E=E_{\phi}$ and $a:=\Omega(\phi)$. 
Then we write $a=\sum_{i=1}^ra_i\delta_i$ and 
$a_i=\sum_{k=0}^{q-1}a_{i,k}\epsilon^k$. Let $K^*$ be the splitting field 
of $\prod_{i\in[1,r]}(x^n-a_{i,0})$, $R^*$ 
the integral closure of $R$ in $K^*$ and $W:=\Spec R^*$. 
Hence $W$ is naturally a reduced irreducible $S=Z$-scheme 
which is finite radical normal and flat over $S$. Since 
$H_W\simeq\mu_{r,W}$, we obtain 
$H^2(H_W,\bG_{m,W})=D_{r,R^*}/D_{r,R^*}^n$. 
Let $\alpha_{i,0}\in R^*$ with $\alpha_{i,0}^n=a_i$ $(i\in[1,r])$ and 
$\alpha=\sum_{i=1}^r\alpha_{i,0}\delta_i\in D_{r,R^*}$. 
  Then 
$E_W=E_{\phi\otimes R^*}$ and $\Omega(\sigma(E_W))
=a\otimes 1_{R^*}=\sum_{i=0}^{r-1}a_{i,0}
\delta_i\in D_{r,R^*}^{\times}$. It follows  
$a\otimes 1_{R^*}=\sum_{i=1}^r\alpha_{i,0}^n\delta_i=\alpha^n\in D_{r,R^*}^n$. 
Hence $E_W$ splits. This completes the proof in Case (ii), and hence of Lemma.
\end{proof}

\subsection{Commutative extensions}
Let $H$ be a finite flat commutative group $S$-scheme acting trivially on $\bG_{m,S}$. Let $\phi\in Z^2(H,\bG_{m,S})$ and $E:=E_{\phi}$. 
\begin{defn}We define a commutator form $e_E$ of $E$ by 
\begin{equation}
e_E(x,y)=\phi(x,y)/\phi(y,x)\quad (x,y\in H(T))
\end{equation}for any $S$-scheme $T$.  We call $\phi$ {\it symmetric} if $e_E$ is trivial. Note that $e_E$ is trivial, that is, $e_E$ is identically equal to 1, iff $E$ is commutative. Let $Z^2_{\symm}(H,\bG_{m,S})$ be the subset of $Z^2(H,\bG_{m,S})$ of all symmetric 2-cocycles.  Since $B^2(H,\bG_{m,S})$ is a submodule of $Z^2_{\symm}(H,\bG_{m,S})$, we set 
$$H^2_{\symm}(H,\bG_{m,S}):=Z^2_{\symm}(H,\bG_{m,S})/B^2(H,\bG_{m,S}).$$ 
Hence we have $H^2_{\symm}(H,\bG_{m,S})\subset H^2(H,\bG_{m,S})$.
\end{defn} 

The following is obvious:
\begin{lemma}\label{lemma:Ext_comm=H2_symm}
Let $H$ be a finite flat commutative group $S$-scheme, and 
$\Ext_{\comm}(H,\bG_{m,S})$ the subset of $\Ext((H,\bG_{m,S})$ of all commutative extensions of $H$ by $\bG_{m,S}$. Then 
$$\Ext_{\comm}(H,\bG_{m,S})\simeq H^2_{\symm}(H,\bG_{m,S})\quad 
\mathrm{(}bijective\mathrm{)}.$$
\end{lemma} 

\begin{lemma}\label{lemma:direct sum of H^2_symm}
Let $H_1,H_2$ be finite flat commutative   
group $S$-schemes.   Then 
$$H^2_{\symm}(\bigoplus_{i\in[1,2]}H_i,\bG_{m,S})\simeq\bigoplus_{i\in[1,2]}H^2_{\symm}(H_i,\bG_{m,S})\ \ \mathrm{(}bijective\mathrm{)}.$$ 
\end{lemma}
\begin{proof} We define a bijection 
$\tau:H^2_{\symm}(H,\bG_{m,S})
\to\bigoplus_{i\in[1,2]}H^2_{\symm}(H_i,\bG_{m,S})$. 
Let $H=H_1\oplus H_2$. Let $T$ be any $S$-scheme, 
$x_1,y_1\in H_1(T)$, $x_2,y_2\in H_2(T)$ and 
$x=(x_1,x_2),y=(y_1,y_2)\in H(T)$.  

Let $\phi\in Z^2(H,\bG_{m,S})$. Then we set $E:=E_{\phi}$, 
$\zeta$ a lifting of $H$ to $E$, and 
\begin{gather*}
\phi_1(x_1,y_1):=\phi(x_1,0;y_1,0),\  
\phi_2(x_2,y_2):=\phi(0,x_2;0,y_2),\\
\phi_{\new}(x_1,x_2;y_1,y_2):=\phi_1(x_1,y_1)\phi_2(x_2,y_2),\\
\zeta_1(x_1):=\zeta(x_1,0),\ \zeta_2(x_2):=\zeta(0,x_2).
\end{gather*}

We define a new lifting $\zeta_{\new}$ of $H$ to $E$ by 
\begin{gather*}\zeta_{\new}(x_1,0):=\zeta(x_1,0),\ 
\zeta_{\new}(0,x_2):=\zeta(0,x_2),\\
\zeta_{\new}(x_1,x_2):=\zeta(x_1,0)\cdot\zeta(0,x_2).
\end{gather*} 
Then we obtain 
\begin{align*}
\zeta_1(x_1)\zeta_1(y_1)&=\phi(x_1,0;y_1,0)\zeta(x_1+y_1,0)
=\phi_1(x_1,y_1)\zeta_1(x_1+y_1),\\
\zeta_2(x_2)\zeta_2(y_2)&=\phi(0,x_2;0,y_2)\zeta(0,x_2+y_2)
=\phi_2(x_2,y_2)\zeta_2(x_2+y_2).
\end{align*}

Moreover since $E$ is commutative, we have 
\begin{align*}
\zeta_{\new}(x_1,x_2)&\zeta_{\new}(y_1,y_2)
=\zeta_1(x_1)\zeta_2(y_1)\zeta_1(x_2)\zeta_2(y_2)\\
&=\phi_1(x_1,y_1)\phi_2(x_2,y_2)\zeta_1(x_1+x_2)\zeta_2(y_1+y_2)\\
&=\phi_1(x_1,y_1)\phi_2(x_2,y_2)\zeta_{\new}(x_1+x_2,y_1+y_2)
\end{align*}
The cohomology class $[\phi]$ does not depend 
on the choice of a lifting $\zeta$. Hence  
$[\phi]=[\phi_{\new}]\in H^2(H,\bG_{m,S})$. Now we define  
 $\tau([\phi])=([\phi_1],[\phi_2])$.

Conversely, 
let $\phi_1\in Z^2_{\symm}(H_1,\bG_{m,S})$,  
$\phi_2\in Z^2_{\symm}(H_2,\bG_{m,S})$ and $\phi=\phi_1\phi_2$. 
Then $\phi\in Z^2_{\symm}(H,\bG_{m,S})$ 
by Eq.~(\ref{eq:two cocycles x,y,z in H(T)}) and  the map $([\phi_1],[\phi_2])\mapsto [\phi]$ is the inverse of $\tau$. This completes the proof.
\end{proof}
\begin{thm}\label{thm:comm extension splits when}
Let $H',H''$ be finite commutative groups, 
$H:=H'_S\oplus D(H''_S)$ and $E\in\Ext_{\comm}(H,\bG_{m,S})$.
Then there exists a finite radical normal cover $W$ of $S$ 
such that the pullback $E_W$ splits. 
\end{thm}
\begin{proof}This follows from Lemmas~\ref{lemma:splitting E_in_Cent_Z/nZ_Gm}/\ref{lemma:splitting E_in_Cent_mu_n_Gm}/\ref{lemma:direct sum of H^2_symm}. 
\end{proof}

\begin{example}\label{example:H=A, E extension of A by GmS}
Let $A$ be an abelian $S$-scheme and $E\in\Ext(A,\bG_{m,S})$. 
Since $A$ is proper over $S$, $A$ acts trivially on $\bG_{m,S}$.
 Then there exists an open covering 
$(U_{\lambda};\lambda\in\Lambda)$ of $A$ such that 
$E_{U_{\lambda}}\simeq U_{\lambda}\times_S\bG_{m,S}$. 
Let $U:=\coprod_{\lambda\in\Lambda}U_{\lambda}$. 
Let $x,y\in U(T)$ for an $S$-scheme $T$. 
Then we have liftings $\zeta_U(x)$ and $\zeta(y)$ of $x$ and $y$ to $E_U$, so that we obtain $\phi_U\in Z^2(A_U,\bG_{m,U})$ such that $E_U=E_{\phi_U}$, and 
$\zeta_U(x)\zeta_U(y)=\phi_U(x,y)\zeta_U(x+y)$. We define 
$e_E(x,y)=\phi_U(x,y)/\phi_U(y,x)$. This is well-defined 
as a morphism of fppf sheaves $e_E:A\times A\to \bG_{m,S}$, which is 
therefore a morphism of $S$-schemes. 
See \S~\ref{subsec:fppf presheaf hW is sheaf} for fppf sheaves. 
 Since $A$ is proper over $S$, $e_E$ is a constant morphism, and hence identically equal to $1$.  This implies that $\phi$ is symmetric, and $E$ 
is a commutative extension of $A$ by $\bG_{m,S}$. 
Therefore the pushout $\bF_x^{\times}$ $(x\in\Hom(T,\bG_{m,S}))$ 
in \S~\ref{subsec:Raynaud extensions split case} 
is a commutative extension. 
\end{example}

\subsection{The $S$-scheme $\Aut(W/S)$}
\label{subsc:Gal(W/S)}
We return to Definition~4.15. 
Let $R$ be a CDVR, $k(\eta)$ its fraction field and $S=\Spec R$. 
Let $K^*$ be a finite radical normal extension of $k(\eta)$,
$R^*$ the integral closure of $R$ in $K^*$ and $W:=\Spec R^*$. 
We set $B:=\Aut(K^*/k(\eta))=\Aut(W/S)$, 
which is a group $S$-scheme of finite type with $B_{\red}$ finite over $S$. 
Let $B^0$ be the identity component of $B$. 
Then $B^0$ is a closed normal subgroup $S$-scheme of $B$.
If $\chara k(\eta)=0$, then $B=\Gal(W/S)$ and $B^0$ is trivial.  
If $\chara k(\eta)=p>0$, the quotient $B/B^0$ is a reduced 
finite group $S$-scheme with $B/B^0\simeq B_{\red}$. 
Since $B_{\red}$ is a closed subscheme of $B$, $B$ is 
a semi-direct product $B^0\rtimes B_{\red}$. 

Let $n\in\bN$, $\alpha\in R^{\times}$ and $F(X):=X^n-\alpha$. 
We assume that $F$ is irreducible in $R[X]$. Let  
$R^*:=R[X]/(F)$, 
$K^*$ the fraction field of $R^*$, $W:=\Spec R^*$ 
and $B=\Aut(W/S)$. 
Let $A=R[Z]/(Z^n-1)$ and  $\mu_{n,S}:=\Spec A$. 
We denote 
the class of $X$ (resp. $Z$)  by $x$ (resp. $z$).  
By definition, $\mu_{n,S}(\Omega)$ is 
the set of all $n$-th roots of unity in $\Omega=\overline{k(\eta)}$, while 
$\mu_{n,S}(S)=\{\zeta\in R; \zeta^n=1\}$.  
$\mu_{n,S}$ is a group $S$-scheme 
which acts on $W$ by 
$\sigma:\mu_{n,S}\times_SW\to W$ defined by      
$\sigma^*(w)=zw\in A\otimes_RR^*$. 

\begin{lemma}\label{lemma:semi-direct Gal}
If $X^n-\alpha$ is irreducible in $R[X]$,  
then $\Aut(W/S)\simeq   
\mu_{n,S}\rtimes\Aut(\mu_{n,S}/S),$ where 
$\Aut(\mu_{n,S}/S)=\{\rho\in(\bZ/n\bZ)_S^{\times}; 
\rho_{|\mu_{n,S}(S)}=\id_{\mu_{n,S}(S)}\}$.
\end{lemma}
\begin{proof}Let $B=\Aut(W/S)$. 
Since $W$ is finite affine over $S$, it suffices to determine 
$B(S)$ in order to determine $B$. 
Let $\zeta\in\mu_{n,S}(S)=\Hom_R(A,R)$, namely, 
$\zeta^*z=\zeta\in R$.  
We define $\phi_{\zeta}\in\End_R(R^*)$ by 
$\phi_{\zeta}^*(x)=\zeta x$ and $\phi_{\zeta}^*(z)=z$. 
Then $\phi_{\zeta}\in B(S)$.  
Similarly we define $\rho_i\in\End(\mu_{n,S})$ by $\rho_i(z)=z^i$ 
for $i\in (\bZ/n\bZ)$.
Let $\phi\in B(S)$. Then $\phi^*(x)^q=\alpha$, so that 
$\zeta:=(\phi^*x)/x\in A^{\times}$. 
Let $\rho:=\phi_{\zeta}^{-1}\circ\phi$. Then 
$\rho^*(x)=x$. Hence $\rho^*\in\End_R(A)$, so that 
\begin{equation}\label{eq:rho* homom} 
\rho^*(ab)=\rho^*(a)\rho^*(b)\ (a,b\in A).
\end{equation} 
This implies that $\rho$ is an automorphism 
of $\mu_{n,S}$ as a group $S$-scheme. 
Let $\Delta$ be comultiplication of $\mu_{n,S}$. 
Then Eq.~(\ref{eq:rho* homom}) is just the relation 
$$\Delta(\rho^*(z))=(\rho^*\otimes\rho^*)(\Delta(z))=\rho^*(z)\otimes\rho^*(z)
$$where $\Delta(z^i)=z^i\otimes z^i$. 
Let $\rho^*(z):=\sum_{i=0}^{n-1}a_iz^i$ $(a_i\in R)$. 
Then we obtain
$$\sum_{i=0}^{n-1}a_iz^i\otimes z^i=\sum_{i,j=0}^{n-1}a_ia_jz^i\otimes z^j,
$$from which we derive $a_i\delta_{i,j}=a_ia_j$. Hence we have 
$\rho^*(z)=z^i$ for a unique $i\in[0,n-1]$, so that $\rho=\rho_i$. 
Since $\rho$ is an automorphism, 
we have $i\in(\bZ/n\bZ)^{\times}$. It follows that 
$\phi=\phi_{\zeta}\circ\rho_i$ and 
\begin{equation*}
\Aut(\mu_{n,S}/S)=\{\rho_i; i\in(\bZ/n\bZ)^{\times},\ 
\rho_{i|\mu_{n,S}(S)}=\id_{\mu_{n,S}(S)}\}.
\end{equation*} 

Since $\mu_{n,S}$ is a normal subgroup $S$-scheme of $B$,
$B$ is the semi-direct product $\mu_{n,S}\rtimes\Aut(\mu_{n,S}/S)$.  
 Note that $B\simeq\mu_{n,S}$ by Kummer theory \cite[IV, \S~3]{Neu99} 
if $\chara k(\eta)=0$ and $\mu_{n,S}(\Omega)\subset k(\eta)$.  
\end{proof}

Next we study $B^0$  
in the case where $X^n-\alpha$ is reducible 
and $\chara k(\eta)=p>0$. 
Let $q$ be the maximal power of $p$ dividing $n$, and $m:=n/q$. 
Since $Y^m-\alpha$ is a separable polynomial, 
the splitting field $K'$ of 
$Y^m-\alpha$ is unramified over $k(\eta)$ 
and $B^0$ for the extension $K'/k(\eta)$ is trivial.  
Hence it suffices to consider $F(X)=X^q-\alpha$ when $n=q$. 

\begin{cor}\label{cor:B^0 when reducible}
Let $F(X)=X^q-\alpha$, $B=\Aut(W/S)$ and $K^*$ the splitting field of $F$. 
If $\chara k(\eta)=p>0$, then 
there exist a divisor $q'$ of $q$ and $\gamma\in R^{\times}$ 
such that $[K^*;k(\eta)]=q'$, $\alpha=\gamma^{q/q'}$,  
$B^0\simeq\mu_{q',S}$ and $B_{\red}\simeq(\bZ/q'\bZ)_S^{\times}$.
\end{cor}
\begin{proof}If $F$ is irreducible in $R[X]$, 
then by Lemma~\ref{lemma:semi-direct Gal}, 
$B^0=\mu_{q,S}$ and $B_{\red}=(\bZ/q\bZ)_S^{\times}$ 
because $\mu_{q,S}(S)=\{1\}$. 
 Next suppose that $F$ is reducible. 
Let $\beta$ be a root of $F(X)=0$ in $\Omega$. Then $K^*=k(\eta)(\beta)$.   
Since $F$ is reducible, there exist monic polynomials 
$F',F''\in R[X]$ such that 
$F=F'F''$. Hence there exist $a,b,q'\in\bN$ such that 
$(a,p)=(b,p)=1$, $a+b=q''\geq 2$, $q'q''=q$, $q'<q$, $\deg F'=aq'$ 
and $\deg F'=bq'$, so that 
$F'=(X-\beta)^{aq'}$ and $F''=(X-\beta)^{bq'}$. 
It follows $\beta^{aq'},\beta^{bq'}\in R$, so that $\beta^{q'}\in R^{\times}$. 
We choose the minimal $q'\in\bN$ with the above property. 
Let $H(X)=X^{q'}-\beta^{q'}\in R[X]$ and $\gamma:=\beta^{q'}\in R^{\times}$. 
Then $F(X)=(X-\beta)^q=(X^{q'}-\gamma)^{q''}=H^{q''}$, and  
$\alpha=\beta^q=\gamma^{q/q'}$. 
If $H$ is reducible, then we can find $q''\in\bN$ such that    
$q''|q'$, $q''<q'$ 
and $\beta^{q''}\in R^{\times}$, 
which contradicts the minimality of $q'$. 
Hence $H$ is irreducible in $R[X]$ and $q'=[K^*:k(\eta)]$. 
By Lemma~\ref{lemma:semi-direct Gal}, 
$B^0\simeq\mu_{q',S}$ and 
$B_{\red}=\Aut(\mu_{q',S}/S)\simeq(\bZ/q'\bZ)_S^{\times}$ 
because $\mu_{q',S}(S)=\{1\}$. 
\end{proof}

\section{Degeneration data of the N\'eron model}
\label{sec:Deg data of Neron model}

\subsection{Extension of $\zeta$}
\label{subsec:extension of zeta}
We return to \S~\ref{subsec:split obj in DDample}. 
Let $S=\Spec R$ with $t$ uniformizer of $R$. 
Let $(G,\cL)$ be a split semiabelian $S$-scheme,   
$\zeta_{l}:=\FC(G,\cL^{\otimes l})$  and $\zeta=\zeta_1$. 
The split object $\zeta_{l}$ is 
assumed to be symmetric.

\begin{defn}\label{defn:miscellany}
Let $v_t$ be a valuation of $R$ with $v_t(t)=1$ and  
$B_{\tau}:=v_t(\tau)$. We define   
$\beta:=\beta_{\tau}\in\Hom_{\bZ}(Y,X^{\vee})$ 
 by $\beta_{\tau}(y)(x)=B_{\tau}(y,x)$ $(y\in Y,x\in X)$, 
$N:=N_{\tau}:=|X^{\vee}/\beta(Y)|$   
 and $\mu:=\phi\circ\beta^{-1}\circ(N\id_{X^{\vee}})
\in\Hom_{\bZ}(X^{\vee},X)$.   Let $\Omega:=\overline{k(\eta)}$ be 
an algebraic closure of $k(\eta)$, $K$ 
an algebraic extension of $k(\eta)$ and 
$R_{K}$ the integral closure of $R$ in $K$. 
\end{defn}

The following is a pd-counterpart of \cite[5.5]{MN24}.
\begin{lemma}
\label{lemma:minimal Galois with extensions cte iotae psie taue}
By taking a finite radical normal extension $K$ of $k(\eta)$ possibly ramified 
over the maximal ideal $I$ of $R$, we can find  
\begin{enumerate}
\item[(i)]$c^{t,e}\in\Hom_{\bZ}(X^{\vee},A(R_K))$ such that
$c^t=c^{t,e}\circ\beta$ and $c\circ\mu=N\lambda\circ c^{t,e}$;
\item[(ii)] a trivialization of biextension 
$\tau^e:1_{X^{\vee}\times X}\simeq (c^{t,e}\times c)^*\cP_{K}^{\otimes (-1)}$, which is a bilinear 
form $\tau^e\in\Hom_{\bZ}(X^{\vee}\times X, K^{\times})$ 
such that $\tau=\tau^e\circ(\beta\times\id_X)$, $v_t\tau^e(u,x)=u(x)$\ 
and $\tau^e(u,\mu(v))=\tau^e(v,\mu(u))$ 
$(\forall u,v\in X^{\vee},\forall x\in X)$;
\item[(iii)] $\iota^e\in\Hom_{\bZ}(X^{\vee},\tG(K))$ 
such that $\iota=\iota^e\circ\beta$ and $c^{t,e}=\pi\circ\iota^e$, which is given by $\iota^e(u)=\bigoplus_{x\in X}\tau^e(u,x)$ $(u\in X^{\vee})$;
\item[(iv)] a cubical trivialization $\psi^e:1_{X^{\vee}}\simeq (c^{t,e})^*\cM_{K}^{\otimes (-2)}$  
compatible with  
$\tau^e\circ(\id_{X^{\vee}}\times (2\mu)):1_{X^{\vee}\times X^{\vee}}
\simeq (c^{t,e}\times(c\circ(2\mu)))^*\cP_K^{\otimes (-1)}$ in the sense that 
\begin{equation}\label{eq:defn of psie}\psi^e(u)=\tau^e(u,\mu(u)) 
\ \ (\forall u\in X^{\vee})
\end{equation} is a $K^{\times}$-valued function on $X^{\vee}$ such that 
$\psi^{2N}=\psi^e\circ\beta$ and 
\begin{equation}\label{eq:cubic triv psie}
\psi^e(u+v)=\psi^e(u)\psi^e(v)\tau^e(u,2\mu(v))\  
(\forall u,v\in X^{\vee}).
\end{equation} 
\end{enumerate}  

\end{lemma}
\begin{proof}Let $r=\rank_{\bZ}X$. 
We choose elementary divisors $e_j$ of $X^{\vee}/\beta(Y)$:
$$X^{\vee}/\beta(Y)\simeq \bigoplus_{i=1}^rH_i,\  
H_i=\bZ/e_i\bZ,\ e_i|e_{i+1}\ , e_i\geq 1\ \ (i\in[1,r-1]). 
$$Note that $N=\prod_{i=1}^re_i$. 
Let $x\in X^{\vee}$ and $\barx$ the class of $x$ in $X^{\vee}/\beta(Y)$. 
Let $(u_i\in X;i\in[1,r])$ (resp. $(y_i\in Y;i\in[1,r])$) 
be a $\bZ$-basis of $X^{\vee}$ 
(resp. $Y$) such that $\baru_i\in H_i$ 
is a generator of $H_i$ and $e_iu_i=\beta(y_i)$ $(\forall i\in[1,r])$. 

Since $e_i\id_A:A\to A$ is surjective, we have 
$\alpha_i\in A(R_{\Omega})$ such that 
$e_i\alpha_i=c^t(y_i)\in A(R_{\Omega})$. Since $e_iu_i=\beta(y_i)$, we have 
$\mu(u_i)=(N/e_i)\phi(y_i)$. Hence  
$$N\lambda(\alpha_i)=(N/e_i)\lambda(c^t(y_i))=(N/e_i)c(\phi(y_i))
=c((N/e_i)\phi(y_i))=c(\mu(u_i)). 
$$

We define $c^{t,e}(u_i)=\alpha_i$ and extend 
$c^{t,e}$ to $X^{\vee}$ additively. Thus we have 
$c^{t,e}\in\Hom_{\bZ}(X^{\vee},A(R_{\Omega}))$ 
with $N\lambda\circ c^{t,e}=c\circ\mu$. Moreover 
$$c^{t,e}(\beta(y_i))=c^{t,e}(e_iu_i)=e_i\alpha_i=c^t(y_i),
$$which shows $c^{t,e}\circ\beta=c^t$ on $Y$. This proves (i). 

Let $\bar\tau^e(u,x):=t^{-u(x)}\tau^e(u,x)$ and 
$\bar\tau(y,x):=t^{-B_{\tau}(y,x)}\tau(y,x)$. 
To construct $\tau^e$ satisfying (ii), 
we first construct 
$\bar\tau^e\in\Hom_{\bZ}(X^{\vee}\times X,\Omega^{\times})$ 
such that 
\begin{equation}\label{eq:bar_taue_betayx=bartau_yx}
\bar\tau^e(\beta(y), x)=\bar\tau(y, x)\ \ 
(\forall x\in X, \forall y\in Y).
\end{equation}
Since $e_iu_i=\beta(y_i)$ $(\forall i\in[1,r])$,  
Eq.~(\ref{eq:bar_taue_betayx=bartau_yx}) is equivalent to the following:
\begin{equation}\label{eq:bar_taue_beta eiui xk}
\bar\tau^e(u_i, x_k)^{e_i}=\bar\tau(y_i, x_k)\ \ 
(\forall k\in[1,r]) 
\end{equation}for each fixed $i\in[1,r]$. It is clear that Eq.~(\ref{eq:bar_taue_beta eiui xk}) has a solution $\bar\tau^e(u_i, x_k)\in R_{\Omega}^{\times}$. 
Then we extend $\bar\tau^e$ (and hence $\tau^e$) to $X^{\vee}\times X$ 
additively. If $r=1$, this proves (ii). 
Next we assume $r\geq 2$. For $i\neq j$, we have
\begin{align*}&\bar\tau^e(u_i,\mu(u_j))=\bar\tau^e(u_i,\phi(y_j))^{N/e_j}
=\bar\tau^e(e_iu_i,\phi(y_j))^{N/e_ie_j}\\
&=\bar\tau(y_i,\phi(y_j))^{N/e_ie_j}
=\bar\tau^e(y_j,\phi(y_i))^{N/e_ie_j}=\bar\tau^e(u_j,\mu(u_i))
\end{align*}   
by Theorem~\ref{thm:degeneration data2}~(\ref{item:psi/tau}). 
This proves  (ii). 

Next we define $\iota^e\in\Hom_{\bZ}(X^{\vee},\tG(\Omega))$ by 
$\iota^e(u)^*:=\bigoplus_{x\in X}\tau^e(u,x)$ 
 $(u\in X^{\vee})$.  
This implies $\iota^e(\beta(y))^*=\bigoplus_{x\in X}\tau(y,x)=\iota(y)^*$ $(y\in Y)$ by Eq.~(\ref{eq:iotay}). 
Since $\tau^e(u,x)$ is additive in $u$, 
so is $\iota^e$, which lifts  
$c^{t,e}\in\Hom_{\bZ}(X^{\vee},A(\Omega))$ 
to $\Hom_{\bZ}(X^{\vee},\tG(\Omega))$. 
This proves (iii).  

Finally we define $\psi^e(u):=\tau^e(u,\mu(u))$
$(u\in X^{\vee})$. 
By (ii),  
$\psi^e(\beta(y))=\tau(y,\mu(\beta(y)))=\psi(y)^{2N}$. 
Eq.~(\ref{eq:cubic triv psie}) is clear, \footnote{Caution.  
Any function $\psi':X^{\vee}\to K$ 
satisfies the condition (iii) iff there exists  
$\chi\in\Hom(X^{\vee}/\beta(Y),\{\pm 1\})$ such that 
$\psi'=\chi\psi^e$. }
which shows that $\psi^e$ is 
(viewed as) a cubical trivialization 
$\psi^e:1_{X^{\vee}}\simeq (\iota^e)^*\tcL_{\Omega}^{\otimes (-2)}
\simeq (c^{t,e})^*\cM_{\Omega}^{\otimes (-2)}$ compatible with $\tau^e\circ(\id_{X^{\vee}}\times(2\mu))$. 
This proves (iv). The rest is clear.
\end{proof}

\begin{defn}\label{defn:minimal K with extensions cte iotae psie taue} 
Let $k(\eta^e)$ be the minimal normal extension of $k(\eta)$ 
among those fields $K$ for which we can find (i)-(iv) of 
Lemma~\ref{lemma:minimal Galois with extensions cte iotae psie taue}. 
Let $R^e$ be the integral closure of $R$ in $k(\eta^e)$, $S^*:=\Spec R^e$,  
$\eta^e$ the generic point of $S^e$ and $0^e$ the closed point of $S^e$. 
By Eq.~(\ref{eq:bar_taue_beta eiui xk}), 
a primitive $e_r$-th root $\zeta_{e_r}$ of unity belongs to $R^e$ where 
$e_r$ is the maximal elementary divisor of $X/\beta(Y)$.  
We denote $S^e$ (resp. $R^e$, $k(\eta^e)$) by $B$ (resp. $R_K$, $K$) if no confusion is possible. We sometimes write $R=R(B)$ and $K=k(B)$ if necessary.  
\end{defn}

\begin{lemma}\label{lemma:closed subgr Xvee/betaY}
Let $(G,\cL)$ be a split symmetric 
semiabelian $B$-scheme and  
$\cG$ the N\'eron model of $G_{\eta}$. 
The $K$-valued point $\iota^e(u)\in \tG(K)$ $(u\in X^{\vee})$ induces  
the $R_K$-valued point $\bar\iota^e(u)\in\cG(R_K)$. 
Let $H:=\{\bar{\iota}^e(u);u\in X^{\vee}/\beta(Y)\}$. 
Then  $H$ is a subgroup  
of $\cG(R_K)$ such that 
$\Phi(\cG/G)=H\simeq X^{\vee}/\beta(Y)$.
\end{lemma}
\begin{proof}Let $\cG^{\flat}$ be
the N\'eron model of $G_{K}$, $\tP$ 
a relatively complete model of $G$ and 
$P$ the algebraization of $\tP^{\wedge}/\iota(Y)$.  
There do exist $\tP$ and $P$, for instance by \cite[3.24]{AN99},    
in which case $\Sigma=\{0\}$. 
The condition \cite[(*), p.~62, line~6]{FC90}
hereby turns out to be trivially true:  
$2v_t(I_{y,\alpha}I_y)
\geq 2v_t(\tau(y,\alpha)\psi(y))
=2B_{\tau}(y,\alpha)+B_{\tau}(y,\phi(y))
=B_{\tau}(y,\phi(y))\geq 0$ $(\forall y\in Y, \forall\alpha\in\Sigma=\{0\})$
where $I_{y,\alpha}=\tau(y,\alpha)R$, $I_y=\psi(y)R$ and  
$\psi(y)^2=\tau(y,\phi(y))$. 
Hence $\tP$ is indeed a relatively complete model 
by \cite[III, 3.3]{FC90}. 

By Lemma~\ref{lemma:minimal Galois with extensions cte iotae psie taue} 
and \cite[III, 3.3]{FC90},  
$\iota^e(X^{\vee})\subset\tG(K)\subset\tP(R_K)$.  
Each $\iota^e(u)$ $(u\in X^{\vee})$ 
gives rise to an $R_K$-valued point 
$\bar{\iota}^e(u)\in\cG^{\flat}(R_K)$ 
because $\tP(R_K)/\iota(Y)
=\tP^{\wedge}(R_K)/\iota(Y)=P^{\wedge}(R_K)=P(R_K)=P(K)
=G(K)=\cG^{\flat}(R_K)$ 
by  \cite[III, 8.1]{FC90}. 
Since 
$v_t(\tau^e(u,x))=u(x)\in\bZ$ $(\forall x\in X)$, we have 
$\bar{\iota}^e(u)\in\cG(R_K)$, so that 
$H$ is a subgroup of $\cG(R_K)$. 

We define a homomorphism 
$\delta_t:\tG(\Omega)\to X^{\vee}_{\bQ}$ 
by $\delta_t(Q)(x)=v_t(w^x(Q))$ $(Q\in
\tG(\Omega), x\in X)$. 
Let $\tG(\Omega)^0:=\{Q\in\tG(\Omega);\delta_t(Q)\in X^{\vee}\}$
and $\tG(R_{\Omega})^0:=\tG(R_{\Omega})\cap\tG(\Omega)^0$. 
 By \cite[III, line 11, p.78]{FC90}, 
$\Phi(\cG/G)=\cG(R^{\unr})/G(R^{\unr})
=\tG(k(\eta)^{\unr})/\tG(R^{\unr})\cdot\iota(Y)
\overset{\delta_t}{\simeq} X^{\vee}/\beta(Y)$. 
Similarly  
$\Phi(\cG/G)=\cG(R_{\Omega})/G(R_{\Omega})
=\tG(\Omega)^0/\tG(R_{\Omega})^0\cdot\iota(Y)
\overset{\delta_t}{\simeq} X^{\vee}/\beta(Y)$. 
Since $\delta_t(\iota^e(u))
=v_t(\tau^e(u,\bullet))=u$\ $(\forall u\in X^{\vee})$, 
$H\ (\subset\cG(R_{\Omega}))$ is mapped onto $X^{\vee}/\beta(Y)$.  
Hence $\Phi(\cG/G)\simeq H$.
\end{proof}

\begin{lemma}\label{lemma:cN Phi-inv}
Let $\cN$ be the unique cubical invertible sheaf 
on $\cG$ extending the sheaf $\cL^{\otimes 2N}_{\eta}$ on $G_{\eta}$. 
Then $T_a^*(\cN_B^{\otimes N})
\simeq\cN_B^{\otimes N}$\ $(\forall a\in H)$.
\end{lemma}
\begin{proof}Since $\cN$ is cubical, 
we have an isomorphism on $\cG^3$:
\begin{align}\label{eq:cubic isom of Ndagger}
\Theta(\cN)&=
\bigotimes_{\emptyset\neq I\subset \{1,2,3\}}m_I^*(\cN)^{(-1)^{|I|+1}}\simeq 
\cO_{\cG^3}.
\end{align}

We define a map $\phi:H\to\Pic(\cG_B)$ by 
$\phi(a):=T_a^*(\cN_B)\otimes_{\cO_{\cG_B}}\cN_B^{\otimes (-1)}$. 
Let $a,b\in H_B(B)$  
and $c=\id_{\cG_B}\in\cG_B(\cG_B)$. 
By the pullback of (\ref{eq:cubic isom of Ndagger}) 
by  the morphism $(a,b,c):\cG_B\to\cG^3_B$,   
we have $T_{a+b}^*(\cN_B)\otimes_{\cO_{\cG_B}}\cN_B\simeq 
T_{a}^*(\cN_B)\otimes_{\cO_{\cG_B}}T_{b}^*(\cN_B)$ on $\cG_B$. 
See \cite[\S~6, Cor.~4, p.~57]{Mumford12}. 
Hence $\phi\in\Hom_{\bZ}(H, \Pic(\cG_B))$, so that 
$N\phi(a)=\phi(Na)=0$, {\it i.e.},    
$T_a^*(\cN^{\otimes N}_B)\simeq \cN^{\otimes N}_B$.
\end{proof}

Let $\cN^{\dagger}:=\cN_{B}^{\otimes N}$.
Since $\cN^{\dagger}_{|G_B}$ is cubical on $G_B$, it is 
the unique cubical extension of $\cL^{\otimes 2N^2}_{K}$ to $G_B$ 
by \cite[2.13~(2)]{MN24}, 
so that $\cN^{\dagger}_{|G_B}=\cL^{\otimes 2N^2}_B$. Hence   
$(\cN^{\dagger}_{|G_B})^{\wedge}=(\cL^{\otimes 2N^2})_B^{\wedge}\simeq(\pi^*\cM^{\otimes 2N^2})_B^{\wedge}$ where $\pi:\tG\to A$ is 
the morphism in Eq.~(\ref{eq:exact seq of tG}). 
By Lemma~\ref{lemma:cN Phi-inv}, $(\cN^{\dagger}_{|T_a(G_B)})^{\wedge}
=(T_a^*\cN^{\dagger}_{|G_B})^{\wedge}
\simeq(\pi^*\cM^{\otimes 2N^2})_B^{\wedge}$ $(\forall a\in H)$. 
Since $(\cG_B)^{\wedge}=\coprod_{a\in H}(a\times G_B^{\wedge})$, 
we have  
\begin{gather}\label{eq:Ndaggerwedge}
\cN^{\dagger,\wedge}:=(\cN^{\dagger})^{\wedge} 
=H\times 
\left(\bigoplus_{x\in X}
(\cM^{\otimes 2N^2})_{x,B}\right)^{I\op{-adic}}.
\end{gather} 
Hence we have a natural homomorphisms 
\begin{equation*}
\Gamma(\cG_B,\cN^{\dagger})\hookrightarrow 
\Gamma(\cG_B^{\wedge},\cN^{\dagger,\wedge})
\hookrightarrow H\times\prod_{x\in X}
\Gamma(A_B,(\cM^{\otimes 2N^2})_{x,B}))
\end{equation*}because $\Gamma(A_B,(\cM^{\otimes 2N^2})_{x,B})^{\wedge}\simeq\Gamma(A_B,(\cM^{\otimes 2N^2})_{x,B})$.  Let 
\begin{gather*}
V_x:=\Gamma(A_B,(\cM^{\otimes 2N^2})_{x,B}).
\end{gather*}

\begin{defn}\label{defn:2nd sigma_x}
Let $\sigma_x:\Gamma(\cG,\cN^{\dagger})
\to V_x$ be the natural projection.
\end{defn}

\subsection{Theta groups of $\cN^{\dagger}$}
\label{subsec:closed subsch of KcNdagger_f}
\begin{defn}
Let $\cG^t$ be the N\'eron model of $G^t$, 
$\lambda(\cN_{\eta}):G_{\eta}\to G^t_{\eta}$ 
the polarization 
morphism of $\cN_{\eta}$ and $\lambda(\cN):\cG\to\cG^t$ 
the extension of $\lambda(\cN_{\eta})$ to $\cG$ 
by the universal property of the N\'eron model $\cG^t$ and 
$\lambda(\cN^{\dagger}):=\lambda(\cN^{\otimes N})_B$. 
Hence $\lambda(\cN^{\dagger})$ extends  
$\lambda(\cN^{\dagger}_K)=\lambda(\cL_{\eta}^{\otimes 2N^2})_K$ 
to $\cG_B$. Let  
\begin{gather*}K^{\sharp}(\cN^{\dagger}):=
\ker(\lambda(\cN^{\dagger})),\ \  
K(\cN^{\dagger}):=K^{\sharp}(\cN^{\dagger})\cap G_{B}.
\end{gather*}   
Let $K^{\sharp}(\cN^{\dagger})^f$ (resp. $K(\cN^{\dagger})^f$) 
be the finite part of $K^{\sharp}(\cN^{\dagger})$ 
(resp. $K(\cN^{\dagger})$) and 
$e^{\sharp}$ the bilinear form 
on $K^{\sharp}(\cN^{\dagger})$ extending the Weil pairing 
$e^{\cN^{\dagger}_K}$ of $K^{\sharp}(\cN^{\dagger}_K)$. See 
\cite[IV, 2.4~(ii)]{MB85}. Finally we set 
$K(\cN^{\dagger})^m:=K(\cN^{\dagger})^f\cap T$ by abuse of notation.
To be more precise, $K(\cN^{\dagger})^m$ is defined as follows: 
Since $K(\cN^{\dagger})^f$ is a $B$-finite subscheme of $G$, 
so is $K(\cN^{\dagger})^{f,\wedge}$ of $G^{\wedge}=\tG^{\wedge}$. 
Hence $K(\cN^{\dagger})^{f,\wedge}\cap T^{\wedge}$ is viewed as 
a closed $B^{\wedge}\ (=B)$-finite formal subscheme of $G^{\wedge}$. 
By \cite[III, 5.4.5]{EGA}, there exists a closed $B$-finite subscheme 
$K(\cN^{\dagger})^m$ of $G$ 
such that $K(\cN^{\dagger})^{m,\wedge}\simeq 
K(\cN^{\dagger})^{f,\wedge}\cap T^{\wedge}$. This is what we mean by  
$K(\cN^{\dagger})^m:=K(\cN^{\dagger})^f\cap T$. Let 
$K^{\sharp}(\cN^{\dagger})^m:=K(\cN^{\dagger})^m$.
\end{defn}
\begin{rem}\label{rem:Weil pairing Ksharp_cNdagger_f}
By \cite[IV, 3.4]{MB85}
$K^{\sharp}(\cN^{\dagger})^f$, $K(\cN^{\dagger})^f$ and 
$K(\cN^{\dagger})^m$ are flat finite over $B$, 
while  $K^{\sharp}(\cN^{\dagger})$ 
and $K(\cN^{\dagger})$ are flat quasi-finite over $B$. 
By \cite[IV, 2.4~(iv)-(v)]{MB85}, \footnote{The \'etale base change mentioned in \cite[iV, 2.4~(v)]{MB85} 
is unnecessary here because there (already) exists $\cM$ on $A$ such that 
$\cL^{\wedge}=\pi^{\wedge,*}(\cM^{\wedge})$ in 
Eq.~(\ref{eq:exact seq of Gwedge}). 
See \S~\ref{subsec:Raynaud extensions split case}.} 
$K(\cN^{\dagger})^m$ is 
a closed subgroup $B$-scheme of $K(\cN^{\dagger})^f$ such that 
\begin{enumerate}
\item[(a)]
$K(\cN^{\dagger})^f$ is $e^{\sharp}$-orthogonal to 
$K(\cN^{\dagger})^m$: $e^{\sharp}(K(\cN^{\dagger})^f,K(\cN^{\dagger})^m)=1$;
\item[(b)] $e^{\sharp}$ induces a perfect dual pairing $e^{\ab}$ on 
$K^{\ab}(\cN^{\dagger}):=K(\cN^{\dagger})^f/K(\cN^{\dagger})^m$;
\item[(c)] 
$(K^{\ab}(\cN^{\dagger}),e^{\ab})\simeq (K(\cM^{\otimes 2N^2}),
e^{\cM^{\otimes 2N^2}})$ as $R_K$-modules with bilinear form  
where $e^{\cM^{\otimes 2N^2}}$ is 
the Weil pairing of $\cM^{\otimes 2N^2}$.
\end{enumerate}  
\end{rem}

\begin{defn}
We define non-commutative group $B$-schemes by 
\begin{gather*}
\cG^{\sharp}(\cN^{\dagger}):
=(\cN^{\dagger})^{\times}_{K^{\sharp}(\cN^{\dagger})},\ \  
\cG(\cN^{\dagger}):=(\cN^{\dagger})^{\times}_{K(\cN^{\dagger})},\\
\cG^{\sharp}(\cN^{\dagger})^f:=(\cN^{\dagger})^{\times}_{K^{\sharp}
(\cN^{\dagger})^f},\ \  
\cG(\cN^{\dagger})^f:=(\cN^{\dagger})^{\times}_{K(\cN^{\dagger})^f},\\ 
\cG(\cN^{\dagger})^m:=(\cN^{\dagger})^{\times}_{K(\cN^{\dagger})^m}.
\end{gather*} 
   Each non-commutative group $B$-scheme 
 is a central extension of a finite group scheme (say $\cK$) 
by the center $\bG_{m,B}$ whose commutator 
form equals (the restriction of) $e^{\sharp}$ by \cite[IV, 2.4~(iii)]{MB85},  
which we call the {\it theta group of $\cK$} where 
$\cK$ is one of 
$K^{\sharp}(\cN^{\dagger})^{\alpha}$ or $K(\cN^{\dagger})^{\alpha}$ 
$(\alpha\in\{\emptyset,f,m\})$ 
with $K^{\sharp}(\cN^{\dagger})^m:=K(\cN^{\dagger})^m$.
Let  
$\varpi^{\sharp}:\cG^{\sharp}(\cN^{\dagger})\to K^{\sharp}(\cN^{\dagger})$ 
(resp. $\varpi:\cG(\cN^{\dagger})\to K(\cN^{\dagger})$) be 
the natural projection. \par
Let $\tG$ be the Raynaud extension of $G$, and 
$\pi:\tG\to A$ the morphism in Eq.~(\ref{eq:exact seq of tG}). 
Let $\cM'$ be an ample invertible sheaf on $A$, 
$\lambda(\cM'):=\lambda_A(\cM'):A\to A^t$ the polarization morphism 
$K(\cM'):=K(A,\cM'):=\ker(\lambda(\cM'):A\to A^t)$ and 
$\cG(\cM'):=\cG(A,\cM'):=(\cM')^{\times}_{K(\cM')}$. Hence 
$\cG(\cM')$ is a central extension of $K(\cM')$ 
by $\bG_{m,B}$, which we call the {\it theta group of} $K(\cM')$,  or 
the {\it theta group of} of $\cM'$ if no confusion is possible. 
See \cite[\S~23]{Mumford12}. 
\end{defn}

\begin{defn}A {\it level subgroup $B$-scheme} $H'$
 of $\cG^{\sharp}(\cN^{\dagger})^f$ is defined to be  
a closed subgroup $B$-scheme of 
$\cG^{\sharp}(\cN^{\dagger})^f$ {\it finite flat over $B$}  
such that  $H'\cap\bG_{m,B}=\{1\}$. By definition, $\varpi(H')$ is a 
closed subgroup $B$-scheme of $K^{\sharp}(\cN^{\dagger})^f$ such that 
$H'\simeq\varpi^{\sharp}(H')$, in which case we say that  
$H'$ lifts $\varpi^{\sharp}(H')$, or 
$\varpi^{\sharp}(H')$ is liftable to $H'$. 
See \cite[I, p.291]{Mumford66}. 
\end{defn}

\begin{lemma}\label{lemma:KflatNdagger_m}
By taking a finite radical normal cover $B$ of $S$
 (Definition~\ref{defn:radical normal base change}) if necessary, 
 we have the following:
\begin{enumerate}
\item\label{item:dual of KNdagger_m} 
$K(\cN^{\dagger})^m
=D(X/2N^2\phi(Y))_B$;
\item\label{item:KflatNdagger_m} $K(\cN^{\dagger})^m$ is liftable 
to a level subgroup $B$-scheme $K^{\flat}(\cN^{\dagger})^m$
of $\cG(\cN^{\dagger})^m$;
\item\label{item:cG_N_dagger_m comm gr sch}$\cG(\cN^{\dagger})^m\simeq\bG_{m,B}\oplus K^{\flat}(\cN^{\dagger})^m$ as commutative group $B$-schemes;
\item\label{item:center of cG_N_dagger_f}$K^{\flat}(\cN^{\dagger})^m$ is 
contained in the center of $\cG(N^{\dagger})^f$.
\end{enumerate}
\end{lemma}
\begin{proof}Let $K^m:=K(\cN^{\dagger})^m$ and $K^f:=K(\cN^{\dagger})^f$.
Since $K^m$ is $B$-finite flat, so is 
$K^{m,\wedge}$. Since $G^{\wedge}\simeq\tG^{\wedge}$ and 
$\tG$ is of finite type, 
there exists by \cite[III, 5.4.5]{EGA} 
a closed $B$-subscheme 
$K^{\natural}(\cN^{\dagger})^m$ 
of $\tG$ $B$-finite flat such that 
$K^{m,\wedge}\simeq (K^{\natural}(\cN^{\dagger})^m)^{\wedge}$. Hence 
by \cite[III, 5.4.1]{EGA}
$$D(K^m)\simeq D((K^m)^{\wedge})
\simeq D((K^{\natural}(\cN^{\dagger})^m)^{\wedge})
\simeq D(K^{\natural}(\cN^{\dagger})^m).$$ 
  It follows from $(\cN^{\dagger})^{\wedge}_{|\tG^{\wedge}}
=(\cL^{\otimes 2N^2})^{\wedge}_{|\tG^{\wedge}}$ 
that $K^m\simeq K^{\natural}(\cN^{\dagger})^m
=\ker(2N^2\lambda_T)=D(X/2N^2\phi(Y))_B$.  
 This proves (\ref{item:dual of KNdagger_m}). 
See \cite[p.~667]{Nakamura99}. 

Since $K^m\subset K(\cN^{\dagger})^f$,  
$K^m$ is totally $e^{\sharp}$-isotropic 
by Remark~\ref{rem:Weil pairing Ksharp_cNdagger_f}~(a).     It follows that $\cG(N^{\dagger})^m$ is 
a commutative extension of $K^m$ by $\bG_{m,S}$.
By Theorem~\ref{thm:comm extension splits when}, 
there exists a finite radical normal cover $\pi:W\to S$ of $S$ 
such that $\cG(N^{\dagger})^m_W$ splits, that is, 
$\cG(N^{\dagger})^m_W=K^m_W\oplus\bG_{m,W}$. 
Therefore by taking $W$ as $B$ if necessary, 
we have a level subgroup scheme 
$K^{\flat}(\cN^{\dagger})^m$ of $\cG(\cN^{\dagger})^m$ 
lifting $K^m$.  
This proves (\ref{item:KflatNdagger_m}).  By (\ref{item:KflatNdagger_m}), 
$\cG(\cN^{\dagger})^m\simeq\bG_{m,B}\oplus K^{\flat}(\cN^{\dagger})^m$ as commutative group schemes because $e^{\sharp}$ is trivial on $K^m\times K^m$. 
This proves (\ref{item:cG_N_dagger_m comm gr sch}).  
Since $\cG(\cN^{\dagger})^m$ is a subgroup of $\cG(\cN^{\dagger})^f$, so is $K^{\flat}(\cN^{\dagger})^m$. Since $K^m$ is $e^{\sharp}$-orthogonal to $K^f$ by Remark~\ref{rem:Weil pairing Ksharp_cNdagger_f}~(a), $K^{\flat}(\cN^{\dagger})^m$ 
is in the center of $\cG(N^{\dagger})^f$. This proves (\ref{item:center of cG_N_dagger_f}).  
\end{proof}

By
Remark~\ref{rem:Weil pairing Ksharp_cNdagger_f}~(c)/\cite[IV, 2.4 (iv)]{MB85}, 
\begin{equation}\label{eq:AV part}
\cG(\cN^{\dagger})^f/\cG(\cN^{\dagger})^m\simeq
K(\cN^{\dagger})^f/K(\cN^{\dagger})^m\simeq K(\cM^{\otimes 2N^2}).
\end{equation}
By Lemma~\ref{lemma:cN Phi-inv}, $H_B$ is a closed 
subgroup $B$-scheme of $K^{\sharp}(\cN^{\dagger})^f$ such that 
$K^{\sharp}(\cN^{\dagger})^f=H_B\oplus K(\cN^{\dagger})^f$, so that 
\begin{equation}\label{eq:isom to HB}
\cG^{\sharp}(\cN^{\dagger})^f/\cG(\cN^{\dagger})^f
\simeq K^{\sharp}(\cN^{\dagger})^f/K(\cN^{\dagger})^f\simeq H_B.
\end{equation}

\begin{lemma}\label{lemma:pairing esharp}
The pairing $e^{\sharp}$ 
is given on $H_B\times K(\cN^{\dagger})^m$ by 
$$e^{\sharp}(\bar\iota^e(u),\alpha)=\mu(u)(\alpha)^{-2N}
$$where $u\in X^{\vee}$ and $\alpha$ is any functorial point of 
$K(\cN^{\dagger})^m$. In particular, $e^{\sharp}$ induces the  
monomorphism $2N\mu:H\to D(K(\cN^{\dagger})^m)(B)\ (=X/2N^2\phi(Y))$ 
and hence the epimorphism $D(2N\mu):K(\cN^{\dagger})^m(B)\to D(H)$.
\end{lemma}
\begin{proof}Let $u\in X^{\vee}$. 
By Lemma~\ref{lemma:cN Phi-inv}, there exists a $B$-isomorphism 
$d_u:\cN^{\dagger}\to T_{\bar\iota^e(u)}^*(\cN^{\dagger})$ 
whose $I$-adic completion 
lifts the translation $T_{\bar\iota^e(u)}^{\wedge}$ of $\cG^{\wedge}$.  
Then 
$d_u^{\wedge}$ is decomposed via Eq.~(\ref{eq:Ndaggerwedge})
into  $d_u^{\wedge}=\coprod_{b\in H}
\left(\prod_{x\in X}d^{b,\wedge}_{u,x}\right)$:
\begin{align*}
d^{b,\wedge}_{u,x}:b\times (\cM^{\otimes 2N^2})_x
&\to (\bar\iota^e(u)+b)\times 
T^*_{c^{t,e}(u)}((\cM^{\otimes 2N^2})_x)\\
&\overset{\simeq}{\to} (\bar\iota^e(u)+b)\times
(\cM^{\otimes 2N^2})_{x+2N\mu(u)}.
\end{align*}sending $f^b_x\in(\cM^{\otimes 2N^2})_x$ to 
$T_{c^{t,e}(u)}^*f^{b}_x\in(\cM^{\otimes 2N^2})_{x+2N\mu(u)}$ by $c\circ(2N\mu)=2N^2\lambda\circ c^{t,e}$. 
Since any functorial point 
$\alpha\in K(\cN^{\dagger})^m$ 
acts on $f\in\cO_x$ by $t\cdot f=x(\alpha)f$ by definition,     
it acts on $\cN^{\dagger,\wedge}$ by: 
\begin{align*}
\rho(\alpha)\cdot\left(b\times \prod_{x\in X}f^b_x\right)
&:=\left(b\times \prod_{x\in X}x(\alpha)f^b_x\right).
\end{align*} 
  By $T_{-c^{t,e}(u)}^*f^b_x\in(\cM^{\otimes 2N^2})_{x-2N\mu(u)}$,   
we see 
\begin{align*}
e^{\sharp}(\bar\iota^e(u),\alpha)(b\times f^b_x)
&=d_u^{\wedge}\circ\rho(\alpha)\circ 
(d_u^{\wedge})^{-1}\circ(\rho(\alpha))^{-1}
(b\times f^b_x)\\
&=x(\alpha)^{-1}d_u^{\wedge}\circ\rho(\alpha)((b-\bar\iota^e(u))\times T_{-c^{t,e}(u)}^*f^{b}_x)\\
&=x(\alpha)^{-1}(x-2N\mu(u))(\alpha)(b\times f^b_x)\\
&=\mu(u)(\alpha)^{-2N}(b\times f^b_x),\end{align*}
whence $e^{\sharp}(c^{t,e}(u),\alpha)=\mu(u)(\alpha)^{-2N}$. 
\end{proof}

\begin{lemma}\label{lemma:lift of HB to GshcNdag}
There exists a totally $e^{\sharp}$-isotropic subgroup $H^{\sharp}$ 
of  $K^{\sharp}(\cN^{\dagger})^f(B)$ 
such that $H^{\sharp}\simeq H$ and  
$K^{\sharp}(\cN^{\dagger})^f=H^{\sharp}_B\oplus K(\cN^{\dagger})^f$.
\end{lemma}
\begin{proof}Let $K^m:=K(\cN^{\dagger})^m$, $K^f:=K(\cN^{\dagger})^f$ 
and $K^{\sharp,f}:=K^{\sharp}(\cN^{\dagger})^f$ for now. We have an exact sequence of commutative group schemes: $0\to K^f\to K^{\sharp,f}\to H_B\to 0$, 
and hence an exact sequence of fppf cohomology groups:
$$0\to K^f(B)\to K^{\sharp,f}(B)\overset{j(B)}{\to} H_B(B)
=H\overset{\delta(B)}{\to} H^1_{\fppf}(\Spec B,K^f).$$  
Since $K^f$ is finite flat over $B$, $H^1_{\Zar}(\Spec B,K^f)=0$ by 
\cite[3.2]{Bhatt12}, and {\it a fortiori} $H^1_{\fppf}(\Spec B,K^f)=0$. 
Hence $j(B)$ is surjective.  

Let $(h_i;i\in[1,r])$ be a minimal set 
of generators of $H$.  Since $j(B)$ is surjective, 
there exist elements $\tildeh_i\in K^{\sharp,f}(B)$ 
$(i\in[1,r])$ such that $j(B)(\tildeh_i)=h_i$. 
We shall prove that we can rechoose $\tildeh_i$ by induction 
so that $e^{\sharp}(\tildeh_i,\tildeh_j)=1$ $(\forall i,j\in[1,r])$. 
First we set $\tildeh'_1:=\tildeh_1$ because $e^{\sharp}$ is alternating: 
$e^{\sharp}(\tildeh'_1,\tildeh'_1)=1$. 

Next 
assume that we have chosen $\tildeh'_i$ $(i\in[1,s-1])$ such that 
$\tildeh'_i\in\tildeh_i+K^m(B)$ 
and $e^{\sharp}(\tildeh'_i,\tildeh'_j)=1$ $(\forall i,j\in[1,s-1])$. 
Let $H^{\sharp}_{s-1}$ be the subgroup of $K^{\sharp,f}(B)$ 
generated by  $\tildeh'_i$ and $\tildeh_j$ $(i\in[1,s-1],j\in[s,r])$. 
Then $H^{\sharp}_{s-1}\simeq H$, whence $D(H^{\sharp}_{s-1})\simeq D(H)$. 
By Lemma~\ref{lemma:pairing esharp},   
we can choose  $\epsilon_s\in K^m(B)$ such that $e^{\sharp}(\epsilon_s,h)=e^{\sharp}(\tildeh_s,h)$ $(\forall h\in H^{\sharp}_{s-1})$. We 
set $\tildeh'_s:=\tildeh_s-\epsilon_s$. Then  
$e^{\sharp}(\tildeh'_s,\tildeh'_j)=e^{\sharp}(\tildeh'_s,\tildeh'_s)=1$ 
$(\forall j\in[1,s-1])$ by the anti-symmetry of $e^{\sharp}$. 
Let $H^{\sharp}:=H^{\sharp}_r\simeq H$. 
Then $H^{\sharp}$ is totally $e^{\sharp}$-isotropic 
and  $K^{\sharp,f}=H^{\sharp}_B\oplus K^f$. 
\end{proof}

\begin{lemma}\label{lemma:level subgroups} 
$H^{\sharp}_B$  
is liftable to a level subgroup $B$-scheme 
$H^{\flat}_B$ of $\cG^{\sharp}(\cN^{\dagger})^f$ such that 
$H^{\flat}_B\cap\cG(\cN^{\dagger})^f=\{1\}$ and 
$\cG^{\sharp}(\cN^{\dagger})^f=H^{\flat}_B\cdot\cG(\cN^{\dagger})^f$. 
\end{lemma}
\begin{proof} 
Let $\cH^{\sharp}:=
\cG^{\sharp}(\cN^{\dagger})^f\times_{K^{\sharp,f}}H^{\sharp}_B$. 
Since $H^{\sharp}_B$ is totally $e^{\sharp}$-isotropic 
by Lemma~\ref{lemma:lift of HB to GshcNdag}, 
$\cH^{\sharp}$ is a commutative central extension of 
$H^{\sharp}_B$ by $\bG_{m,B}$. 
By Lemma~\ref{lemma:splitting E_in_Cent_Z/nZ_Gm},  
there exists a finite radical normal 
cover $W$ of $B$ such that 
$\cH^{\sharp}_W$  splits.  
By taking $W$ as $B$ if necessary, $H^{\sharp}_B$  
is liftable to a level subgroup $B$-scheme 
$H^{\flat}_B$ of $\cG^{\sharp}(\cN^{\dagger})^f$. 
The rest is clear from Lemma~\ref{lemma:lift of HB to GshcNdag}. 
\end{proof}

\begin{rem}
$H^{\flat}_B$ does not commute with $K(\cN^{\dagger})^f$ 
by Lemma~\ref{lemma:pairing esharp}.
Hence $\cG^{\sharp}(\cN^{\dagger})^f\neq 
H^{\flat}_B\times_B\cG(\cN^{\dagger})^f$ as group schemes.
\end{rem}

\subsection{$\cG(\cN^{\dagger})^f$-modules}
In this subsection we basically follow \cite[p.47]{FC90}. 
In general, let $M$ be a $B$-module 
with $\mu:=K(\cN^{\dagger})^m\ (=D((X/2N^2\phi(Y))_B))$-action $\rho_M$. 
Then $M=\bigoplus_{\chi\in D(\mu)}M^{\chi}$ where $M^{\chi}:=\{m\in M; 
\rho(\alpha)(m)=\chi(\alpha)m\ 
\text{for any functorial point $\alpha$ of $\mu$}\}$. 
By Lemma~\ref{lemma:KflatNdagger_m}, 
 $K(\cN^{\dagger})^m$ acts on $\cN^{\dagger}$ via 
$K^{\flat}(\cN^{\dagger})^m$, so that  
$\Gamma(\cG_B,\cN^{\dagger})$ and $\Gamma(G_B,\cN^{\dagger}_{|G_B})
=\Gamma(G_B,\cL^{\otimes 2N^2}_B)$ 
are $K(\cN^{\dagger})^m$-modules, which are 
 decomposed respectively into 
\begin{align*}
\Gamma(\cG_B,\cN^{\dagger})&=\bigoplus_{\barx\in X/2N^2\phi(Y)}
\Gamma(\cG_B,\cN^{\dagger})^{\barx},\\
\Gamma(G_B,\cL^{\otimes 2N^2}_B)&=\bigoplus_{\barx\in X/2N^2\phi(Y)}
\Gamma(G_B,\cL^{\otimes 2N^2}_B)^{\barx}.
\end{align*} 
 Note that 
$\Gamma(\cG_K,\cN^{\dagger}_K)^{\barx}
=\Gamma(G_K,\cL^{\otimes 2N^2}_K)^{\barx}.$  

By Lemma~\ref{lemma:KflatNdagger_m}~(\ref{item:KflatNdagger_m}), 
every $\barx\in X/2N^2\phi(Y)$ defines a homomorphism  
$h_{\barx}\in\Hom(\cG(\cN^{\dagger})^m,\bG_{m,B})$ 
such that $h_{\barx}=\id$ on $\bG_{m,B}$ and $\barx$ on 
$K^{\flat}(\cN^{\dagger})^m$. Since  
$\ker(h_{\barx})$ $(\simeq K^{\flat}(\cN^{\dagger})^m)$ 
is a closed  subgroup $B$-scheme 
of $\cG(\cN^{\dagger})^f$ $B$-finite flat, so that we have a quotient 
$\cG(\cL)^f_{\barx}:=\cG(\cN^{\dagger})^f/\ker(h_{\barx})$ 
flat over $B$,  which is a central extension of $K((\cM^{\otimes 2N^2})_x)_B$ by $\bG_{m,B}$. Let 
$\eta_{\barx}:\cG(\cN^{\dagger})^f\to\cG(\cL)^f_{\barx}$ be 
the natural morphism. 
Since $K((\cM^{\otimes 2N^2})_x)=K(\cM^{\otimes 2N^2})$, 
in view of Eq.~(\ref{eq:AV part}), 
we obtain a commutative diagram of exact sequences:
\begin{equation}\label{diagram:cGNdagger_barx}
\smallskip
\begin{diagram}
0&\rTo &\cG(\cN^{\dagger})^m&\rTo^i &\cG(\cN^{\dagger})^f&\rTo& K((\cM^{\otimes 2N^2})_x)_B&\rTo& 0\\
&&\dTo^{h_{\barx}}&&\dTo^{\eta_{\barx}}&&\dTo_{\id}&&\\
0&\rTo &\bG_{m,B}&\rTo &\cG(\cN^{\dagger})^f_{\barx}&\rTo&
K((\cM^{\otimes 2N^2})_x)_B&\rTo& 0.
\end{diagram}
\smallskip
\end{equation} 
Meanwhile we have an exact sequence:
\begin{equation}\label{exact seq for AV}
0\to \bG_{m,B}\to \cG((\cM^{\otimes 2N^2})_{x})_B
\to K((\cM^{\otimes 2N^2})_{x})_B\to 0\quad (x\in X).
\end{equation}
Since $\cN^{\dagger}$ is cubical on $\cG$,  
the biextension 
$\Lambda(\cN^{\dagger})_{K(\cN^{\dagger})^f\times \cG}$ 
is trivial by \cite[IV, 3.3.2-3.3.3]{MB85}. Therefore we have an 
isomorphism $\cN^{\dagger}_x\otimes \cN^{\dagger}_y\simeq\cN^{\dagger}_{x+y}$ 
where $x\in K(\cN^{\dagger})^f, y\in\cG$ are functorial points, so that  
by \cite[I, 3.1, p.~23]{MB85} we have  
an action $\rho$ of 
$\cG(\cN^{\dagger})^f=(\cN^{\dagger})^{\times}_{K(\cN^{\dagger})^f}$ 
on $\cN^{\dagger}$:
\begin{equation}\label{eq:action rho}
(\cN^{\dagger})^{\times}_{K(\cN^{\dagger})^f}\times_B \cN^{\dagger}\overset{\rho}{\to} \cN^{\dagger},
\end{equation}which induces an action of $\cG(\cN^{\dagger})^f$:
\begin{equation}\label{eq:Gamma(rho)}
\Gamma(\rho):
\cG(\cN^{\dagger})^f\times_B 
\Gamma(\cG_B,\cN^{\dagger})\to\Gamma(\cG_B,\cN^{\dagger}).
\end{equation}  
Since $K^{\flat}(\cN^{\dagger})^m$ is contained in the center 
of $\cG(\cN^{\dagger})^f$ 
by Lemma~\ref{lemma:KflatNdagger_m}~(\ref{item:center of cG_N_dagger_f}), 
$\Gamma(\rho)$ is decomposed into 
$\bigoplus_{\barx\in X/2N^2\phi(Y)}\Gamma(\rho)^{\barx}$ where
\begin{equation}\label{eq:Gamma(rho)_barx}
\Gamma(\rho)^{\barx}:
\cG(\cN^{\dagger})_{\barx}^f\times_B 
\Gamma(\cG_B,\cN^{\dagger})^{\barx}\to\Gamma(\cG_B,\cN^{\dagger})^{\barx}.
\end{equation}

Thus we obtain the following compatible data below: 
\begin{enumerate}
\item[(i)]the group structure 
of $\cG(\cN^{\dagger})^f=(\cN^{\dagger})^{\times}_{K(\cN^{\dagger})^f}$
by restricting $\rho$ to 
$$\left((\cN^{\dagger})^{\times}_{K(\cN^{\dagger})^f}\right)^2\ 
(\subset (\cN^{\dagger})^{\times}_{K(\cN^{\dagger})^f}
\times_B \cN^{\dagger});$$
\item[(ii)](in view of Eq.~(\ref{eq:Ndaggerwedge}))
a formal morphism \\ 
from $\left((\cN^{\dagger})^{\times}_{K(\cN^{\dagger})^f}\right)^{\wedge}$ to 
$\left(\bar\iota^e(0)\times((\cM^{\otimes 2N^2})_x)^{\times}_{K((\cM^{\otimes 2N^2})_x)}\right)_B$; 
\item[(iii)] a morphism of formal group $B$-schemes induced from (ii) 
via Eq.~(\ref{diagram:cGNdagger_barx}):
$$\gamma_x^{\wedge}:\left(\cG(\cN^{\dagger})^f_{\barx}\right)^{\wedge}\to 
\left(\cG((\cM^{\otimes 2N^2})_x)_B\right)^{\wedge}$$
which is an isomorphism by Eq.~(\ref{eq:AV part});
\item[(iv)] an action
$\rho_{\barx}:=\Gamma(\rho)^{\barx}$ defined by Eq.~(\ref{eq:Gamma(rho)_barx});
\item[(v)] a formal action 
$\rho_{A,x}^{\wedge}$ of 
$\left(((\cM^{\otimes 2N^2})^{\times}_x)_{K((\cM^{\otimes 2N^2})_x)}\right)^{\wedge}$\\ on $\Gamma(A^{\wedge},((\cM^{\otimes 2N^2})_x)^{\wedge})$, which is 
algebraizable because $A$ is projective.
\end{enumerate}

Since $K(\cN^{\dagger})^f$ and $K((\cM^{\otimes 2N^2})_x)$ 
are $B$-finite, the formal morphism in (iii) is algebraizable. 
Therefore we have a commutative diagram of algebraic/formal actions 
and morphisms:
\begin{equation}\label{diagram:Big diagram of actions}
\begin{diagram}
\cG(\cN^{\dagger})^f_{\barx}&\times&\Gamma(G_B,\cN^{\dagger})^{\barx}&\rTo^{\rho_{\barx}}
&\Gamma(G_B,\cN^{\dagger})^{\barx}\\
\dTo^{\gamma_x^{\wedge}}&&\dTo^{s_x^{\wedge}}&&\dTo^{s_x^{\wedge}}\\
\cG((\cM^{\otimes 2N^2})_x)_B^{\wedge}&\times&\Gamma(A^{\wedge},(\cM^{\otimes 2N^2})_x)^{\wedge}_B&\rTo^{\rho_{A,x}^{\wedge}}&\Gamma(A^{\wedge},(\cM^{\otimes 2N^2})^{\wedge}_x)_B\\ 
\uTo&&\uTo^{\simeq}&&\uTo^{\simeq}\\
\cG((\cM^{\otimes 2N^2})_x)_B&\times&\Gamma(A,(\cM^{\otimes 2N^2})_x)_B&\rTo
^{\rho_{A,x}}&\Gamma(A,(\cM^{\otimes 2N^2})_x)_B, 
\end{diagram}
\end{equation}where $s_x:=(\sigma_x)_{|\Gamma(G,\cN^{\dagger})^{\barx}}$ 
by Definition~\ref{defn:2nd sigma_x}.
\smallskip

Let $\gamma_x$ be the algebraization of $\gamma_x^{\wedge}$. By  
Eqs.~(\ref{diagram:cGNdagger_barx})-(\ref{diagram:Big diagram of actions}), 
we have 
a commutative diagram of exact sequences:
\begin{equation}\label{diag:isom of cGNdaggerf and cGM2N2}
\begin{diagram}
0&\rTo&\bG_{m,B}&\rTo&
\cG(\cN^{\dagger})^f_{\barx}&\rTo&K((\cM^{\otimes 2N^2})_x)_B&\rTo&0\\
&&\dTo^{\|}&&\dTo^{\gamma_x}&&\dTo^{\|}\\
0&\rTo&\bG_{m,B}&\rTo&
\cG((\cM^{\otimes 2N^2})_x)_B&\rTo&K((\cM^{\otimes 2N^2})_x)_B&\rTo&0.
\end{diagram}
\end{equation}Since  
the central extensions in Eq.~(\ref{diag:isom of cGNdaggerf and cGM2N2}) 
are Zariski-locally products of $\bG_{m,B}$ and 
$K((\cM^{\otimes 2N^2})_x)_B$, $\gamma_x$ is an isomorphism. 
\medskip

 Let $d'_a:=\rank_R\Gamma((\cM^{\otimes 2N^2})_x)$ and  
$d'_t:=|X/2N^2\phi(Y)|$.  
By \cite[V, 2.4.2~(iii)]{MB85}, the rank of  
any nonzero $\cG((\cM^{\otimes 2N^2})_{x,K})$-module of weight one 
is equal to an integral multiple of $d'_a$, so at least $d'_a$. 
Since $\gamma_x$ is an isomorphism 
by Eq.~(\ref{diag:isom of cGNdaggerf and cGM2N2}), 
any nonzero $\cG(\cN^{\dagger})^f_{\barx}$-module of weight one 
has rank at least $d'_a$. 
Since $\sigma_x$ is nonzero 
by Theorem~\ref{thm:degeneration data2}~(\ref{item:nonzero sigma}), 
we have $\Gamma(G_B,\cN^{\dagger})^{\barx}\neq 0$, whence    
$\rank_B\Gamma(G_B,\cN^{\dagger})^{\barx}\geq d'_a$. 
Hence 
$$d'_ad'_t\leq\rank_B \Gamma(G_B,\cN^{\dagger})
=\rank_K \Gamma(G_K,\cL^{\otimes 2N^2}_K)=d'_ad'_t.$$ It follows that  
$\Gamma(G_B,\cN^{\dagger})^{\barx}$ is an irreducible 
$\cG(\cN^{\dagger})^f_{\barx}$-module of rank $d'_a$ 
$(\forall x\in X/2N^2\phi(Y))$. 
Moreover since 
$s_x=(\sigma_x)_{|\Gamma(G_B,\cN^{\dagger})^{\barx}}$ is nonzero, 
 $s_x\otimes_BK$ is an isomorphism 
of $\cG(\cN^{\dagger})^f_{\barx,K}$-modules 
via $\gamma_x\otimes_BK$ 
by \cite[V, 2.4.2~(iii)]{MB85}. 

\begin{defn}\label{defn:irreducible cH-module over R}
Let $\cH$ be a group $B$-scheme and $M$ an $B$-module 
on which $\cH$ acts. We call $M$ an {\it irreducible $\cH$-module} if  
any $\cH$ and $B$-submodule $N$ of $M$ is of 
the form $N=JM$ for some ideal $J$ of $B$. 
See \cite[V, 2.3]{MB85}.
\end{defn}

Summarizing the above, we obtain by \cite[V, 2.4.2~(ii)]{MB85}: 
\begin{lemma}\label{lemma:equivariant isom}The following are true: 
\begin{enumerate}
\item\label{item:gammax isom}
 $\gamma_x:\cG(\cN^{\dagger})^f_{\barx}\simeq\cG((\cM^{\otimes 2N^2})_x)$;
\item\label{item:irred cG modules}
$\Gamma(G_B,\cN^{\dagger})^{\barx}$ 
(resp. $\Gamma(A,(\cM^{\otimes 2N^2})_x)$) is an irreducible 
$\cG(\cN^{\dagger})^f_{\barx}$-module (resp. 
 an irreducible $\cG((\cM^{\otimes 2N^2})_x)$-module) of weight one;
\item\label{item:isom phixK} $s_x:\Gamma(G,\cN^{\dagger})^{\barx}
\simeq\Gamma(A,(\cM^{\otimes 2N^2})_{x})$  
as $\cG(\cN^{\dagger})^f_{\barx}$-modules 
via $\gamma_x$;
\item\label{item:extra isom} $s_x$ induces an isomorphism of 
$\cG(\cN^{\dagger})^f_{\barx}$-modules:\\
$\Gamma(A,(\cM^{\otimes 2N^2})_{x})
\simeq\Gamma(A,(\cM^{\otimes 2N^2})_{x+2N^2\phi(y)})$\   
$(\forall y\in  Y)$. 
\end{enumerate}
\end{lemma} 
\begin{proof}(\ref{item:gammax isom})-(\ref{item:isom phixK}) 
have been proved above. 
(\ref{item:extra isom}) follows 
from (\ref{item:gammax isom})-(\ref{item:isom phixK}). 
\end{proof}

\subsection{Degeneration data of $\cG_B$}
\label{subsec:deg data of Neron cGB}
Now we shall construct degeneration data of the pullback $\cG_B$ 
of the N\'eron model $\cG$. Note $\cG_K=G_K$. 
\begin{lemma}\label{lemma:deg data nonzero constants}
There exist nonzero constants $m^{\sharp}(u,x)\in K$ such that 
\begin{gather*}
\sigma_{x+2N\mu(u)}\circ T_{\bar\iota^e(u)}^*
=m^{\sharp}(u,x)T_{c^{t,e}(u)}^*\circ\sigma_{x}, \\
m^{\sharp}(\beta(y),x)=m(y,x):=\psi(y)^{2N^2}\tau(y,x)\ \  
(\forall x\in X,\forall y\in Y, \forall u\in X^{\vee}).
\end{gather*}
\end{lemma}
\begin{proof}
Since $T_{\bar\iota^{e}(u)}^*\cN^{\dagger}\simeq\cN^{\dagger}$\ 
$(\forall u\in X^{\vee}/\beta(Y))$, we have  
by Lemma~\ref{lemma:pairing esharp},  
$T_{\bar\iota^{e}(u)}^*\Gamma(G_K,\cN^{\dagger}_K)^{\barx}
\simeq\Gamma(G_K,\cN^{\dagger}_K)^{\overline{x+2N\mu(u)}}$. By Lemma~\ref{lemma:equivariant isom}~(\ref{item:isom phixK})
\begin{align*}
T_{c^{t,e}(u)}^*\Gamma(A_K,(\cM^{\otimes 2N^2})_{x,K})
&\simeq T_{\bar\iota^{e}(u)}^*\Gamma(G_K,\cN^{\dagger}_K)^{\barx}
\simeq \Gamma(G_K,\cN^{\dagger}_K)^{\overline{x+2N\mu(u)}}\\
&\simeq \Gamma(A_K,(\cM^{\otimes 2N^2})_{x+2N\mu(u),K})
\end{align*}as $\cG(\cN^{\dagger})^f_{\overline{x+2N\mu(u)},K}$-modules. 
Hence  $\sigma_{x+2N\mu(u)}\circ T_{\bar\iota^e(u)}^*$ 
and $T_{c^{t,e}(u)}^*\circ\sigma_{x}$ 
are equivalent isomorphisms from 
$\Gamma(G_K,\cN^{\dagger}_K)^{\barx}$ to 
$\Gamma(A_K,(\cM^{\otimes 2N^2})_{x+2N\mu(u),K})$. 
It follows from Lemma~\ref{lemma:equivariant isom}~
(\ref{item:irred cG modules}) 
that there exists a nonzero constant $m^{\sharp}(u,x)\in K$ such that 
$\sigma_{x+2N\mu(u)}=m^{\sharp}(u,x)T_{c^{t,e}(u)}^*\sigma_{x}$\ $(\forall x\in X, \forall u\in X^{\vee})$, where 
$m^{\sharp}(\beta(y),x)=m(y,x):=\psi(y)^{2N^2}\tau(y,x)$ $(\forall y\in Y)$
by Theorem~\ref{thm:degeneration data2}.
\end{proof}

\begin{cor}\label{cor:msharp l}
Let $\sigma_x^{(l)}:\Gamma(\cG_B,(\cN^{\dagger})^{\otimes l})
\to \Gamma(A_B,(\cM^{\otimes 2N^2l})_{x,B})$  $(x\in X)$ 
be the natural projection. 
Then there exists   
$m^{\sharp,l}(u,x)\in K^{\times}$ such that 
\begin{gather*}
\sigma^{(l)}_{x+2Nl\mu(u)}\circ T_{\bar\iota^e(u)}^*
=m^{\sharp,l}(u,x)T_{c^{t,e}(u)}^*\circ\sigma_{x}^{(l)}, \\
m^{\sharp,l}(\beta(y),x)=\psi(y)^{2N^2l}\tau(y,x)\ \  
(\forall x\in X,\forall y\in Y, \forall u\in X^{\vee}).
\end{gather*}
\end{cor}

\begin{defn}\label{defn:me(u,x)}
We define 
$$\psi^{\sharp,l}(u):=m^{\sharp,l}(u,0),\ 
\tau^{\sharp,l}(u,x):=m^{\sharp,l}(u,x)/\psi^{\sharp,l}(u)\ 
(x\in X,u\in X^{\vee}).$$
 
By definition, $\psi^{\sharp,l}(0)
=\tau^{\sharp,l}(u,0)=\tau^{\sharp,l}(0,x)=1$. 
\end{defn}

\begin{thm}
\label{thm:degeneration data Neron} 
The $B$-module homomorphisms 
$$\sigma_{x}:\Gamma(\cG_B,\cN^{\dagger})\to\Gamma(A_B,(\cM^{\otimes 2N^2})_{x,B})\ \ (x\in X)$$ 
have the following properties:  
\begin{enumerate}
\item\label{item:nonzero sigma Neron} 
$\sigma_{x}\neq 0$\quad $(\forall x\in X)$;
\item\label{item:sigma_psi_tau Neron}there exist a cubic 
trivialization {\small$\psi^{\sharp} : 1_{X^{\vee}} \to 
(c^{t,e})^*\cM^{\otimes (-2N^2)}_K$} 
and a trivialization of biextension 
{\small$\tau^{\sharp} :1_{X^{\vee}\times X}\to 
(c^{t,e}\times c)^*\cP_K^{\otimes(-1)}$}  
such that 
$\sigma_{x+2N\mu(u)}(T^*_{\bar\iota^e(u)}\theta)=
\psi^{\sharp}(u)\tau^{\sharp}(u,x)T^*_{c^{t,e}(u)}(\sigma_{x}(\theta))$\quad 
$(\forall \theta\in \Gamma(\cG_B,\cN^{\dagger}), 
\forall x\in X,\forall u\in X^{\vee})$; 
in particular, $\psi^{\sharp}(0)=\tau^{\sharp}(0,x)=\tau^{\sharp}(u,0)=1$;  
\item\label{item:tau_sharp Neron} 
$\tau^{\sharp}(u,x)$ is bilinear in $u$ and $x$; 
\item\label{item:psi_sharp/tau_sharp Neron} 
$\tau^{\sharp}(u,2N\mu(v))  
=\psi^{\sharp}(u+v)/\psi^{\sharp}(u)\psi^{\sharp}(v)$\ 
$(\forall u,v\in X^{\vee})$;
\item\label{item restriction to beta of tausharp psisharp}
 $\tau^{\sharp}(\beta(y),x)=\tau(y,x)$,   
$\psi^{\sharp}(\beta(y))=\psi(y)^{2N^2}$\ 
$(\forall x\in X,y\in Y)$;
\item\label{item:positive tau_sharp Neron} $\tau^{\sharp}(u,2N\mu(u))\in I^2$ 
 $(\forall u\in X^{\vee}\setminus\{0\})$;
\item\label{item:positive psi_sharp Neron} for every $n\geq 0$, 
$\psi^{\sharp}(u)\in I^n$ for all but finitely many $u\in X^{\vee}$;  
\item\label{eq:Gamma GK Ndagger Neron}
$\Gamma(G_K,\cN^{\dagger}_K)$ is identified with 
the $K$-vector subspace of Fourier series 
$\theta=\sum_{x\in X}\sigma_x(\theta)$ satisfying 
$\sigma_{x+2N^2\phi(y)}(\theta)=
\psi(y)^{2N^2}\tau(y,x)T^*_{c^t(y)}(\sigma_{x}(\theta))$ 
with $\sigma_{x}(\theta)\in \Gamma(\cA_K,(\cM^{\otimes 2N^2})_{x,K})$ $(\forall x\in X,\forall y\in Y)$.
\end{enumerate}
\end{thm}
\begin{proof}See \cite[II, 5.4]{FC90}. 
First note that (\ref{item:nonzero sigma Neron}) 
follows from 
Theorem~\ref{thm:degeneration data2}~(\ref{item:nonzero sigma}), 
while  (\ref{item:sigma_psi_tau Neron}) 
follows from Definition~\ref{defn:me(u,x)}, (\ref{item:tau_sharp Neron}) and 
(\ref{item:psi_sharp/tau_sharp Neron}) which we prove below.  

Now we shall prove (\ref{item:tau_sharp Neron}) and 
(\ref{item:psi_sharp/tau_sharp Neron}). 
First we prove that $\tau^{\sharp}(u,x)$ is linear in $x$. 
Let $\Phi:G_K\times_KG_K\to G_K\times_KG_K$ 
(resp. $\Phi_A:A_K\times_KA_K\to A_K\times_KA_K$) 
be the isogeny $(z,w)\mapsto (z+w, z-w)$. 
Let $u,u'\in X^{\vee}$ and $x\in X$.
By \cite[(AF'), p.~48]{FC90} and 
Theorem~\ref{thm:degeneration data2} 
(applied to $\cN^{\dagger}$ symmetric) 
\begin{gather}
\Phi_A^*\circ (\sigma_{x+2N\mu(u)}\otimes_{K}\sigma_{x'+2N\mu(u)})
=(\sigma^{(2)}_{x+x'+4N\mu(u)}\otimes_{K}\sigma^{(2)}_{x-x'})\circ\Phi^*,
\label{eq:cubic rel 1st}\\
\Phi_A^*\circ (\sigma_{x}\otimes_{K}\sigma_{x'})
=(\sigma^{(2)}_{x+x'}\otimes_{K}\sigma^{(2)}_{x-x'})\circ\Phi^*.
\label{eq:cubic rel 2nd}
\end{gather} 

We compute Eqs.~(\ref{eq:cubic rel 1st}) 
and (\ref{eq:cubic rel 2nd}):
\begin{align*}
&\Phi_A^*\circ (\sigma_{x+2N\mu(u)}\otimes_{K}\sigma_{x'+2N\mu(u)})\\
&=\psi^{\sharp}(u)^2\tau^{\sharp}(u,x)\tau^{\sharp}(u,x')
\left(\Phi_A^*\circ(T_{c^{t,e}(u)}^*\sigma_x\otimes_KT_{c^{t,e}(u)}^*\sigma_{x'})\right)\\
&=\psi^{\sharp}(u)^2\tau^{\sharp}(u,x)\tau^{\sharp}(u,x')
(T_{c^{t,e}(2u)}^*\sigma^{(2)}_{x+x'}
\otimes_K\sigma^{(2)}_{x-x'})\circ\Phi^*,\\
&(\sigma^{(2)}_{x+x'+4N\mu(u)}\otimes_{K}\sigma^{(2)}_{x-x'})
\circ\Phi^*\\
&=m^{\sharp,2}(u,x+x')(T_{c^{t,e}(2u)}^*\sigma^{(2)}_{x+x'}\otimes_K\sigma^{(2)}_{x-x'})\circ\Phi^*,
\end{align*}whence 
we have $\psi^{\sharp}(u)^2\tau^{\sharp}(u,x)\tau^{\sharp}(u,x')=
m^{\sharp,2}(u,x+x')$.  By substituting $x+x'$ for $x$ and $0$ for $x'$, 
we have 
\begin{equation}\label{eq:msharp 2}
\psi^{\sharp}(u)^2\tau^{\sharp}(u,x+x')\tau^{\sharp}(u,0)=
m^{\sharp,2}(u,x+x').
\end{equation}
By $\tau^{\sharp}(u,0)=1$, 
$\tau^{\sharp}(u,x)\tau^{\sharp}(u,x')=
\tau^{\sharp}(u,x+x')$: $\tau^{\sharp}(u,x)$ is linear in $x$.

By definition,
\begin{equation}\label{eq:psi_sharp/tau_sharp} 
\begin{aligned}
m^{\sharp}(u+u',x)&=m^{\sharp}(u,x)m^{\sharp}(u',x+2N\mu(u)),
\\ 
\psi^{\sharp}(u+u')&=m^{\sharp}(u+u',0)
=\psi^{\sharp}(u)m^{\sharp}(u',2N\mu(u))\\
&=\psi^{\sharp}(u)\psi^{\sharp}(u')\tau^{\sharp}(u',2N\mu(u)),
\end{aligned}
\end{equation}
which proves (\ref{item:psi_sharp/tau_sharp Neron}). 
Since $\tau^{\sharp}(u,x)$ is linear in $x$, 
we obtain by Eq.~(\ref{eq:psi_sharp/tau_sharp}) 
\begin{align*}
&\tau^{\sharp}(u+u',x)
=m^{\sharp}(u,x)m^{\sharp}(u',x+2N\mu(u))/\psi^{\sharp}(u+u')\\
&=\psi^{\sharp}(u)\psi^{\sharp}(u')\tau^{\sharp}(u,x)
\tau^{\sharp}(u',x+2N\mu(u))/\psi^{\sharp}(u+u')
=\tau^{\sharp}(u,x)\tau^{\sharp}(u',x),
\end{align*}whence $\tau^{\sharp}(u,x)$ is linear in $u$. This 
proves  (\ref{item:tau_sharp Neron}) 
and hence (\ref{item:sigma_psi_tau Neron}) by 
 (\ref{item:psi_sharp/tau_sharp Neron}).  

For $u\in X^{\vee}$, there exists $y\in Y$ with 
$\beta(y)=Nu$. Then  
(\ref{item:positive tau_sharp Neron}) 
follows from (\ref{item:tau_sharp Neron}) and 
Theorem~\ref{thm:degeneration data2}~(\ref{item:positive tau}) because  
$\tau^{\sharp}(u,2N\mu(u))^{N}=\tau(y,2N\phi(y)))\in I^{2N}$.  
 Since  
$\psi^{\sharp}(u)^{2N^2}
=\tau^{\sharp}(u,2N\mu(u))^{N^2}=\tau(y,2N^2\phi(y))
=\psi(y)^{4N^2}$, (\ref{item:positive psi_sharp Neron}) 
follows from 
Theorem~\ref{thm:degeneration data2}~(\ref{item:positive psi}). 
Finally (\ref{eq:Gamma GK Ndagger Neron}) 
follows from  
Theorem~\ref{thm:degeneration data2}~(\ref{item:Gamma Geta L general}) and 
Lemma~\ref{lemma:deg data nonzero constants}. 
This completes the proof.
\end{proof}

\begin{cor}\label{cor:tae_e=tau_sharp}The following are true:
\begin{enumerate}
\item\label{item:tau_sharp=tau_e} $\tau^{\sharp}=\tau^e$ on $X^{\vee}\times X$;
\item\label{item:psi_sharp(u)=psi_sharp(-u)} $\psi^{\sharp}(u)=\psi^{\sharp}(-u)$\ $(\forall u\in X^{\vee})$;
\item\label{item:psi_sharp_2=psi_e_2N} $(\psi^{\sharp})^2(u)=(\psi^e)^{2N}(u)=\tau^{\sharp}(u,2N\mu(u))$ 
$(\forall u\in X^{\vee})$.
\end{enumerate}
\end{cor}
\begin{proof}
As in Eq.~(\ref{eq:transf tSy on fx}), the 
translation $T_{\iota^e(u)}$ of $\tG_K$ $(u\in X^{\vee})$
induces a homomorphism 
$\tdelta_u:T_{c^{t,e}(u)}^*\cO_{x}\otimes_RK\to\cO_{x}\otimes_RK$ by 
$$\tdelta_u\left(T_{c^{t,e}(u)}^*(f_x)\right)
=\tau^e(u,x)T_{c^{t,e}(u)}^*f_x.$$ 
Moreover $T_{\iota^e(u)}$ also induces 
a homomorphism 
$$\tdelta_u:\tT_{\bar\iota^e(u)}^*\Gamma(G_K,\cN^{\dagger}_{|G_K})
\to \Gamma(G_K,\cN^{\dagger}_{|G_K})$$ 
as follows. Let $s=\sum_{x\in X}\sigma_x(s)\in\Gamma(G_K,\cN^{\dagger}_{|G_K})$. 
By Theorem~\ref{thm:degeneration data Neron}~(\ref{item:sigma_psi_tau Neron})
\begin{gather*}
\tdelta_{u}\left(\tT_{\bar\iota^e(u)}^*(s)\right)
:=\sum_{x\in X}\sigma_{x+2N\mu(u)}\left(\tT_{\bar\iota^e(u)}^*(s)\right),\label{eq:action delta_u on sigma_x(s)}\\ 
\sigma_{x+2N\mu(u)}\left(\tT_{\bar\iota^e(u)}^*(s)\right)=\psi^{\sharp}(u)\tau^{\sharp}(u,x)T_{c^{t,e}(u)}^*(\sigma_{x}(s)). 
\label{eq:relation of sigma_x(s) and sigma_x+2nmuu(s)}
\end{gather*}
 
Now choose $s\in\Gamma(G_K,\cN^{\dagger}_{|G_K})$ with 
$\sigma_x(s)\neq 0$ and $\sigma_0(s)\neq 0$. 
Then 
$$\sigma_{z+2N\mu(u)}(\tT^*_{\bar\iota^e(u)}s)=
\psi^{\sharp}(u)\tau^{\sharp}(u,x)T^*_{c^{t,e}(u)}(\sigma_{z}(s))
\neq 0$$ for $z=0,x$. It follows 
\begin{equation}\label{eq:delta_u on sigma_x sigma0}
\begin{aligned}\tau^e(u,&\,x)T_{c^{t,e}(u)}^*
\left(\frac{\sigma_x(s)}{\sigma_0(s)}\right)=
\tdelta_uT_{c^{t,e}(u)}^*\left(\frac{\sigma_x(s)}{\sigma_0(s)}\right):=
\frac{\tdelta_uT_{c^{t,e}(u)}^*\sigma_x(s)}{\tdelta_uT_{c^{t,e}(u)}^*\sigma_0(s)}\\
&=\frac{\sigma_{x+2N\mu(u))}(\tT_{\bar\iota^e(u)}^*(s))}{\sigma_{2N\mu(u)}(\tT_{\bar\iota^e(u)}^*(s))}
=\tau^{\sharp}(u,x)T_{c^{t,e}(u)}^*\left(\frac{\sigma_x(s)}{\sigma_0(s)}
\right). 
\end{aligned}
\end{equation} 

Hence $\tau^{\sharp}(u,x)=\tau^e(u,x)$. 
This proves (\ref{item:tau_sharp=tau_e}). 
Let $s\in\Gamma(G_K,\cN^{\dagger}_{|G_K})$. 
Since $\cM$ and $\cN^{\dagger}$ are symmetric, we have $\sigma_{x+2N\mu(u)}\circ T^*_{\bar\iota^e(u)}(s)=\sigma_{-x-2N\mu(u)}\circ T^*_{\bar\iota^e(-u)}(s)$ and $\psi^{\sharp}(u)\tau^{\sharp}(u,x)
=\psi^{\sharp}(-u)\tau^{\sharp}(-u,-x)=
\psi^{\sharp}(-u)\tau^{\sharp}(u,x)$. 
This proves (\ref{item:psi_sharp(u)=psi_sharp(-u)}).
By 
Theorem~\ref{thm:degeneration data Neron}~(\ref{item:psi_sharp/tau_sharp Neron}), $\psi^{\sharp}(u)^2=\tau^{\sharp}(u,2N\mu(u))$, while 
$\psi^e(u)^2=\tau^e(u,2\mu(u))$ by
Lemma~\ref{lemma:minimal Galois with extensions cte iotae psie taue}~(iv).  
Hence we have (\ref{item:psi_sharp_2=psi_e_2N}) 
by (\ref{item:tau_sharp=tau_e}).  
\end{proof}

By Theorem~\ref{thm:degeneration data Neron},   
we obtain a similar result for $(\cG_B,(\cN^{\dagger})^{\otimes l})$. 
Morover in view of Corollary~\ref{cor:tae_e=tau_sharp}, 
we can derive a similar result for 
$\cN^{\otimes 2Nl}$ over $S^e$ which can be expressed  
in terms of $\psi^e$ and $\tau^e$.

\begin{thm}
\label{thm:degeneration data Neron 2l case} 
Let $l\in\bN$. Let $R^e$ the integral closure of $R$ in $k(\eta^e)$, 
$S^e=\Spec R^e$\ 
(Definition~\ref{defn:minimal K with extensions cte iotae psie taue}) 
and $\cN^{\natural}:=\cN^{\otimes N}_{S^e}$. 
The $R^e$-module homomorphisms 
$$\sigma^{(2l)}_{x}:\Gamma(\cG_{S^e},(\cN^{\natural})^{\otimes 2l})\to \Gamma(A_{S^e},(\cM^{\otimes 4N^2l})_{x,S^e})\ (x\in X)$$ 
have the following properties:  
\begin{enumerate}
\item\label{item:nonzero sigma Neron 2l} 
$\sigma^{(2l)}_{x}\neq 0$\quad $(\forall x\in X)$ and 
\begin{equation}\label{eq:Phi acts on Gamma(cG,cNdagger 2l)}
\sigma^{(2l)}_{x+4Nl\mu(u)}(T^*_{\bar\iota^e(u)}\theta)=
\psi^e(u)^{2Nl}\tau^e(u,x)
T^*_{c^{t,e}(u)}(\sigma^{(2l)}_{x}(\theta))
\end{equation}
$(\forall \theta\in\Gamma(G,(\cN^{\natural})^{\otimes 2l}), 
\forall x\in X,\forall u\in X^{\vee})$;
\item\label{eq:Gamma GK Ndagger 2l}
$\Gamma(G_{\eta^e},(\cN^{\natural}_{\eta^e})^{\otimes 2l})$ 
is identified with 
the $k(\eta^e)$-vector subspace of Fourier series 
$\theta=\sum_{x\in X}\sigma^{(2l)}_x(\theta)$ such that 
\begin{equation}\label{eq:Sy inv theta}
\sigma^{(2l)}_{x+4N^2l\phi(y)}(\theta)=
\psi(y)^{4N^2l}\tau(y,x)T^*_{c^t(y)}(\sigma^{(2l)}_{x}(\theta))
\end{equation} 
with $\sigma^{(2l)}_{x}(\theta)\in \Gamma(\cA_{\eta^e},(\cM^{\otimes 4N^2l})_{x,\eta^e})$ $(\forall x\in X,\forall y\in Y)$.
\end{enumerate}
\end{thm}
\begin{proof}
By Lemma~\ref{lemma:minimal Galois with extensions cte iotae psie taue} 
and Corollary~\ref{cor:tae_e=tau_sharp}, all the extensions 
$c^{t,e}$, $\iota^e$, $\psi^e$ and $\tau^e$ are defined over $k(\eta^e)$. 
The finite cover $B\to S^e$ is radical normal such that   
$\Aut(B/S^e)=\mu:=D(F)$: the Cartier dual of a finite 
commutative constant group $S^e$-scheme $F$. 
Since $\cN^{\dagger}=\cN^{\natural}_B$, 
by \cite[11.3/11.9]{MN24}
\begin{align*}
\Gamma(G_{Z},(\cN^{\natural}_Z)^{\otimes 2l})
&=\Gamma(G_{W},(\cN^{\dagger}_W)^{\otimes 2l})^{\mu},\\
\Gamma(A_{Z},(\cM^{\otimes 4N^2l})_{x,Z})
&=\Gamma(A_{W},(\cM^{\otimes 4N^2l})_{x,W})^{\mu} 
\end{align*}where $(Z,W)=(S^e,B), (\eta^e,K)$.   
Hence we obtain Corollary.  
\end{proof}

\begin{rem}\label{rem:Phi acts on includes Sy inv}
Eq.~(\ref{eq:Sy inv theta})  
is a particular case 
of Eq.~(\ref{eq:Phi acts on Gamma(cG,cNdagger 2l)}) 
under $S_y^*(\theta)=\theta$. 
\end{rem}

\section{Twisted Mumford families}
\label{sec:twisted Mumford families}

\subsection{The eFC data of a semiabelian scheme}
\label{subsec:eFC data of semiabelian}
Let $(G,\cL)$ be a split semiabelian $S$-scheme 
and $\zeta:=\FC(G,\cL)$. 

\begin{defn}\label{defn:minimal Galois with extensions cte iotae psie taue}
By Lemma~\ref{lemma:minimal Galois with extensions cte iotae psie taue}, 
we define $K:=K_{\min}(\zeta)$ to 
be the minimal normal extension of $k(\eta)$ 
such that all the extensions $c^{t,e},\iota^e, \psi^e$ 
and $\tau^e$ are defined over $K$. 
Let $R_{\min}(\zeta)$ be the integral closure of $R$ 
in $K_{\min}(\zeta)$,  
$S_{\min}(\zeta):=\Spec R_{\min}(\zeta)$ 
and $\eta_{\min}(\zeta)$ 
the generic point of $S_{\min}(\zeta)$. 
 Note that 
\begin{equation}\label{eq:extension of cte iotae taue psie by beta}
c^{t,e}\circ\beta=c^t,\ \iota^e\circ\beta=\iota,\  
\tau^e\circ(\beta\times\id)=\tau,\ \psi^e\circ\beta=\psi^{2N}.
\end{equation} 
\end{defn}

\begin{notation}
\label{notation:Notation Rinit Sinit keta=Kmin(xi) for zeta}
In what follows, let $R_{\init}$ be a CDVR with $t:=s_{\init}$ uniformizer, 
$v_t$ the valuation of $k(\eta_{\init})$ with $v_t(t)=1$, 
$S_{\init}:=\Spec R_{\init}$, 
$\eta_{\init}$ (resp. $0_{\init}$) 
the generic point (resp. the closed point) of $S_{\init}$, 
$k(\eta_{\init})$ (resp. $k(0_{\init})$)
the fraction field (resp. the residue field) of $R_{\init}$, 
$(G,\cL)$ a split semiabelian $S_{\init}$-scheme 
and $\zeta=\FC(G,\cL)$ the symmetric split object over $S_{\init}$ 
given by Eq.~(\ref{eq:split obj zeta in DDample}).  
Let $k(\eta):=K_{\min}(\zeta)$, $R:=R_{\min}(\zeta)$, 
$S:=S_{\min}(\zeta)=\Spec R$, $\eta$ the generic point of $S$, 
$s$ a uniformizer of $R$ 
and $e(\zeta)$ the ramification index of 
$k(\eta)/k(\eta_{\init})$, {\it i.e.}, $v_s(t)=e(\zeta)$. 
We extend $v_t$ to $R$ by $v_t=v_s/e(\zeta)$. 
\end{notation}

\begin{rem}\label{rem:zeta and zetam}
Let $\zeta_{l}:=\FC(G,\cL^{\otimes l})$. Then 
$K_{\min}(\zeta_{l})=K_{\min}(\zeta)$,\ 
$R_{\min}(\zeta_{l})=R_{\min}(\zeta)$ and 
 $\eta_{\min}(\zeta_{l})=\eta_{\min}(\zeta)$\ $(\forall l\in\bN)$
by Theorem~\ref{thm:degeneration data2} and Lemma~\ref{lemma:minimal Galois with extensions cte iotae psie taue}. Note that $K_{\min}(\zeta)=k(\eta^e)$ and Theorem~5.23 is true over $k(\eta^e)$. Hence we can set $B=S_{\min}(\zeta)=S=S^e$ 
and $\cN^{\natural}=\cN^{\dagger}$ in Theorem~\ref{thm:degeneration data Neron 2l case}.  
\end{rem}

\begin{defn}\label{defn:2nl th eFC}  We define the  
{\it $4N^2l$-th eFC datum of $(G,\cL)$} extending $\zeta$ 
(or the {\it $4N^2l$-th  extension of $\zeta$})  by  
{\small\begin{equation*}
\zeta^e_{4N^2l}=(\tG,A,T,X,Y,c,c^{t,e},\iota^e,
4N^2l\lambda,4N^2l\phi,\tau^e,\tcL^{\otimes 4N^2l},
(\psi^e)^{2Nl},\cM^{\otimes 4N^2l}). 
\end{equation*}}

  Recall that $\psi^e(u):
=\tau^e(u,\mu(u))$ $(u\in X^{\vee})$. Then $\psi^e$ 
is compatible with $\tau^e$ 
in the sense that Eq.~(\ref{eq:cubic triv psie}) is true.

Under Notation~\ref{notation:Notation Rinit Sinit keta=Kmin(xi) for zeta}, 
every eFC datum of $(G,\cL)$ extending $\zeta$ 
is conjugate to each other under the action of 
$\Aut(k(\eta)/k(\eta_{\init}))$ by 
Lemma~\ref{lemma:minimal Galois with extensions cte iotae psie taue}.
\end{defn}

\begin{defn}\label{defn:Btau Sigma_dagger_l0}
Let $B_{\tau}$, $\beta_{\tau}$ and $\mu$ 
be the same as Definition~\ref{defn:miscellany}.
Let $B_{\tau^e}:=v_t(\tau^e)$. Note that 
$B_{\tau^e}(u,x)=u(x)$ 
$(\forall u\in X^{\vee}, \forall x\in X)$
by \cite[5.10]{MN24}.  
Let  $\epsilon^{\dagger}(u):=\tau^e(u,2N\mu(u))=\psi^e(u)^{2N}$ 
and $E^{\dagger}(u):=v_t\epsilon^{\dagger}(u)=2Nu(\mu(u))$\ 
$(u\in X^{\vee})$.  
For $l\in\bN$, we define 
$\epsilon^{\dagger}_{l}=(\epsilon^{\dagger})^{l}=(\psi^e)^{2Nl}$, 
$E^{\dagger}_{l}=lE^{\dagger}$ and 
\begin{gather*}
\Sigma^{\dagger}_{l}(0)=\{x\in X_{\bR}; E^{\dagger}_{l}(u)+u(x)\geq 0\ 
(\forall u\in X^{\vee})\},\  
\Sigma^{\dagger}_{l}=\Sigma^{\dagger}_{l}(0)\cap X,\\ 
\Sigma_{l}^{\dagger,\alpha}=(\alpha+4Nl\mu(X^{\vee}))
\cap\Sigma^{\dagger}_{l}\ \ (\alpha\in\Sigma^{\dagger}_{l}).
\end{gather*}
\end{defn}

\begin{rem}\label{rem:identify E_dagger_l with E_2Nl}
By \cite[6.3]{MN24}, $\Sigma^{\dagger}_{l}+4Nl\mu(X^{\vee})=X$.
For $c\in X^{\vee}$ we define  
 $\Sigma^{\dagger}_{l}(c):=\Sigma^{\dagger}_{l}(0)+4Nl\mu(c)$. 
The set $\Vor^{\dagger}_{l}$  consisting of 
all $\Sigma^{\dagger}_{l}(c)$ $(c\in X^{\vee})$ and their closed faces 
 forms a polyhedral decomposition of $X_{\bR}$:  
$$X_{\bR}=\bigcup_{c\in X^{\vee}}\Sigma^{\dagger}_{l}(c),\  
\Sigma^{\dagger}_{l}(c)^0\cap\Sigma^{\dagger}_{l}(c')^0=\emptyset\ 
(c\neq c').
$$ 

We can identify under the notation of \cite[5.21]{MN24}:
\begin{equation}\label{eq:E_dagger_l=E_2Nl}
E^{\dagger}_{l}=E_{2Nl},\ 
\Sigma^{\dagger}_{l}(0)=\Sigma_{2Nl}(0),\ 
\Sigma^{\dagger}_{l}=\Sigma_{2Nl},\ \Vor^{\dagger}_{l}=\Vor_{2Nl}.
\end{equation} 
\end{rem}

\begin{lemma}
\label{lemma:lth basic formulae general}
The following are true: 
\begin{enumerate}
\item\label{item:lth epsilon beta= psi2N} 
$\epsilon^{\dagger}_{l}(\beta(y))=\psi(y)^{4N^2l}$ and 
$\tau^e(\beta(y),x)=\tau(y,x)$;
\item\label{item:lth be} $\tau^e(u,\mu(v))=\tau^e(v,\mu(u))$;
\item\label{item:lth epsilon u+v} $\epsilon^{\dagger}_{l}(u+v)
=\epsilon^{\dagger}_{l}(u)\epsilon^{\dagger}_{l}(v)\tau^e(u,\mu(v))^{4Nl}$;
\item\label{item:lth epsilon u+betay} $\epsilon^{\dagger}_{l}(u+\beta(y))=\epsilon^{\dagger}_{l}(u)\psi(y)^{4N^2l}\tau(y,\mu(u))^{4Nl}$;
\item\label{item:geq 0}
$\epsilon^{\dagger}_{l}(u)\tau^e(u,\alpha)\in R$ 
\end{enumerate}where $x\in X$, $y\in Y$,  $u,v\in X^{\vee}$ and 
$\alpha\in\Sigma^{\dagger}_{l}$.  
\end{lemma}
\begin{proof}(\ref{item:lth epsilon beta= psi2N}) follows from 
$\epsilon^{\dagger}_{l}(\beta(y))
=\tau(y,2N^2l\phi(y))=\psi(y)^{4N^2l}$. 
(\ref{item:lth be})-(\ref{item:geq 0}) are clear from 
Lemma~\ref{lemma:minimal Galois with extensions cte iotae psie taue} and Definitions~\ref{defn:minimal Galois with extensions cte iotae psie taue}/\ref{defn:Btau Sigma_dagger_l0}.
\end{proof}

\subsection{The Mumford construction associated with a N\'eFC kit}
\label{subsec:twisted Mumford construction}
\begin{defn}\label{defn:NeFC kit}
Let $l\in\bN$. We define an {\it $l$-th N\'eFC kit of $\zeta^e_{4N^2}$} by
\begin{equation}\label{deq:simple kit of zetae}
\xi^{\dagger}_{l}=\xi^{\dagger}(\zeta^e_{4N^2l})
=(X,Y,\epsilon^{\dagger}_{l},\tau^e,E^{\dagger}_{l},\Sigma^{\dagger}_{l}),
\ \ \xi^{\dagger}:=\xi^{\dagger}_1. 
\end{equation} 
We define  $\NeFC(G,\cL^{\otimes 4N^2l})=(\zeta^e_{4N^2l},\xi^{\dagger}_{l})$ which we call the {\it N\'eFC pair of $(G,\cL^{\otimes 4N^2l})$}.
Note that $\zeta^e_{4N^2l}$ and $\xi^{\dagger}_{l}$ are defined 
over $S_{\min}(\zeta)$. 
\end{defn}

We start the Mumford construction associated with 
$\xi^{\dagger}_{l}$: 
\begin{defn}For $l\in\bN$, 
we define fractional ideals of $k(\eta)$ as follows: 
\begin{gather*}
J_u:=R\epsilon^{\dagger}(u),\ J_{u,x}:=R\tau^e(u,x),\   
J_{l,u}:=J_u^{l}\ 
(x\in X, u\in X^{\vee}).
\end{gather*} 

By Lemma~\ref{lemma:lth basic formulae general}, 
$$J_{l,u+v}=J_{l,u}\cdot J_{l,v}\cdot  J_{u,\mu(v)}^{4Nl},\  
J_{u,\mu(v)}=J_{v,\mu(u)}.$$ 
\end{defn}

\begin{defn}\label{defn:R1Sigma,R2Sigma}(\cite[III, 3.1]{FC90})\ \ 
We define two graded algebras 
$\cS_{l,i}(\xi^{\dagger}):=\cS_{l,i}$ 
and their graded subalgebras 
$\cR_{i,\Sigma^{\dagger}_{l}}(\xi^{\dagger}):=\cR_{i,\Sigma^{\dagger}_{l}}$ $(i=1,2)$ by:
\begin{equation}\label{eq:R1Sigma,R2Sigma}
\begin{aligned}
\cS_{l,1}&=\left(\bigoplus_{n\geq 0}\pi_*(\cO_{\tG})\theta_{l}^n\right)\otimes_Rk(\eta)=\left(\bigoplus_{n\geq 0}(\bigoplus_{x\in X}\cO_{x})\theta_{l}^n\right)\otimes_Rk(\eta),\\
\cS_{l,2}&
=\left(\bigoplus_{n\geq 0}\left(\bigoplus_{x\in X}\cM^{\otimes 4N^2ln}\otimes_{\cO_A}\cO_{x}\right)\theta_{l}^n\right)\otimes_Rk(\eta),\\
\cR_{1,\Sigma^{\dagger}_{l}}&=\cO_A[(J_{l,v}\cdot J_{v,\beta}\cdot\cO_{\beta+4Nl\mu(v)})\theta_{l};v\in X^{\vee},\beta\in\Sigma^{\dagger}_{l}],\\
\cR_{2,\Sigma^{\dagger}_{l}}&=\cO_A[(J_{l,v}\cdot J_{v,\beta}
\cdot\cM^{\otimes 4N^2l}\otimes_{\cO_A}
\cO_{\beta+4Nl\mu(v)})\theta_{l};
v\in X^{\vee},\beta\in\Sigma^{\dagger}_{l}]
\end{aligned}
\end{equation}where 
$\beta+4Nl\mu(v)\in\Sigma^{\dagger}_{l}+4Nl\mu(X^{\vee})=X$ by 
\cite[6.3]{MN24}.  
Note that 
\begin{equation}\label{eq:cSl2}
\cS_{l,2}=\left(\bigoplus_{n\geq 0}\pi_*(\tcL^{\otimes 4N^2ln})\theta_{l}^n
\right)\otimes_Rk(\eta).
\end{equation}

Since $\cM$ is an invertible $\cO_A$-module, we have
$\Proj \cR_{1,\Sigma^{\dagger}_{l}}=\Proj\cR_{2,\Sigma^{\dagger}_{l}}$, 
which we denote by $\tQ_{l}$.  The $S$-scheme $\tQ_{l}$ is endowed with a $\pi^Q_{l}$-ample tautological invertible sheaf $\cN^Q_{l}$ as $\Proj \cR_{1,\Sigma^{\dagger}_{l}}$, while it is endowed with a $\pi^Q_{l}$-ample tautological invertible sheaf $\tcL^Q_{l}:=\cN^Q_{l}\otimes_{\cO_A}(\pi^Q_{l})^*\cM^{\otimes 4N^2l}$ as $\Proj \cR_{2,\Sigma^{\dagger}_{l}}$ where $\pi^Q_{l}:\tQ_{l}\to A$ 
is the natural projection.   
\end{defn}

\begin{defn}\label{defn:action tildeSy}
Let $\tT_{\iota^e(u)}$ be the translation of $\tG$ by $\iota^e(u)$ 
 $(u\in X^{\vee})$ and 
set $\delta_{l,u}:=\tT_{\iota^e(u)}$. 
 Let 
$\sigma_{x},\sigma'_{x}\in\cM^{\otimes 4N^2l}\otimes_{\cO_A}\cO_{x,\eta}$. 
Then  
we define a $k(\eta)$-isomorphism 
$\tdelta_{l,u}:\delta_{l,u}^*\cS_{l,2}\otimes_Rk(\eta)\to\cS_{l,2}\otimes_Rk(\eta)$  by
\begin{align}\label{eq:action delta_lv}
\tdelta_{l,u}\delta_{l,u}^*(\bigoplus_{x\in X}\sigma_{x}\theta_{l})
&=\bigoplus_{x\in X}\sigma'_{x}\theta_{l},\ 
\sigma'_{x+4Nl\mu(u)}=\epsilon^{\dagger}_{l}(u)\tau^e(u,x)
T_{c^{t,e}(u)}^*\sigma_{x}.
\end{align} The sheaf automorphism $\tdelta_{l,u}$ lies 
over the identity of the base $A$. 

Let $S_{l,y}:=\delta_{l,\beta(y)}=\tT_{\iota(y)}$ and 
$\tS_{l,y}:=\tdelta_{l,\beta(y)}$\ $(y\in Y)$. Hence  
\begin{equation}\label{eq:action S_ly}
\begin{aligned}
\tS_{l,y}S_{l,y}^*(\bigoplus_{x\in X}\sigma_{x}\theta_{l})
&=\bigoplus_{x\in X}\sigma''_{x}\theta_{l},\ 
\sigma''_{x+4N^2l\phi(y)}=\psi(y)^{4N^2l}\tau(y,x)
T_{c^t(y)}^*\sigma_{x}.
\end{aligned}
\end{equation}
  
Since $\theta_{l}\in\cO_0\theta_{l}=\cO_A\theta_{l}$, we have  
\begin{equation}\label{eq:tdelta_lu delta_lu*} 
\begin{aligned}
\tdelta_{l,u}\delta_{l,u}^*(\theta_{l})&
=\epsilon^{\dagger}_{l}(u)\theta_{l}\in
\cO_{4Nl\mu(u)}\theta_{l},\\ 
\tS_{l,y}S_{l,y}^*(\theta_{l})&=\psi(y)^{4N^2l}\theta_{l}\in\cO_{4N^2l\phi(y)}\theta_{l}.
\end{aligned}
\end{equation}
\end{defn}

\begin{lemma}\label{lemma:tdelta_u stable}
$\tdelta_{l,u}$ (resp. $\tS_{l,y}$) is an isomorphism   
$\delta_u^*\cS_{l,2}\otimes_Rk(\eta)\simeq\cS_{l,2}\otimes_Rk(\eta)$ 
(resp. $S_y^*\cS_{l,2}\otimes_Rk(\eta)\simeq\cS_{l,2}\otimes_Rk(\eta)$), 
which keeps the $R$-subalgebra $\cR_{2,\Sigma^{\dagger}_{l}}$ 
of $\cS_{l,2}$ stable: 
$\tdelta_u\delta_u^*\cR_{2,\Sigma^{\dagger}_{l}}
=\cR_{2,\Sigma^{\dagger}_{l}}$,\ (resp. 
$\tS_yS_y^*\cR_{2,\Sigma^{\dagger}_{l}}=\cR_{2,\Sigma^{\dagger}_{l}}).$  
\end{lemma}
\begin{proof} 
 By Lemma~\ref{lemma:lth basic formulae general}~(\ref{item:lth epsilon u+v}), 
\begin{align*}
&\tdelta_{l,u}\delta_{u,\eta}^*
(J_{l,v}\cdot J_{v,\alpha}\cdot\cM^{\otimes 4N^2l}_{\alpha+4Nl\mu(v)}
\theta_{l})\\
&=\epsilon_{l}(u)\tau^e(u,\alpha+4Nl\mu(v))\cdot 
T_{c^{t,e}(u)_{\eta}}^*(J_{l,v}\cdot J_{v,\alpha}\cdot\cM^{\otimes 4N^2l}_{\alpha+4Nl\mu(v)}\theta_{l})\\
&=J_{l,u}\cdot J_{l,v}\cdot J_{v,\alpha}\cdot\tau^e(u,\alpha+4Nl\mu(v))\cdot  T_{c^{t,e}(u)_{\eta}}^*(\cM^{\otimes 4N^2l}_{\alpha+4Nl\mu(v)})\theta_{l}\\
&=J_{l,u+v}\cdot J_{u+v,\alpha}\cdot 
T_{c^{t,e}(u)_{\eta}}^*(\cM^{\otimes 4N^2l}_{\alpha+4Nl\mu(v)})\theta_{l}\\
&=J_{l,u+v}\cdot J_{u+v,\alpha}\cdot \cM^{\otimes 4N^2l}_{\alpha+4Nl\mu(u)+4Nl\mu(v)}\theta_{l}.
\end{align*}
 This proves Lemma.  See \cite[III, 3.3]{FC90}.
\end{proof}  

\begin{rem}\label{rem:rel complete model}
The condition \cite[(*), p.~62]{FC90} is Lemma~\ref{lemma:lth basic formulae general}~(\ref{item:geq 0}), which 
ensures that $\tG$ is an open subscheme of $\tP_{l}$. 
See Lemma~\ref{lemma:tG}.  By \cite[III, 3.3]{FC90}, 
$\tQ_{l}$ is a relatively complete model 
$\tP_{(4N^2l\phi,\psi^{4N^2l}),\Sigma^{\dagger}_{l}}$ 
associated with a symmetric split object 
$\zeta^e_{4N^2l}$. 
\end{rem}

\begin{defn}\label{defn:tPl tcLl}
Let $\tP_{l}$ be the normalization of the $S$-scheme $\tQ_{l}$ locally of finite type, and $\tcL_{l}$ (resp. $\tcN_{l}$) 
the pullback to $\tP_{l}$ of $\tcL^Q_{l}$ (resp. $\tcN^Q_{l}$) on 
$\tQ_{l}$. Then $\tcL_{l}=\tcN_{l}\otimes_{\cO_{\tP_{l}}}\tpi^*_{l}\cM^{\otimes 4N^2l}$ where $\tpi_{l}:\tP_{l}\to A$ is the natural projection.
The actions $(S_{y},\tS_{l,y})$ and 
$(\delta_{u},\tdelta_{l,u})$ on $(\tQ_{l},\tcL^Q_{l})$ lift 
to those on $(\tP_{l},\tcL_{l})$, which we denote by the same letters 
$(S_{l,y},\tS_{l,y})$ and 
$(\delta_{l,u},\tdelta_{l,u})$ respectively.  
\end{defn}

 Let  $(\tP_{l}^{\wedge},\tcL_{l}^{\wedge})$ be the $I$-adic completion 
of $(\tP_{l},\tcL_{l})$ and $(S^{\wedge}_{l,y},\tS^{\wedge}_{l,y})$ 
the $I$-adic completion of $(S_{l,y},\tS_{l,y})$ $(y\in Y)$. 
Since the action of $Y$ on $(\tP_{l}^{\wedge},\tcL_{l}^{\wedge})$ via $((S_{l,y}^{\wedge},\tS_{l,y}^{\wedge}); y\in Y)$ is properly discontinuous, there exists a formal quotient $(\tP_{l}^{\wedge},\tcL_{l}^{\wedge})/Y$ 
of $(\tP_{l}^{\wedge},\tcL_{l}^{\wedge})$ 
by $Y$, which is a formal projective $S^{\wedge}$-scheme 
with $\tcL_{l}^{\wedge}$ relatively ample.  
By \cite[III$_1$, 5.4.5]{EGA}, 
there exists an algebraization $(P_{l},\cL_{l})$ of 
$(\tP_{l}^{\wedge},\tcL_{l}^{\wedge})/Y$. 

\begin{defn}\label{defn:lth Mumford fam Pl cLl ass to zetae xidagger_l}
We call $(P_{l},\cL_{l})$ the {\it $l$-th twisted Mumford family associated 
with} $\xi^{\dagger}$ (or the {\it twisted Mumford family associated 
with} $\xi^{\dagger}_{l}$),  which we denote 
by $(P_{l}(\xi^{\dagger}),\cL_{l}(\xi^{\dagger}))$ if necessary. 
We also denote $\cL_{l}$ by $\cL^P_{l}$ if necessary.  
 Let $\varpi_{l}:P_{l}\to S$  and $\pi_{l}:\tP_{l,0}\to P_{l,0}$ 
be the natural projections respectively. 
\end{defn}

\begin{notation}\label{notation:cD(0)}
For any (general) graded algebra $\cD$, 
we denote the degree zero part of $\cD$ by $\cD(0)$.
\end{notation}
\begin{defn}\label{defn:W_l,alpha,u}
Let $u\in X^{\vee}$, $\alpha\in\Sigma^{\dagger}_{l}$ and let 
$\xi_{l,\alpha,u,\cM}$ be a (local) generator of the $\cO_A$-module 
$J_{l,u}\cdot J_{u,\alpha}\cdot\cM^{\otimes 4N^2l}\otimes_{\cO_A}
\cO_{\alpha+4Nl\mu(u)}$ and 
\begin{equation*}
\cR_{2,\Sigma^{\dagger}_{l},\alpha, u}:
=\left(\cR_{2,\Sigma^{\dagger}_{l}}[1/(\xi_{l,\alpha,u,\cM}\theta_{l})]\right)(0).
\end{equation*}
This definition makes sense because 
$\cM$ and $\cO_x$ $(x\in X)$ are invertible and $R$ is a UFD. 
Then we define 
\begin{gather*}
U_{l,\alpha,u}=\Spec\cR_{2,\Sigma^{\dagger}_{l},\alpha,u},\\
W_{l,\alpha,u}=\text{the normalization of $U_{l,\alpha,u}$}. 
\end{gather*}    
For any affine open subset $V$ of $A$, $U_{l,u,\alpha}\times_AV$ (resp. 
$W_{l,\alpha,u}\times_AV$) is an affine open $S$-subscheme of $\tQ_{l}$ 
(resp. $\tP_{l}$) and  
$\cW:=(W_{l,\alpha,u}; \alpha\in\Sigma^{\dagger}_{l}, u\in X^{\vee})$ 
is an open covering of $\tP_{l}$. 
\end{defn}

\begin{rem}
Since $A$ is proper over $S$, 
 by \cite[7.6]{MN24},  
 $\tP_{l}(K')=\tP_{l}(R')$ for 
any finite extension $K'$ of $k(\eta)$ and 
the integral closure $R'$ of $R$ in $K'$. 
\end{rem}

\subsection{Relation with torus embeddings}
\label{subsec:relation with torus embeddings}
Let $\bZ m_0\ (\simeq \bZ)$ be a rank one $\bZ$-module, 
$(\bZ m_0)^{\vee}=\bZ f_0$ with $f_0(m_0)=1$, $\tX=\bZ m_0\oplus X$, 
$\tX^{\vee}=\bZ f_0\oplus X^{\vee}$ and $t$ a uniformizer 
of $R$ identified with $w^{m_0}$. \par
\begin{defn}\label{defn:ep_dagger,E_dagger,D_dagger}
Let $\epsilon^{\dagger}_{l}$, $E^{\dagger}_{l}$ and 
$\Sigma^{\dagger}_{l}$ be 
the same as Definition~\ref{defn:Btau Sigma_dagger_l0}. 
Let $\gamma\in\Sigma^{\dagger}_{l}$, 
$z\in X^{\vee}$ and $x=\gamma+4Nl\mu(z)$. Then we define 
$D^{\dagger}_{l}(x)=E^{\dagger}_{l}(z)+z(\gamma)$.  
This is well-defined as is shown in the same manner as in   
\cite[6.6/6.7]{MN24}.
\end{defn} 
\begin{defn}\label{defn:fan ass with xi_dagger_l}
Let $\xi^{\dagger}_{l}$ 
be the $\NeFC$ kit of $\zeta^e_{2Nl}$ in Definition~\ref{defn:NeFC kit}. 
We define a {\it fan (in $X^{\vee}$) associated with $\xi^{\dagger}_{l}$} by:
\begin{gather*}
\Fan(\xi^{\dagger}_{l}):=
\left\{\tau^{\dagger}_{l,\alpha,u}\ \text{and their faces}; 
\alpha\in\Sigma^{\dagger}_{l}, u\in X^{\vee}\right\},\\
\tau^{\dagger}_{l,\alpha,u}:=
\left\{\begin{matrix}
v_0f_0+v\in \bR f_0\oplus X_{\bR}^{\vee};\ 
v_0\geq 0, v\in X_{\bR}^{\vee}, \forall y\in X\\ 
e(\zeta)(D^{\dagger}_{l}(y)-D^{\dagger}_{l}(\alpha+4Nl\mu(u)))v_0
+v(y-(\alpha+4Nl\mu(u)))\geq 0
\end{matrix}
\right\}.
\end{gather*}
\end{defn}

\begin{rem}\label{rem:comparison with MN24 xi_natural}
We set $\xi^{\natural}_{2Nl}:=
(X,Y,t^{E^{\dagger}_{l}},t^{B_{\tau^e}},E^{\dagger}_{l},
\Sigma^{\dagger}_{l})$.
Then  $\xi^{\natural}_{2Nl}$ is an eFC datum, 
which is the $l$-th N\'eFC datum induced from $\xi^{\natural}_{2N}$ 
\cite[5.21]{MN24}.  
Note that $\Fan(\xi^{\dagger}_{l})
=\Fan(\xi^{\natural}_{2Nl})$ by \cite[7.4]{MN24} because 
$E^{\dagger}_{l}(u)=2Nl u\mu(u)=E_{2Nl}(u)$ and hence 
$D^{\dagger}_{l}=D_{2Nl}$
by \cite[5.24]{MN24}.  
\end{rem}

Let $U(\Fan(\xi^{\dagger}_{l}))$ be the torus $T_X$-embedding over $S$ 
associated with $\Fan(\xi^{\dagger}_{l})$. 
Let $l_1:T\times_S\tG\to \tG$ 
(resp. $l_2:T\times_S U(\tau^{\dagger}_{l,\alpha,u})
\to U(\tau^{\dagger}_{l,\alpha,u})$)
 be the left action of $T$ via $\iota_T$ (resp.  
via the isomorphism $T\simeq T_{X,S}$). 
See Eq.~(\ref{eq:exact seq of tG}) for $\iota_T$. 
Let $\varphi:=(l_1\times_Sl_2)\circ(\id_T\times\inv_T)$  
where $\inv_T$ is the inverse morphism of $T$.  
We define the contraction product of $U(\tau^{\dagger}_{l,\alpha,u})$ and 
$\tG$ by $T$ by  
\begin{equation}\label{eq:contraction product}
U(\tau^{\dagger}_{l,\alpha,u})\times_{S}\tG/\varphi(T).
\end{equation} 
Since $\tG$ is a Zariski $T$-bundle, 
the quotient (\ref{eq:contraction product}) 
is an $S$-scheme. 
 As is suggested by \cite[3.5]{K98}, we have an isomorphism  
\begin{equation}\label{eq:isom with contraction prod}
W_{l,\alpha,u}\simeq U(\tau^{\dagger}_{l,\alpha,u})\times^{T}_{S}\tG
:=U(\tau^{\dagger}_{l,\alpha,u})\times_{S}\tG/\varphi(T).
\end{equation}

Let $\xi_{l,\alpha,v}:
=\epsilon^{\dagger}_{l}(v)\tau^e(v,\alpha)w^{\alpha+4Nl\mu(v)}$. 
Let $\tcL'_{l}$ be the invertible sheaf 
on $U(\Fan(\xi^{\dagger}_{l}))$ defined by a one cocyle 
$f_{l,\alpha,u;\beta,v}:=\xi_{l,\alpha,u}/\xi_{l,\beta,v}$ 
with respect to the covering $\cW$. 
By \cite[7.5]{MN24} and Eq.~(\ref{eq:isom with contraction prod}), we see 
\begin{equation}
\begin{aligned}
\tP_{l}&\simeq U(\Fan(\xi^{\dagger}_{l}))\times^{T}_{S}\tG,\ \ 
\tcL_{l}\simeq \tcL'_{l}\otimes_{S}\pi^*\cM^{\otimes 4N^2l}.
\end{aligned}
\end{equation} 
Since both $U(\Fan(\xi^{\dagger}_{l}))$ and $\tG$ are 
irreducible, so is $\tP_{l}$.

\subsection{A down-to-earth description of $\tP_{l}$}
\label{subsec:down_to_earth tPl}
This subsection is a presheaf version of  
\S~\ref{subsec:twisted Mumford construction}. 
In terms of local torus coordinates $w^x_U$, 
the description of $\tP_{l}$ becomes 
very similar to \cite[7.1-7.3]{MN24}. \par
First we choose an affine covering $\cF^0$ of $A$ such that 
$\tG_U\simeq T_{U}$ 
and $\cM_{|U}\simeq \cO_U$ $(\forall U\in\cF^0)$.  
Let $v\in X^{\vee}$, $U\in\cF^0$. 
$U^v:=T_{-c^{t,e}(v)}(U)$,  
$\cF^v:=\{U^v; U\in \cF^0\}$ and  
$\cF:=\coprod_{v\in X^{\vee}}\cF^v$.
By the assumption on $\cF$, there exists a generator 
$w^x_U$ of $(\cO^x)_U$
 $(\forall x\in X,\forall U\in\cF)$ such that 
 $\tG_U\simeq\Spec\Gamma(\cO_U)[w^x_U;x\in X]$ with  
$w^x_U\cdot w^{x'}_U=w^{x+x'}_U$ $(\forall x,x'\in X)$. 
We set $g^x_{UV}=w^x_V/w^x_U$. By Eq.~(\ref{eq:-Fx=c(x)=Ox}), 
 $(g^x_{UV};U,V\in\cF)$ is a one-cocycle representing 
$\cO^{x}\ (\simeq \cO_{x})$. 
Let $(\mu_{UV};U,V\in\cF)$ be 
a one-cocycle representing 
 $\cM\in H^1(A,\cO^*_{A})$ with respect to $\cF$.

Let $l\in\bN$ and $(\pi^{-1}\cM)^{x}:=(\pi^{-1}\cM)\otimes_{\cO_A}\cO^x$. 
Then as in \S~\ref{subsec:Fourier series general}, 
we obtain a homomorphism
\begin{align*}
&\Gamma(\cG,(\cN^{\dagger})^{\otimes 2l})
\hookrightarrow\Gamma(G^{\wedge},(\cN^{\dagger})^{\otimes 2l,\wedge})
\simeq\Gamma(\tG^{\wedge},(\pi^{\wedge})^*\cM^{\otimes 4N^2l,\wedge})\\
&\simeq\left(\bigoplus_{x\in X}\Gamma(\tG,(\pi^{-1}\cM^{\otimes 4N^2l})^{x})\right)^{I\text{-adic}}
\hookrightarrow 
\prod_{x\in X}\Gamma(\tG,(\pi^{-1}\cM^{\otimes 4N^2l})^{x}), 
\end{align*}
whose $x$-th projection sending $\Gamma(\cG,(\cN^{\dagger})^{\otimes 2l})$ to $\Gamma(\tG,(\pi^{-1}\cM^{\otimes 4N^2l})^{x})$ 
we denote by $s^{(2l)}_x$.  Let 
$\nu_x:\Gamma(A,(\cM^{\otimes 4N^2l})_{x,B})\to
\Gamma(\tG,(\pi^{-1}\cM^{\otimes 4N^2l})^{x})$ 
be the natural isomorphism.  
By Theorem~\ref{thm:degeneration data Neron 2l case}, 
we have $s^{(2l)}_x=\nu_x\circ\sigma^{(2l)}_x$, so that  
$\theta\in\Gamma(\cG,(\cN^{\dagger})^{\otimes(2l)})$ can be expressed as 
$$\theta=\sum_{x\in X}s^{(2l)}_x(\theta),\ 
s^{(2l)}_x(\theta)\in \Gamma(\tG,(\pi^{-1}\cM^{\otimes 4N^2l})^{x}).
$$

\begin{defn}\label{defn:rho_x_2l_theta theta_l_U}
 Let $u,v\in X^{\vee}$ and $U,V\in\cF$. 
We write $s^{(2l)}_x(\theta)=(\rho^{(2l)}_{x,U}(\theta)w^x_U;U\in\cF)$ 
with $\rho^{(2l)}_{x,U}(\theta)\in\Gamma(U,\cM^{\otimes 4N^2l})$ 
 subject to $\rho^{(2l)}_{x,U}(\theta)=\mu^{4N^2l}_{UV}\rho^{(2l)}_{x,V}(\theta)$, and let $\rho^{(2l)}_x(\theta):=(\rho^{(2l)}_{x,U}(\theta);U\in\cF)\in
\Gamma(A,\cM^{\otimes 4N^2l})$. 
Let $\theta_{l,U}$ be an indeterminate of degree one 
subject to  
$\theta_{l,U}=\mu^{-4N^2l}_{UV}\theta_{l,V}$. By definition, $\rho^{(2l)}_{x,U}(\theta)\theta_{l,U}=\rho^{(2l)}_{x,V}(\theta)\theta_{l,V}$.
\end{defn}
\begin{defn}For $U\in\cF$, we define 
\begin{gather*}
\cR^{\dagger}_{2,\Sigma^{\dagger}_{l},U}=
\Gamma(\cO_U)\left[
\rho_U\xi_{l,\alpha,v,U}\theta_{l,U};
\begin{matrix}
\rho=(\rho_V;V\in\cF)\in\Gamma(A,\cM^{\otimes 4N^2l})\\
v\in X^{\vee},\ \alpha\in\Sigma^{\dagger}_{l}
\end{matrix}
\right],\\
\xi_{l,\alpha,v,U}=\epsilon^{\dagger}_{l}(v)
\tau^e(v,\alpha)w_U^{\alpha+4Nl\mu(v)}.
\end{gather*}  
Since $\Gamma(A,\cM^{\otimes 4N^2l})$ is base point free, 
$\cR^{\dagger}_{2,\Sigma^{\dagger}_{l},U}\otimes_{\Gamma(\cO_U)}\cO_U
=\cR^{\dagger}_{2,\Sigma^{\dagger}_{l}}\otimes_{\cO_A}\cO_U.$
\end{defn}

\begin{defn}\label{defn:U_alpha_u_U W_alpha_u_U}
Let $\alpha\in\Sigma^{\dagger}_{l}$, $u\in X^{\vee}$ and $U\in\cF$. 
We define 
\begin{align*}
A_{l,\alpha,u, U}&:=\Gamma(\cO_U)[\xi_{l,\gamma,v,U}/\xi_{l,\alpha,u,U}; \gamma\in\Sigma^{\dagger}_{l}, v\in X^{\vee}],\\
B_{l,\alpha,u,U}&:=\text{the normalization of $A_{l,\alpha,u, U}$,}\\
U_{l,\alpha,u, U}&:=\Spec A_{l,\alpha,u, U},\ 
W_{l,\alpha,u, U}:=\Spec B_{l,\alpha,u, U}.
\end{align*}

Note that $U_{l,\alpha,u, U}\simeq U_{l,\alpha,u}\times_{A}U$ and 
$W_{l,\alpha,u, U}\simeq W_{l,\alpha,u}\times_AU$. 
Let  $\tQ_{l,U}:=\Proj \cR^{\dagger}_{2,\Sigma_{l},U}$ and 
$\tP_{l,U}$ the normalization of $\tQ_{l,U}$. Then  
$\tQ_{l,U}$ (resp. $\tP_{l,U}$) 
is covered by $U_{l,\alpha,u,U}$ (resp. 
$W_{l,\alpha,u,U}$) $(\alpha\in\Sigma^{\dagger}_{l}, u\in X^{\vee})$. 
Note that $\tQ_{l,U}\simeq \tQ_{l}\times_{A}U$ and 
$\tP_{l,U}\simeq \tP_{l}\times_{A}U$ 
with the notation of Definition~\ref{defn:tPl tcLl}.
\end{defn}

\begin{defn}\label{defn:delta_l_u on U}
The sheaf automorphisms $\tdelta_{l,u}$ and $\tS_{l,y}$ given in 
Definition~\ref{defn:action tildeSy} 
induce automorphisms of $\cR^{\dagger}_{2,\Sigma^{\dagger}_{l},U}$
in the natural manner:
\begin{align*}
\tdelta_{l,u}\delta_{l,u}^*(\rho_U\xi_{l,\alpha,v,U}\theta_{l,U})&:
=(T_{c^{t,e}(u)}^*\rho)_U\xi_{l,\alpha,v+u,U}\theta_{l,U},\\
\tS_{l,y}S_{l,y}^*(\rho_U\xi_{l,\alpha,v,U}\theta_{l,U})&:
=(T_{c^{t}(y)}^*\rho)_U\xi_{l,\alpha,v+\beta(y),U}\theta_{l,U},
\end{align*}where $\rho:=(\rho_U;U\in \cF)\in\Gamma(A,\cM^{\otimes 4N^2l})$. 
By Definition~\ref{defn:rho_x_2l_theta theta_l_U},
\begin{align*}
\tdelta_{l,u}
\delta_{l,u}^*(\rho^{(2l)}_x(\theta)_U\theta_{l,U})
&=(T_{c^{t,e}(u)}^*\rho^{(2l)}_x(\theta))_U\xi_{l,0,u,U}\theta_{l,U}\\
&=(T_{c^{t,e}(u)}^*\rho^{(2l)}_x(\theta))_U\epsilon^{\dagger}_{l}(u)w_U^{4Nl\mu(u)}\theta_{l,U}.
\end{align*} 
 Note that 
$\tdelta_{l,u}$ (resp. $\tS_y$) lifts $\delta_{l,u}$ (resp. $S_y$) 
from $\tP_{l}$ to $\tcL_{l}$.
\end{defn}

\begin{lemma}\label{lemma:delta_u_delta_u well defined}
$\tdelta_{l,u}$ is well-defined.
\end{lemma}
\begin{proof}This is clear from Eq.~(\ref{eq:action delta_lv}). 
However here we prove Lemma based on Definition~\ref{defn:delta_l_u on U}. 
 Let $\rho:=(\rho_U;U\in\cF)\in\Gamma(A,\cM^{\otimes 4N^2l})$. Then  
\begin{equation}\label{eq:T_cteu_rho in cM_4N2l_4Nl_mu(u)}
T_{c^{t,e}(u)}^*\rho\in\Gamma(A,(\cM^{\otimes 4N^2l})_{4Nl\mu(u)}).
\end{equation}   
By $c\circ\mu=N\lambda\circ c^{t,e}$, we have $T_{c^{t,e}(u)}^*\cM^{\otimes N}\simeq\cM^{\otimes N}\otimes \cO_{\mu(u)}$, so that we may assume $T_{c^{t,e}(u)}^*\mu^{N}_{UV}=\mu^{N}_{UV}g^{\mu(u)}_{UV}$. Also we may assume 
$T_{c^{t,e}(u)}^*g^x_{UV}=g^x_{UV}$ by 
$T_{c^{t,e}(u)}^*\cO_x\simeq\cO_x$. Hence  
$(T_{c^{t,e}(u)}^*\rho)_U=\mu^{4N^2l}_{UV}
\cdot g^{4Nl\mu(u)}_{UV}\cdot (T_{c^{t,e}(u)}^*\rho)_V$. 
Let 
\begin{gather*}
\Theta_U:=\rho_U\cdot\xi_{l,\alpha,v,U}\cdot\theta_{l,U}\in\cR_{2,\Sigma^{\dagger}_{l},U},\\ 
\Theta_V:=\rho_V\cdot\xi_{l,\alpha,v,V}\cdot\theta_{l,V}\in\cR_{2,\Sigma^{\dagger}_{l},V}. 
\end{gather*} Therefore by $\xi_{l,\alpha,v,U}=g_{UV}^{-\alpha-4Nl\mu(v)}
\xi_{l,\alpha,v,V}$  and $\rho_U\cdot\theta_{l,U}=\rho_V\cdot\theta_{l,V}$. 
\begin{center}
$\Theta_U$ and $\Theta_V$ are 
 identified iff $\Theta_U=g^{-\alpha-4Nl\mu(v)}_{UV}\Theta_V$. 
\end{center}
Then we obtain $\tdelta_{l,u}\delta_{l,u}^*(\Theta_U)=g^{-\alpha-4Nl\mu(v)}_{UV}\tdelta_{l,u}\delta_{l,u}^*(\Theta_V)$:
\begin{align*}
&\tdelta_{l,u}\delta_{l,u}^*(\Theta_U)=(T_{c^{t,e}(u)}^*\rho)_U\cdot\xi_{l,\alpha,u+v,U}\theta_{l,U}\\
&=\mu^{4N^2l}_{UV}\cdot g^{4Nl\mu(u)}_{UV}\cdot 
(T_{c^{t,e}(u)}^*\rho)_V\cdot g^{-\alpha-4Nl\mu(u+v)}_{UV}\cdot
\xi_{l,\alpha,u+v,V}\cdot\mu^{-4N^2l}_{UV}\cdot\theta_{l,V}\\
&=g^{-\alpha-4Nl\mu(v)}_{UV}(T_{c^{t,e}(u)}^*\rho)_V
\cdot\xi_{l,\alpha,u+v,V}\theta_{l,V}
=g^{-\alpha-4Nl\mu(v)}_{UV}\tdelta_{l,u}\delta_{l,u}^*(\Theta_V).
\end{align*}  This proves that $\tdelta_{l,u}\delta_{l,u}^*$, or rather, $\tdelta_{l,u}$ is well-defined. \end{proof}

\begin{notation}\label{notation:theta expansion in pd case}Our notation in 
Definitions~\ref{defn:R1Sigma,R2Sigma}-\ref{defn:action tildeSy} can be 
sometimes confusing because it is vague 
to which direct factor $\sigma_x$ may belong. 
To remedy this, we introduce the following notation 
by taking \S~\ref{subsec:down_to_earth tPl} into  
account. Let $\theta_{l}$ be an indeterminate and $\theta\in\Gamma(P_{l}^{\wedge},\cL_{l}^{\wedge})$. Then we write 
\begin{align}\label{eq:defn sum rho2l_xw_x}
\theta&=\sum_{x\in X}\rho_x^{(2l)}(\theta)w^x\theta_{l}
\end{align}for $\rho_x^{(2l)}(\theta)\in\Gamma(A^{\wedge},(\cM^{\otimes 4N^2l})^{\wedge}_x)$. 
By Eq.~(\ref{eq:defn sum rho2l_xw_x}),  we mean the following: 
\begin{align*}
\theta_{U}&=\sum_{x\in X}(\rho_x^{(2l)}(\theta))_Uw^x_U
\theta_{l,U}\ \  (U\in\cF_0).
\end{align*}
Under Eq.~(\ref{eq:defn sum rho2l_xw_x}), 
we understand Eq.~(\ref{eq:tdelta_lu delta_lu*}) as
\begin{equation}\label{eq:ideltadelta tSySy w weight}
\begin{aligned}
\tdelta_u\delta_u^*(\theta_{l})
&=\epsilon^{\dagger}_{l}(u)w^{4Nl\mu(u)}\theta_{l},\\
\tS_yS_y^*(\theta_{l})&=\psi(y)^{4N^2l}
w^{4N^2l\phi(y)}\theta_{l}.
\end{aligned}
\end{equation}
\end{notation}

\section{The relative compactification $P_{l}$}
\label{sec:str of Pl in pd case}

\subsection{The N\'eFC kit $\xi^{\dagger}$ over $S_{\min}$}
Let $(G,\cL)$ a split {\it partially degenerate} 
semiabelian $R_{\init}$-scheme 
 with $\cL$ ample symmetric cubic (and hence rigidified), 
and $\cG$ the N\'eron model of $G_{\eta_{\init}}$. 
Let $\zeta:=\FC(G,\cL)$ 
and $(\zeta^e_{4N^2l},\xi^{\dagger}_{l}):=\NeFC(G,\cL^{\otimes 4N^2l})$ 
(Definitions~\ref{defn:2nl th eFC}/\ref{defn:NeFC kit}). 
Then $\zeta$ is a split object over $S$ in $\DD_{\ample}$. The data  
$\zeta$, $\zeta^e_{4N^2l}$ and $\xi^{\dagger}_{l}$ are given as follows:
\smallskip
{\small\begin{gather*}\label{eq:split obj/kit in DDample as deg data}
\zeta=(\tG,A,T,X,Y,c,c^t,\iota,\lambda,\phi,\tau,\tc,\psi,\cM),\\
\zeta^e_{4N^2l}=(\tG,A,T,X,Y,c,c^{t,e},\iota^e,
4N^2l\lambda,4N^2l\phi,\tau^e,\tcL^{\otimes 4N^2l},
(\psi^e)^{2Nl},\cM^{\otimes 4N^2l}),\\
\xi^{\dagger}_{l}=\xi^{\dagger}(\zeta^e_{4N^2l})
=(X,Y,\epsilon^{\dagger}_{l},\tau^e,E^{\dagger}_{l},\Sigma^{\dagger}_{l}),
\ \ \xi^{\dagger}=\xi^{\dagger}_1 . 
\end{gather*} 
}

\vspace*{-0.2cm} 
Let $(\tP_{l},\tcL_{l})$ 
the normalization of a relatively complete model 
$(\tQ_{l},\tcL_{l})$ associated with $\zeta^e_{4N^2l}$ 
(Remark~\ref{rem:rel complete model}) and $(P_{l},\cL_{l})$   
the $l$-th twisted Mumford family associated with $\xi^{\dagger}$ 
(Definition~\ref{defn:lth Mumford fam Pl cLl ass to zetae xidagger_l}).  
We study the structures of 
$(\tP_{l},\tcL_{l})$ and $(P_{l},\cL_{l})$. 
In \S~\ref{sec:str of Pl in pd case}, 
we follow Notation~\ref{notation:Notation Rinit Sinit keta=Kmin(xi) for zeta}, 
and hence $k(\eta)=K_{\min}(\zeta)$, $R=R_{\min}(\zeta)$ and 
$S=\Spec R$.  

\subsection{The structure of $\tP_{l}$ and $\tP_{l,0}$}
\label{subsec:structure of tPl tPl0}
Let $E^{\dagger}_{l}(u)=v_t(\epsilon^{\dagger}_{l})$ and 
let $\Sigma^{\dagger}_{l}(0)$ and  $\Sigma^{\dagger}_{l}$ 
(resp. $\Fan(\xi^{\dagger}_{l})$ and $\tau_{l,\alpha,u}$) be
the same as in Definition~\ref{defn:Btau Sigma_dagger_l0} (resp. 
Definition~\ref{defn:fan ass with xi_dagger_l}). 
By \cite[6.5]{MN24}, 
there exists a positive integer $l_0$ such that 
$\Sigma^{\dagger}_{l_0}(0)$ is integral, 
and so is $\Sigma^{\dagger}_{l_0l}(0)$ 
$(\forall l\in\bN)$. In what follows we assume
\begin{equation}
\label{assump:Sigma_dagger_l(0) integral}
\text{$\Sigma^{\dagger}_{l}(0)$ is integral.}
\end{equation}

\begin{lemma}\label{lemma:theta on Pl general}
Let $N:=|B_{\tau}|$. There exists $m_0$ such that for any $m\geq m_0$, 
$\Gamma(P_{l,\eta},\cL_{l,\eta}^{\otimes m})$ is identified with  
the following $k(\eta)$-vector space:
\begin{equation*}\left\{
\theta=\sum_{x\in X}\sigma_x;
\begin{matrix}
\sigma_{x+4N^2l m\phi(y)}(\theta)
=\psi(y)^{4N^2l m}\tau(y,x)T_{c^t(y)}^*\sigma_x(\theta)\\
\sigma_x(\theta)\in\Gamma(A_{\eta},
(\cM^{\otimes 4N^2l m}\otimes_{\cO_A}\cO_x)_{\eta})\\
(\forall x\in X, \forall y\in Y)\end{matrix}
\right\}
\end{equation*}
\end{lemma}
\begin{proof}Lemma is proved in parallel to \cite[7.12]{MN24} by 
Theorem~\ref{thm:G' isom G'' iff FC same} and Lemma~\ref{lemma:H0(G,Lm) as pd Fourier} instead of \cite[3.13/4.21]{MN24}.  The pd-counterpart 
of \cite[4.16~(3)]{MN24} applied in \cite[4.21]{MN24} is true 
by \cite[4.2]{Mumford72}/\cite[III, \S~4]{FC90}.  
\end{proof}

\begin{thm}\label{thm:scheme normalization isom general}
The following are true:
\begin{enumerate}
\item there exists an $S$-isomorphism
$\tpsi_{l,ll'}:(\tP_{l},\tcL_{l}^{\otimes l'})
\simeq (\tP_{ll'},\tcL_{ll'})$
such that $\tpsi_{l,ll'}(W_{l,\alpha,u})=W_{ll',l'\alpha,u}$ 
$(\forall \alpha\in\Sk^0(\Vor^{\dagger}_{l}),\forall u\in X^{\vee})$;
\item $\tpsi_{l,ll'}$ induces an $S$-isomorphism 
$\psi_{l,ll'}:(P_{l},\cL_{l}^{\otimes l'})
\simeq (P_{ll'},\cL_{ll'}).$
\end{enumerate} 
\end{thm}
\begin{proof}This is proved in parallel to 
\cite[7.25]{MN24}. 
\end{proof}

The following describes the local structure of $\tP_{l}$, 
which is a corollary of 
Eq.~(\ref{eq:isom with contraction prod}) and \cite[8.12]{MN24}:
\begin{thm}\label{thm:W_l,alpha,0 general}
Let $\alpha\in\Sigma^{\dagger}_{l}$, $u\in X^{\vee}$ and $U\in\cF$. 
Then $B_{l,\alpha,u,U}\otimes_Rk(0)$ is a finitely generated 
$\cO_{A\cap U}$-algebra  
given by 
\begin{align*}
&\begin{cases}
\cO_{A_0\cap U}[s^{e(\xi)u(x)}w^x_U;x\in X]
&\text{if $\alpha\in\Sigma^{\dagger}_{l}(0)^0$,}\\
\cO_{A_0\cap U}
\begin{bmatrix}s^{e(\xi)(v_{\gamma}+u)(x)}w^x_U;\\ 
\gamma\in(\Sigma^{\dagger}_{l})^{\alpha}, 
x\in \Cone(C_{l}^{\gamma})\cap X
\end{bmatrix}
&\text{if $\alpha\in\partial\Sigma^{\dagger}_{l}(0)$,}
\end{cases}
\end{align*}
with 
fundamental relations  
in the second case 
\begin{align*}
&s^{e(\xi)(v_{\gamma'}+u)(x)}w^{x}_U\cdot 
s^{e(\xi)(v_{\gamma''}+u)(y)}w^{y}_U\\
&=\begin{cases}s^{e(\xi)(v_{\gamma}+u)(x+y)}w^{x+y}_U
&\text{if 
$x,y\in\Cone(C_{l}^{\gamma})\cap X$\ 
$(\exists \gamma\in\Sigma^{\dagger,\alpha}_{l}$)},\\
0&\text{otherwise}
\end{cases}
\end{align*}
where $\gamma',\gamma''\in\Sigma^{\dagger,\alpha}_{l}$, 
$x\in\Cone(C_{l}^{\gamma'})\cap X$ and $y\in\Cone(C_{l}^{\gamma''})\cap X$.  
\end{thm}

The following is a corollary 
of Eq.~(\ref{eq:isom with contraction prod}) and 
\cite[8.12~(1)]{MN24}:
\begin{lemma}\label{lemma:tG}
$\tG\simeq W_{l,\alpha,u}$\ 
$(\forall\alpha\in\Sigma^{\dagger}_{l}(0)^0)$. 
To be more precise, 
\begin{gather*}
\tG_U\simeq W_{l,\alpha,u,U}
=\Spec \Gamma(\cO_U)[\tau^e(u,x)w^x_U; x\in X],\\ 
\tG_{U^v}\simeq W_{l,\alpha,u+v,U^v}
=\Spec \Gamma(\cO_{U^v})[\tau^e(u+v,x)w^x_{U^v}; x\in X],
\end{gather*}where $v\in  X^{\vee}$ and $U\in\cF$. Moreover  
$\delta_{l,v}(W_{l,\alpha,u+v,U^v})=W_{l,\alpha,u,U}$. 
\end{lemma}

\begin{rem}\label{rem:same as thm:W_l,alpha,0}Since $v_s(\tau^e(u,x))=B^e(u,x)=e(\xi)u(x)$ $(u\in X^{\vee}, x\in X)$, Theorem~\ref{thm:W_l,alpha,0 general} is 
the same as \cite[8.12]{MN24} if $A_0=0$. Since $v_t(\tau^e(u,x))=u(x)$, we have in Lemma~\ref{lemma:tG}, 
\begin{gather*}
W_{l,\alpha,u,U}=\Spec \Gamma(\cO_U)[t^{u(x)}w^x_U; x\in X],\\ 
W_{l,\alpha,u+v,U^v}=\Spec \Gamma(\cO_U)[t^{(u+v)(x)}w^x_{U^v}; x\in X].
\end{gather*}
\end{rem}

\begin{cor}\label{item:Pl0 reduced CM} $\tP_{l,0}:=\tP_{l}\times_S0$ 
is reduced and Cohen-Macaulay.
\end{cor}

\begin{defn}\label{defn:Vorldagger}
Let $\overline{\Vor}^{\dagger}_{l}
:=\Vor^{\dagger}_{l}/4Nl\mu(X^{\vee})$. 
Let    
$\overline{\vor}^{\dagger}_{l,ll'}:\overline{\Vor}^{\dagger}_{l}\to \overline{\Vor}^{\dagger}_{ll'}$ be multiplication by $l'$. 
\end{defn}

\begin{defn}\label{defn:barZdaggerl_Delta}
Let   
$W_{l,\Delta}:=\bigcap_{a\in\Sk^0(\Delta)}W_{l,a}$ 
for $\Delta\in\Vor^{\dagger}_{l}$. 
See \cite[\S\S~8.2-8.3]{MN24} for more details. 
Let $Z^{\dagger}_{l}(\Delta)$ be the unique reduced closed $G_0$-orbit 
of $(W_{l,\Delta})_0$ and $\barZ^{\dagger}_{l}(\Delta)$ 
the closure of $Z^{\dagger}_{l}(\Delta)$ in 
$\tP_{l,0}$ (with reduced structure). 
\end{defn}

The following is proved in parallel 
to \cite[8.31]{MN24}:  
\begin{thm}
\label{thm:str of P_l general}  
Let $\pi_{l}:\tP_{l,0}\to P_{l,0}$ be the natural morphism and 
$\overline{\Orb}(P_{l,0})$ 
the set of closures in $P_{l,0}$ of $G_0$-orbits of $P_{l,0}$. 
We define $\omega_{l}$ by 
$\omega_{l}(\Delta)=\pi_{l}\barZ^{\dagger}_{l}(\Delta)$ for 
$\Delta\in\overline{\Vor}^{\dagger}_{l}$. Then 
\begin{enumerate}
\item\label{item:omega_l general} 
$\omega_{l}$  
is a bijective correspondence between $\overline{\Vor}^{\dagger}_{l}$ and 
$\overline{\Orb}(P_{l,0})$;
\item\label{item:indep l} $\omega_{l}$ is independent of $l$ in the sense that $\psi_{l,ll'}\circ\omega_{l}=\omega_{ll'}\circ\overline{\vor}^{\dagger}_{l,ll'}$ $(\forall l'\in\bN)$;
\item\label{item:Irred P_l_0} the set of irreducible components of $P_{l,0}$ 
consists of $\pi_{l}(\barZ^{\dagger}_{l}(\Sigma^{\dagger}_{l}(v)))$\ 
$(v\in X^{\vee}/\beta(Y))$ where $Z^{\dagger}_{l}(\Sigma^{\dagger}_{l}(v))=(W_{l,0,v})_0$.
\end{enumerate}
\end{thm}

\subsection{The structure of $P_{l}$}
\label{subsec:str of Pl in pd case}
  Throughout this section,  
we assume that $\Sigma^{\dagger}_{l}(0)$ is integral. The normalization $\tP_{l}$ is covered with open affine sets $W_{l,\alpha,u}$\ 
$(\alpha\in\Sigma^{\dagger}, u\in X^{\vee})$. 
Let $\tG_{l}:=W_{l,0,0}$. Since  
$\Proj \cR_{1,\Sigma^{\dagger}_{l}}\simeq
\Proj \cR_{2,\Sigma^{\dagger}_{l}}$, 
we have
$$\tG_{l}=W_{l,0,0}=\Spec \bigoplus_{x\in X}\cO_x,$$ 
which is independent of $l$.
Since $0\in\Sigma^{\dagger}_{l}\cap\Sigma^{\dagger}_{l}(0)^0$ 
by \cite[6.2]{MN24}, $\tG_{l}$ is described also by Lemma~\ref{lemma:tG}. 
Let $\delta_u:=\tT_{c^{t,e}(u)}$\ $(u\in X^{\vee})$ 
and 
$$\tC_{l}:=(\tP_{l}\setminus \bigcup_{u\in X^{\vee}}\delta_u(\tG_{l}))_{\red},\ \tD_{l}:=(\tP_{l}\setminus \bigcup_{y\in Y}S_y(\tG_{l}))_{\red}. 
$$ 

Let $(P_{l},\cL_{l})$ be the algebraization 
of the (formal) quotient $(\tP_{l}^{\wedge},\tcL_{l}^{\wedge})/Y$,    
$C_{l}$ (resp. $D_{l}$) the algebraization 
of the quotient $\tC_{l}^{\wedge}/Y$ (resp. 
$\tD_{l}^{\wedge}/Y$) and  
\begin{equation*} 
G_{l}
:=P_{l}\setminus D_{l},\quad 
\cG_{l}:=P_{l}\setminus C_{l}
=\bigcup_{u\in X^{\vee}/\beta(Y)}\delta_u(G_{l}),
\end{equation*}where   
$G_{l}$ is semiabelian with $G_{l}^{\wedge}\simeq\tG_{l}^{\wedge}$ 
by \cite[III, 5.7-5.8]{FC90}.   

\begin{thm}
\label{thm:pd summary of Pl Gl}
The following are true:
\begin{enumerate}
\item\label{item:Pl}$P_{l}$ is an irreducible flat 
projective normal Cohen-Macaulay 
$S$-scheme;
\item\label{item:open} 
$G_{l}$ and $\cG_{l}$ 
are open subschemes of $P_{l}$ and semiabelian group $S$-schemes
such that $P_{l,\eta}=G_{l,\eta}=\cG_{l,\eta}$ is an abelian variety  
and $G_{l}=\cG_{l}^0$;
\item\label{item:conn neron model}
$G_{l}$ is the connected N\'eron model of $G_{l,\eta}$;
\item\label{item:component gp acts on Pl}
 $(\cG_{l,0}/G_{l,0})(\overline{k(0)})\ 
(\simeq X^{\vee}/\beta(Y))$ acts on $P_{l}$;  
\item\label{item:Phi acts}
$(\cG_{l,0}/G_{l,0})(\overline{k(0)})$   
is naturally a subgroup 
$e(\xi)(X^{\vee}/\beta(Y))$ of 
the component group $X^{\vee}/e(\xi)\beta(Y)$ 
of the N\'eron model of $G_{l,\eta}$; 
\item\label{item:codim two}
 $\codim_{P_{l,0}} (P_{l,0}\setminus\cG_{l,0})=1$,  
and any irreducible component of $P_{l,0}$ 
is the closure in $P_{l,0}$ of an irreducible component of $\cG_{l,0}$;
\item\label{item:independence of l}
$P_{l}$, $G_{l}$ and $\cG_{l}$ are independent of the choice of $l$;
\item\label{item:pullback of neron over Rinit} 
$\cG_{l}\simeq\cG_S$ : the pullback to $S$ of the N\'eron model $\cG$ of 
$G_{\eta_{\init}}$; 
the quotient of $\cG_{l}$ 
by $\Aut(k(\eta)/k(\eta_{\init}))$ is the N\'eron model of 
$G_{\eta_{\init}}$;
\item\label{item:ample cubical cLl}
$\cL_{l|\cG_{l}}$ 
(resp. $\cL_{l|G_{l}}$) is a unique ample cubical invertible sheaf 
on $\cG_{l}$ (resp. $G_{l}$) with 
$(G_{l},\cL_{l|G_{l}})\simeq (G_S,\cL^{\otimes 4N^2l}_S)$. 
\end{enumerate}
\end{thm}
\begin{proof} See \cite[9.1]{MN24}.
 First we note the following. 
By Lemmas~\ref{lemma:H0(G,Lm) as pd Fourier}/\ref{lemma:theta on Pl general},
there exists a positive integer $m_0$ such that 
$\Gamma(P_{l,\eta},\cL_{l,\eta}^{\otimes m})=\Gamma(G_{\eta},\cL_{\eta}^{\otimes 4N^2lm})$ if $m\geq m_0$. Hence 
$\FC(G_{l},\cL_{l}^{\otimes m})=\FC(G_S,\cL_S^{\otimes 4N^2lm})$. 
Therefore $(P_{l,\eta},\cL_{l,\eta}^{\otimes m})\simeq 
(G_{\eta},\cL_{\eta}^{\otimes 4N^2lm})$, 
which is a polarized abelian variety. Moreover $(G_{l},\cL_{l}^{\otimes m})
=(G_S,\cL_S^{\otimes 4N^2lm})$ by Theorem~\ref{thm:G' isom G'' iff FC same}. 
By this fact, (\ref{item:Pl})-(\ref{item:ample cubical cLl})  
are proved verbatim in the same manner 
as \cite[9.1]{MN24}.   
\end{proof}

\begin{cor}\label{cor:pd neron if exi=1}
If $e(\xi)=1$, then $\cG_{l}$ is the N\'eron model of 
$G_{\eta}=G_{\eta_{\init}}$. 
\end{cor}
\begin{proof}
Clear from Theorem~\ref{thm:pd summary of Pl Gl}~
(\ref{item:pullback of neron over Rinit}).  
\end{proof}

\begin{cor}\label{cor:pd rel compact of G}
Let $\Gamma:=\Aut(k(\eta)/k(\eta_{\init}))$. 
The N\'eron model $\cG$ of $G_{\eta_{\init}}$ 
is the quotient $\cG_{l}/\Gamma$,
 which is relatively compactified by the quotient $P_{l}/\Gamma$.
\end{cor}
\begin{proof}
Clear from Theorem~\ref{thm:pd summary of Pl Gl}~(\ref{item:pullback of neron over Rinit}).  
\end{proof}

\begin{cor}\label{cor:pd m0=1}
We can choose $m_0=1$ in Lemma~\ref{lemma:theta on Pl general}.
\end{cor}
\begin{proof}This follows from 
Lemma~\ref{lemma:H0(G,Lm) as pd Fourier}/Theorem~\ref{thm:pd summary of Pl Gl}~(\ref{item:ample cubical cLl}).
\end{proof}

\subsection{The action of $G_{l}$ on $P_{l}$}
\label{sec:action of Gl on Pl}
See \cite[\S~10]{MN24} for a quick survey on rigid geometry. 
Let $\tP_{l,\eta}^{\an}$ be 
the rigid analytification of the $k(\eta)$-scheme 
$\tP_{l,\eta}$ and $\tP^{\wedge}_{l,\rig}$ the generic fiber 
(in the sense of rigid geometry \cite[10.5]{MN24}) 
of the formal $R$-scheme $\tP^{\wedge}_{l}$. 
Since $\tP_{l}$ satisfies the completeness condition 
\cite[7.6]{MN24}, we have
$\tP^{\wedge}_{l,\rig}\simeq\tP_{l,\eta}^{\an}$ 
by \cite[10.9~(3)-(4)]{MN24}.  
Similarly let $P_{l,\eta}^{\an}$ be 
the rigid analytification of the $k(\eta)$-scheme 
$P_{l,\eta}$ and $P^{\wedge}_{l,\rig}$ the generic fiber 
(in the sense of rigid geometry) 
of the formal $R$-scheme $P^{\wedge}_{l}$. Then 
$P^{\wedge}_{l,\rig}\simeq P_{l,\eta}^{\an}$ as rigid $k(\eta)$-spaces.
 Moreover $Y$ acts on $\tP_{l,\eta}^{\an}$ via $(S_y:y\in Y)$, so that 
\begin{equation}\label{eq:P_eta_an}
\tP^{\an}_{l,\eta}/Y\simeq P^{\an}_{l,\eta}\simeq G^{\an}_{l,\eta}.
\end{equation}

\begin{defn}\label{defn:abs value |cdot|}
Let $\Omega:=\overline{k(\eta)}$ be an algebraic closure of $k(\eta)$.
We define an 
{\it absolute value $|\cdot|$ of $k(\eta)$} 
by $|a|=e^{-v_s(a)}$ for $a\in k(\eta)$. Then  
$v_s$ and $|\cdot|$ can be uniquely extended to $\Omega$ 
\cite[A, Th.~3]{Bosch14}. Let $Q\in |\tP^{\wedge}_{l,\rig}|=\tP_{l}(\Omega)/G_{\abs}$ and $\tQ\in\tP_{l}(\Omega)$ 
 a representative of $Q$. Then 
$\log|w^{x}(Q)|:=\log|\tQ^*(w^x_U)|\in\bQ$ $(x\in\tX, U\in\cF^0_B)$ 
is well-defined in the sense that it is independent 
of the choice of $\tQ$ and $U$ because 
$w^x_U(Q)=\tilde\pi^*g^x_{UV}(\tQ)w^x_V(Q)$ and 
$\tilde\pi^*g^x_{UV}(\tQ)\in R^{\times}_{\Omega}$ 
where $\tilde\pi:\tP_{l}\to A$ is the natural projection. 
\end{defn}
 
\begin{thm}\label{thm:pd action of Gl on Pl}
Assume $\Sigma^{\dagger}_{l}(0)$ is integral. Then  
there exists an $S$-morphism $G_{l}\times_SP_{l}\to P_{l}$ 
which extends the group law $G_{l,\eta}\times_\eta G_{l,\eta}\to G_{l,\eta}$ of $G_{l,\eta}$.
\end{thm}

This is proved in parallel to that of \cite[10.14]{MN24} only by replacing 
\cite[10.12]{MN24} with Definition~\ref{defn:abs value |cdot|}.

\section{Proof of Theorem~\ref{thm:main thm}}
\label{sec:proof of main thm}
The purpose of this section is to prove 
Theorem~\ref{thm:main thm} in the {\it partially degenerate} case. 
Since we have proved Theorem~\ref{thm:pd summary of Pl Gl}, 
it suffices to descend Theorem~\ref{thm:pd summary of Pl Gl} to $S_{\init}$ 
 following \cite[\S\S~11-12]{MN24}. This is what we do in \S~\ref{sec:proof of main thm}.

\subsection{Plan of descent to $S_{\init}$}
\label{subsec:Plan}

First we present a precise plan of descent from $S_{\min}$ to $S_{\init}$. 
Let $R$ be a CDVR, $k(\eta)$ 
its fraction field, $S:=\Spec R$ and 
$\eta$ (resp. $0$) the generic (resp. closed) point of $S$. 
We basically follow the notation and 
the argument in \cite[\S\S~11-12]{MN24}.
Therefore $k(\eta)$ (resp. $R$, $S$, $0$) 
in Theorem~\ref{thm:main thm} is denoted here  
by $k(\eta_{\init})$ (resp. $R_{\init}$, $S_{\init}$, $0_{\init}$). 
Let $(G,\cL)$ be a semiabelian $S_{\init}$-scheme 
with $\cL$ symmetric ample cubical invertible and $\cG$ 
the N\'eron model of $G_{\eta_{\init}}$. 
We can choose such an $\cL$ by 
\cite[2.13~(1)]{MN24}.
By taking a finite \'etale Galois cover $S_{\spl}$ of $S_{\init}$ 
by \cite[X, \S~1]{SGA3}, we have a split  
semiabelian $S_{\spl}$-scheme $(G_{S_{\spl}},\cL_{S_{\spl}})$, that is, 
whose closed fiber $G_{0_{\spl}}$ is 
an extension of an abelian $k(0_{\spl})$-variety $A$ 
by a split $k(0_{\spl})$-torus $T$ where 
$R_{\spl}:=\Gamma(\cO_{S_{\spl}})$ is a CDVR and 
$0_{\spl}$ is the closed point of $S_{\spl}$. 
Let $k(\eta_{\spl})$ be the fraction field of $R_{\spl}$. 
By Definition~\ref{defn:FC ample datum}, we have the FC data 
$\zeta_m=\FC(G_{S_{\spl}},\cL_{S_{\spl}}^{\otimes m})$ 
 $(m\geq 1)$ and $\zeta=\zeta_1$ by 
Eq.~(\ref{eq:split obj zetam in DDample}) and Corollary~\ref{cor:psim,taum}:
\begin{equation*}
\zeta_{m}=(\tG,A,T,X,Y,c,c^t,\iota,\lambda^{(m)},
m\phi,\tau,\tcL^{\otimes m},
\psi^m,\cM^{\otimes m}).
\end{equation*}
Each $\zeta_m$ is a split object over $S_{\spl}$
in $\DD_{\ample}$ 
(Definitions~\ref{defn:split obj zeta in DDample}/\ref{defn:FC ample datum}). 

Then we choose a finite normal extension $k(\eta_{\AB})$ of 
the composite field 
$k(\eta_{\spl})\cdot k(\eta_{\init})$ 
so that $c^{t,e}:X\to A(k(\eta_{\AB}))$ is a homomorphism 
satisfying the condition (i) of 
Lemma~\ref{lemma:minimal Galois with extensions cte iotae psie taue} is 
true. See \S~\ref{subsec:pd Smin to Sspl} Step~3 for more detail. 
 Since $A$ is proper over $S_{\init}$, we have 
a homomorphism $c^{t,e}:X_{S_{\AB}}\to A$ such that (i) is true 
over $S_{\AB}$ where $S_{\AB}=\Spec R_{\AB}$ 
and $R_{\AB}$ is the integral closure 
of $R_{\spl}$ in $k(\eta_{\AB})$.  Then 
we choose a finite normal extension $k(\eta_{\min})$ of the composite field 
$k(\eta_{\AB})\cdot k(\eta_{\spl})\cdot k(\eta_{\init})$ such that 
Lemma~\ref{lemma:minimal Galois with extensions cte iotae psie taue}
~(ii)-(iv), so that we can find the eFC datum  
$\zeta^e_{4N^2l}$ (Definition~\ref{defn:2nl th eFC}): 
$$
\zeta^e_{4N^2l}=(\tG,A,T,X,Y,c,c^{t,e},\iota^e,
4N^2l\lambda,4N^2l\phi,\tau^e,\tcL^{\otimes 4N^2l},
(\psi^e)^{2Nl},\cM^{\otimes 4N^2l})$$  
by enlarging $k(\eta_{\AB})$ 
to a finite radical normal extension $k(\eta^e)$ 
of $k(\eta_{\init})$, where we choose the minimal such $k(\eta^e)$. 
Let $R^e$ be the integral closure of $R$ in $k(\eta^e)$ and $S^e:=\Spec R^e$. 
By Theorem~\ref{thm:degeneration data Neron 2l case} and 
Definition~\ref{defn:NeFC kit}, we have  
the degeneration data of the pullback $\cG_{S^e}$ of $\cG$ 
\begin{equation*}
\NeFC(G_{S^e},\cL^{\otimes 4N^2l}_{S^e})
=(\zeta^e_{4N^2l},\xi^{\dagger}_{l}).
\end{equation*}
 Note that the $\NeFC$ kit 
$\xi^{\dagger}_{l}$ plays the same role 
as $\xi^{\natural}_{l}$ in \cite[5.21]{MN24}.

The rest is the same as \cite[\S~11.1]{MN24}. 
By Theorem~\ref{thm:degeneration data Neron 2l case} 
and Definition~\ref{defn:NeFC kit}, 
let $\eta_{\min}:=\eta_{\min}(\zeta)=\eta^e$. 
By Remark~\ref{rem:zeta and zetam}, 
$k(\eta_{\min}(\zeta_m))=k(\eta_{\min}(\zeta))$ 
$(\forall m\in\bN)$.  Hence we have a sequence of extensions:
\begin{equation}\label{eq:field extensions}
k(\eta_{\min})\supset k(\eta_{\AB})\supset k(\eta_{\spl})
\supset k(\eta_{\init}).
\end{equation}

Theorem~\ref{thm:pd summary of Pl Gl} constructs a quadruple 
$(P_{l},\cL_{l},G_{l},\cG_{l})$ over $R^e$. 
To prove Theorem~\ref{thm:main thm}, 
it suffices to descend the quadruple to $R_{\init}$ as follows:
\begin{enumerate}
\item[]Case 1.\ descent from $S_{\min}=S^e$ to $S_{\spl}$\ 
(\S~\ref{subsec:pd Smin to Sspl});
\item[]Case 2.\ descent from $S_{\spl}$ to $S_{\init}$\ 
(\S~\ref{subsec:pd Sspl to Sinit}).
\end{enumerate}

\subsection{Descent from $S_{\min}$ to $S_{\spl}$}
\label{subsec:pd Smin to Sspl}
We basically follow \cite[\S~11.4]{MN24}. 
We set $k(\eta)=K_{\min}(\zeta)$, 
$R=R_{\min}(\zeta)$ and $S=S_{\min}(\zeta)=S^e$. 
We discuss descent from $S_{\min}$ to $S_{\ab}$ 
(resp. $S_{\ab}$ to $S_{\spl}$) 
in Steps 1-2 (resp. in Step 3). Let $s$ (resp. $t$, $t_{\spl}$) be a uniformizing parameter of $R_{\min}$ (resp. $R_{\ab}$, $R_{\spl}$), 
$e_{\ma}:=v_s(t)$ and $e_{\as}:=v_t(t_{\spl})$. 
Hence $e(\zeta):=v_s(t_{\spl})=e_{\ma}e_{\as}$. 
\par
\m 
\noindent{\it Step 1.}\ \ 
Let $(G,\cL)$ be a semiabelian $S_{\AB}$-scheme with $\cL$ 
symmetric ample cubical invertible and $\cG$ 
the N\'eron model of $G_{\eta_{\AB}}$ with $G=\cG^0$.   
By $\zeta=\zeta_1=\FC(G,\cL)$, we are given 
a split $S_{\AB}$-torus $T$ and 
a $T$-bundle $\pi:\tG\to A$ over an abelian $S_{\AB}$-scheme $A$ with 
$\tG^{\wedge}\simeq G^{\wedge}$. 
By Eq.~(\ref{eq:isom of Lwedge and sum of Mxwedge}) 
\begin{equation}\label{eq:pi_*cL_wedge}
(\pi^{\wedge})_*\cL^{\otimes m,\wedge}
=\left(\bigoplus_{x\in X}(\cM^{\otimes m})_x\right)^{I\op{-adic}}.
\end{equation} 
 \par
By Lemma~\ref{lemma:minimal Galois with extensions cte iotae psie taue}, 
there exist an N\'eFC pair $(\zeta^e_{4N^2l},\xi^{\dagger}_{l})$ 
over $S=S_{\min}=S^e$. 
Let 
$\tilde\Gamma:=\Aut(k(\eta)/k(\eta_{\AB}))
=\Aut(S/S_{\AB})$. If $\chara k(\eta)=0$, then 
$\tilde\Gamma=\Gal(S/S_{\AB})$. 
If $\chara k(\eta)=p>0$, then $\tilde\Gamma$ 
is a finite group $S_{\AB}$-scheme. 
The identity 
component $\mu$ of $\tilde\Gamma$ is a normal subgroup $S_{\ab}$-scheme 
of $\tilde\Gamma$ such that 
$\tilde\Gamma/\mu\simeq\tilde\Gamma_{\red}$ : 
 the reduced part of $\tilde\Gamma$, which is 
a subgroup scheme of $\tilde\Gamma$. 
By Lemma~\ref{lemma:semi-direct Gal}/Corollary~\ref{cor:B^0 when reducible} 
and the proof of 
Lemma~\ref{lemma:minimal Galois with extensions cte iotae psie taue}, 
$\mu\simeq\prod_{i=1}^g\mu_{q_i,S}$ where 
$q_i$ is some power of $p$ dividing $e_i$ $(i\in[1,r])$, 
and  there exists a finite group $\Gamma$ 
with $\tilde\Gamma_{\red}\simeq\Gamma_S$. 
 Regardless of $\chara k(\eta)$, 
$S/\tilde\Gamma\simeq S_{\AB}$ in any case. 

By Definition~\ref{defn:NeFC kit}, 
the $l$-th N\'eFC kit $\xi^{\dagger}_{l}$ of $\zeta^e_{4N^2l}$ 
is given by:
$$\xi^{\dagger}_{l}=(X,Y,\epsilon^{\dagger}_{l},\tau^e_{l},E^{\dagger}_{l},\Sigma^{\dagger}_{l}),\ \  \xi^{\dagger}:=\xi^{\dagger}_1.
$$

By \cite[6.5]{MN24}, 
there exists $l_0\in\bN$ such that 
\begin{equation}
\label{eq:integral Sigmal(0)}
\text{$\Sigma^{\dagger}_{l_0l}(0)$ is integral for any integer $l\in\bN$.}
\end{equation}

By Theorems~\ref{thm:pd summary of Pl Gl}/\ref{thm:pd action of Gl on Pl}, 
we obtain
\begin{claim}\label{claim:pd Pl0l Gl0l cGl0l GS cGS}
Let $S=S_{\min}(\xi)$ and 
$(P_{l_0l},\cL_{l_0l}):
=(P_{l_0l}(\xi^{\dagger}),\cL_{l_0l}(\xi^{\dagger}))$.
Then 
\begin{enumerate}
\item[(i')]\label{item:Cohen-Macaulay} 
$(P_{l_0l},\cL_{l_0l})$ is an irreducible flat 
projective Cohen-Macaulay $S$-scheme;
\item[(ii')]\label{item:Xvee/betaY}there exist 
open $S$-subschemes $G_{l_0l}$ and $\cG_{l_0l}$ of $P_{l_0l}$ 
such that $(G_{l_0l},\cL_{l_0l})\simeq (G_S,\cL^{4N^2l_0l}_S)$, 
$\cG_{l_0l}\simeq\cG_S$ and  
 $(\cG_{l_0l,0}/G_{l_0l,0})(\overline{k(0)})\simeq X^{\vee}/\beta(Y)$;
\item[(iii')]\label{item:codim 2} $P_{l_0l}\setminus \cG_{l_0l}$ is 
of codimension two in $P_{l_0l}$;
\item[(iv')]\label{item:action of cGl on Pl}
 $\cG_{l_0l}$ acts on $P_{l_0l}$ extending 
multiplication of  $\cG_{l_0l}$;
\item[(v')]\label{item:cubical sheaves cLl}
 $(\cL_{l_0l})_{|\cG_{l_0l}}$ is 
ample cubical with $\cL_{l_0l,\eta}=\cL^{4N^2l_0l}_{\eta}$.
\end{enumerate}
\end{claim}

{\it Step 2.}\ \  We shall descend 
Claim~\ref{claim:pd Pl0l Gl0l cGl0l GS cGS} from $S$ to $S_{\ab}$. 
Let $s$ (resp. $t$) 
be a uniformizer of $R$ (resp. $R_{\ab}$), $e_{\ma}:=v_s(t)$ 
and $I:=tR_{\ab}$. 
Since $A$ is an $S_{\ab}$-scheme, $\tilde\Gamma$ 
acts on $\cO_A$, $\cO_x$ $(x\in X)$ and 
$\cM$ trivially. By \cite[11.8/11.9]{MN24}, $\tilde\Gamma$ 
acts on $M_{m,x}:=R\otimes_{R_{\ab}}(\cM^{\otimes m})_x$ by 
\begin{gather*}
\rho_{M_{m,x}}(cf_x)=\rho_R(c)\cdot (1\otimes f_x)
\in\Gamma(\cO_{\bf\mu})\otimes_{R_{\ab}} (\cM^{\otimes m})_x,\\ 
\sigma(cf_x)=\sigma(c)f_x\ \  
(\forall \sigma\in\Gamma, \forall c\in R, 
\forall f_x\in(\cM^{\otimes m})_x, \forall x\in X)
\end{gather*}where $\rho_{M_{m,x}}$ (resp. $\rho_R$) is the coaction 
of $\mu$ on $M_{m,x}$ (resp. $R$).

Let $V_{m,\ab}:=\Gamma(G_{\eta_{\ab}},\cL_{\eta_{\ab}}^{\otimes m})$ 
and $V_{m}:=\Gamma(G_{\eta},\cL_{\eta}^{\otimes m})=V_{m,\ab}\otimes_{R_{\ab}}R$\ $(m\geq 1)$. Then $\tilde\Gamma$ acts 
on $V_{m}$, so that $V_{m}^{\tilde\Gamma}=V_{m,\ab}$ 
by Lemma~\ref{lemma:H0(G,Lm) as pd Fourier}.
Next we define the (co)action $\rho_{\cS_{l_0l,2}}$ of $\mu$ 
(resp. $\psi_{\bullet}$ of $\Gamma$) on  
$\cS_{l_0l,2}=\cS_{l_0l,2}(\xi^{\dagger})$  by 
\begin{gather*}
\rho_{\cS_{l_0l,2}}(cm_x\theta_{l_0l})
=\rho_{M_{4N^2l_0l}}(cm_x)(1\otimes\theta_{l_0l}),\\
\psi_{\sigma}(cm_x\theta_{l_0l})=\sigma(c)m_x\theta_{l_0l}\\
(\forall\sigma\in\Gamma, \forall c\in R, \forall x\in X, \forall m_x\in(\cM^{\otimes 4N^2l_0l})_x).
\end{gather*}
 
Now we recall 
\begin{align*}
\cR_{2,\Sigma^{\dagger}_{l_0l}}(\xi^{\dagger})
&=\cO_{A_S}\left[
\begin{matrix}\left(J_{l_0l,v}\cdot J_{v,\beta}
\cdot(\cM^{\otimes 4N^2l_0l}\otimes_{\cO_A}
\cO_{\beta+4Nl_0l\mu(v)})\right)\theta_{l_0l}\\
v\in X^{\vee},\beta\in\Sigma^{\dagger}_{l_0l}
\end{matrix}\right]
\end{align*}where $J_{l_0l,v}=s^{e_{\ma}E^{\dagger}_{l_0l}(v)}R$ 
and $J_{v,\beta}=s^{e_{\ma}v(\beta)}R$.
Then $\cR_{2,\Sigma^{\dagger}_{l_0l}}(\xi^{\dagger})$ 
is an $R$-subalgebra of $\cS_{l_0l,2}$ 
stable under $\tilde\Gamma$ 
because $v_s\circ\psi_{\sigma}=v_s$ on $R$. 

Let $J_{l_0l,v,\ab}:=t^{E^{\dagger}_{l_0l}(v)}R_{\ab}$, 
$J_{v,\beta,\ab}:=t^{v(\beta)}R_{\ab}$ and 
\begin{align*}\label{eq:R form of cR2_l_xi_dagger}
\cS_{l_0l,2,\ab}&:=\left(\bigoplus_{n\geq 0}\left(\bigoplus_{x\in X}\cM^{\otimes 4N^2l_0l n}\otimes_{\cO_A}\cO_{x}\right)\theta_{l_0l}^n\right)\otimes_{R_{\ab}}k(\eta_{\ab}),\\
\cR_{2,\Sigma^{\dagger}_{l_0l},\ab}(\xi^{\dagger})
&:=\cO_A\left[
\begin{matrix}\left(J_{l_0l,v,\ab}\cdot J_{v,\beta,\ab}
\cdot(\cM^{\otimes 4N^2l_0l}\otimes_{\cO_A}
\cO_{\beta+4Nl_0l\mu(v)})\right)\theta_{l_0l}\\
v\in X^{\vee},\beta\in\Sigma^{\dagger}_{l_0l}
\end{matrix}\right].
\end{align*}  
Then $\cR_{2,\Sigma^{\dagger}_{l_0l},\ab}(\xi^{\dagger})
$ is an $R_{\ab}$-subalgebra of $\cS_{l_0l,2,\ab}$ such that 
\begin{gather*}
\cR_{2,\Sigma^{\dagger}_{l_0l}}(\xi^{\dagger})
=\cR_{2,\Sigma^{\dagger}_{l_0l},\ab}(\xi^{\dagger})\otimes_{R_{\ab}}R,\\  
\left(\cR_{2,\Sigma^{\dagger}_{l_0l}}(\xi^{\dagger})\right)^{\tilde\Gamma}=\cR_{2,\Sigma^{\dagger}_{l_0l},\ab}(\xi^{\dagger}).
\end{gather*}
 
We see that $S_y^*$ acts 
on $\cR_{2,\Sigma^{\dagger}_{l_0l},\ab}(\xi^{\dagger})$ 
by Definition~\ref{defn:action tildeSy}, while 
$\delta_u^*$ does not keep $\cR_{2,\Sigma^{\dagger}_{l_0l},\ab}(\xi^{\dagger})$ stable in general because $\tau^e(u,x)\not\in k(\eta_{\ab})$ can happen. 
Let $\tQ^*_{l_0l}:=\Proj \cR_{2,\Sigma^{\dagger}_{l_0l},\ab}(\xi^{\dagger})$  
and 
$(\tP^*_{l_0l},\tcL^*_{l_0l})$ be 
the normalization of $\tQ^*_{l_0l}$ with  
$\tcL^*_{l_0l}$ the pullback of $\cO_{\tQ^*_{l_0l}}(1)$ 
to  $\tP^*_{l_0l}$.  Let $(P^*_{l_0l},\cL^*_{l_0l})$ 
be the algebraization 
of the (formal) quotient 
$(\tP_{l_0l}^{*\wedge},\tcL_{l_0l}^{*\wedge})/Y$. 
By our construction we have a natural 
finite surjective $S_{\ab}$-morphism  
$f:(P_{l_0l},\cL_{l_0l})\to (P^*_{l_0l},\cL^*_{l_0l})$. 
Let $C^*_{l_0l}=f(C_{l_0l})_{\red}$,  
$D^*_{l_0l}=f(D_{l_0l})_{\red}$, 
$\cG^*_{l_0l}=P^*_{l_0l}\setminus C^*_{l_0l}$ and  
$G^*_{l_0l}=P^*_{l_0l}\setminus D^*_{l_0l}$.

\begin{claim}\label{claim:pd P*l_ab}
The following are true over $S_{\ab}$:
\begin{enumerate}
\item[(i$^*$)] $(P^*_{l_0l},\cL^*_{l_0l})$ is an irreducible flat 
projective Cohen-Macaulay $S_{\ab}$-scheme;
\item[(ii$^*$)]  
$(G^*_{l_0l},\cL^*_{l_0l})\simeq 
(G,\cL^{\otimes 4N^2l_0l})$ and $\cG^*_{l_0l}\simeq\cG$; 
\item[(iii$^*$)] $P^*_{l_0l}\setminus\cG^*_{l_0l}$ is of codimension two;
\item[(iv$^*$)] $\cG^*_{l_0l}$ acts on $P^*_{l_0l}$ 
extending multiplication of $\cG^*_{l_0l}$;
\item[(v$^*$)]\label{item:cubical sheaves cL*l}
 $(\cL^*_{l_0l})_{|\cG^*_{l_0l}}$ is 
ample cubical with $\cL^*_{l_0l,\eta_{\ab}}=\cL^{4N^2l_0l}_{\eta_{\ab}}$.
\end{enumerate}
\end{claim}
\begin{proof}
This is proved in the same manner as \cite[11.14]{MN24} here by using 
Claim~\ref{claim:pd Pl0l Gl0l cGl0l GS cGS}, Theorem~\ref{thm:pd summary of Pl Gl}~(\ref{item:conn neron model})/(\ref{item:codim two})/(\ref{item:pullback of neron over Rinit})/(\ref{item:ample cubical cLl}) and Theorem~\ref{thm:pd action of Gl on Pl}.  
\end{proof}

{\it Step 3.}\ \  Next we descend Claim~\ref{claim:pd P*l}
from $S_{\ab}$ to $S_{\spl}$. We set $\tilde\Gamma:=\Aut(S_{\ab}/S_{\spl})$. 
Since  $S_{\ab}/\tilde\Gamma\simeq S_{\spl}$, 
we can descend $(P^*_{l_0l},\cL^*_{l_0l},\cG^*_{l_0l},G^*_{l_0l})$ 
from $S_{\ab}$ to $S_{\spl}$ 
in the same manner as Step~2 because
\begin{enumerate}
\item[(a)] $(G,\cL)$, $\cM$, $A$, $T$ and 
$\cO_x$ $(x\in X)$ are defined over $S_{\spl}$; 
\item[(b)] by (a), we can 
define the $R_{\spl}$-algebras $\cS_{l_0l,2,\spl}$ and 
$\cR_{\Sigma^{\dagger}_{l_0l,\spl}}$ 
in the same manner as Step~2;
\item[(c)] $S_y^*$ can be defined over $S_{\spl}$ 
without using $\alpha_i$ $(i\in[1,r])$.
\end{enumerate} 

Although the identity component $\tilde\Gamma^0$ of $\tilde\Gamma$ can be 
a nonreduced additive group $S_{\spl}$-scheme if $\chara k(\eta)=p>0$, 
this is harmless for descent. Thus 
we can descend Claim~\ref{claim:pd P*l_ab} to its $S_{\spl}$-version:
\begin{claim}\label{claim:pd P*l}
$(P^*_{l_0l},\cL^*_{l_0l})$, $\cG^*_{l_0l}$ and 
$G^*_{l_0l}$ descend from $S_{\ab}$ 
to $S_{\spl}$ and satisfy the conditions 
Claim~\ref{claim:pd P*l_ab}~(i$^*$)-(v$^*$) over $S_{\spl}$.
\end{claim}

This completes the proof of Theorem~\ref{thm:main thm} in Case 1. 
In what follows we denote the descents to $S_{\spl}$ of 
$(P^*_{l_0l},\cL^*_{l_0l})$, $\cG^*_{l_0l}$ and 
$G^*_{l_0l}$ by the same letters if no confusion is possible.

\subsection{Descent from $S_{\spl}$ to $S_{\init}$}
\label{subsec:pd Sspl to Sinit}
The descent of $(P^*_{l_0l},\cL^*_{l_0l},\cG^*_{l_0l},G^*_{l_0l})$ 
 from $S_{\spl}$ to $S_{\init}$ is 
a standard finite Galois descent, which is carried out 
in the same manner as \cite[\S~12]{MN24}. 
To do so, we prove Lemmas~\ref{lemma:nomal form of Gamma action}\,-\,\ref{lemma:a b lambda_sigma}, which substitute \cite[12.2-12.3]{MN24}. \par
We have started in \S~\ref{subsec:Plan} with a semiabelian $S$-scheme $(G,\cL)$
such that $G_{0_{\init}}$ is a torus $T_{0_{\init}}$-bundle over an abelian 
$k(0_{\init})$-variety. Then there exists a finite \'etale Galois cover $\varpi:S_{\spl}\to S_{\init}$ of $S_{\init}$ such that $T_{0_{\spl}}$ is a split torus. 
For $Z\in\{S,R,\eta,0\}$, we denote $Z_{\spl}$ and $Z_{\init}$ 
by $Z^*$ and $Z$ respectively,   
{\it e.g.}, $S^*=S_{\spl}$ and $S=S_{\init}$. 
Let $\Gamma=\Gal(S^*/S)$ and $s$ a uniformizer of $R=R_{\init}$. 
By Hensel's lemma, 
\begin{equation}\label{eq:Hensel}
\Gal(k(0^*)/k(0))=\Gal(S^*/S)=\Gal(k(\eta^*)/k(\eta)). 
\end{equation} 

\begin{defn}\label{defn:pd phi_sigma f_sigma} 
Let $\sigma\in\Gamma$. 
We define an $R$-homomorphism $s_{\sigma}^*:R^*\to R^*$ by  
$s_{\sigma}^*(a)=\sigma(a)$ $(\forall a\in R^*)$, 
and let $s_{\sigma}:S^*\to S^*$  the $S$-morphism 
induced from $s_{\sigma}^*$. 
Since $\sigma$ acts on $S^*$, so does $\sigma$ on $G_{S^*}$ and the $I^*$-adic completion $G_{S^*}^{\wedge}$ of $G_{S^*}$ where
$I^*=IR^*$.  We define 
$\psi_{\sigma}\in\Aut_R(\Gamma(\cO_{G_{S^*}^{\wedge}}))$ 
by $\psi_{\sigma}(z\otimes a)=z\otimes\sigma(a)$\ 
$(\forall a\in R^*,\forall z\in\Gamma(\cO_{G^{\wedge}}))$ via the 
identification $\Gamma(\cO_{G_{S^*}^{\wedge}})
=\Gamma(\cO_{G^{\wedge}})\otimes_RR^*$ \cite[11.1]{MN24}. 
Then $\psi_{\sigma\tau}=\psi_{\sigma}\circ\psi_{\tau}$ $(\forall \sigma,\tau\in\Gamma)$. Thus  $(\psi_{\sigma};\sigma\in\Gamma)$ is a $\Gamma$-structure of 
$\Gamma(\cO_{G_{S^*}^{\wedge}})$ and $\Gamma(\cO_{G_{S^*}^{\wedge}})^{\Gamma}
=\Gamma(\cO_{G^{\wedge}})$ (\cite[1.2]{MN24}). 
Let $\cL_{S^*}$ be the pullback of $\cL$ to $G_{S^*}$. Then  
 $\Gamma$ 
acts on $\Gamma(G_{S^*},\cL_{S^*}^{\otimes m})$ $(m\geq 1)$, where 
$\Gamma$ acts on  $\Gamma(G,\cL^{\otimes m})$ trivially 
because $(G,\cL)$ is defined over $S$. By \cite[11.3]{MN24},  
\begin{gather*}
\Gamma(G_{S^*},\cL_{S^*}^{\otimes m})=\Gamma(G,\cL^{\otimes m})\otimes_RR^*,
\ \ 
\Gamma(G_{S^*},\cL_{S^*}^{\otimes m})^{\Gamma}=\Gamma(G,\cL^{\otimes m}).
\end{gather*} 
\end{defn}

\begin{lemma}
\label{lemma:nomal form of Gamma action}
The following are true: 
\begin{enumerate}
\item\label{item:Hom(Gamma,Aut)} 
there exists $\lambda\in\Hom(\Gamma,\Aut_{\bZ}(X))$ 
such that $\lambda_{\sigma}(Y)=Y$ and 
$\psi_{\sigma}(w^x)=w^{\lambda_{\sigma}(x)}$\ 
$(\forall x\in X, \forall \sigma\in\Gamma)$  
where $\lambda_{\sigma}:=\lambda(\sigma)$;
\item\label{item:Gamma effective}
$\Gamma$ acts effectively on $X$ via $(\lambda_{\sigma};\sigma\in\Gamma)$.   
\end{enumerate}
\end{lemma}
\begin{proof}
Let $0\to T\to \tG\overset{\pi}{\to} A\to 0$ be the Raynaud extension 
Eq.~(\ref{eq:exact seq of tG}) associated to $G$. 
There exists an open affine covering $\cF^0$ 
of $A$ such that $\tG_U\simeq T\times_S U$\ $(\forall U\in\cF^0)$. 
Since $T_{0^*}$ is a split $k(0^*)$-torus, 
by \S~\ref{subsec:Raynaud extensions split case},
$\tG_{U_{S^*}}\simeq\Spec \Gamma(\cO_U)[w^x_U;x\in X]\otimes_RR^*$. 
Let  $\cB_U:=\Gamma(\cO_{G_U})^{\wedge}$ and 
$\cB^*_U:=\Gamma(\cO_{G_{U_{S^*}}})^{\wedge}=\cB_U\otimes_RR^*$. 
By the proof of \cite[12.2]{MN24}, there exists 
$\lambda_{\sigma}\in\Aut_{\bZ}(X)$ $(\forall\sigma\in\Gamma)$ such that 
$\psi_{\sigma}(w^x)=w^{\lambda_{\sigma}(x)}$ $(\forall x\in X)$ and 
$\lambda_{\sigma\tau}=\lambda_{\sigma}\lambda_{\tau}$\ 
$(\forall\sigma,\tau\in\Gamma)$.  Hence 
$\lambda\in\Hom(\Gamma,\Aut_{\bZ}(X))$. 
Since the natural morphism $G\times_US^*\to G^t\times_US^*$ is $\Gamma$-equivariant,  $\Gamma(\cO_{G^t_{U}})$ is a $\Gamma$-submodule of 
$\Gamma(\cO_{G_{U_{S^*}}})$, so that 
$\lambda_{\sigma}(Y)=Y$.   
This proves (\ref{item:Hom(Gamma,Aut)}).  
It is obvious that $\lambda_{\sigma}$ 
is independent of the choice of $U\in\cF^0$. 
(\ref{item:Gamma effective}) is proved in the same manner as 
\cite[12.2~(2)]{MN24}.
\end{proof}

\begin{lemma}
\label{lemma:a b lambda_sigma}Under the notation 
in Lemma~\ref{lemma:nomal form of Gamma action}, we have  
\begin{gather*}
\sigma(\psi(y))=\psi(\lambda_{\sigma}(y)),\ \sigma(\tau(y,x))=\tau(\lambda_{\sigma}(y),\lambda_{\sigma}(x))\\ 
(\forall y\in Y, \forall x\in X, \forall\sigma\in\Gamma).
\end{gather*}
\end{lemma}
\begin{proof}We follow the proof of \cite[12.3]{MN24}. 
 By using Lemma~\ref{lemma:H0(G,Lm) as pd Fourier} with $m=1$, 
we obtain $\psi(\lambda_{\sigma}(y))\tau(\lambda_{\sigma}(y),\lambda_{\sigma}(x))=\psi(y)\tau(y,x)$\ $(\forall y\in Y, \forall x\in X)$, 
from which Lemma follows.
\end{proof}

Since $G_{0^*}$ is a split $k(0^*)$-torus, we obtain 
a split object $\zeta$ and  
N\'eFC pairs $(\zeta^2_{4N^2l},\xi^{\dagger}_{l})$. 
Let $\xi^{\dagger}=\xi^{\dagger}_1$. By \S~\ref{subsec:pd Smin to Sspl}, 
by choosing $l_0\in\bN$ such that 
$\Sigma_{l_0}(0)$ is integral, 
$(P^*_{l_0l},\cL^*_{l_0l})$ 
\ $(l\in\bN)$ is 
the twisted Mumford family obtained from 
\begin{align*}
\cR^{\dagger}_{l_0l}(\xi^{\dagger})&:=
\cR_{2,\Sigma^{\dagger}_{l_0l},\spl}(\xi^{\dagger}). 
\end{align*} 

Let $(U;U\in\cF^0)$ be an affine covering of $A$.  
Recall that 
$D^{\dagger}_{l_0l}(x)$ is well-defined for
$x\in X$ by Definition~\ref{defn:ep_dagger,E_dagger,D_dagger}. 
For $U\in\cF^0$, we define 
\begin{align*}
\cR^{\dagger}_{l_0l, U}(\xi^{\dagger})&
=\Gamma(\cO_{U_{S^*}})
\left[
\begin{matrix}
\zeta_{l_0l,x, U}\theta_{l_0l,U}; x\in X
\end{matrix}
\right],
\ \ \zeta_{l_0l,x,U}=t^{D^{\dagger}_{l_0l}(x)}w_U^{x}.
\end{align*}
Then $\cR^{\dagger}_{l_0l, U}(\xi^{\dagger})
=\cR^{\dagger}_{l_0l}(\xi^{\dagger})\otimes_{R}\Gamma(\cO_U)$, which  
will play the same role as $R^*_{l_0l}(\xi^{\natural})$ 
in \cite[12.3]{MN24}. 
 We can descend $(P^*_{l_0l},\cL^*_{l_0l},\cG^*_{l_0l},G^*_{l_0l})$ from $S^*$ to $S$ in the same manner as \cite[12.5-12.9]{MN24}.
Now we define 
\begin{equation}\label{eq:descents}
\begin{aligned}
P^{\flat}_{l_0l}:=P^*_{l_0l}/\Gamma,&\ \ \cL^{\flat}_{l_0l}
:=\text{the descent of $\cL^*_{l_0l}$ to 
$P^{\flat}_{l_0l}$},\\   
G^{\flat}_{l_0l}&:=G^*_{l_0l}/\Gamma,\ 
\cG^{\flat}_{l_0l}:=\cG^*_{l_0l}/\Gamma.
\end{aligned}
\end{equation}

Thus we obtain
\begin{claim}\label{claim:action T val points}
The group $S^*$-scheme $G^*$ (resp. $\cG^*)$ descends to a group $S$-scheme 
$G^{\flat}$ (resp. $\cG^{\flat})$, and the following are true:
\begin{enumerate}
\item[(i$^{\flat}$)] $P^{\flat}$ is an irreducible 
flat projective Cohen-Macaulay $S$-scheme;
\item[(ii$^{\flat}$)] $G^{\flat}$ (resp. $\cG^{\flat}$) 
is a group $S$-scheme 
such that $G^{\flat}\simeq G$ (resp. $\cG^{\flat}\simeq\cG$);
\item[(iii$^{\flat}$)] $P^{\flat}\setminus\cG^{\flat}$ is of codimension two;
\item[(iv$^{\flat}$)] there exists a morphism $G^{\flat}\times P^{\flat}\to P^{\flat}$ 
(resp. $\cG^{\flat}\times P^{\flat}\to P^{\flat}$) extending 
multiplication of $G^{\flat}$ (resp. $\cG^{\flat})$;
\item[(v$^{\flat}$)] $(\cL^{\flat}_{l_0l})_{|\cG^{\flat}}$ is ample cubical with $\cL^{\flat}_{l_0l,\eta}=\cL^{4N^2l_0l}_{\eta}$.
\end{enumerate}
\end{claim}

This completes the proof of Theorem~\ref{thm:main thm} in Case 2, and hence in 
any case.

\section{The totally degenerate case revisited}
\label{sec:totally deg case revisited}
In this section, we revisit the totally degenerate case \cite[\S~8]{MN24}, 
where $G$ is a semiabelian $S$-scheme such that $G_0=T_0$ is 
a split $k(0)$-torus and $A_0=\{0\}$.
We use the notation of \cite[\S~8]{MN24} freely. 
We assume that $\Sigma_{l}(0)$ is integral. 
We prove a projective embedding theorem of toric closed subschemes of 
the closed fiber $P_{l,0}$  \cite[\S~9]{MN24}. 

\subsection{The Voronoi polytopes $\Sigma_{l}(v)$} 
\label{subsec:piecewise linear}

 Let $\xi:=\FC(G,\cL)=(X,Y,a,b,A,B)$, 
$\xi^{\natural}$ an $\NeFC$ kit over $S$ 
extending $\xi^e_{2N}$, and  
$\xi^{\natural}_{l}$\ $(l\in\bN)$ the $l$-th $\NeFC$ kit over $S$
induced from $\xi^{\natural}$ in \cite[5.21]{MN24}  
with $\xi^{\natural}_1=\xi^{\natural}$:
$$\xi^{\natural}_{l}=(X,Y,\epsilon_{l},b^e,E_{l},\Sigma_{l}).
$$

To $\xi^{\natural}$, we associate a fan $\Fan(\xi^{\natural}_{l})$\ 
\cite[5.18/5.21]{MN24} and a decomposition 
$\Vor_{l}:=\Vor_{l}(\xi^{\natural})$ of 
$X_{\bR}$ \cite[6.1/6.2]{MN24}: 
\begin{gather*}
\Vor_{l}=\{\text{$\Sigma_{l}(v)$ $(v\in X^{\vee})$ 
and their faces}\},\\
\Sigma_{l}(v)=\{x\in X_{\bR};  E_{l}(u)+u(x-2l\mu(v))\geq 0
\ (\forall u\in X^{\vee})\}
\end{gather*}where $E_{l}(u)=lu(\mu(u))$. 
For $\alpha\in\Sigma_{l}(0)$, let  
$\Star_{l}(\alpha):=\bigcup_{\alpha\in\Sigma_{l}(v)}\Sigma_{l}(v)$.

\begin{lemma}\label{lemma:torus embedding Zl_Sigma(0)}
Let $\alpha\in\Sk^0(\Sigma_{l}(0))$. 
Then
\begin{enumerate}
\item\label{item:tau_alpha}
 $\Cone(\Cut(\tau_{l,\alpha,0}))
=\{u\in X^{\vee}_{\bR};
u(\beta)\geq u(\alpha)\ 
(\forall \beta\in\Sigma_{l})\}=\Cone(C_{l}^{\alpha})^{\vee}$;
\item\label{item Xvee ConeCuttau}
$X^{\vee}_{\bR}=\bigcup_{\alpha\in\Sk^0(\Sigma_{l}(0))}
\Cone(\Cut(\tau_{l,\alpha,0}))$;
\item\label{item:Sigma_0} 
$\Sigma_{l}(0)
=\left\{x\in X_{\bR}; u(x)\geq u(\beta)\ 
\left(\begin{matrix}\forall \beta\in\Sk^0(\Sigma_{l}(0))\\ 
\forall u\in X^{\vee}\cap\Cone(\Cut(\tau_{l,\beta,0}))
\end{matrix}\right)
\right\}$.
\end{enumerate}
\end{lemma}
\begin{proof}Lemma is  \cite[8.28]{MN24}. 
\end{proof}

\begin{defn}\label{defn:torus emb Zl_Sigma(0)}
We define a finite fan in $X^{\vee}$ by
\begin{equation*}
\Fan(\Sigma_{l}(0))=\left\{
\Cone(\Cut(\tau_{l,\Delta}));  
\Delta\in\Sk(\Sigma_{l}(0))\right\}.
\end{equation*}  Recall that 
$Z_{l}(\Sigma_{l}(0))=(W_{l,\Sigma_{l}(0)})_0\simeq G_0$ 
by \cite[8.20]{MN24} and the closure 
$\tZ_{l}(\Sigma_{l}(0))$ of $Z_{l}(\Sigma_{l}(0))$ is 
an irreducible component of $\tP_{l,0}$. 
By Lemma~\ref{lemma:torus embedding Zl_Sigma(0)}, 
$\tZ_{l}(\Sigma_{l}(0))$ is 
a torus embedding associated with $\Fan(\Sigma_{l}(0))$. 
It is covered with affine open subsets 
$\tZ_{l}(\Sigma_{l}(0))\cap W_{l,\alpha}$ 
$(\alpha\in\Sk^0(\Sigma_{l}(0))$.
Since $\tZ_{l}(\Sigma_{l}(v))\simeq\tZ_{l}(\Sigma_{l}(0))$\ $(\forall v\in X^{\vee})$, $\tZ_{l}(\Sigma_{l}(v))$ is covered 
with $\tZ_{l}(\Sigma_{l}(v))\cap W_{l,a}$\ $(a\in\Sk^0(\Sigma_{l}(v))$ 
and $Z_{l}(\Sigma_{l}(v))=(W_{l,\Sigma_{l}(v)})_0
=\Spec k(0)[t^{v(x)}w^x;x\in X]\simeq G_0$. 
\end{defn}

\begin{lemma}\label{lemma:Sigma0 Star_alpha mgeq4}
Let $w\in X^{\vee}$, $\alpha\in\Sigma_{l}(0)$ 
and $\mu^*:=2l\mu$. Then
\begin{enumerate}
\item\label{item:2Sigma0}
 if $\Sigma_{l}(w)\cap\Sigma_{l}(0)\neq\emptyset$, 
then $\mu^*(w)\in 2\Sigma_{l}(0)$;
\item\label{item:Star_alpha} 
$\Star_{l}(\alpha):=\bigcup_{\alpha\in\Sigma_{l}(v)}\Sigma_{l}(v)
\subset 3\Sigma_{l}(0)$;
\item\label{item:4Sigma0} 
if $\Sigma_{l}(w)\cap\Star_{l}(\alpha)\neq\emptyset$, then 
$\mu^*(w)\in 4\Sigma_{l}(0)$;
\item\label{item:Star_alpha 4Sigma0} 
if $(\mu^*(w)+\Star_{l}(\alpha))\cap\Star_{l}(\alpha)\neq\emptyset$, then
$\mu^*(w)\in 4\Sigma_{l}(0)$; 
\item\label{item:m_geq 4} 
$(\mu^*(mw)+\Star_{l}(\alpha))\cap\Star_{l}(\alpha)\neq\emptyset$ 
for some $m\geq 4$, then $w=0$.  
\end{enumerate}
\end{lemma}
\begin{proof}
Recall $\Sigma_{l}(w)=\Sigma_{l}(0)+\mu^*(w)$. 
 We set $\Delta:=\Sigma_{l}(0)\cap\Sigma_{l}(w)\neq\emptyset$. 
Since $\Delta':=\Delta-\mu^*(w)$ is a face of $\Sigma_{l}(0)$,  
there exists $\beta\in\Sk^0(\Delta)$ such that 
$\beta':=\beta-\mu^*(w)\in\Sk^0(\Delta')$. This implies 
$\mu^*(w)=\beta-\beta'\in\Sigma_{l}(0)-\Sigma_{l}(0)=2\Sigma_{l}(0)$. Ths proves (\ref{item:2Sigma0}). 
If $\Sigma_{l}(w)\cap\Sigma_{l}(0)\neq\emptyset$, then 
$\Sigma_{l}(w)\subset 3\Sigma_{l}(0)$ because 
$\Sigma_{l}(w)=\Sigma_{l}(0)+\mu^*(w)\subset
\Sigma_{l}(0)+2\Sigma_{l}(0)=3\Sigma_{l}(0)$ by (\ref{item:2Sigma0}). 
Hence (\ref{item:Star_alpha}) follows. 

Next we shall prove (\ref{item:4Sigma0}). If $\Sigma_{l}(w)\cap\Star_{l}(\alpha)\neq\emptyset$, then $\Sigma_{l}(w)\cap\Sigma_{l}(w')\neq\emptyset$ for some $w'\in X^{\vee}$ 
such that $\alpha\in\Sigma_{l}(w')$, whence 
$\mu^*(w')\in 2\Sigma_{l}(0)$ by (\ref{item:2Sigma0}).
Meanwhile $\Sigma_{l}(w-w')\cap\Sigma_{l}(0)\neq\emptyset$ 
by $\Sigma_{l}(w)\cap\Sigma_{l}(w')\neq\emptyset$. 
By (\ref{item:2Sigma0}) 
$\mu^*(w-w')\in 2\Sigma_{l}(0)$, so that 
$\mu^*(w)\in\mu^*(w')+2\Sigma_{l}(0)\subset 4\Sigma_{l}(0)$. 
This proves (\ref{item:4Sigma0}).

Next we shall prove (\ref{item:Star_alpha 4Sigma0}). 
If $(\mu^*(w)+\Star_{l}(\alpha))\cap\Star_{l}(\alpha)\neq\emptyset$, then 
$\Sigma_{l}(w+w')\cap\Sigma_{l}(w'')\neq\emptyset$ 
for some $w',w''\in X^{\vee}$ 
such that $\alpha\in\Sigma_{l}(w')\cap\Sigma_{l}(w'')$. Hence 
$\Sigma_{l}(w'-w'')\cap\Sigma_{l}(0)\neq\emptyset$, 
whence $\mu^*(w'-w'')\in2\Sigma_{l}(0)$ by (\ref{item:2Sigma0}). 
Meanwhile $\Sigma_{l}(w+w'-w'')\cap\Sigma_{l}(0)\neq\emptyset$, so that 
$\mu^*(w+w'-w'')\in 2\Sigma_{l}(0)$ by  (\ref{item:2Sigma0}). It follows 
that $\mu^*(w)=\mu^*(-w'+w'')+\mu^*(w+w'-w'')\in 2\Sigma_{l}(0)+2\Sigma_{l}(0)\subset 4\Sigma_{l}(0)$. 
This proves (\ref{item:Star_alpha 4Sigma0}).
Finally we shall prove (\ref{item:m_geq 4}). If 
$(\mu^*(mw)+\Star_{l}(\alpha))\cap \Star_{l}(\alpha)\neq\emptyset$, 
then we have $\mu^*(mw)\in 4\Sigma_{l}(0)$ 
by (\ref{item:Star_alpha 4Sigma0}), 
whence $\mu^*(mw/4)\in\Sigma_{l}(0)$. 
If $m\geq 4$, then $\mu^*(w)\in\Sigma_{l}(0)$ 
because $\Sigma_{l}(0)$ is convex with $0\in\Sigma_{l}(0)$. 
By $\mu^*(w)\in\Sigma_{l}(w)^0$,  
we have $\Sigma_{l}(w)=\Sigma_{l}(0)$ and $w=0$. 
This proves (\ref{item:m_geq 4}).\end{proof}

\begin{cor}\label{cor:theta single term}
Let $m,n\in\bN$ such that  $m\geq 4n$. 
If $x\in n\Star_{l}(\alpha)$ and 
$x+\mu^*(mw)\in n\Star_{l}(\alpha)$, then $w=0$.
\end{cor}
\begin{proof}Let $z=x+\mu^*(mw)\in n\Star_{l}(\alpha)$. Then 
$z=x+\mu^*(mw)\in(\mu^*(mw)+n\Star_{l}(\alpha))
\cap n\Star_{l}(\alpha)$. 
By Lemma~\ref{lemma:Sigma0 Star_alpha mgeq4}~(\ref{item:m_geq 4}), 
$w=0$. 
\end{proof}

\subsection{The closed subscheme $\tZ_{l}(\Delta)$}
\label{subsec:barZ_l_Delta}

\begin{defn}\label{defn:Bla,BlDelta}
For any $\Delta\in\Vor_{l}$, we define following \cite[8.17]{MN24}
\begin{gather*}
\tau_{l,\Delta}=
\bigcap_{a\in\Sk^0(\Delta)}\tau_{l,a},\ \  
B_{l,\Delta}=R[\tau_{l,\Delta}^{\vee}\cap\tX],\ \ 
W_{l,\Delta}=\Spec B_{l,\Delta}.
\end{gather*} 
\end{defn}

\begin{defn}\label{defn:Zl_Delta}
For $\Delta\in\Vor_{l}$, 
let $Z_{l}(\Delta)$ be the unique reduced closed $\tG_0$-orbit of 
$(W_{l,\Delta})_0$, $\tZ_{l}(\Delta)$ 
the closure of $Z_{l}(\Delta)$ in $\tP_{l,0}$ and 
$\barZ_{l}(\Delta)$ 
the closure of $\pi_{l}(Z_{l}(\Delta))$ in $P_{l,0}$
with reduced structure. 
\end{defn}

\begin{lemma}
\label{lemma:strata of tP_l0} 
For $\Delta,\Delta'\in\Vor_{l}$,
\begin{enumerate}
\item\label{item:dim Delta} $\dim_k\tZ_{l}(\Delta)=\dim_{\bR}\Delta$;
\item\label{item:Delta subset Delta'}  
$\tZ_{l}(\Delta)\subset \tZ_{l}(\Delta')$ iff $\Delta\subset\Delta'$;
\item\label{item:Delta and Delta'} 
$\tZ_{l}(\Delta\cap\Delta')
=\tZ_{l}(\Delta)\cap \tZ_{l}(\Delta')$.
\end{enumerate} 
\end{lemma}
\begin{proof}
(\ref{item:dim Delta})-(\ref{item:Delta and Delta'}) are standard 
in torus embedding.  See \cite[1.6, p.~10]{Oda85}. 
\end{proof}

\begin{defn}\label{defn:fan Fan(Delta)}
Let $\Delta\in\Sk(\Sigma_{l}(0))$ and let   
$X(\Delta)$ be the saturation in $X$ of 
the $\bZ$-submodule of $X$ generated by $\Delta-\Delta$. 
There is a natural epimorphism 
$j_{\Delta}:X^{\vee}_{\bR}\to X(\Delta)^{\vee}_{\bR}$.   
For $\Delta'\in\Sk(\Delta)$, let  
$D_{\Delta}(\Delta'):=j_{\Delta}(\Cone(\Cut(\tau_{l,\Delta'})))$ and 
let $\Fan_{l}(\Delta)$ be 
a fan in $X(\Delta)^{\vee}$ consisting of  
$D_{\Delta}(\Delta')$\ 
$(\Delta'\in\Sk(\Delta))$.   
See also \cite[\S~5.7]{D78}.
\end{defn}

\begin{lemma}\label{lemma:fan_Delta}Let 
$\Delta\in\Sk(\Sigma_{l}(0))$ and $\Delta'\in\Sk(\Delta)$. 
Then 
\begin{enumerate}
\item\label{item:tau_lDelta orthogonal X(Delta)} 
$X(\Delta)_{\bR}=\{x\in X_{\bR}; F(x)=0\ 
(\forall  F\in\tau_{l,\Delta})\}$;
\item\label{item:D_Delta_Delta} $D_{\Delta}(\Delta)=0$;
\item\label{item:each cone of fan_Delta}
$D_{\Delta}(\Delta') 
=\Cone(j_{\Delta}(e(\xi)v);  e(\xi)v\in\Cut(\tau_{l,\Delta'})\setminus\Cut(\tau_{l,\Delta}))$;
\item\label{item:barZ_Delta_simeq_barZ_Delta+2lmuv} 
$\tZ_{l}(\Delta)\simeq\tZ_{l}(2l\mu(v)+\Delta)$\ 
$(\forall v\in X^{\vee})$. 
\end{enumerate}
\end{lemma}
\begin{proof}By applying \cite[8.21~(3)]{MN24} for $f=\id_{X_{\bR}}$, 
$X(\Delta)_{\bR}$ is the {\it maximal} linear subspace of 
$\Cone(B_{l,\Delta})=\tau_{l,\Delta}^{\vee}$ 
 Definition~\ref{defn:Bla,BlDelta}.  
It follows 
$X(\Delta)_{\bR}=\{x\in X_{\bR}; F(x)=0\ (\forall  F\in\tau_{l,\Delta})\}$.
This proves (\ref{item:tau_lDelta orthogonal X(Delta)}). 
(\ref{item:D_Delta_Delta}) follows 
from (\ref{item:tau_lDelta orthogonal X(Delta)}). 
(\ref{item:each cone of fan_Delta}) follows from 
$\Sk^0\Cut(\tau_{l,\Delta'})=\{e(\xi)v; \Delta'\subset\Sigma_{l}(-v), 
v\in X^{\vee}\}$ \cite[8.27~(2)]{MN24}. (\ref{item:barZ_Delta_simeq_barZ_Delta+2lmuv}) is clear. 
\end{proof}

Now we describe $\tZ_{l}(\Delta)$ 
as a toric variety through the polytope $\Delta$. 
\begin{defn}\label{defn:Fan*_Delta}
Let $\Delta\in\Vor_{l}$ and 
$\Delta'\in\Sk(\Delta)$. We define a cone $\sigma_{\Delta'}$  by 
\begin{equation}
\sigma_{\Delta'}=\Cone(\Delta-a; a\in\Delta')\subset X(\Delta)_{\bR}
\end{equation}and set 
$\Fan^*_{l}(\Delta):=\{\sigma^{\vee}_{\Delta'};\Delta'\in\Sk(\Delta)\}$.
\end{defn}

\begin{lemma}\label{lemma:Fan*_Delta=Fan_Delta}
Let $\Delta\in\Vor_{l}$ and $\Delta'\in\Sk(\Delta)$, 
Then $D_{\Delta}(\Delta')=\sigma_{\Delta'}^{\vee}$ and 
$\Fan^*_{l}(\Delta)=\Fan_{l}(\Delta)$. 
\end{lemma}
\begin{proof}Without loss of generality, we 
may assume $\Delta\subset\Sigma_{l}(0)$.
Let $\Delta':=\alpha\in\Sk^0(\Delta)$. Then  
$D_{\Delta}(\Delta')=
\Cone(e(\xi)v_{|X(\Delta)_{\bR}};\alpha\in\Sigma_{l}(-v),v\in X^{\vee})$. 
Let $e(\xi)v\in D_{\Delta}(\Delta')$. We shall prove 
$e(\xi)v\in\sigma_{\alpha}^{\vee}$.  Let 
$x\in\Delta$. Since $\Delta\subset\Sigma_{l}(0)$, 
we have $E_{l}(v)+v(x)\geq 0$. 
Since $\alpha\in\Sigma_{l}(-v)$, we have 
$E_{l}(-v)-v(\alpha+2l\mu(v))\geq 0$. It follows 
$0\leq 2E_{l}(v)+v(x-\alpha)-2lv\mu(v)=v(x-\alpha)$. Hence 
$v(x-\alpha)\geq 0$ $(\forall x\in\Delta)$, {\it i.e.}, 
$v\in\sigma_{\alpha}^{\vee}$. Since  
$\tau_{l,\Delta'}=\bigcap_{\alpha\in\Sk^0(\Delta')}\tau_{l,\alpha}$, 
we have \begin{align*}
D_{\Delta}(\Delta')&=\bigcap_{\alpha\in\Sk^0(\Delta')}D_{\Delta}(\alpha)
=\bigcap_{\alpha\in\Sk^0(\Delta')}\sigma_{\alpha}^{\vee}
=\sigma_{\Delta'}^{\vee}.
\end{align*}  
Since $|\Fan^*_{l}(\Delta)|=|\Fan_{l}(\Delta)|=X(\Delta)_{\bR}$,  
we have $D_{\Delta}(\Delta'))
=\sigma^{\vee}_{\Delta'}$\ $(\forall \Delta'\in\Sk(\Delta))$. 
Hence $\Fan^*_{l}(\Delta)=\Fan_{l}(\Delta)$.
\end{proof}

\begin{lemma}\label{lemma:Z_Delta=torus emb Fan_Delta}
The closed subscheme 
$\tZ_{l}(\Delta)$ of $\tP_{l,0}$ is a complete 
torus embedding associated with $\Fan_{l}(\Delta)$ and hence 
$\Fan^*_{l}(\Delta)$.
\end{lemma}
\begin{proof} 
This is standard in torus embedding. See 
\cite[\S\S~5.7-5.8]{D78}. 
\end{proof}

\begin{cor}\label{cor:ZDelta' subset WDelta}
Let $\Delta$ and $\Delta'\in\Vor_{l}$. Then  
$\Delta'\subset\Delta\Leftrightarrow W_{l,\Delta'}\supset W_{l,\Delta}\Leftrightarrow Z_{l}(\Delta')\subset \tZ_{l}(\Delta)\cap W_{l,\Delta'}
\Leftrightarrow Z_{l}(\Delta)\subset (W_{l,\Delta'})_0$.
\end{cor}

\begin{proof}By 
Lemma~\ref{lemma:strata of tP_l0}~(\ref{item:Delta subset Delta'}), 
$\tZ_{l}(\Delta')\subset \tZ_{l}(\Delta)\Rightarrow \Delta'\subset\Delta\Rightarrow W_{l,\Delta'}\supset W_{l,\Delta}\Rightarrow Z_{l}(\Delta)\ (\subset W_{l,\Delta})\subset W_{l,\Delta'}\Rightarrow Z_{l}(\Delta')\subset\tZ_{l}(\Delta)\cap W_{l,\Delta'}\Rightarrow Z_{l}(\Delta')\subset\tZ_{l}(\Delta)\Rightarrow \tZ_{l}(\Delta')\subset\tZ_{l}(\Delta)$.  
Meanwhile $Z_{l}(\Delta)\subset W_{l,\Delta'}\Leftrightarrow Z_{l}(\Delta)\subset (W_{l,\Delta'})_0$ by $Z_{l}(\Delta)\subset(\tP_{l})_0$.
This proves Corollary.
\end{proof}

\subsection{Support functions}

\begin{defn}\label{defn:convex function F}
By Lemma~\ref{lemma:torus embedding Zl_Sigma(0)}~(\ref{item Xvee ConeCuttau}), 
$X^{\vee}_{\bR}$ is covered with 
$\Cone(\Cut(\tau_{l,\alpha,0}))$ 
$(\alpha\in\Sk^0(\Sigma_{l}(0)))$. 
For each $v\in X^{\vee}$, we define a function $h_v$ 
on $X^{\vee}_{\bR}$ by 
\begin{align}\label{eq:defn h_v}
h_v(u)&=u(\alpha+2l\mu(v))
\quad\text{if $u\in\Cone(\Cut(\tau_{l,\alpha,0}))$.}
\end{align}
Then $h_v$ is a {\it $\Fan(\Sigma_{l}(0))$-linear 
support function} 
 in the sense of \cite[\S~2.1]{Oda85}, {\it i.e.}, $h_v$ is 
$\bZ$-valued strictly upper convex on 
$X\cap|\Fan(\Sigma_{l}(0))|$ and linear on each 
$\sigma\in\Fan(\Sigma_{l}(0))$. 
\end{defn}

\begin{lemma}\label{lemma:F convex}The function 
$h_v$\ $(v\in X^{\vee})$ is continuous strictly upper convex on 
$X^{\vee}_{\bR}$ with respect to $\Fan(\Sigma_{l}(0))$.
\end{lemma}
\begin{proof}Let 
$\Star^{\vee}(0):=\bigcup_{\alpha\in\Sk^0(\Sigma_{l}(0))}
\Cut(\tau_{l,\alpha,0})$. 
First we prove that  $h_v$ is continuous on $\Star^{\vee}(0)$.  Let  
$\beta\in\Sk^0(\Sigma_{l}(0))$. 
Any one-codimensional face $\sigma$ 
of the cone $\tau_{l,\beta,0}$ 
is, in view of \cite[8.23~(5)]{MN24}, 
of the form $\sigma=\tau_{l,\beta,0}\cap \tau_{l,\gamma,0}$ 
where $\gamma\in\Sk^0(\Sigma_{l}(0))$ 
with $\Conv(\beta,\gamma)\in\Sk^1(\Sigma_{l}(0))$.
In other words,
$\sigma=\{u\in \Star^{\vee}(0); u(\beta)=u(\gamma)\}$,  
which shows that $h_v$ is continuous on $\Star^{\vee}(0)$ and hence so is 
$h_v$ on $X^{\vee}_{\bR}$ by 
Lemma~\ref{lemma:torus embedding Zl_Sigma(0)}
~(\ref{item Xvee ConeCuttau}). 
It follows from Lemma~\ref{lemma:torus embedding Zl_Sigma(0)}~(\ref{item:tau_alpha}) and Definition~\ref{defn:convex function F} 
that $h_v$ is strictly upper convex.  
\end{proof}

\begin{lemma}\label{lemma:hv as order of tcL_l on W(v)}
Let $v\in X^{\vee}$, 
$\barZ_{l}(v):=\barZ_{l}(\Sigma_{l}(v))$ and 
$W_{l,a}(v):=W_{l,a}\cap\barZ_{l}(v)$. 
We define an $\cO_{\barZ_{l}(v)}$-module $\cD_{l}(v)$ by 
$\cD_{l}(v)_{|W_{l,a}(v)}:=t^{v(a)}w^a\cO_{W_{l,a}(v)}.$  
Then 
\begin{enumerate}
\item\label{item:isom of tcL_l and cD_l_(v)}
$\tcL_{l}\otimes_{\cO_{\tP_{l}}}\cO_{\barZ_{l}(v)}
\simeq\cD_{l}(v)$ as $\cO_{\barZ_{l}(v)}$-modules;
\item\label{item:support function of D(v)}$h_v$ 
is the support function of 
$\cD_{l}(v)$ in the sense of 
\cite[\S~2.1, p.~65]{Oda85}.
\end{enumerate}    
\end{lemma}
\begin{proof}A one cocycle associated with the sheaf 
$\tcL_{l}\otimes_{\tP_{l}}\cO_{\tZ_{l}(v)}$ is given by 
$\left(f^*_{a,b}=\xi_{l,b}/\xi_{l,a};a,b\in\Sk^0(\Sigma_{l}(v))\right)$. 
By \cite[\S~4.1/6.7]{MN24}, 
$\wt_{\tX}(\xi_{l,c})=(D_{l}(c)+v(c))m_0+v(c)$ 
if $c\in\Sk^0(\Sigma_{l}(v))$. 
Since $\tZ_{l}(v)$ is a complete torus embedding over $k(0)$, we have 
$H^1(\tZ_{l}(v),k(0)^*)=\{1\}$. Hence we may replace $f^*_{a,b}$ by 
$$f_{a,b}=t^{D_{l}(b)+v(b)}w^b/t^{D_{l}(a)+v(a)}w^a$$ 
by \cite[6.7]{MN24}, that is, the one cocycle 
$(f^*_{a,b}/f_{a,b};a,b\in\Sk^0(\Sigma_{l}(v)))\in H^1(\tZ_{l}(v),k(0)^*)$ 
is trivial.  
Since $a,b\in\Sigma_{l}(v)$, we have $D_{l}(a)=D_{l}(b)$ 
by \cite[6.7]{MN24}, so that $f_{a,b}=t^{v(b)}w^b/t^{v(a)}w^a$, 
which proves (\ref{item:isom of tcL_l and cD_l_(v)}). 
Since $Z_{l}(v)\simeq\Spec k(0)[t^{v(x)}w^x;x\in X]$ 
by Definition~\ref{defn:torus emb Zl_Sigma(0)}, 
by \cite[2.1, pp.~68-69]{Oda85}, 
$h_v$ is the support function of the sheaf 
$\cD_{l}(v)$, that is, $h_v(u)=u(a)$ if $u\in\Cone(\Cut(\tau_{l,\alpha,0}))$ 
where $a=\alpha+2l\mu(v)\in\Sigma_{l}(v)$. 
This proves (\ref{item:support function of D(v)}).  
\end{proof}

\subsection{Projective embedding of $\tZ_{l}(\Delta)$}
\label{subsec:projective embedding of tZlDelta}

\begin{defn}\label{defn:cone and semi Delta-alpha}
Let $\Delta\in\Vor_{l}$, 
$\alpha\in\Sk^0(\Delta)$ and 
$\Cone(\Delta-\alpha)$ the convex cone in $X(\Delta)_{\bR}$ 
generated by $\beta-\alpha$ $(\forall \beta\in\Sk^0(\Delta))$. 
Let $T$ be a closed subset of $X_{\bR}$, $T^0$ 
the relative interior of\/ $T$ and $mT:=\{\sum_{i=1}^mv_i; v_i\in T\}$. 
\end{defn}

\begin{defn}\label{defn:primitive}
Let $C$ be a convex cone in $X_{\bR}$ 
generated by finitely many integral vectors.  
We say that $w\in C\cap X$ is a {\it primitive element of $C\cap X$} 
or $w$ is {\it primitive in $C\cap X$}   
if $w-v\not\in C\cap X$ 
for any $v\in C\cap X\setminus\{0,w\}$.
\end{defn}

\begin{lemma}\label{lemma:SmDelta_alpha generate C}
Let $\Delta\in\Vor_{l}$, $r=\dim_{\bR}\Delta$,  
$\alpha\in\Sk^0(\Delta)$ and 
$S_{m,\Delta-\alpha}:=\Semi(m(\Delta-\alpha)\cap X)$. Then 
\begin{enumerate}
\item\label{item:m_geq r+1}$(m\Delta)^0\cap X\neq\emptyset$ if $m\geq r+1$;
\item\label{item:primitive element}any primitive element belongs to $S_{r,\Delta-\alpha}$;any element of $S_{m,\Delta-\alpha}\setminus S_{r,\Delta-\alpha}$ is 
a sum of two or more primitive elements if $m\geq r+1$;
\item\label{item:generated by mDelta-malpha}
 $S_{m,\Delta-\alpha}=\Cone(\Delta-\alpha)\cap X$ if $m\geq r$. 
\end{enumerate}
\end{lemma}
\begin{proof}Since $\Sigma_{l}(0)$ is integral by assumption, 
 $\Delta$ is the convex closure of $\Sk^0(\Delta)$. 
Since $\dim_{\bR} \Delta=r$, there exist  
$\alpha_i\in\Sk^0(\Delta)$ $(i\in[0,r])$  
such that $D:=\Conv(\alpha_i,i\in[0,r])$ 
is a subset of $\Delta$ with $\dim_{\bR}D=r$. 
Then $\sum_{i=0}^r\alpha_i\in((r+1)D)^0\cap X\subset((r+1)\Delta)^0\cap X$. 
Hence there exist $\epsilon>0$ such that 
$\sum_{i=0}^ra_i\alpha_i\in(r+1)D$ if $|a_i-1|<\epsilon$\ 
$(\forall i\in[0,r])$. Therefore 
$(m-r-1+a_0)\alpha_0+\sum_{i=1}^ra_i\alpha_i\in mD$, whence 
$(m-r)\alpha_0+\sum_{i=1}^r\alpha_i\in(mD)^0\cap X\subset(m\Delta)^0\cap X$. 
This proves (\ref{item:m_geq r+1}).  
Let $C:=\Cone(\Delta-\alpha)$. Then $C$ is a cone 
of dimension $r$ with $0$ vertex. 
Let $w\in C\cap X$ be any {\it primitive} element of $C\cap X$.
Let $v_i\in(\Delta-\alpha)\cap X$ $(i\in[1,r])$ 
such that $w$ is contained in $\Cone(v_i; i\in [1,r])$.  
Since $w$ is primitive, we can write 
$w=\sum_{i=1}^ra_iv_i$ for some $0\leq a_i\leq 1$. 
Hence $w\in r(\Delta-\alpha)\cap X$, so that 
$w\in S_{r,\Delta-\alpha}$. This proves (\ref{item:primitive element}).
Since any $w'\in C\cap X$ 
is written as $w=\sum_{i=1}^sw_i$ 
with $w_i$ primitive in $C\cap X$ and $s\in\bN$, we have 
$C\cap X\subset S_{r,\Delta-\alpha}\subset 
S_{m,\Delta-\alpha}$ $(\forall m\geq r)$.  
This proves (\ref{item:generated by mDelta-malpha}) because 
 $S_{m,\Delta-\alpha}\subset C$.
\end{proof}

\begin{thm}\label{thm:cohomology of Zl_Delta}
Let $m\geq 1$, $\Delta\in\Sk(\Sigma_{l}(v))$,   
$\tcL_{l}(\Delta):=\tcL_{l}\otimes_{\cO_{\tP_{l,0}}}
\cO_{\tZ_{l}(\Delta)}$\ 
and $\square_{h_v}(\Delta):=\{x\in\Delta+X(\Delta)_{\bR}; 
u(x)\geq h_v(u)\ (\forall u\in X^{\vee}_{\bR})\}$. 
Then 
\begin{enumerate}
\item\label{item: SquareL_h_v Delta} $\square_{h_v}(\Delta)=\Delta$; 
\item\label{item:vanishing cL_Delta} 
$H^q(\tZ_{l}(\Delta),\tcL_{l}(\Delta)^{\otimes m})=0$ $(q>0)$; 
\item\label{item:base points free} $\tcL_{l}(\Delta)^{\otimes m}$ is 
generated by global sections; 
\item\label{item:Gamma Zl_Delta cL_Delta}
$(t^{v(x)}w^x;x\in m\Delta\cap X)$ 
is a $k(0)$-basis of 
$H^0(\tZ_{l}(\Delta),\tcL_{l}(\Delta)^{\otimes m})$;
\item\label{item:very ample if m_geq_r}
 $\tcL_{l}(\Delta)^{\otimes m}$ is very ample if 
$m\geq\dim_{\bR}\Delta$.
\end{enumerate}
\end{thm}
\begin{proof}Let $r=\dim_{\bR}\Delta$. 
By Lemma~\ref{lemma:torus embedding Zl_Sigma(0)}
~(\ref{item:Sigma_0}), $\square_{h_0}(\Sigma_{l}(0))=\Sigma_{l}(0)$ and  
$\square_{h_v}(\Sigma_{l}(v))=\Sigma_{l}(v)$.  
By \cite[8.21~(1)]{MN24}, 
$\square_{h_v}(\Delta)=(\Delta+X(\Delta)_{\bR})
\cap\square_{h_v}(\Sigma_{l}(v))=(\Delta+X(\Delta)_{\bR})
\cap\Sigma_{l}(v)=\Delta$. 
This proves (\ref{item: SquareL_h_v Delta}). 
(\ref{item:vanishing cL_Delta}) and  (\ref{item:base points free}) 
follow from Lemma~\ref{lemma:F convex} and  \cite[2.7/2.9]{Oda85}.
Since $mh_v$ is the support function of 
$\tcL^{\otimes m}_{l}\otimes_{\cO_{\tP_{l}}}\cO_{\tZ_{l}(v)}$ 
by Lemma~\ref{lemma:hv as order of tcL_l on W(v)} 
and $\vor_{l,ml}(\Delta)=m\Delta$ 
by \cite[8.15]{MN24}, we have 
$\tZ_{l}(\Delta)=\tZ_{lm}(m\Delta)$ and  
$H^0(\tZ_{l}(\Delta),\tcL_{l}(\Delta)^{\otimes m})=
H^0(\tZ_{lm}(m\Delta),\tcL_{lm}(m\Delta)).$ 
Hence  
(\ref{item:Gamma Zl_Delta cL_Delta}) follows from 
(\ref{item: SquareL_h_v Delta})   
and \cite[2.3]{Oda85}.  
The $k(0)$-scheme $\tZ_{l}(\Delta)$ is 
naturally a closed subscheme of not only $\tP_{l,0}$ but also 
$\tP_{l,0}/nY\simeq(\tP_{l}^{\wedge}/nY)_0$ 
for a very large $n\in\bN$ such that 
$(n(Y\setminus\{0\})+\tZ_{l}(\Delta))\cap\tZ_{l}(\Delta)
=\emptyset$.  
Since $\tcL_{l}^{\wedge}/nY$ is ample on $\tP_{l}^{\wedge}/nY$, 
so is $\tcL_{l,0}/nY$ on $\tP_{l,0}/nY$. 
Hence $\tcL_{l}(\Delta)$ is ample. 
Let $\alpha\in\Sk^0(\Delta)$ and  
$D_{\Delta}(\alpha):=j_{\Delta}(\Cone(\Cut(\tau_{l,\alpha})))$. 
Then $D_{\Delta}(\alpha)=\sigma^{\vee}_{\alpha}=\Cone(\Delta-\alpha)^{\vee}$. 
by Lemma~\ref{lemma:Fan*_Delta=Fan_Delta}. 
By Lemma~\ref{lemma:SmDelta_alpha generate C}~(\ref{item:generated by mDelta-malpha}), we obtain $S_{m,\Delta-\alpha}=\Cone(\Delta-\alpha)\cap X=
D_{\Delta}(\alpha)^{\vee}\cap X$ if $m\geq r$.
This proves (\ref{item:very ample if m_geq_r}) 
by \cite[2.13~(b)]{Oda85}. 
\end{proof}

\begin{rem}\label{rem:very ampleness}Here is another proof of 
Theorem~\ref{thm:cohomology of Zl_Delta}~(\ref{item:very ample if m_geq_r}).
Let $\alpha\in\Delta\subset\Sigma_{l}(v)$. 
The equality 
$D(\Sigma_{l}(v))=\Cone(\Sigma_{l}(v)-\alpha)^{\vee}$ 
follows from Lemma~\ref{lemma:torus embedding Zl_Sigma(0)}~(\ref{item:tau_alpha}). By Lemma~\ref{lemma:SmDelta_alpha generate C}~(\ref{item:generated by mDelta-malpha}), $S_{m,\Sigma_{l}(v)-\alpha}=\Cone(\Sigma_{l}(v)-\alpha)\cap X$ if $m\geq r$. Hence $\tcL_{l}(\Sigma_{l}(v))^{\otimes m}$ is very ample. 
Since the restriction map 
$\Gamma(\tZ_{l}(v),\tcL_{l}(\Sigma_{l}(v))^{\otimes m})
\to \Gamma(\tZ_{l}(\Delta),\tcL_{l}(\Delta)^{\otimes m})$ 
is surjective by Theorem~\ref{thm:cohomology of Zl_Delta}~(\ref{item:Gamma Zl_Delta cL_Delta}),  
$\tcL_{l}(\Delta)^{\otimes m}$ is very ample. This proves Theorem~\ref{thm:cohomology of Zl_Delta}~(\ref{item:very ample if m_geq_r}).
\end{rem}

\subsection{The subscheme $Z_{l}(\Delta)$}
\label{subsec:Zl_Delta}

\begin{defn}\label{defn:Sup_vee}
Let $Z$ be a subset of $\tP_{l,0}$. We define 
\begin{align*}
\Sup^{\vee}(Z)&
=\{e(\xi)v\in X^{\vee}; Z\subset\overline{(W_{l,\Sigma_{l}(-v)})_0}, 
v\in X^{\vee}\},\\
\Sup(Z)&=\{a\in X; Z\subset (W_{l,a})_0\}. 
\end{align*}
\end{defn}

\begin{lemma}\label{lemma:Sup vee}Let $\Delta\in\Vor_{l}$.
If $q\in Z_{l}(\Delta)$, then 
\begin{enumerate}
\item\label{item:Sup_wedge} 
$\Sup^{\vee}(q)=\Sup^{\vee}(Z_{l}(\Delta))
=\Cut(\tau_{l,\Delta})\cap e(\xi)X^{\vee}$;
\item\label{item:Sup}    
$\Sup(q)=\Sup(Z_{l}(\Delta))=\Delta\cap X$.
\end{enumerate}
\end{lemma}
\begin{proof}For $v\in X^{\vee}$, 
$Z_{l}(\Delta)\subset\tZ_{l}(\Sigma_{l}(-v))$ iff 
 $\Delta\subset\Sigma_{l}(-v)$. 
Since $\Cut(\tau_{l,\Delta})=
\Conv(e(\xi)v\in X^{\vee}; \Delta\subset\Sigma_{l}(-v), 
v\in X^{\vee})$ by \cite[8.27~(2)]{MN24},
we have $\Cut(\tau_{l,\Delta})\cap e(\xi)X^{\vee}
=\Sup^{\vee}(Z_{l}(\Delta))$.  Let $q\in Z_{l}(\Delta)$. Then 
$\Sup^{\vee}(Z_{l}(\Delta))\subset\Sup^{\vee}(q)$. 
If $e(\xi)v\in\Sup^{\vee}(q)$, then $q\in Z_{l}(\Delta)\cap\tZ_{l}(\Sigma_{l}(-v))\neq\emptyset$ because $\overline{(W_{l,\Sigma_{l}(-v)})_0}=\tZ_{l}(\Sigma_{l}(-v))$. Since $\tZ_{l}(\Sigma_{l}(-v))$ 
is a torus embedding, $Z_{l}(\Delta)$ is a stratum of 
$\tZ_{l}(\Sigma_{l}(-v))$, 
so that $Z_{l}(\Delta)\subset\tZ_{l}(\Sigma_{l}(-v))$. 
Hence $e(\xi)v\in\Sup^{\vee}(Z_{l}(\Delta))$, so that  
$\Sup^{\vee}(q)\subset \Sup^{\vee}(Z_{l}(\Delta))$. 
This proves (\ref{item:Sup_wedge}).\par 

Let $q\in Z_{l}(\Delta)$. 
For $a\in X$, there exists $\Delta'\in\Vor_{l}$ 
such that $a\in(\Delta')^0=$
the relative interior of $\Delta'$.
\footnote{by taking for $l$ its suitable multiple $ll'$ if necessary.} 
 Then $W_{l,a}=W_{l,\Delta'}$ by \cite[8.22~(3)]{MN24}.  
Since $W_{l,\Delta'}$ is a torus embedding, we see that
$q\in(W_{l,\Delta'})_0$ iff 
$Z_{l}(\Delta)\subset(W_{l,\Delta'})_0$  
iff $\Delta'\subset\Delta$ 
by Corollary~\ref{cor:ZDelta' subset WDelta}.  In particular, 
$a\in\Sup(q)$ iff  
$q\in (W_{l,a})_0$ iff $Z_{l}(\Delta)\subset(W_{l,a})_0$ 
($\Leftrightarrow a\in\Sup(Z_{l}(\Delta))$) 
iff $a\in\Delta$. 
This proves (\ref{item:Sup}). 
\end{proof}

\subsection{Specialization of points of $\tP_{l}(\Omega)$}
\label{subsec:appendix specialization}
\begin{defn}\label{defn:specialization of product}
Let $Q\in\tP_{l}(\Omega)$. 
By \cite[7.6]{MN24}, 
$\tP_{l}(\Omega)=\tP_{l}(R_{\Omega})$, so that 
$\tP_{l}(R_{\Omega})\ni Q$, which is 
 an $R$-flat finite 
formal $R$-scheme.\footnote{As is equivalent, it is 
an admissible finite formal $R$-scheme. 
See \cite[10.8]{MN24}.}
The closed fiber $Q_0$ of $Q$ is 
the specialization $\specialization(Q)$ of $Q$ 
in the sense of \cite[p.~200]{Bosch14}.  
\end{defn}

Let $Q\in\tP_{l}(\Omega)$ and 
$\log|w^{x}(Q)|:=\log|Q^*(w^x)|\in\bQ$ $(x\in \tX)$.
Since $w^{m_0}(Q)=s$, 
$\log(Q):=-\sum_{i=0}^g(\log|w^{m_i}(Q)|)f_i\in\tX_{\bQ}^{\vee}$ and 
$\cutlog(Q):=-f_0+\log(Q)=-\sum_{i=1}^g(\log|w^{m_i}(Q)|)f_i
\in X_{\bQ}^{\vee}$. See \cite[10.14]{MN24}. 

We also recall from \cite[10.15]{MN24}: 
\begin{lemma}\label{lemma:cutlog}Let 
$a\in\Sk^0(\Vor_{l})$, 
$u\in X^{\vee}$ and $Q\in\tP_{l}(R_{\Omega})$. 
Then 
\begin{enumerate}
\item\label{item:Q0} if $Q_0\in (W_{l,a})_0$, 
then $Q\in W_{l,a}(R_{\Omega})$;
\item\label{item:cutlog}
$Q\in W_{l,a}(R_{\Omega})\Leftrightarrow 
\log(Q)\in\tau_{l,a}\Leftrightarrow \cutlog(Q)\in\Cut(\tau_{l,a})$.
\end{enumerate}
\end{lemma}

\begin{thm}\label{thm:condition Q0 in ZDelta}
Let $\Delta\in\Vor_{l}$   
and $Q\in P_{l}(R_{\Omega})$. 
Then $Q_0\in Z_{l}(\Delta)$ iff  
$\cutlog(Q)\in\Cut(\tau_{l,\Delta})^0$. 
\end{thm}
\begin{proof}First we prove the if-part. 
Recall that $Z_{l}(\Delta)$ is 
the unique reduced closed $G_0$-orbit of $W_{l,\Delta}$, 
and $W_{l,\Delta}$ is an affine open subscheme of $\tP_{l}$ 
with $\Gamma(\cO_{W_{l,\Delta}})=B_{l,\Delta}
=R[(\tau_{l,\Delta})^{\vee}\cap \tX]
=R[w^{\tx};\tx\in(\tau_{l,\Delta})^{\vee}\cap \tX]$. Let $X(\Delta)$ be the saturation in $X$ of the $\bZ$-submodule of $X$ generated by $\Delta-\Delta$, and $\cI_{l,\Delta}$ the ideal of $\cO_{W_{l,\Delta}}$ 
generated by all $w^{\tx}\in B_{l,\Delta}$ such that 
$\tx=x_0m_0+x\in\tau_{l,\Delta}^{\vee}\cap \tX$ 
with $x\in X\setminus X(\Delta)$. Then $Z_{l}(\Delta)$ is a 
closed subscheme of $W_{l,\Delta}$ defined by $\cI_{l,\Delta}$. 
Since $\cutlog(Q)\in\Cut(\tau_{l,\Delta})^0$, we have 
$v_sw^{\tx}(Q)=\langle\log(Q),x\rangle= \cutlog(Q)(x)>0$ by Lemma~\ref{lemma:fan_Delta}~(\ref{item:tau_lDelta orthogonal X(Delta)}) if  
$\tx=x_0m_0+x\in\tau_{l,\Delta}^{\vee}\cap\tX$ 
with $x\in X\setminus X(\Delta)$. Hence $w^{\tx}(Q)\in sR$,  
$w^{\tx}\in \cI_{l,\Delta}$ and $Q_0\in Z_{l}(\Delta)$.\par

Next we prove the only-if-part. Let $Q_0\in Z_{l}(\Delta)$. 
Then $\cutlog(Q)\in\Cut(\tau_{l,\Delta})$ by 
Lemmas~\ref{lemma:cutlog}/\ref{lemma:Sup vee}~(\ref{item:Sup}). 
Assume $\cutlog(Q)\not\in\Cut(\tau_{l,\Delta})^0$. Then there exists a face $\Delta'$ of $\Delta$ such that $Q\in\Cut(\tau_{l,\Delta'})^0$. By what we have proved first, $Q_0\in Z_{l}(\Delta')$, which contradicts $Q_0\in Z_{l}(\Delta)$ because $Z_{l}(\Delta)\cap Z_{l}(\Delta')=\emptyset$. It follows 
$\cutlog(Q)\in\Cut(\tau_{l,\Delta})^0$.
This completes the proof.
 \end{proof}
\begin{example}\label{example:Fan_l_Delta}
We recall the scheme $\tP^{\natural}=\Proj R^{\natural}$ 
from \cite[\S~13.3]{MN24}. Let $R$ be a CDVR with $s$ uniformizer, $I=sR$, 
$k:=k(0)$, $e(\xi)=1$, $l=1$, $X=\bZ^2=\bZ m_1+\bZ m_2$, $X^{\vee}=\bZ^2=\bZ f_1+\bZ f_2$ with $f_i(m_j)=\delta_{i,j}$ $(i,j\in[1,2])$,  
$m_3:=-(m_1+m_2)$ and $f_3:=-(f_1+f_2)$.  
Let $u=u_1f_1+u_2f_2$ and $y=y_1m_1+y_2m_2$. 
Let $\theta_{+}$ be an indeterminate, and we define 
$R^{\natural}$, $\Sigma(0)$ and $\Sigma$  as follows:
\begin{gather*}
R^{\natural}=R[\epsilon_{+}(u)b^e(u,\alpha)w^{\alpha+\mu(u)}
\theta_{+}; u\in X^{\vee},\alpha\in\Sigma],\\
\Sigma(0)=\{x\in X_{\bR};E_{+}(u)+u(x)\geq 0\ (\forall u\in X^{\vee})\},\ 
\Sigma=\Sigma((0)\cap X 
\end{gather*}
where 
\begin{gather*}
\epsilon_{+}(u):=s^{u_1^2+u_1u_2+u_2^2},\ 
b^e(u,y)=s^{u_1y_1+u_2y_2},\\
\mu(u):=(2u_1+u_2)m_1+(u_1+2u_2)m_2,\\ 
E_{+}(u):=v_s\epsilon_{+}(u)=u(\mu(u))/2=u_1^2+u_1u_2+u_2^2. 
\end{gather*} 
 We see   
\begin{gather*}
\Sigma(0)=\Conv(\pm m_i,i\in[1,3]),\ 
\Sk^0(\Sigma(0))=\{\pm m_i\ (i\in[1,3])\},\\ 
\Sk^1(\Sigma(0))=\{\pm[m_i,-m_j],\ 
(i,j\in[1,3]\},\ \Sk^2((\Sigma(0)))=\{(\Sigma(0))\}.
\end{gather*}

\begin{figure}[ht]
\label{fig:Sigma0}
   \vspace*{-0.5cm}
   \begin{picture}(110,90)(56.5,-10)
   \multiput(90,10)(20,0){2}{\circle*{5}}
   \multiput(90,30)(20,0){3}{\circle*{5}}
   \multiput(110,50)(20,0){2}{\circle*{5}}
   \put(70,0){$m_3$}
   \put(140,35){$m_1$}
   \put(90,60){$m_2$}
   \put(90,10){\line(0,1){20}}
   \put(110,-10){\line(0,1){80}}
   \put(130,30){\line(0,1){20}}
   \put(90,10){\line(1,0){20}}
   \put(70,30){\line(1,0){80}}
   \put(110,50){\line(1,0){20}}
   \put(90,30){\line(1,1){20}}
   \put(110,10){\line(1,1){20}}
   \put(90,10){\line(1,1){20}}
   \put(110,30){\line(1,1){20}} 
   \end{picture} 
\caption{$\Sigma$ and $\Sigma(0)$}
\end{figure}

Let $\Delta=[m_1,-m_3]$. 
Then $\Sk(\Delta)=\{\Delta, \Delta', \Delta''\}$ and 
$X(\Delta)=\bR f_2$, where 
$\Delta':=\{m_1\}$ and $\Delta''=\{-m_3\}$. By \cite[8.27]{MN24}, 
$\Cut(\tau_{\Delta})=\Conv(0,-f_1)$,
$\Cut(\tau_{\Delta'})=\Conv(0,-f_1,-f_2)$ and 
$\Cut(\tau_{\Delta''})=\Conv(0,-f_1,-f_1+f_2)$, Hence 
$\Fan(\Delta)=\{0,\bR{\geq 0}f_2,\bR_{\geq 0}(-f_2)\}$, 
which is a complete fan.
\end{example}

\begin{example}
\label{examle:specialization}
Let $W_{-m_3}:=\Spec R[w_1^{-1},w_2^{-1},sw_3^{-1}]$. Then 
$W_{-m_3}$ is an open affine subset of $\tP^{\natural}$ and 
all the strata $Z(\Delta)$ of $(W_{-m_3})_0$ 
are listed with relevant data in Table~\ref{table:Q and Q0}, 
where we denote $\Conv(A)$ by $\langle A\rangle$. 

Now let $\alpha,\beta\in R_{\Omega}^{\times}$ and $\gamma=(\alpha\beta)^{-1}$. 
Let $\Delta=\langle m_1,-m_3\rangle^0$. 
Let us choose $Q\in W_{-m_3}(R_{\Omega})$ 
with $\cutlog(Q)=-\epsilon f_1$\ $(0<\epsilon <1)$\ 
(Table~\ref{table:Q and Q0}-$\text{sp}_2$):
$$Q:(w_1^{-1},w_2^{-1},sw_3^{-1})
=(s^{\epsilon}\alpha,\beta,s^{1-\epsilon}\gamma),$$  
then $Q_0:(w_1^{-1},w_2^{-1},sw_3^{-1})
=(0,\beta\op{mod} I,0)\in Z(\Delta)$ 
as $s$ approaches $0$. By Theorem~\ref{thm:condition Q0 in ZDelta}, 
if we assume $Q_0\in Z(\Delta)$, then  
$\cutlog(Q)\in\Cut(\tau_{l,\Delta})^0$, {\it i.e.}, 
$Q$ is of the above form.  
\par
 If $\Delta'=\langle -m_3,-2m_3\rangle^0$ and  
$\cutlog(Q)=-(\epsilon f_1+(1-\epsilon)f_2)$\ $(0<\epsilon <1)$ 
\ (Table~\ref{table:Q and Q0}-$\text{sp}_3$), then 
$Q:(w_1^{-1},w_2^{-1},sw_3^{-1})
=(s^{\epsilon}\alpha,s^{1-\epsilon}\beta,\gamma),$
so that $Q_0:(w_1^{-1},w_2^{-1},sw_3^{-1})
=(0,0,\gamma\op{mod} I)\in Z(\Delta')$.

Finally if $\Delta''=\langle -m_3\rangle$ and  
$\cutlog(Q)=-\epsilon_1f_1-\epsilon_2f_2$\ 
(Table~\ref{table:Q and Q0}-$\text{sp}_1$), then  
$Q:(w_1^{-1},w_2^{-1},sw_3^{-1})
=(s^{\epsilon_1}\alpha,s^{\epsilon_2}\beta,
s^{1-\epsilon_1-\epsilon_2}\gamma),$
so that $Q_0\in Z(\Delta'')$.
\end{example}

\begin{table}[ht]
\caption{$Q_0\in Z(\Delta)$ and $\cutlog(Q)\in\Cut(\tau_{l,\Delta})^0$}
\centering
\renewcommand{\arraystretch}{1.5}
 \renewcommand{\arraycolsep}{0.7mm}
 $\begin{array}{|c|c|c|c|c|}
\hline
&\Delta&\Gamma(\cO_{Z(\Delta)})&\Cut(\tau_{l,\Delta})
&\cutlog(Q)\\
\hline
\text{sp}_1&\langle -m_3\rangle&k&\langle 0,-f_1,-f_2\rangle
&\begin{matrix}-\epsilon_1f_1-\epsilon_2f_2\\
(0<\epsilon_i, 0<\epsilon_1+\epsilon_2<1)
\end{matrix}\\
\text{sp}_2&\langle m_1,-m_3\rangle^0&k[w_2^{\pm 1}]
&\langle 0,-f_1\rangle&
\begin{matrix}-\epsilon f_1\ (0<\epsilon<1)\end{matrix}\\
\text{sp}_3&\langle -m_3,-2m_3\rangle^0&k[(sw_3^{-1})^{\pm 1}]
&\langle -f_1,-f_2\rangle&
\begin{matrix}-\epsilon f_1-(1-\epsilon)f_2\\
 (0<\epsilon<1)\
\end{matrix}\\
\text{sp}_4&\langle -m_3,m_2\rangle^0&k[w_1^{\pm 1}]
&\langle 0,-f_2\rangle&\begin{matrix}-\epsilon f_2\ (0<\epsilon<1)
\end{matrix}\\
\text{sp}_5&\Sigma(0)&k[w_1^{\pm 1},w_2^{\pm 1}]&\langle 0\rangle&\begin{matrix}0\end{matrix}\\
\text{sp}_6&\Sigma(f_1)&k[w_2^{\pm 1},(sw_3^{-1})^{\pm 1}]&\langle -f_1\rangle&\begin{matrix}-f_1\end{matrix}\\ 
\text{sp}_7&\Sigma(f_2)&k[w_1^{\pm 1},(sw_3^{-1})^{\pm 1}]&\langle -f_2\rangle&\begin{matrix}-f_2\end{matrix}\\
\hline
\end{array}$
\label{table:Q and Q0}
\end{table} 

\section{The cohomology groups of $P$}
\label{sec:cohomology of Pl}
Let $R_{\init}$ be a CDVR and $S_{\init}:=\Spec R_{\init}$.  
We have  $(P_{l,l_0l},\cL_{l,l_0l})$ or $(P_{l},\cL_{l})$ 
 over $S_{\min}=\Spec R_{\min}$
in Theorem~\ref{thm:pd summary of Pl Gl} by assuming that 
$\Sigma^{\dagger}_{l_0}(0)$ is integral. Then $(P,\cN)$ is the descent 
of $(P_{l_0l},\cL_{l_0l})$ to $S_{\init}$. 
 In this section we study $H^q(P_{l},\cL_{l}^{\otimes m})$ and 
$H^q(P,\cN^{\otimes m})$\ $(m\geq 0)$.

\subsection{$H^q(P,\cO_P)$}
\label{subsec:Hq(P,OP)}
First we consider the case $m=0$.
\begin{thm}
\label{thm:cohom of P}Let $G$ be a semiabelian $S$-scheme, $\cG$ 
the N\'eron model of $G_{\eta}$ and 
$(P,i,\cN)$ be a cubical compactification of $\cG$ over $S$. Then 
\begin{enumerate}
\item\label{item:omegaPl0} $P_0$ is Gorenstein with dualizing sheaf 
$\omega_{P_0}$ trivial:
\item\label{item:omegaPl} 
$P$ is Gorenstein with dualizing sheaf $\omega_{P}$ trivial:
\item\label{item:hi_P0_cO} $h^q(P_0,\cO_{P_0})=\binom{g}{q}$;
\item\label{item:hi_P_cO} $H^q(P,\cO_{P})$ is an $R$-free module of rank $\binom{g}{q}$.
\end{enumerate}
\end{thm}
\begin{proof}First we assume $P=P_{l}$. We can prove 
$\omega_{P_{l,0}}\simeq
\cO_{P_{l,0}}$ 
in the same manner as \cite[4.2]{AN99}, 
which is (\ref{item:omegaPl0}). 
Since $P_{l}$ is Cohen-Macaulay 
by Theorem~\ref{thm:pd summary of Pl Gl}~(\ref{item:Pl}), 
we have a dualizing sheaf $\omega_{P_{l}}$ of $P_{l}$, so that  
$\omega_{P_{l}}(P_{l,0})\otimes_{\cO_{P_{l}}}
\cO_{P_{l,0}}\simeq\omega_{P_{l,0}}\simeq
\cO_{P_{l,0}}$. 
 Since $P_{l,0}$ is a Cartier divisor 
defined by the ideal $s\cO_{P_{l}}$, 
we have $\cO_{P_{l}}(P_{l,0})\simeq\cO_{P_{l}}$.  
Hence $\omega_{P_{l}}\simeq\omega_{P_{l}}(P_{l,0})\simeq\cO_{P_{l}}$ 
on a  Zariski open subset $U$ of $P_{l}$ containing $P_{l,0}$. 
Meanwhile $\omega_{P_{l,\eta}}
=\omega_{G_{\eta}}\simeq\cO_{G_{\eta}}=\cO_{P_{l,\eta}}$. 
Therefore $\pi_*(\omega_{P_{l}})$ is a 
torsion free $\cO_S$-module of rank one. It follows that  
$\Gamma(\pi_*(\omega_{P_{l}}))=R\ni 1$, which gives,  by the isomorphism 
$\omega_{P_{l}}\otimes_{\cO_{P_{l}}}\cO_{P_{l,0}}\simeq\cO_{P_{l,0}}$, 
 a nonvanishing section of $\omega_{P_{l}}$ and hence an isomorphism 
$\omega_{P_{l}}\simeq\cO_{P_{l}}$. 
This proves (\ref{item:omegaPl}). 
(\ref{item:hi_P0_cO}) is proved in 
the same manner as \cite[4.3]{AN99}. 
By \cite[Cor.~2, p.~121]{Mumford12},  
$H^q(P_{l,\eta},\cO_{P_{l,\eta}})=H^q(G_{\eta},\cO_{G_{\eta}})$ is 
a $k(\eta)$-module of rank $\binom{g}{q}$.  
Therefore (\ref{item:hi_P_cO}) follows from 
(\ref{item:hi_P0_cO}) and \cite[Cor.~2, p.~48]{Mumford12}. \par
Now we consider the general case. The process 
of constructing $P$ in \cite[\S\S~11-12]{MN24} 
and \S~\ref{sec:proof of main thm} 
of this article is as follows. First we construct 
the descent of $(P_{l_0l},\cL_{l_0l})$ by the action of $\mu$
in \cite[\S~11.3]{MN24}, which is $(P^*_{l_0l},\cL^*_{l_0l})$ 
over $S^*=S_{\spl}$. See \cite[11.13]{MN24}. Next we construct the Galois descent of $(P^*_{l_0l},\cL^*_{l_0l})$ by $\Gamma:=\Gal(k(\eta^*)/k(\eta))$ where $\eta^*=\eta_{\spl}$, which is $(P^{\flat}_{l_0l},\cL^{\flat}_{l_0l})$ denoted by $(P,\cN)$ here.  See \cite[12.9]{MN24}.

For $(P,\cN)$ there exists an open immersion $i:\cG\hookrightarrow P$ such that $i^*\cN$ is ample cubical and 
$\cN_{\eta}\simeq\cL_{\eta}^{\otimes 4N^2l_0l}$ by 
Theorem~\ref{thm:pd summary of Pl Gl}/\cite[12.9]{MN24}. Thus $(P,i,\cN)$ is a 
cubical compactification of $\cG$. 
Moreover 
\begin{enumerate}
\item[($\alpha$)] $(P^*_{l_0l},\cL^*_{l_0l})$ is the quotient of 
$(P_{l_0l},\cL_{l_0l})$ (over $R_{\min}$ with uniformizer $s$) by the group scheme $\mu$.  $(P^*_{l_0l},\cL^*_{l_0l})$ is obtained from the relatively complete model $(\tP^*_{l_0l},\tcL^*_{l_0l})$ over $S^*$ 
with uniformizer $t$, whence Theorem~\ref{thm:cohom of P} 
is true as well for $(P^*_{l_0l},\cL^*_{l_0l})$; 
\item[($\beta$)] $(P^{\flat}_{l_0l},\cL^{\flat}_{l_0l})$ is obtained from 
$(P^*_{l_0l},\cL^*_{l_0l})$  by Galois descent, where $\Gamma:=\Gal(S^*/S)$ and $S^*$ is finite \'etale (hence flat) over $S$, 
\end{enumerate}where we can set $l_0=1$ 
under the assumption (\ref{assump:Sigma_dagger_l(0) integral}). 

It remains to consider the case ($\beta$). Let 
$Z=P^*_{l}$. Note that by ($\alpha$), Theorem~\ref{thm:cohom of P} 
is true for $(P^*_{l},\cL^*_{l})$. 
 By Theorem~\ref{thm:pd summary of Pl Gl}/\cite[12.9]{MN24}, 
(\ref{item:omegaPl}) and ($\alpha$), 
$\omega_Z\simeq\cO_Z$, $Z=P\times_SS^*$ and 
$P=Z^{\Gamma}$. Since $S^*$ is flat over $S$, $\omega_Z\ (\simeq \cO_S)$ is $R$-torsion free. 
Since $P$ is projective over $S$, there is a closed immersion 
$\iota:P\to \bP_S$ into some projective space $\bP_S$ over $S$. 
The descent $(\beta)$ from $Z$ to $P$ shows that 
$\iota\times_SS^*:Z=P\times_SS^*\hookrightarrow\bP_{S^*}$ is 
a closed immersion.   

Since $\omega_S=\cO_S$ and $P$ is Cohen-Macaulay, 
the dualizing sheaf of $P$ is given by  
$\omega_P:={\cal Ext}^r_{\cO_{\bP_S}}(\cO_P,\omega_{\bP_S})$ 
(by identifying $P$ with $\iota(P)$) in view of 
\cite[7.2]{Hartshorne66} 
where $r=\codim_{\bP_S} P$. Since $S^*$ is flat over $S$, 
$\omega_Z:={\cal Ext}^r_{\cO_{\bP_{S^*}}}(\cO_{P_{S^*}},\omega_{\bP_{S^*}})
=\omega_P\otimes_RR^*$ and $\omega_P=(\omega_P\otimes_RR^*)^{\Gamma}=
(\omega_Z)^{\Gamma}$. It follows that $\Gamma(P,\omega_P)$ is 
$R$-torsion free, and hence $R$-free.
Since $P_{\eta}$ is an abelian variety, 
$\Gamma(P_{\eta},\omega_{P_{\eta}})\simeq k(\eta)$. 
Hence $\Gamma(P,\omega_P)\simeq R$, whose generator we denote by $\tau$. 
Since $\omega_Z=\omega_P\otimes_RR^*$, $R^*\simeq\Gamma(Z,\omega_Z)=\Gamma(P,\omega_P)\otimes_RR^*$ is generated by the pullback of $\tau$, which is nonvanishing everywhere on $Z$, and hence $\tau$ is nonvanishing everywhere on $P$. 
Hence $\omega_P=\cO_P$. 
 This proves (\ref{item:omegaPl}). 
Since $P_0$ is a Cartier divisor defined by $t=0$, we have 
$\omega_{P_0}=\omega_P(P_0)\otimes_{\cO_{P}}\cO_{P_0}=\cO_{P_0}$. 
This proves (\ref{item:omegaPl0}). 
Since $S^*$ is finite flat over $S$,  
by \cite[Cor.~5, p.~51]{Mumford12}, 
\begin{gather*}
H^q(Z,\cO_Z)=H^q(P\times_SS^*,\cO_{P\times_SS^*})=
H^q(P,\cO_P)\otimes_RR^*.\end{gather*}
Hence  $H^q(P,\cO_P)=H^q(Z,\cO_Z)^{\Gamma}$ 
is a free $R$-submodule of a free $R$-module $H^q(Z,\cO_Z)$.
Hence $\rank_RH^q(P,\cO_P)
=\rank_{k(\eta)}H^q(P_{\eta},\cO_{P_{\eta}})=\binom{g}{q}$. 
This proves (\ref{item:hi_P_cO}), 
from which (\ref{item:hi_P0_cO}) follows.   
\end{proof}

\subsection{The closed subscheme $\tZ_{l}^{\dagger}(\Delta)$ of $\tP_{l,0}$}
\label{subsec:closed subsch barZ_dagger_Delta}
Let $\Delta\in\Vor^{\dagger}_{l}$ 
and let $W_{l,\Delta}$ be as in Definition~\ref{defn:Bla,BlDelta}. 
Let $Z^{\dagger}_{l}(\Delta)$ be the unique reduced closed $G_0$-orbit 
of $(W_{l,\Delta})_0$ and $\tZ^{\dagger}_{l}(\Delta)$ 
the closure of $Z^{\dagger}_{l}(\Delta)$ in 
$\tP_{l,0}$ (with reduced structure) and 
$\tcL_{l}(\Delta):=\tcL_{l}\otimes_{\cO_{\tP_{l,0}}}
\cO_{\tZ^{\dagger}_{l}(\Delta)}$. Then  $\tcL_{l}(\Delta)$ is ample on $\tZ^{\dagger}_{l}(\Delta)$ by 
Theorem~\ref{thm:cohomology of Zl_Delta}~(\ref{item:very ample if m_geq_r}). 
We denote the restriction of 
the natural morphism $\tpi_{l}:\tP_{l,0}\to A_0$ 
to $\tZ^{\dagger}_{l}(\Delta)$ by 
$\tpi_{l,\Delta}:\tZ^{\dagger}_{l}(\Delta)\to A_0$. 
Each $\tpi_{l,\Delta}$ is a {\it fiber bundle} over $A_0$ 
by Theorem~\ref{thm:W_l,alpha,0 general}, 
whose fiber is a torus embedding associated 
with $\Fan_{2Nl}(\Delta)$ by Remark~\ref{rem:identify E_dagger_l with E_2Nl} and Lemma~\ref{lemma:Z_Delta=torus emb Fan_Delta}. 
See Definition~\ref{defn:fan Fan(Delta)}. 
The $G_0$-orbit $Z^{\dagger}_{l}(\Sigma^{\dagger}_{l}(0))$ 
is a $T_0$-bundle over $A_0$ 
associated with extension class $c\times_S0$ 
by \S~\ref{subsec:Raynaud extensions split case}. 
Since $\Delta\in\Vor^{\dagger}_{l}=\Vor_{2Nl}$, 
$Z_{2Nl}(\Delta)$ is the pushforward of 
the $T_0$-bundle $G_0$ over $A_0$ by the inclusion 
$T_{X(\Delta),k(0)}\hookrightarrow T_0=T_{X,k(0)}$, which is compactified by 
$\tZ_{2Nl}(\Delta)$. 
  Therefore $\tZ^{\dagger}_{l}(\Delta)$ is 
a $\tZ_{2Nl}(\Delta)$-bundle over $A_0$, 
which  compactifies $Z^{\dagger}_{l}(\Delta)$ $G_0$-equivariantly. 
By Definition~\ref{defn:split obj zeta in DDample}~(\ref{item:tcL and cM}), 
let $\pi:\tG\to A$ be the natural morphism of the Raynaud extension,
$\tcL=\pi^*\cM$ and $\pi_*(\tcL)=\bigoplus_{x\in X}\cM_x$. 
We denote $\cM_x\otimes_{\cO_A}\cO_{A_0}$ by $\cM_{x,0}$.

\begin{lemma}\label{lemma:cohomology of tZ_dagger_l_Delta}
Let $m\geq 1$, $\Delta\in\Sk(\Sigma^{\dagger}_{l}(v))$ 
and $x\in m\Delta\cap X$. 
Then 
\begin{enumerate}
\item\label{item:vanishing tZ_dagger_cL_Delta} 
$H^p(\tZ^{\dagger}_{l}(\Delta),\tcL_{l}(\Delta)^{\otimes m})=0$ $(p>0)$; 
\item\label{item:Gamma tZdagger cL_Delta} let 
$(\theta_{x,\lambda};\lambda\in\Lambda_x)$ 
be a $k(0)$-basis of $\Gamma(A_0,(\cM^{\otimes 4N^2lm})_{x,0})$; 
then \linebreak
$(\theta_{x,\lambda};x\in m\Delta\cap X,\lambda\in\Lambda_x)$ 
is a $k(0)$-basis of 
$\Gamma(\tZ^{\dagger}_{l}(\Delta),\tcL_{l}(\Delta)^{\otimes m})$, 
where $|\Lambda_x|=(4N^2lm)^g
\dim_{k(0)} \Gamma(A_0,\cM\otimes_{\cO_A}\cO_{A_0})$. 
\end{enumerate}
\end{lemma}
\begin{proof}Let $\pi:=\pi_{l,\Delta}$. 
Then $\pi:\tZ^{\dagger}_{l}(\Delta)\to A_0$ 
is a fiber bundle 
with $\tZ_{2Nl}(\Delta)$ fiber.  
By Theorem~\ref{thm:cohomology of Zl_Delta}
~(\ref{item:vanishing cL_Delta})/(\ref{item:Gamma Zl_Delta cL_Delta}), 
$R^q\pi_*\tcL_{l}(\Delta)^{\otimes m}=0$ $(q>0)$ and 
\begin{align*}
R^0(\pi_{A_0\cap U})_*\tcL_{l}(\Delta)^{\otimes m}&
=\bigoplus_{x\in m\Delta\cap X} (t^{v(x)}w^x_{A_0\cap U})
\cM^{\otimes 4N^2lm}_{A_0\cap U}\ \ (U\in\cF_B),
\end{align*}whence globally over $A_0$ we have an isomorphism 
by Eqs.~(\ref{eq:cO^x})-(\ref{eq:-Fx=c(x)=Ox})
\begin{equation}\label{eq:isom R0pi*tcL(Delta)}
R^0\pi_*\tcL_{l}(\Delta)^{\otimes m}\simeq 
\bigoplus_{x\in m\Delta\cap X}
(\cM^{\otimes 4N^2lm})_{x,0}.
\end{equation}
Since the following spectral sequence is degenerate:
$$E_{2}^{p,q}=H^p(A_0,R^q\pi_*\tcL_{l}(\Delta)^{\otimes m})
\Rightarrow H^{p+q}(\tZ^{\dagger}_{l}(\Delta),\tcL_{l}(\Delta)^{\otimes m})\ \ (m\geq 1),$$
we obtain $H^p(A_0,R^0\pi_*\tcL_{l}(\Delta)^{\otimes m})
\simeq H^{p}(\tZ^{\dagger}_{l}(\Delta),\tcL_{l}(\Delta)^{\otimes m})$.
By Eq.~(\ref{eq:isom R0pi*tcL(Delta)}) and \cite[\S~16]{Mumford12},
$H^p(\tZ^{\dagger}_{l}(\Delta),\tcL_{l}(\Delta)^{\otimes m})=0\ (p>0)
$ because $\cM^{\otimes 4N^2lm}\otimes_{\cO_{A_0}}\cO_x$ is ample $(\forall x\in X)$.  
 Moreover by Eq.~(\ref{eq:isom R0pi*tcL(Delta)}), 
\begin{align*}
\Gamma(\tZ^{\dagger}_{l}(\Delta),\tcL_{l}(\Delta)^{\otimes m})&\simeq 
\bigoplus_{x\in m\Delta\cap X}
\Gamma(A_0,(\cM^{\otimes 4N^2lm})_{x,0}).
\end{align*}   
\end{proof}

\subsection{The cohomology groups of $P_{l}$ and $P_{l,0}$}
\label{subsec:module of global sections}

\begin{defn}\label{defn:Map(Y,k)}
Let $Y$ be a lattice, and $\Map(Y,k(0))$ 
the $k(0)$-module of of all $k(0)$-valued functions on $Y$, $y\in Y$ and 
$\delta_y$ an element of $\Map(Y,k(0))$ such that $\delta_y(y)=1$ and 
$\delta_y(z)=0$ $(\forall z\in Y\setminus\{y\})$. 
Any $\theta\in\Map(Y,k(0))$ is 
a formal sum $\sum_{y\in Y}\theta(y)\delta_y$. 
We define an action $\rho:Y\times \Map(Y,k(0))\to \Map(Y,k(0))$ of $Y$ on $V$ 
by $\rho(y,\theta)(z)=\theta(z-y)$ $(y,z\in Y,\theta\in\Map(Y,k(0)))$. 
Let $Z$ be a subgroup of $Y$ with $[Y:Z]<\infty$, and $\Map(Y,k(0))^Z$ 
the submodule of $\Map(Y,k(0))$ of all $Z$-invariants.
\end{defn}  
\begin{lemma}\label{lemma:Hp_Y_Vmax=0} 
$H^p(Y,\Map(Y,k(0)))=0$ $(\forall p>0)$.
\end{lemma}
\begin{proof}This is proved in the same manner as    
\cite[1.3.7]{NSW99}. 
\end{proof}

\begin{defn}\label{defn:xi_m_lalpha,u}
Let 
$\alpha\in\Sigma^{\dagger}_{lm}$, $u\in X^{\vee}$ and 
$x=\alpha+4Nlm\mu(u)\in X$. Then we define 
$\xi^{(m)}_{l,\alpha,u}:=\epsilon^{\dagger}_{l}(u)^{m}
\tau^e(u,\alpha)w^{\alpha+4Nlm\mu(u)}$, which we also denote 
by $\xi^{(m)}_{l,x}$ if no confusion is possible. We note 
$\xi_{l,\alpha_i,u}=\xi^{(1)}_{l,\alpha_i,u}$.
\end{defn}

\begin{lemma}\label{lemma:theta_l_x_varthetax}
Let $\alpha\in\Sigma^{\dagger}_{lm}$, $u\in X^{\vee}$,
$x=\alpha+4Nlm\mu(u)\in X$ and 
$\vartheta_{x}\in\Gamma(A,(\cM^{\otimes 4N^2lm})_{x})$. 
We define 
{\small\begin{align}\label{eq:defn of theta_x_varthetax}
\Theta_{l,x,\vartheta_x}&
=\sum_{y\in Y}\epsilon^{\dagger}_{l}(u+\beta(y))^{m}
\tau^e(u+\beta(y),\alpha)T^*_{c^t(y)}(\vartheta_{x})
w^{\alpha+4Nlm\mu(u+\beta(y))}\theta^m_{l}.
\end{align}}
\hskip-0.15cm Then  $\Theta_{l,x,\vartheta_x}\in\Gamma(P^{\wedge}_{l},\cL_{l}^{\otimes lm,\wedge})$ and 
$\Theta_{l,x,\vartheta_x}\otimes_R k(0)\neq 0$ 
if $\vartheta_x\otimes_R k(0)\neq 0$. 
\end{lemma}
\begin{proof}Since $\alpha\in\Sigma^{\dagger}_{lm}$, there exist 
$\alpha_i$ $(i\in[1,m])$ such that $\alpha_i\in\Sigma^{\dagger}_{l}$ $(\forall i\in[1,m])$ such that $\alpha=\sum_{i=1}^m\alpha_i$. Let $x_i:=\alpha_i+4Nl\mu(u)$. By Definition~\ref{defn:delta_l_u on U} 
and Notation~\ref{notation:theta expansion in pd case}, we have  
$\xi_{l,\alpha_i,u}=\epsilon^{\dagger}_{l}(u)\tau^e(u,\alpha_i)
w^{\alpha_i+4Nl\mu(u)}$, so that 
\begin{align*}
&\tS_zS_z^*(\vartheta_{x}\prod_{i=1}^m\xi_{l,\alpha_i,u}\theta_{l})
=T^*_{c^t(z)}(\vartheta_{x})\prod_{i=1}^m
\xi_{l,\alpha_i,u+\beta(z)}\theta_{l}
=T^*_{c^t(z)}(\vartheta_{x})\xi^{(m)}_{l,\alpha,u+\beta(z)}\theta_{l}^m\\
&=T^*_{c^t(z)}(\vartheta_{x})
\epsilon^{\dagger}_{l}(u+\beta(z))^m\tau^e(u+\beta(z),\alpha)
w^{\alpha+4Nlm\mu(u+\beta(z))}\theta_{l}^m. 
\end{align*} Hence $\Theta_{l,x,\vartheta_x}
=\sum_{z\in Y}\tS_zS_z^*(\vartheta_x\prod_{i=1}^m\xi_{l,\alpha_i,u}\theta_{l})$.  It follows  
$\tS_yS_y^*(\Theta_{l,x,\vartheta_x})
=\Theta_{l,x,\vartheta_x}$ $(\forall y\in Y)$.  
It is clear that $\Theta_{l,x,\vartheta_x}/\xi_{l,\gamma,v}^m\in A_{l,\gamma,v}^{\wedge}$ $(\forall \gamma\in\Sigma^{\dagger}_{l}, \forall v\in X^{\vee})$, whence $\Theta_{l,x,\vartheta_x}\in\Gamma(P^{\wedge}_{l},\cL_{l}^{\otimes lm,\wedge})=\Gamma(P_{l},\cL_{l}^{\otimes lm})$. 
Since $v_t(\tau^e(u,x))=u(x)$,  
$$(W_{l,0,u})_0=\Spec \cO_{A_0}[t^{u(x)}w^x;x\in X]\simeq 
\Spec \cO_{A_0}[\tau^e(u,x)w^x;x\in X]. 
$$ 

The restriction of $\Theta_{l,x,\vartheta_x}$ 
to $(W_{l,0,u})_0$ is given by:
\begin{align*}
&\Theta_{l,x,\vartheta_x}/\xi_{l,0,u}^m\theta_{l}^m
=\sum_{z\in Y}T^*_{c^t(z)}(\vartheta_x)\xi^{(m)}_{l,\alpha,u+\beta(z)}/\xi_{l,0,u}^m\\
&=\sum_{y\in Y}\epsilon^{\dagger}_{l}(u+\beta(y))^{m}
\tau^e(u+\beta(y),\alpha)
T^*_{c^t(y)}(\vartheta_{x})
w^{\alpha+4Nlm\mu(u+\beta(y))}\theta^m_{l}. 
\end{align*}
 Its initial term is   
$\vartheta_x\tau^e(u,\alpha)w^{\alpha}\neq 0$ (for $y=0$). This proves Lemma. 
\end{proof}

\begin{thm}\label{thm:cohomology of Pl Pl0}
Let $g=\dim_{k(0)}P_{l,0}$. Then  
\begin{enumerate}
\item\label{item:vanishing Hq=0}
 $H^q(P_{l},\cL^{\otimes m}_{l})
=H^q(P_{l,0},\cL^{\otimes m}_{l,0})=0$ $(\forall q>0, \forall m>0)$;
\item\label{item:H0_Pl0}
$\Gamma(P_{l,0},\cL^{\otimes m}_{l,0})
=\Gamma(P_{l},\cL^{\otimes m}_{l})\otimes_Rk(0)$;
\item\label{item:rank H0_Pl}
$\Gamma(P_{l},\cL^{\otimes m}_{l})$ is a free $R$-module of rank 
$(4N^2lm)^g|X/Y|\rank_R \Gamma(A,\cM)$; 
\item\label{item:H0_Pl} 
$\Gamma(P_{l},\cL^{\otimes m}_{l})$ is 
identified via $I$-adic completion 
with the finite free $R$-module spanned by $\Theta_{l,x,\vartheta_x}$ $(x=\alpha+4Nlm\mu(u), \alpha\in\Sigma^{\dagger}_{lm}, u\in X^{\vee}, \vartheta_x\in\Gamma(A,(\cM^{\otimes 4N^2lm})_{x}))$. 
\end{enumerate}
\end{thm}
\begin{proof} See \cite[4.4]{AN99}/\cite[3.9/4.10]{Nakamura99}. 
To explicitly know the $k(0)$-generators 
of $\Gamma(P_{l,0},\cL^{\otimes m}_{l,0})$ $(m\geq 1)$, we basically 
follow the proof of \cite[3.9/4.10]{Nakamura99}. 
Let $\varpi_{l}:\tP_{l,0}\to P_{l,0}$ be the natural morphism. Note that  
$\tcL_{l,0}=\varpi_{l}^*\cL_{l,0}$. 
It is the quotient of $\tP_{l,0}$ by a free 
and proper action of $Y$. \footnote{
In contrast with \cite{Nakamura99}, here we do not assume 
that $S_y(\tZ^{\dagger}_{l}(\sigma))\cap\tZ^{\dagger}_{l}(\sigma))
=\emptyset$ $(\forall y\neq 0)$.} 

First we prove $H^q(\tP_{l,0},\tcL_{l,0}^{\otimes m})=0$ $(q>0)$. 
Let $r=\rank X$. 
Let  $\Sk^i(\Vor^{\dagger}_{l})$ be the set of all $i$-dimensional closed 
polytopes of $\Vor^{\dagger}_{l}$ 
and $\sigma\in\Vor^{\dagger}_{l}$. 
Let $\varpi_{l,\sigma}:=(\varpi_{l})_{|\tZ^{\dagger}_{l}(\sigma)}$. 
We define 
\begin{gather*}
\cS:=\tcL_{l}^{\otimes m}\otimes_{\cO_{\tP_{l}}}{\cO_{\tP_{l,0}}},\\ 
\cS(\sigma):=\cS\otimes_{\cO_{P_{l,0}}}\cO_{\tZ^{\dagger}_{l}(\sigma)},\ 
\cS^q(\sigma):=R^q(\varpi_{l,\sigma})_*\cS(\sigma).
\end{gather*}

By Theorem~\ref{thm:W_l,alpha,0 general}, 
there is an 
exact sequence of $\cO_{P_{l,0}}$-modules:
\begin{equation}\label{eq:exact seq OPl0} 
0 \to \cO_{\tP_{l,0}}\to
\bigoplus_{\sigma_r} \cO_{\tZ^{\dagger}_{l}(\sigma_r)}
\overset{\partial_{r}}\to \cdots 
\overset{\partial_{2}}\to 
\bigoplus_{\sigma_1} \cO_{\tZ^{\dagger}_{l}(\sigma_1)}
\overset{\partial_1}\to 
\bigoplus_{\sigma_0} \cO_{\tZ^{\dagger}_{l}(\sigma_0)}\to 0
\end{equation}
where the direct sum is taken over all 
$\sigma_i\in \Sk^i(\Vor^{\dagger}_{l})$. The homomorphism 
$\partial_i : \bigoplus \cO_{\tZ^{\dagger}_{l}(\sigma_i)} 
\to \bigoplus \cO_{\tZ^{\dagger}_{l}(\sigma_{i-1})}$ is 
defined by 
$\partial_i(\bigoplus_\sigma \phi_{\sigma})=
\bigoplus_\tau 
\sum_{\tau\subset\sigma} [\sigma:\tau]\phi_{\sigma}$
where the sum $\sum_{\tau\subset\sigma}$ is taken
over the set of all $\sigma\in\Vor^{\dagger}_{l}$ 
containing $\tau$ as a face of codimension one, 
while the direct sum $\bigoplus_{\tau}$ is taken over 
the set $\Sk^{i-1}(\Vor^{\dagger}_{l})$.  
By tensoring Eq.~(\ref{eq:exact seq OPl0}) 
with $\tcL^{\otimes m}_{l}$, we have an exact sequence with $Y$-action 
 $(S_y;y\in Y)$:
\begin{equation}\label{eq:exact seq cS}
0 \to \tcL^{\otimes m}_{l,0}\to
\bigoplus_{\sigma_r} \cS(\sigma_r)\overset{\partial_{r}}\to \cdots 
\overset{\partial_{2}}\to 
\bigoplus_{\sigma_1} \cS(\sigma_1)\overset{\partial_1}\to 
\bigoplus_{\sigma_0} \cS(\sigma_0)\to 0.
\end{equation}by abuse of the notation $\partial_{\bullet}$.  

Let 
$\wtW_i:=\bigoplus_{\sigma_i\in \Sk^i(\Vor_{l}^{\dagger})}
\Gamma(\tZ^{\dagger}_{l}(\sigma_i),\cS(\sigma_i)).$  
Let $\widetilde\bW$ be a complex: 
$$0\to \wtW_r\overset{\partial_{r-1}}{\to}\wtW_{r-1}
\overset{\partial_{r-2}}{\to} \cdots 
\overset{\partial_{1}}{\to} \wtW_1\overset{\partial_{0}}{\to} 
\wtW_0\to 0. 
$$

By Eq.~(\ref{eq:exact seq cS}), we have a spectral sequence 
$$E_2^{p,q}=H^p(\bigoplus_{\sigma_{\bullet}}H^q(\tZ^{\dagger}_{l}(\sigma_{\bullet}),\cS(\sigma_{\bullet})),\partial_{\bullet})\Rightarrow 
H^{p+q}(\tP_{l,0},\tcL^{\otimes m}_{l,0}).$$ 
This spectral sequence degenerates because  
$H^q(\tZ^{\dagger}_{l}(\sigma),\cS(\sigma))=0$ $(\forall \sigma\in\Vor^{\dagger}_{l},\forall q\geq 1)$ by Lemma~\ref{lemma:cohomology of tZ_dagger_l_Delta}
~(\ref{item:vanishing tZ_dagger_cL_Delta}). 
Hence $H^p(\tP_{l,0},\cL^{\otimes m}_{l,0})
\simeq H^p(\widetilde\bW)\ (\forall p)$.

Let $\Star(x)$ be the union of all $\sigma\in\Vor^{\dagger}_{l}$ 
such that $x\in\sigma$ and $\Star(x)^0$ the relative interior of $\Star(x)$. 
Then $H^p(\widetilde\bW)$ is given by 
$$\bigoplus_{x\in X}H^p(\Star(\frac{x}{m})^0,k(0)))\otimes_{k(0)}H^0((\cM^{\otimes 4N^2lm})_{x,0}).$$

Since $\Star(\frac{x}{m})^0$ is contractible and connected, 
$H^p(\Star(\frac{x}{m})^0,k(0)))=0$ $(\forall p>0, \forall x\in X)$, so that 
$H^p(\widetilde\bW)=0$\ $(\forall p>0)$. 
By \cite[Appendix to \S~2, pp.~21-23]{Mumford12}, we obtain 
\begin{equation}\label{eq:isomorphism theorem}
H^p(P_{l,0},\cL_{l,0}^{\otimes m})\simeq 
H^p(Y,\Gamma(\tP_{l,0},\tcL_{l,0}^{\otimes m}))
\simeq H^p(Y,\Gamma(\widetilde{\bW})),
\end{equation}where 
$\Gamma(\tP_{l,0},\tcL_{l,0}^{\otimes m})=\Gamma(\widetilde{\bW})
=\ker(\partial_{r-1}:\wtW_r\to\wtW_{r-1})\subset\wtW_r$. 

Let $\tZ^{\dagger}_{l}(v):=\tZ^{\dagger}_{l}(\Sigma_{l}(v))$ and 
$\tcL^{\dagger}_{l}(v):=\tcL^{\dagger}_{l}(\Sigma_{l}(v))$..  
By Lemma~10.2, $\Gamma(\tZ^{\dagger}_{l}(v),\tcL^{\dagger}_{l}(v))$ 
is spanned by 
$(\theta_{x,\lambda_x};x\in m\Sigma_{l}(v), \lambda_x\in\Lambda_x)$.
Hence we have 
\begin{equation*}
\Gamma(\widetilde{\bW})\simeq\bigoplus_{x\in X}
\Gamma(A_0,(\cM^{\otimes 4N^2lm})_{x,0}).
\end{equation*} 
Let $\Delta$ be a subset of $X$ of all  
representatives of $X/4N^2lm\phi(Y)$. We set 
\begin{gather*}
\Gamma(\widetilde{W})=\bigoplus_{\alpha\in\Delta}
\Gamma(\widetilde{W})(\alpha),\\
\Gamma(\widetilde{W})(\alpha)
:=\bigoplus_{x\in \alpha+4N^2lm\phi(Y)}
\Gamma(A_0,(\cM^{\otimes 4N^2lm})_{x,0}). 
\end{gather*}
For each $\alpha\in\Delta$, 
there is a $Y$-isomorphism: 
$$\Map(Y,k(0))\otimes_{k(0)}\Gamma(A_0,(\cM^{\otimes 4N^2lm})_{\alpha}\otimes_{\cO_A}\cO_{A_0})\simeq\Gamma(\widetilde{W})(\alpha)$$
sending 
$\delta_y\otimes\lambda_{\alpha}\mapsto\theta_{\alpha+4N^2lm\phi(y),T^*_{c^t(y)}\lambda_{\alpha}}\theta_{l}^m$. By Lemma~\ref{lemma:Hp_Y_Vmax=0}
we obtain 
\begin{gather*}H^p(P_{l,0},\cL_{l,0}^{\otimes m})
=H^p(Y,\Gamma(\widetilde\bW))=\bigoplus_{\alpha\in\Delta}H^p(Y,\Gamma(\widetilde{W})(\alpha))=0\  
(\forall p>0),\\
\Gamma(P_{l,0},\cL_{l,0}^{\otimes m})
=\Gamma(\widetilde\bW)^{Y\op{-inv}}
\simeq\bigoplus_{\alpha\in\Delta}\Gamma(\widetilde\bW(\alpha))^{Y\op{-inv}}. 
\end{gather*}This proves (\ref{item:vanishing Hq=0}).  
By (\ref{item:vanishing Hq=0}) and \cite[Cor.~3, p.~50]{Mumford12}, 
we have 
$$\Gamma(P_{l,0},\cL^{\otimes m}_{l,0})
=\Gamma(P_{l},\cL^{\otimes m}_{l})\otimes_Rk(0).$$ 
This proves (\ref{item:H0_Pl0}). 
Since $P_{l}$ is flat and proper over $S$, 
$h^0(P_{l,0},\cL^{\otimes m}_{l,0})$ equals 
\begin{equation}\label{h0 Pl0 cL_lm_0}
\chi(P_{l,0},\cL^{\otimes m}_{l,0})
=\chi(P_{l,\eta},\cL^{\otimes m}_{l,\eta})
=(4N^2lm)^g|X/Y|\rank\Gamma(A,\cM).
\end{equation} 
Since $\Gamma(P_{l},\cL^{\otimes m}_{l})$ is $R$-torsion free, 
it is $R$-free, which  proves (\ref{item:rank H0_Pl}).  

Let $\alpha\in\Delta$. If $\vartheta_{\alpha}\otimes_Rk(0)\neq 0$, 
$\Theta_{l,\alpha,\vartheta_{\alpha}}\otimes_Rk(0)$ 
is {\it nonzero} by 
Lemma~\ref{lemma:theta_l_x_varthetax}, 
which is  a $Y$-invariant element of $\Gamma(\widetilde{W}(\alpha))$. 
This proves (\ref{item:H0_Pl}) by (\ref{item:H0_Pl0}) and Nakayama's lemma. 
\end{proof}

\subsection{The cohomology of $\barZ^{\dagger}_{l}(\Delta)$}
\label{subsec:closed subscm barZlDelta}

Let $\Delta\in\Vor^{\dagger}_{l}$ and let 
$\varpi_{l}:\tP_{l,0}\to P_{l,0}$ be the natural morphism. 
Then $\varpi_{l}(Z^{\dagger}_{l}(\Delta))\simeq Z^{\dagger}_{l}(\Delta)$.
Let $\barZ^{\dagger}_{l}(\Delta)$ be the closure 
of $\varpi_{l}(Z^{\dagger}_{l}(\Delta))$ in $P_{l,0}$ and 
$\cL_{l}(\Delta):=\cL_{l}\otimes_{\cO_{P_{l}}}\cO_{\barZ^{\dagger}_{l}(\Delta)}$.  

\begin{lemma}\label{lemma:cohomology of barZldagger_Delta}
Let $r=\dim_{k(0)}\barZ^{\dagger}_{l}(\Delta)$. Then  
\begin{enumerate}
\item\label{item:vanishing Hq(cLdagger_Delta)=0}
 $H^q(\barZ^{\dagger}_{l}(\Delta),\cL^{\otimes m}_{l}(\Delta))
=0$ $(\forall q>0, \forall m>0)$;
\item\label{item:H0_cLdagger_Delta} 
$\Gamma(\barZ^{\dagger}_{l}(\Delta),\cL^{\otimes m}_{l}(\Delta))$ is 
spanned by the restriction $\Theta^{\Delta}_{l,x,\vartheta_x}$  
of $\Theta_{l,x,\vartheta_x}$ to $\barZ^{\dagger}_{l}(\Delta)$ 
$(x\in m\Delta\cap X, 
\vartheta_x\in\Gamma(A_0,(\cM^{\otimes 4N^2lm})_{x,0}))$, 
where  
\begin{equation*}
\Theta^{\Delta}_{l,x,\vartheta_x}:=\sum_{\begin{matrix}
y\in Y\\
 x+4N^2lm\phi(y)\in m\Delta\cap X
\end{matrix}}T^*_{c^t(y)}(\vartheta_x)\xi^{(m)}_{l,x+4N^2lm\phi(y)}\theta_{l}^m;
\end{equation*}
\item\label{item:x in mDelta0}$\Theta^{\Delta}_{l,x,\vartheta_x}=
\vartheta_x\xi^{(m)}_{l,x}\theta_{l}^m$ if $x\in m\Delta^0\cap X$.
\end{enumerate}
\end{lemma}
\begin{proof}
Let $\barZ^{\dagger,\su}_{l}:=\barZ^{\dagger}_{l}(\Delta)\times_{P_{l,0}}\tP_{l,0}$ and $\tcL^{\su,\otimes m}_{l}:=\tcL^{\otimes m}_{l}
\otimes_{\cO_{\tP_{l}}}\cO_{\barZ^{\dagger,\su}_{l}}$.  
Note that $Y$ acts freely on $\barZ^{\dagger,\su}_{l}$, so that 
$\barZ^{\dagger}_{l}(\Delta)\simeq\barZ^{\dagger,\su}_{l}/Y$. 
We have an exact sequence: 
\begin{equation*}
0 \to \tcL^{\su,\otimes m}_{l}\to
\bigoplus_{\sigma_r\subset\Delta} \cS(\sigma_r)
\overset{\partial_{r}}\to \cdots 
\overset{\partial_{2}}\to 
\bigoplus_{\sigma_1\subset\Delta} \cS(\sigma_1)\overset{\partial_1}\to 
\bigoplus_{\sigma_0\subset\Delta} \cS(\sigma_0)\to 0.
\end{equation*}which is obtained in the same manner as Eq.~(\ref{eq:exact seq cS}).  

By Lemma~\ref{lemma:cohomology of tZ_dagger_l_Delta}~(\ref{item:vanishing tZ_dagger_cL_Delta}), $H^q(\tZ^{\dagger}_{l}(\Delta),\tcL^{\otimes m}_{l}(\Delta))
=0$, so that $H^p(\barZ^{\dagger,\su}_{l},\tcL^{\su,\otimes m}_{l})\simeq 
H^p(\widetilde\bW(\Delta))$ $(\forall p)$ where $\widetilde\bW(\Delta)$ is the complex 
$$0\to \wtW_r(\Delta)\overset{\partial_{r-1}}{\to}\wtW_{r-1}(\Delta)
\overset{\partial_{r-2}}{\to} \cdots 
\overset{\partial_{1}}{\to} \wtW_1(\Delta)\overset{\partial_{0}}{\to} 
\wtW_0(\Delta)\to 0
$$
with $\wtW_i(\Delta):=\bigoplus_{\sigma_i\in \Sk^i(\Delta+4N^2lm\phi(Y))}
\Gamma(\tZ^{\dagger}_{l}(\sigma_i),\cS(\sigma_i))$. 
In the same manner as Theorem~\ref{thm:cohomology of Pl Pl0}, 
$H^p(\widetilde\bW(\Delta))$ is given by 
$$\bigoplus_{x\in m\Delta\cap X+4N^2lm\phi(Y)}H^p(\Star(\frac{x}{m})^0,k(0)))\otimes_{k(0)}H^0((\cM^{\otimes 4N^2lm})_{x,0}).$$
Hence $H^p(\tZ^{\dagger,\su}_{l},\tcL^{\su,\otimes m})\simeq H^p(\widetilde\bW(\Delta))=0$\ $(p>0)$. By \cite[Appendix to \S~2, pp.~21-23]{Mumford12}, 
we obtain in the same manner as Theorem~\ref{thm:cohomology of Pl Pl0},
\begin{equation*}
\begin{aligned}H^p(\barZ^{\dagger}_{l}(\Delta),\cL_{l}^{\otimes m}(\Delta))
&\simeq H^p(Y,\Gamma(\tZ^{\dagger,\su}_{l},\tcL^{\su,\otimes m}_{l}))
\simeq H^p(Y,\Gamma(\widetilde{\bW}(\Delta))). 
\end{aligned}
\end{equation*}  This proves (\ref{item:vanishing Hq(cLdagger_Delta)=0}) 
because  $H^p(Y,\Gamma(\widetilde{\bW}(\Delta)))=0$\ 
$(p>0)$ by Lemma~\ref{lemma:Hp_Y_Vmax=0}. Moreover  
\begin{align*}
\Gamma(\barZ^{\dagger}_{l}(\Delta),\cL_{l}^{\otimes m}(\Delta))
&\simeq\left(\bigoplus_{y\in m\Delta+4N^2lm\phi(Y)}\Gamma(A_0,(\cM^{\otimes 4N^2lm})_{y,0})\right)^{Y\op{-inv}}.
\end{align*}
This proves (\ref{item:H0_cLdagger_Delta})  
by Theorem~\ref{thm:cohomology of Zl_Delta}
~(\ref{item:Gamma Zl_Delta cL_Delta}). If $x\in m\Delta^0\cap X$, 
(\ref{item:x in mDelta0}) follows from  
$(x+4N^2lm\phi(Y))\cap m\Delta^0=\{x\}$. 
\end{proof}

\begin{rem}\label{rem:xth coordinates}
Let $\alpha_0\in\Sigma^{\dagger}_{l}$ and 
$x_0:=\alpha_0+4Nlm\mu(u)\in m\Delta^0\cap X$. 
Then there exists $\alpha_{0,i}$ $(i\in[1,m])$ such that 
$\alpha_0=\sum_{i=1}^m\alpha_{0,i}$. For $U\in \cF_B$, 
we define $W_{l,\alpha_0,u,U}:=\bigcap_{i=1}^mW_{l,\alpha_{0,i},u,U}$ by Definition~\ref{defn:U_alpha_u_U W_alpha_u_U}.  
If $\alpha\in\Sigma^{\dagger}_{lm}$ 
and $x=\alpha+4Nlm\mu(u)\in m\Delta^0\cap X$ 
then by (applying the proof of) Theorem~\ref{thm:W_l,alpha,0 general}, 
we see that 
$\tau^e(u,x)w^x$ is the $x$-th coordinate of an irreducible 
component $\barZ_{l}(\Sigma_{l}(u))$ of $(W_{l,\alpha_0,u,U})_0$
containing $\barZ_{l}(\Delta)$. 
By Lemma~\ref{lemma:cohomology of barZldagger_Delta}~(\ref{item:x in mDelta0}),
$$\Theta^{\Delta}_{l,x,\vartheta_x}/\Theta^{\Delta}_{l,x_0,\vartheta_{x_0}}
=\tau^e(u,x)w^{x}/\tau^e(u,x_0)w^{x_0}
=\tau^e(u,\alpha)w^{\alpha}/\tau^e(u,\alpha_0)w^{\alpha_0}.$$  
\end{rem}

\subsection{Projective embedding of $P_{l}$ by $\cL^{\otimes m}_{l}$}

\begin{thm}$\cL^{\otimes m}_{l,0}$ is very ample on $P_{l,0}$ 
if $m\geq 4g$.
\end{thm}
\begin{proof}Since 
$\varpi_{l}(Z^{\dagger}_{l}(\Delta))\simeq 
Z^{\dagger}_{l}(\Delta)$, we identify them for brevity. 
Let $r=\rank X\ (\leq g)$. 
$\Gamma(P_{l,0},\cL^{\otimes m}_{l,0})$ is spanned by $\Theta_{l,x,\vartheta_x}\otimes_R k(0)$ $(x\in X)$ by    
 Theorem~\ref{thm:cohomology of Pl Pl0}~(\ref{item:H0_Pl}) when 
$\vartheta_x$ ranges over $\Gamma(A_0,(\cM^{\otimes 4N^2lm})_{x,0})$. 
First we prove that the linear system 
$|\cL^{\otimes m}_{l,0}|$ is base point free if $m\geq r+1$.  
Let $p$ be any closed point of $P_{l,0}$. Then there exists $\Delta\in\Vor^{\dagger}_{l}$ such that $p\in Z^{\dagger}_{l}(\Delta)(\overline{k(0)})$. 
 Since $4N^2lm\geq 3$, there exists $\vartheta_x\in\Gamma(A_0,(\cM^{\otimes 4N^2lm})_{x,0})$ such that 
$\vartheta_x(\pi_{l,\Delta}(p))\neq 0$.  If $m\geq r+1$, 
then by Lemma~\ref{lemma:SmDelta_alpha generate C}~(\ref{item:m_geq r+1}), 
there exists $x\in (m\Delta)^0\cap X$, which is written as 
$x=\alpha+4Nlm\mu(u)$ for some $\alpha\in\Sigma^{\dagger}_{lm}$ 
and $u\in X^{\vee}$.
 Hence by 
Lemma~\ref{lemma:cohomology of barZldagger_Delta}~(\ref{item:x in mDelta0}), 
\begin{equation}\label{eq:single term theta_lxvarhetax}
\Theta^{\Delta}_{l,x,\vartheta_x}=\vartheta_x\xi^{(m)}_{l,x}\theta_{l}^m\in\Gamma(\barZ^{\dagger}_{l}(\Delta),\cL^{\otimes m}_{l}(\Delta))
\end{equation} 
is nonvanishing at $p$. Hence 
$|\cL^{\otimes m}_{l,0}|$ is base point free if $m\geq r+1$. \par

Let $\phi:P_{l,0}\to\bP(\Gamma(\cL^{\otimes m}_{l,0}))$ be 
the morphism associated with $|\cL^{\otimes m}_{l,0}|$. 
We shall prove that $\phi$ separates two distinct closed 
points of $P_{l,0}$. 
Let $p\in Z^{\dagger}_{l}(\Delta)(\overline{k(0)})$ 
and $p'\in Z^{\dagger}_{l}(\Delta')(\overline{k(0)})$. 
If $\Delta\neq\Delta'$, 
then $(m\Delta)^0\neq(m\Delta')^0$, and 
by Eq.~(\ref{eq:single term theta_lxvarhetax}), 
$\Theta_{l,x,\vartheta_x}(p)=0$  and $\Theta_{l,x,\vartheta_x}(p')$ 
is nonvanishing for $x\in(m\Delta')^0\cap X$. It follows 
$\phi(Z^{\dagger}_{l}(\Delta))\cap\phi(Z^{\dagger}_{l}(\Delta'))=\emptyset$. 
Hence $\phi(p)\neq \phi(p')$. 
If $\Delta=\Delta'$ and $\pi_{l,\Delta}(p)\neq \pi_{l,\Delta}(p')$, then 
there exists $\vartheta_x\in \Gamma(A_0,(\cM^{\otimes 4N^2lm})_{x,0})$ 
such that $x\in m\Delta\cap X$, $\vartheta_x(\pi_{l,\Delta}(p))=0$ and 
$\vartheta_x(\pi_{l,\Delta}(p'))\neq 0$ because $4N^2lm\geq 3$.
Hence $\phi$ separates $p$ and $p'$. 
Finally if $p\neq p'$, $\Delta=\Delta'$ 
and $\pi_{l,\Delta}(p)=\pi_{l,\Delta}(p')$, 
then $p\neq p'$ in the fiber $Z_{2Nl}(\Delta)$ of $\pi_{l,\Delta}
:Z^{\dagger}_{l}(\Delta)\to A_0$. 
Therefore there exists $x\in X(\Delta)$ such that 
the $x$-th coordinates $w^x$  
of the fiber $Z_{2Nl}(\Delta)$ at $p$ and $p'$ are not equal. 
This means by Remark~\ref{rem:xth coordinates} that 
$\tau^e(u,x)w^x(p)\neq\tau^e(u,x)w^x(p')$.
By Lemma~\ref{lemma:SmDelta_alpha generate C}~
(\ref{item:generated by mDelta-malpha}), 
we can find $y_i,z_i\in (m\Delta)^0\cap X$ $(i\in[1,n])$ 
for some $n\in\bN$ such that 
$x=\sum_{i=1}^n(y_i-z_i)$. 
Then there exists at least one $i$ such that 
$\tau^e(u,y_i-z_i)w^{y_i-z_i}(p)\neq\tau^e(u,y_i-z_i)w^{y_i-z_i}(p')$.  
Since $y_i,z_i\in (m\Delta)^0\cap X$, 
by Lemma~\ref{lemma:cohomology of barZldagger_Delta}
~(\ref{item:x in mDelta0}), 
$\Theta^{\Delta}_{l,y_i,\vartheta_{y_i}}(p)
/\Theta^{\Delta}_{l,z_i,\vartheta_{z_i}}(p)
\neq\Theta^{\Delta}_{l,y_i,\vartheta_{y_i}}(p')
/\Theta^{\Delta}_{l,z_i,\vartheta_{z_i}}(p')$. 
It follows that $\phi$ separates $p$ and $p'$.  \par
It remains to prove that $\phi$ is a closed immersion. 
Since $P_{l,0}$ is the union of $G_0$-orbits, 
it is enough to prove it along 
$(g-r)$-dimensional $G_0$-orbits of $p_{l,0}$ where 
$r=\rank X$. 
The $(g-r)$-dimensional $G_0$-orbits are 
$Z^{\dagger}_{l}(\alpha)\simeq p_{\alpha}\times_{k(0)} A_0$ 
$(\alpha\in\Sk^0(\Vor^{\dagger}_{l}))$ where 
$p_{\alpha}=Z_{2Nl}(\alpha)$. Let $p\in Z^{\dagger}_{l}(\alpha)$ and 
$\fm$ the maximal ideal of $\cO_{P_{l,0},p}$. 
Let $\Vor^{\dagger}_{l}(\alpha)$ be the set 
 of all $\Delta\in\Vor^{\dagger}_{l}$ containing $\alpha$, and 
$\Star(\alpha):=\bigcup_{\Delta\in\Vor^{\dagger}_{l}(\alpha)}\Delta
=\bigcup_{w\in X^{\vee},\alpha\in\Sigma^{\dagger}_{l}(w)}\Sigma^{\dagger}_{l}(w)$.
Let $\Delta\in\Vor^{\dagger}_{l}(\alpha)$.  
By Lemma~\ref{lemma:SmDelta_alpha generate C}, if $m\geq r$, then  
$m(\Delta-\alpha)\cap X$ generates $\Cone(\Delta-\alpha)\cap X$. 
Hence $m(\Star(\alpha)-\alpha)$ generates $\Cone(\Delta-\alpha)\cap X$ 
$(\forall \Delta\in \Vor^{\dagger}_{l}(\alpha))$.  
By Lemma~\ref{lemma:cohomology of tZ_dagger_l_Delta}, if 
$x\in\Star(\alpha)\cap X$, then the restriction 
$\Theta_{l,x,\vartheta_x}$ to $P_{l,0}$ is a finite sum of terms 
with $X$-weights in $x+4N^2lm\phi(Y)$. 
Now assume $m\geq 4r$. 
If $x\in r\Star(\alpha)$, then 
by Corollary~\ref{cor:theta single term}, 
\footnote{In this case, $\Vor^{\dagger}_{l}=\Vor_{4Nl}$, that is, 
$2l$ in Corollary~\ref{cor:theta single term} is replaced by $4Nl$. }
we have 
$(x+4Nlm\mu(X^{\vee}))\cap  r\Star(\alpha)=\{x\}$. 
Therefore every nonzero term in 
$(\Theta_{l,x,\vartheta_x}/\xi^{(m)}_{l,\alpha}\theta_{l}^m)_{|P_{l,0}}$ 
with $X$-weight in $x+4N^2lm\phi(Y\setminus\{0\})$ 
is contained in $\fm^2$ (up to units in $\cO_{P_{l,0,p}}$) by 
Lemma~\ref{lemma:SmDelta_alpha generate C}~(\ref{item:primitive element}). 
Since $4N^2lm\geq 2$, there exists 
 $\vartheta_{\alpha}\in\Gamma(A_0,(\cM^{\otimes 4N^2lm})_{\alpha,0})$ 
 nonvanishing at $\pi_{l}(p_{\alpha})$. Then  
$\Theta_{l,\alpha,\vartheta_{\alpha}}/\theta_{l}^m$ is nonvanishing at $p$ 
by Lemma~\ref{lemma:cohomology of barZldagger_Delta}~(\ref{item:x in mDelta0}). It follows  
$\Theta_{l,x,\vartheta_x}/\Theta_{l,\alpha,\vartheta_{\alpha}}\equiv
\vartheta_x\xi^{(m)}_{l,x}/\vartheta_{\alpha}\xi^{(m)}_{l,\alpha}
\mod\fm^2$ if $x\in r\Star(\alpha)$. 
Since $(\cM^{\otimes 4N^2lm})_{x,0}$ 
is very ample on $A_0$ by $4N^2lm\geq 3$,  
$\Theta_{l,x,\vartheta_x}/\Theta_{l,\alpha,\vartheta_{\alpha}}$ 
generate $\fm/\fm^2$ when $x$ (resp. $\vartheta_x$) ranges over $r\Star(\alpha)$ (resp.  
$\Gamma(A_0,(\cM^{\otimes 4N^2lm})_{x,0})$). 
It follows that $\phi$ is a closed immersion 
if $m\geq 4r$, {\it a fortiori,} 
if $m\geq 4g$. This completes the proof. 
\end{proof}

\begin{cor}\label{cor:very ample of Pl}
$\cL_{l}^{\otimes m}$ is very ample for $m\geq 4g$.
\end{cor}

\subsection{The general case}
\begin{thm}
\label{thm:cohomology of P P0}
Let $G$ be a semiabelian $S$-scheme, $\cG$ 
the N\'eron model of $G_{\eta}$ and 
$(P,i,\cN)$ be a cubical compactification of $\cG$ over $S$. Then  
\begin{enumerate}
\item\label{item:vanishing Hq_PN}
 $H^q(P,\cN^{\otimes m})
=H^q(P_0,\cN^{\otimes m}_0)=0$ $(\forall q>0, \forall m>0)$;
\item\label{item:H0_P0}
$\Gamma(P_0,\cN^{\otimes m}_0)
=\Gamma(P,\cN^{\otimes m})\otimes_Rk(0)$;
\item\label{item:rank H0_PN}
$\Gamma(P,\cN^{\otimes m})$ is a free $R$-module of rank 
$(4N^2lm)^g|X/Y|\rank \Gamma(A,\cM)$;
\item\label{item:very ample m_geq 4g}
$\cL^{\otimes m}$ is very ample for $m\geq 4g$.
\end{enumerate}
\end{thm}
\begin{proof}By \S~\ref{subsec:Hq(P,OP)}~($\alpha$), 
$(P^*_{l},\cL^*_{l})$ is constructed 
from the relatively complete model 
$(\tP^*_{l},\tcL^*_{l})$ over $S_{\min}$ with uniformizer $t$. 
Hence Theorem~\ref{thm:cohomology of Pl Pl0}~(1)-(3) 
are true for $(P^*_{l},\cL^*_{l})$. Meanwhile 
by \S~\ref{subsec:Hq(P,OP)}~($\alpha$)-($\beta$), 
$(P,\cN)=(P^{\flat}_{l},\cL^{\flat}_{l})$ 
is the Galois descent of $(P^*_{l},\cL^*_{l})$ 
by $\Gamma=\Aut(k(\eta^*)/k(\eta))$. 
Hence $H^q(P,\cN^{\otimes m})=H^q(P^*_{l},(\cL^*_{l})^{\otimes m})^{\Gamma}$.
This proves (\ref{item:vanishing Hq_PN}) 
for $q>0$ by Theorem~\ref{thm:cohomology of Pl Pl0}. 
For $q=0$, $\Gamma(P,\cN^{\otimes m})$ is an $R$-submodule of 
a free $R$-module $\Gamma(P^*_{l},(\cL^*_{l})^{\otimes m})$ such that 
$$\rank_{R^*}\Gamma(P^*_{l},(\cL^*_{l})^{\otimes m})
=\rank_R\Gamma(P,\cN^{\otimes m}),$$  
which proves (\ref{item:rank H0_PN}). Let $\Lambda$ be an $R$-basis of 
$\Gamma(P,\cN^{\otimes m})$. Then $\Lambda$ is also an $R^*$-basis of 
$\Gamma(P_{l}^*,(\cL^*_{l})^{\otimes m})$. If $m\geq 4g$, then 
$(\cL^*_{l})^{\otimes m}$ is very ample. Hence  
the morphism $\phi_{\Lambda}$ defined by $\Lambda$ from $P_{l}$ 
to $\bP(\Gamma(P,\cN^{\otimes m}))\times_SS^*$ embeds $P_{l}^*$, 
so that $\phi_{\Lambda}$ embeds 
$(P_{l}^*)^{\Gamma}=P$ into $\bP(\Gamma(P,\cN^{\otimes m}))$. 
This proves (\ref{item:very ample m_geq 4g}). 
\end{proof}

\section{Uniqueness of $P$}
\label{sec:uniqueness of P}
\subsection{Codimension two compactifications}
\label{subsec:codim two comp}

\begin{lemma}\label{lemma:lifting of autom of U to P}
Let $S$ be a scheme and $(U,\cL)$   
an $S$-scheme with $\cL$ invertible. 
Let $\pi:(P,i,\cN)\to S$ be an integral proper flat $S$-scheme 
such that 
\begin{enumerate} 
\item[(c1)] $\pi_*(\cO_P)=\cO_S$;
\item[(c2)] $i: U\hookrightarrow P$ is an open immersion with 
$\codim_P(P\setminus i(U))\geq 2$;
\item[(c3)] $P$ is Cohen-Macaulay with $\cN$ 
ample invertible such that $i^*\cN\simeq\cL$.
\end{enumerate}

Let $(P',\iota,\cM)$ be another integral proper flat $S$-scheme 
satisfying (c1)-(c3). Then any $S$-isomorphism 
$g:(U,i^*\cN)\simeq (U,\iota^*\cM)$ extends to 
an $S$-isomorphism $\phi(g):(P,\cN)\simeq (P',\cM)$.
\end{lemma}
\begin{proof}Let $Z:=P\setminus i(U)$, $V_k:=\Gamma(P,\cN^{\otimes k})$ 
and $W_k:=\Gamma(P',\cM^{\otimes k})$.  Since $P$ is Cohen-Macaulay by (c3), 
$\cH^q_{Z}(\cN^{\otimes m})=0$\ 
$(\forall q\in [0,1], \forall m\in\bZ)$ by (c2) and 
\cite[2.2/3.8]{G67}. Hence $H^q_{Z}(P, \cN^{\otimes m})=0$ 
$(\forall q\in [0,1])$ by
the spectral sequence $E_2^{pq}=H^p(P,\cH^q_{Z}(\cN^{\otimes m}))$ 
with abutment $H^q_{Z}(P, \cN^{\otimes m})$ in \cite[2.8, p.~35]{G67}. 
Let $g^*_k:\Gamma(U,\iota^*\cM^{\otimes k})\to\Gamma(U,i^*\cN^{\otimes k})$ 
be the isomorphism induced from $g$. 
It follows that there are isomorphisms $f_k:=(i^*)^{-1}\circ g^*_k\circ \iota^*:W_k\to V_k$ $(k\geq 0)$ 
of $R$-modules 
such that $f_k(a)\cdot f_{l}(b)=f_{k+l}(a\cdot b)$ 
$(\forall a\in W_k,\forall b\in W_{l})$. Moreover $f_0=\id_R$ because 
$\Gamma(\cO_P)=\Gamma(\cO_{P'})=R$ by (c1).
Hence $g$ induces an isomorphism  of graded $R$-algebras : 
$$A(P',\cM):=\bigoplus_{k=0}^{\infty}W_k
\simeq \bigoplus_{k=0}^{\infty}V_k=:A(P,\cN),
$$each of which is isomorphic to $A(U,i^*\cN):=\bigoplus_{k=0}^{\infty}\Gamma(U,i^*\cN^{\otimes k})$. 
Since $\cN$ is ample, 
$P\simeq\Proj A(P,\cN)$.  Similarly $P'\simeq\Proj A(P',\cM)$, 
 so that we have an  $S$-isomorphism  $\phi(g):P\simeq P'$ via
\begin{equation}
\label{eq:P P' isom}
P\simeq\Proj A(P,\cN)\simeq\Proj A(U,i^*\cN)\simeq
\Proj A(P',\cM)\simeq P'.
\end{equation}

Let $\cF:=\phi(g)^*\cM\otimes_{\cO_P}\cN^{\otimes(-1)}$. 
Since $\phi(g)\circ i=\iota\circ g$, we have 
$i^*\cF\simeq \cO_U\simeq i^*\cO_P$. 
Since $\cH^q_{Z}(\cF)=0$\ 
$(\forall q\in [0,1])$ by 
\cite[2.2/3.8]{G67}, we have 
$\cF\simeq i_*i^*\cF\simeq i_*i^*\cO_P\simeq \cO_P$. 
Hence $\phi(g)^*\cM=\cN$. It follows that 
any isomorphism $g:(U,i^*\cN)\simeq (U,\iota^*\cM)$ 
extends to an isomorphism $\phi(g):(P,\cN)\simeq (P',\cM)$. 
This completes the proof. 
\end{proof}

\subsection{Uniqueness of $P$}
\label{subsec:uniqueness of P}
We return to the situation of Theorem~\ref{thm:main thm}. 
\begin{thm}\label{thm:P=Proj A(cG cN)}
Let $(G,\cL)$ be a semi-abelian $S$-scheme, $\cG$ its N\'eron model  
and $(P,i,\cN)$ a cubical compactification of $\cG$ 
extending $(G_{\eta},\cL_{\eta})$.
Then 
\begin{enumerate}
\item\label{item:P=Proj A(cG cN)}  
$P\simeq\Proj A(\cG,i^*\cN)$ with $i_*i^*\cN=\cN$; 
\item\label{item:P unique}  
$P$ is uniquely determined by $(G_{\eta},\cL_{\eta})$; 
\item\label{item:cN unique} 
$\cN$ is uniquely determined by $(G_{\eta},\cL_{\eta})$ 
up to positive multiples;
\item\label{item:extension of autom} 
any $k(\eta)$-automorphism 
$h_{\eta}$ of $G_{\eta}$ 
with $h_{\eta}^*\cL_{\eta}\simeq\cL_{\eta}$ extends 
to an $S$-automorphism $g$ of 
$P$ with $g(\cG)=\cG$ and $g^*\cN=\cN$.
\end{enumerate}
\end{thm}
\begin{proof} $(P,i,\cN)$ satsifies (rc1)-(rc3) and hence 
(c2)-(c3), while (c1) is true for  $(P,i,\cN)$ 
by \cite[Cor.~2, p.~48]{Mumford12}/Theorem~\ref{thm:cohom of P}~(\ref{item:hi_P_cO}). By (rc3) and Eq.~(\ref{eq:P P' isom}),  
$P\simeq\Proj A(\cG,i^*\cN)$. 
 Since $P$ is Cohen-Macaulay with $\codim_P(P\setminus i(\cG))\geq 2$ 
and $\cN$ is invertible, we have $\cN=i_*i^*\cN$.  
This proves (\ref{item:P=Proj A(cG cN)}). 
 Suppose that $\cM$ is an ample invertible 
$\cO_P$-module such that $i^*\cM$ is cubical and 
$\cM_{\eta}\simeq\cL_{\eta}^{\otimes l'}$ for some $l'\in\bN$. 
Then there exist $m,n\in\bN$ such that 
$\cM^{\otimes n}_{\eta}\simeq\cN^{\otimes m}_{\eta}$. 
 Let $\cF:=\cM^{\otimes n}\otimes\cN^{\otimes (-m)}$.
 Since $i^*\cM$ and $i^*\cN$ are cubical, so is $i^*\cF$. 
Since $\cF_{\eta}\simeq\cO_{G_{\eta}}$, we have $i^*\cF\simeq\cO_{\cG}$ 
by \cite[2.13~(3)]{MN24}. 
Therefore $i^*\cM^{\otimes n}\simeq i^*\cN^{\otimes m}$, so that  
$P\simeq \Proj A(\cG,i^*\cM^{\otimes n})
\simeq \Proj A(\cG,i^*\cN^{\otimes m})$ 
by (\ref{item:P=Proj A(cG cN)}).  Since 
$i^*\cN$ is, by \cite[2.13~(3)]{MN24},  
a unique cubical invertible sheaf on $\cG$ 
with $(i^*\cN)_{\eta}\simeq\cL_{\eta}^{\otimes l}$, 
$P$ is uniquely determined by 
$(G_{\eta},\cL_{\eta})$. This proves (\ref{item:P unique}). 

Since $P$ is Cohen-Macaulay with $\codim_P(P\setminus i(\cG))\geq 2$ 
and $\cF$ is invertible, we have 
$\cF\simeq i_*i^*\cF\simeq i_*\cO_{\cG}\simeq\cO_P$  by 
\cite[2.2/3.8]{G67}. It follows $\cM^{\otimes n}\simeq\cN^{\otimes m}$. 
This proves (\ref{item:cN unique}). 
Next we prove (\ref{item:extension of autom}).
Since $\cG$ is the N\'eron model of $G_{\eta}$,  
$h_{\eta}$ extends to an $S$-automorphism $h$ of $\cG$. 
By Theorem~\ref{thm:main thm}, $\cN$ is 
an invertible sheaf of $P$ such that $i^*\cN$ is cubical 
on $\cG$ and $\cN_{\eta}\simeq \cL^{\otimes l}_{\eta}$ by (rc3). 
Since both $h^*i^*\cN$ and $i^*\cN$ are cubical 
on $\cG$ with $(h^*\cN)_{\eta}\simeq h^*_{\eta}\cL_{\eta}^{\otimes l}
\simeq\cL_{\eta}^{\otimes l}\simeq\cN_{\eta}$, 
we have $h^*i^*\cN\simeq i^*\cN$  
by \cite[2.13~(3)]{MN24}. 
 By Lemma~\ref{lemma:lifting of autom of U to P},    
$h$ extends to an $S$-automorphism $g:=\phi(h)$ of $P$ such that 
$g^*\cN=\cN$.  This proves (\ref{item:extension of autom}). 
\end{proof}
\begin{cor}\label{cor:extension of autom of compactification}
Let $(P,i,\cN)$ and $(P',\iota,\cM)$ 
be cubical compactifications of $\cG$ 
extending $(G_{\eta},\cL_{\eta})$. Then 
there exists an isomorphism $\phi:P\to P'$ 
such that $\phi\circ i=\iota$ and $\cN^{\otimes m}
\simeq\phi^*\cM^{\otimes n}$ for some $m,n\in\bN$;
\end{cor}
\begin{proof} 
Let $U=\cG$ and apply Lemma~\ref{lemma:lifting of autom of U to P}. 
By assumption,  we can choose 
$m, n\in\bN$ such that 
$\cN^{\otimes m}_{\eta}
\simeq\cM^{\otimes n}_{\eta}$.
By applying Lemma~\ref{lemma:lifting of autom of U to P} to $g=\id_U$, 
we obtain an $S$-automorphism $\phi:=\phi(\id_U):P\to P'$ such that 
$\cN^{\otimes m}\simeq\phi^*\cM^{\otimes n}$ and $\phi_{|U}=\id_U$, 
so that $\phi\circ i=\iota$. This completes the proof.  
\end{proof}

\section{Relative compactification over a Dedekind domain}
\label{sec:rel compactification Dedekind} 
\subsection{Goal}
\label{subsec:plan}
The purpose of this section is to prove Theorem~\ref{thm:main thm Dedekind} 
which is the version  over a Dedekind domain 
of Theorem~\ref{thm:main thm}.  We note that 
the same type of uniqueness  
as Theorem~\ref{thm:P=Proj A(cG cN)} 
is true over a Dedekind domain. 
We assume that $\Sigma^{\dagger}_{l}(0)$ is integral.
\begin{thm}\label{thm:main thm Dedekind} \ \ 
Let $D$ be a Dedekind domain, $K$ the fraction field of $D$ and 
$S=\Spec D$. Let $(G,\cL)$ be a semiabelian $S$-scheme   
 over $S$, $\cL$ a 
symmetric ample cubical invertible sheaf on $G$ 
and $\cG$ a N\'eron model of $G$.  
Then there exists a relative compactification 
$(P,i,\cN)$   
of $\cG$ extending $(G_K,\cL_K)$ such that 
\begin{enumerate} 
\item[(a)] $P$ is Cohen-Macaulay;
\item[(b)] $i(\cG)=P\setminus\Sing(P/S)$ with   
$\codim_{P}\Sing(P_S)\geq 2$ where $\Sing(P/S)$ denotes the union 
of the singular loci of closed fibers of $P$ over $S$;
\item[(c)] $i^*\cN$ is ample cubical invertible; 
\item[(d)] 
$\cG$ acts on $P$ so that $i$ is $\cG$-equivariant.
\end{enumerate}

Moreover the relative compactification $(P,i,\cN)$ satisfying (a)-(c) 
is unique up to isomorphism 
in the sense of Theorem~\ref{thm:P=Proj A(cG cN)/intro}.
\end{thm}

The rest of this section is devoted to proving 
Theorem~\ref{thm:main thm Dedekind}.

Let $s$ be a closed point of $S$, 
$D_s$ the localization of $D$ at 
(the prime ideal $\fp_s$ of $D$ corresponding to) $s$, $I_s$ 
the maximal ideal of $D_s$ and $k(s):=D_s/I_s$. Let 
$\cB$ be the set of closed points $s$ of $S$ such that 
$\cG_s$ is not an abelian variety over $k(s)$. 
Let $D^{\wedge}_s:=(D_s)^{\wedge}$ be the $I_s$-adic completion 
of $D_s$, $K_s$ the fraction field of $D_s^{\wedge}$, 
$S_s:=\Spec D_s$, $S^{\wedge}_s:=\Spec D^{\wedge}_s$ and $\eta_s$ 
the generic point of $S^{\wedge}_s$. Since $D$ is a Dedekind domain, 
$D^{\wedge}_s$ is a CDVR. Let $\cG_{S_s}:=\cG\times_SS_s$, 
$\cG^s:=\cG\times_SS^{\wedge}_s$, $G_{S_s}:=G\times_SS_s$, 
$G^s:=G\times_SS^{\wedge}_s$, $\cL^s:=\cL\otimes_{\cO_G}\cO_{G^s}$, 
$\Phi^s:=\cG^s/G^s$ and  $N:=\LCM(|\Phi^s|;s\in\cB)$. 
Then $\cG_{S_s}$ (resp. $\cG^s$) 
is the N\'eron model over $D_{s}$ (resp. $D_s^{\wedge}$) 
of the abelian variety $G_K$ 
(resp. $G_{\eta_s}:=G_K\otimes_KK_s$).

\subsection{Proof of Theorem~\ref{thm:main thm Dedekind} - start}
\label{subsec:proof of main thm Dedekind start}
In what follows, we choose and fix any $s\in\cB$ and 
set $R:=D_s$ and $U:=S_s=\Spec R$.  

\begin{defn}
Let $Q$ be an $S$-scheme and $\cF$ an ample invertible sheaf on $Q$ 
and $\Image\left(S^n\Gamma(Q,\cF)\right)$ the image of the homomorphism $S^k\Gamma(Q,\cF)\to \Gamma(Q,\cF^{\otimes k})$.  
Then we define 
\begin{gather*}
A(Q,\cF):=\bigoplus_{k=0}^{\infty}
\Gamma(Q,\cF^{\otimes k}),\ \ S(Q,\cF):=\bigoplus_{k=0}^{\infty}
\Image\left(S^k\Gamma(Q,\cF)\right).
\end{gather*}
\end{defn}

\begin{lemma}\label{lemma:finite generation}
Let $Z$ be a projective $S$-scheme and $\cL$ 
a (relatively) very ample invertible sheaf on $Z$.
Then there exists $n_1$ such that 
$S^n\Gamma(Z,\cL^{\otimes n_1})\to \Gamma(Z,\cL^{\otimes nn_1})$ 
is surjective $(\forall n\geq 1)$. 
\end{lemma}
\begin{proof}
We regard $Z$ as a closed $S$-subscheme of $\bP_S^M$ with $\cI$ 
its defining ideal and 
$\cL=\cO_{\bP_S^M}(1)\otimes_{\cO_{\bP_S^M}}\cO_Z$ for some $M\in\bN$.
There is an exact sequence 
$0\to \cI\cO_{\bP_S^N}(n)\to\cO_{\bP_S^N}(n)\to \cL^{\otimes n}\to 0$\ 
$(\forall n\in\bN)$.  
There exists $n_1$ such that $H^1(\cI\cO_{\bP_S^N}(n))=0$ 
for any $n\geq n_1$, so that 
 $S^n\Gamma(\cO_{\bP_S^N}(1))=\Gamma(\cO_{\bP_S^N}(n))\to
\Gamma(Z,\cL^{\otimes n})\to 0$ is surjective. It follows 
that $S^n\Gamma(Z,\cL)\to \Gamma(Z,\cL^{\otimes n})$ 
is surjective $(\forall n\geq n_1)$. 
If $n\geq 1$, then as $R$-submodules of $\Gamma(Z,\cL^{\otimes nn_1})$, 
\begin{align*}
&\Image\left(S^n\Gamma(Z,\cL^{\otimes n_1})\right)
=\Image\left(S^n\left(S^{n_1}\Gamma(Z,\cL)\right)\right)\\
&=\Image\left(S^{nn_1}\Gamma(Z,\cL)\right)=\Gamma(Z,\cL^{\otimes nn_1}),
\end{align*}which proves Lemma. 
\end{proof}

By \cite[II, 1.2.1]{MB85}/\cite[2.13]{MN24},  
there exists a unique ample 
cubical invertible sheaf $\cM$ on $\cG$ with 
$\cM_K:=\cM\otimes_RK\simeq\cL^{\otimes 2N}_K$. 
Let $\cM_U:=\cM\otimes_DR$ and let 
$\cM_U^{\wedge}$ be the $I_s$-adic completion of $\cM_U$.

Let $l\geq 1$. By Raynaud \cite[VI,\ 1.4.2]{MB85}, 
$\Gamma(\cG,\cM^{\otimes l})$ is a finite $D$-module. 
Since $\cM$ is locally $\cO_{\cG}$-free and $\cG$ is smooth over $S$, 
$\Gamma(\cG,\cM^{\otimes l})$ is $D$-torsion free. Hence 
$\Gamma(\cG,\cM^{\otimes l})$ is locally $\cO_S$-free because  
$D_s$ is a DVR for any closed point $s$ of $S$. 
It follows that for any $l,n\geq 1$,
\begin{gather*}
\Gamma(\cG_U,\cM_U^{\otimes l})=\Gamma(\cG,\cM^{\otimes l})\otimes_DR,\ \ 
\Gamma(\cG_{U^{\wedge}},\cM_{U^{\wedge}}^{\otimes l})
=\Gamma(\cG,\cM^{\otimes l})
\otimes_DR^{\wedge},\\
A(\cG_U,\cM_U^{\otimes l})=A(\cG,\cM^{\otimes l})\otimes_DR,\ \ 
A(\cG_{U^{\wedge}},\cM_{U^{\wedge}}^{\otimes l})=A(\cG,\cM^{\otimes l})\otimes_DR^{\wedge},\\
S(\cG_U,\cM_U^{\otimes l})=S(\cG,\cM^{\otimes l})\otimes_DR,\ \ 
S(\cG_{U^{\wedge}},\cM_{U^{\wedge}}^{\otimes l})=S(\cG,\cM^{\otimes l})\otimes_DR^{\wedge}.
\end{gather*}

\subsection{Proof continued}
By Theorem~\ref{thm:main thm}, 
there exists a relative compactification 
$(P,i^s,\cM^{\dagger})$ of $\cG_{U^{\wedge}}$ such that (a)-(d) are true where 
$i^s:\cG_{U^{\wedge}}\to P$ is an open immersion, 
$P$ is an irreducible Cohen-Macaulay $U^{\wedge}$-scheme 
and $((i^s)^*\cM^{\dagger})_{k(\eta_s)}=\cL_{k(\eta_s)}^{\otimes n_0}$ 
for some $n_0\in\bN$. 
Since $P$ is Cohen-Macaulay, $\Gamma(P,(\cM^{\dagger})^{\otimes n})
=\Gamma(\cG_{U^{\wedge}},\left((i^s)^*(\cM^{\dagger})\right)^{\otimes n})$ 
$(\forall n\geq 0)$. 
By Lemma~\ref{lemma:finite generation}, there exists 
$n_1\in\bN$ such that 
$A(P,(\cM^{\dagger})^{\otimes n_1})=S(P,(\cM^{\dagger})^{\otimes n_1})$. 

Since $\cM_U$ is ample, there exists $m_2\in\bN$ 
such that $\cM^{\otimes m}_U$ is very ample $(\forall m\geq m_2)$. 
In view of \cite[2.13~(3)]{MN24},
by choosing an integral multiple of $\cM^{\dagger}$ 
for $\cM^{\dagger}$ if necessary, we may assume 
$(i^s)^*\cM^{\dagger}\simeq\cM_{U^{\wedge}}^{\otimes n_2}$ 
for some $n_2\in\bN$ with $n_2\geq m_2$.  
It follows from Theorem~\ref{thm:P=Proj A(cG cN)} 
\begin{align*}
P&=\Proj A(P,(\cM^{\dagger})^{\otimes n_1})
=\Proj S(P,(\cM^{\dagger})^{\otimes n_1})\\
&\simeq\Proj S(\cG_{U^{\wedge}},((i^s)^*\cM^{\dagger})^{\otimes n_1})
\simeq\Proj S(\cG_{U^{\wedge}},\cM_{U^{\wedge}}^{\otimes n_1n_2}). 
\end{align*}

By Theorem~\ref{thm:cohomology of P P0}~(\ref{item:vanishing Hq_PN})-(\ref{item:H0_P0}), for any $l\in\bN$, 
\begin{align*}
H^q(P_s,(\cM^{\dagger})^{\otimes l}\otimes_{R}k(s))
&=0\ \  (q\geq 1),\\
\Gamma(P_s,(\cM^{\dagger})^{\otimes l}\otimes_{R}k(s))
&=\Gamma(P,(\cM^{\dagger})^{\otimes l})\otimes_{R^{\wedge}}k(s)\\
&=\Gamma(\cG_{U^{\wedge}},
\cM_{U^{\wedge}}^{\otimes n_2l})\otimes_{R^{\wedge}}k(s)\\
&=\Gamma(\cG_U,\cM_U^{\otimes n_2l})\otimes_{R}k(s).
\end{align*}
By Lemma~\ref{lemma:finite generation}, there exists $n_3\in\bN$ such that 
$$\Image\left(S^{l}\Gamma(P_s,(\cM^{\dagger})^{\otimes n_3}\otimes_{R^{\wedge}}k(s))\right)=\Gamma(P_s,(\cM^{\dagger})^{\otimes n_3l}
\otimes_{R^{\wedge}}k(s))\ \ 
(\forall l\in\bN).$$ 
Let $N_0:=n_1n_2n_3$.
Then  
$\Image\left(S^{l}\Gamma(\cG_U,\cM_U^{\otimes N_0})\right)=
\Gamma(\cG_U,\cM_U^{\otimes N_0l})$ $(\forall l\in\bN)$ 
by Nakayama's lemma,.   
It follows that 
$A(\cG_U,\cM_U^{\otimes N_0})=S(\cG_U,\cM_U^{\otimes N_0})$, 
which is an $R$-algebra generated by a finite free $R$-module 
$\Gamma(\cG_U,\cM_U^{\otimes N_0})$. 

Let $Q:=\Proj S(\cG_U,\cM_U^{\otimes N_0}).$  By definition, $Q$ is 
a projective $R$-scheme (of finite type), 
so that $Q$ is excellent, {\it i.e.}, 
the local ring $\cO_{Q,x}$ is excellent 
for any closed point $x$ of $Q$ by \cite[2.3]{MN24}.
See \cite[IV$_2$, 7.8.3]{EGA}.
Since $P\simeq Q_{U^{\wedge}}$ and $P$ is Cohen-Macaulay, so is 
$Q_{U^{\wedge}}^{\wedge}=Q_U^{\wedge}$, so that 
$Q=Q_U$ is Cohen-Macaulay by \cite[2.1/2.4]{MN24}. 
Moreover $Q$ is $S$-flat by \cite[III, 9.9]{Hartshorne77} because $\chi(Q_s,\cM_s)
=\chi(P_s,\cM_s)=\chi(P_{K_s},\cM_{K_s})=\chi(Q_{K_s},\cM_{K_s})
=\chi(Q_K,\cM_K)$.

Let $V_{N_0}:=\Gamma(\cG_U,\cM_U^{\otimes N_0})$. 
Since $S(\cG_U,\cM_U^{\otimes N_0})$ is generated by 
$V_{N_0}$, $Q$ is a closed $S$-subscheme of $\bP(V_{N_0})$.  
Associated with the linear system $|V_{N_0}|$, 
we have a natural morphism 
$i:\cG_U\to \bP(V_{N_0})$, which factors through $Q$. 
Then $\cM_U^{\otimes N_0}$ is very ample by $N_0\geq n_2\geq m_2$. 
Hence $i$ is an injective immersion which is an isomorphim over $K$, while    
$i_{U^{\wedge}}=i^s:\cG_{U^{\wedge}}\hookrightarrow Q_{U^{\wedge}}=P$ 
is an open immersion over $U^{\wedge}$. 
It follows that  $i$ is an open immersion of $\cG_U$ into $Q$.  
Now we define $\cN:=i_*\cM_U^{\otimes N_0}$. 
Since $P$ is Cohen-Macaulay and $\cM^{\dagger}$ 
is locally $\cO_P$-free, we obtain 
$\cN_{U^{\wedge}}=(i_{U^{\wedge}})_*\cM_{U^{\wedge}}^{\otimes N_0}
=(i^s)_*(i^s)^*((\cM^{\dagger})^{\otimes n_1n_3})
=(\cM^{\dagger})^{\otimes n_1n_3}$.
Since $(\cM^{\dagger})^{\otimes n_1n_3}=\cN_{U^{\wedge}}$ 
is an invertible $\cO_P$-module, 
$\cN$ is an invertible $\cO_Q$-module. Moreover 
$i^*\cN=i^*i_*\cM_U^{\otimes N_0}=\cM_U^{\otimes N_0}$. Since 
$\cM$ is cubical on $\cG$,  $i^*\cN$ is cubical with 
$(i^*\cN)_K=\cM_K^{\otimes N_0}=\cL_K^{\otimes 2NN_0}$. 
This proves that $(Q,i,\cN)$ is 
a relative compactification of $\cG_U$ such that (a)-(c) are true. 
Uniqueness of $(Q,i,\cN)$ in the sense 
of Theorem~\ref{thm:P=Proj A(cG cN)/intro} 
is clear. Now we denote $(Q,i,\cN)$ by $(Q^s,j^s,\cN^s)$. 
Since $(Q^s,j^s,\cN^s)$ is unique over $S_s$, 
 by gluing together all $(Q^s,j^s,\cN^s)$ $(s\in\cB)$ 
over $S\setminus\cB$, 
we can construct a relative compactification $(P,i,\cN)$ 
of $\cG$ over $S$ such that (a)-(c) are true and 
$(P,i,\cN)\times_S{S_s}\simeq (Q^s,j^s,\cN^s)$ 
in the sense of Theorem~\ref{thm:P=Proj A(cG cN)/intro}.   
Thereby uniqueness of $(P,i,\cN)$  
is clear.

\subsection{Proof of (d)} 
\label{subsec:proof of main thm (d)}
The proof of (d) is very formal.   
We set as follows:
\begin{gather*}
P_{12}:=P\times_SP,\ P_{123}:=P\times_SP\times_SP,\\ 
L_{12}:=G\times_SP,\ 
L_{123}:=G\times_SP\times_SP.
\end{gather*}

Since $\cG_K=P_K$ 
is an abelian variety over $K$,  we have an action $\rho_K:G_K\times P_K\to P_K$, which defines a graph, {\it i.e.}, a $K$-flat closed subscheme $\Gamma^K$ of $P_{123,K}$ such that the restriction of the projection to the first two components 
$p_{12,K}:\Gamma^K\to P_{12,K}$ is an isomorphism.  
Similarly by Theorem~\ref{thm:main thm}~(d), we have an 
action $\rho^s:\cG^s\otimes_{D^{\wedge}_s} P_{S^{\wedge}_s}$ with graph $\Gamma^s$ which is a closed $S^{\wedge}_s$-flat subscheme of $P_{123,S^{\wedge}_s}$. Since $P$ is projective over $S$, there is a $S$-flat closed subscheme $\Gamma$ of $P_{123}$ by \cite[IV$_2$, 2.8.5]{EGA} such that $\Gamma_K:=\Gamma\times_S\Spec K=\Gamma^K$ and $\Gamma\times_{S}S^{\wedge}_s=\Gamma^s$. 
The projection $p_{12}:\Gamma\ (\subset P_{123})\to P_{12}$ is an $S$-morphism which is an isomorphism over $K$ and $S^{\wedge}_s$ $(\forall s\in\cB)$. Hence 
$p_{12}$ is an isomorphism everywhere over $S$. Hence $\Gamma$ defines an action of $G$ on $P$ (denoted by $\rho:L_{12}\to L_{123}$) extending $\rho^K$ and $\rho^s$. This proves (d).

\subsection{Relative compactifications in the general case}
\label{sec:rel compact general}
Let $R$ be a CDVR with $k(\eta)$ its fraction field. 
Let $(G_{\eta},\cL_{\eta})$ be 
a polarized abelian variety over $k(\eta)$, 
$\cG$ the N\'eron model of $G_{\eta}$, $G:=\cG^0$ 
the identity component of $\cG$ and $N:=|\cG_0/G_0|$. 
By \cite[II, 1.1]{MB85}, there exists 
an ample {\it cubical} invertible sheaf $\cN$ on $\cG$ with  
$\cN_{\eta}\simeq\cL^{\otimes 2N}_{\eta}$.  
By Raynaud \cite[VI,\ 1.4.2]{MB85} $\Gamma(\cG,\cN^{\otimes k})$ 
is a finite free $R$-module for any $k\in\bN$. 

\begin{conj}\label{conj:rel comp of general Neron model}
Let $P_n:=\Proj A(\cG,\cN^{\otimes n})$. Then 
\begin{enumerate}
\item[(i)] $P_{n_0}$ is a relative compactification of $\cG$ for some $n_0\in\bN$;
\item[(ii)] $\cG\simeq P_{n_0}\setminus\Sing(P_{n_0}/S)$ with 
$\codim_{P_{n_0}}\hskip -0.05cm\Sing(P_{n_0}/S)\geq 2$.
\end{enumerate}
\end{conj}

\begin{rem}
Conjecture~\ref{conj:rel comp of general Neron model} 
is true for elliptic curves over $k(\eta)$ 
by Kodaira-N\'eron's classification of singular fibers 
of elliptic fibrations. 
\end{rem}

\section{The complex analytic case}
\label{sec:complex case}
\subsection{Notation and Terminology}
In this section we consider 
a complex analogue to \S\S~\ref{sec:twisted Mumford families}-\ref{sec:str of Pl in pd case}. 
Let $P,P'$ be groups, and  
 $\Hom_{\bZ}(P,P')$ the set of homomorphisms  
from $P$ to $P'$. Let $Q$ a ring and $Q^{\times}$ 
the group of {\it invertible} elements in $Q$. 
Let $t$ be a complex parameter, 
$\bD:=\{t\in\bC;|t|<1\}$ the unit disc 
with center $0$, $D^{\times}:=\bD\setminus\{0\}$ and 
$t^{\bZ}:=\{t^n; n\in \bZ\}$. 
Let $Z$ be a complex space and $Z(\bD):=\Hom(\bD,Z)$ be 
the set of complex analytic maps from $\bD$ to $Z$. 
Let $\varpi:G\to \bD$ be a relative complex Lie group over $\bD$. 
If moreover any fiber $(G_t,\cL_t)$ of it 
is a complex semiabelian scheme, then we call $\varpi:G\to \bD$ 
a (complex) semiabelian $\bD$-scheme .

\subsection{Theta functions}
\label{subsec:theta functions}
In the complex case, we will see that Eq.~(\ref{eq:action S_ly}) 
is a consequence of transformation equations of theta functions. 
We start with quoting from \cite[Chap.~12, p.~77]{Siegel62} 
by slightly modifying. 
Let $g=g_1+g_2$ with $g_1,g_2\geq 0$, 
$X=Y=\bigoplus_{i=1}^{g_1}\bZ m_i\simeq\bZ^{g_1}$ and  
define an injective homomorphism 
$\phi:Y\to X$ by $\phi(m_i)=d_im_i$ $(i\in[1,g_1])$.  
Let $Z_1:=X^{\vee}=\bigoplus_{i=1}^{g_1}\bZ f_i$, 
$Z_2:=\bigoplus_{i=g_1+1}^{g}\bZ f_i\simeq\bZ^{g_2}$ 
and $Z=Z_1\oplus Z_2$.  Let $V:=Z\otimes_{\bZ}\bC$ and 
$V_k:=Z_k\otimes_{\bZ}\bC$ $(k\in[1,2])$. 
Let $t_i\in V$ be 
a column vector $(i\in[1,g])$ and define 
\begin{gather*}
T_1=(t_1,\cdots,t_{g_1}),\ T_2=(t_{g_1+1},\cdots,t_{g_1+g_2}),\\
T=(T_1, T_2)=
\begin{pmatrix}T_{11}&T_{12}\\
T_{21}&T_{22}
\end{pmatrix},\ 
C=(I_g,T),\  C_i=(I_{g_i},T_{ii}),\\
D=\diag(d_1,\cdots,d_g)=
\begin{pmatrix}D_1&0\\
0&D_2
\end{pmatrix},\ 
W=DT,\  W_1=D_1T_{11} 
\end{gather*}
where $T_{ij}\in M_{g_i\times g_j}(\bC)$,  
$D_i\in M_{g_i,g_i}(\bZ)$ and $W$ is symmetric with $\Impart(W)>0$. 
Let $W_k(u_k,v_k)=u_k^tW_kv_k$ and  $W_k[v_k]=W_k(v_k,v_k)$ 
$(u_k,v_k\in\bC^{g_i})$.

Let $\Delta:=\{a=(a_i)\in Z;0\leq a_i<d_i\ (\forall i\in [1,g])$,
$a\in\Delta$ and $s_a\in\bC$. Then we define $\theta_a$,
$\theta$ and $\sigma_{r_1}(\theta)$ by
\begin{align*}
\theta_a(z)&=\sum_{m\in Y}\be(W[m+D^{-1}a]/2+(Dm+a)^tz)\\
&=\sum_{m\in Z}\be(W[D^{-1}(\phi(m)+a)]/2+(\phi(m)+a)^tz),\\
\theta(z)&=\sum_{a\in\Delta}s_a\theta_a(z)
=\sum_{r_1\in X}\be(r_1^tz_1)\sigma_{r_1}(\theta)(z_2)
\end{align*} where $z:=(z_1,z_2)$ and $z_i\in Z_i\otimes_{\bZ}\bC$. 
For $b,c\in Z$, we have
\begin{equation}\label{eq:transf of theta}
\begin{aligned}
\theta_a(z+D^{-1}b+TD^{-1}c)&=\be(a^tD^{-1}b-W[D^{-1}c]/2-c^tz)
\theta_{a+c}(z),\\
\theta(z+b+Tc)&=\be(-W[c]/2-(Dc)^tz)\theta(z).
\end{aligned}
\end{equation}

Let $G:=V/Z+TZ$ 
and $A_2=V_2/Z_2+T_{22}Z_2$. 
Then $G$ and $A$ are abelian varieties. 
Eq.~(\ref{eq:transf of theta}) shows that 
if $\theta$ is not identically zero, 
the divisor $\theta=0$ defines an ample line bundle 
on $G$ independent of $a\in\Delta$, so that 
the data $(C,D)$ determines 
a unique ample line bundle $\cL_{C,D}$ of $G$. 
By \cite[Th.\,12A]{Siegel62}, 
$(\theta_a;a\in\Delta)$ is a $\bC$-basis of 
$\Gamma(G,\cL_{C,D})$. Let 
$K(\cL_{C,D}):=\{a\in G(\bC); T_a^*\cL_{C,D}\simeq\cL_{C,D}\}$ 
be the Heisenberg group scheme of $\cL_{C,D}$. Then 
$K(\cL_{C,D})=(D^{-1}Z+TD^{-1}Z)/(Z+TZ)$,  
$|K(\cL_{C,D})|=|D|^2$ and $G^t:=\Pic_{G/\bC}=G/K(\cL_{C,D})
\simeq\bC^{g}/D^{-1}Z+TD^{-1}Z$.  

Let $c_1\in Z_1$ and $c':=(c_1,0_{g_2})^t\in Z$. 
Then by Eq.~(\ref{eq:transf of theta}), we obtain
$$\theta(z+Tc')=\be(-W_1[c_1]/2)\sum_{r_1\in X}\be(r_1^tz_1)
\sigma_{r_1+\phi(c_1)}(\theta)(z_2).
$$From this we derive a relation : 
\begin{equation}\label{eq:transf formula of sigma and sigma}
\sigma_{r_1+\phi(c_1)}(\theta)(z_2)=\be(W_1[c_1]/2)
\be(r_1^tT_{11}c_1)\sigma_{r_1}(\theta)(z_2+T_{21}c_1).
\end{equation} 
This is Theorem~\ref{thm:degeneration data2}~(\ref{item:Gamma Geta L general})/Eq.~(\ref{eq:Gamma(geta,cLeta) psi tau}) if we set $\psi(y):=\be(W_1[y]/2)$ and 
$\tau(y,x):=\be(y^tT_{11}^tx)$ $(y\in Y, x\in X)$. 
Note $\psi(y)=\tau(y,\phi(y))^{1/2}$. Thus 
Eq.~(\ref{eq:transf formula of sigma and sigma}) 
is a classical equation 
behind Eqs.~(\ref{eq:Gamma(geta,cLeta) psi tau})/(\ref{eq:action Sy on sigma_x})/(\ref{eq:action delta_lv})/(\ref{eq:action S_ly}).  

Let $c_2\in Z_2$ and $c''=(0_{g_1},c_2)^t\in Z$. 
By computing $\theta(z+Tc'')$, 
we obtain 
{\small\begin{align*}\sigma_{r_1}(\theta)(z_2+T_{22}c_2)&
=\be(-W_2[c_2]/2-(D_2c_2)^tz_2)\be(-r_1^tT_{12}c_2))
\sigma_{r_1}(\theta)(z_2).
\end{align*}} \hskip -0.15cm Hence the divisor $\sigma_{r_1}(\theta)=0$  defines a unique line bundle $\cL_{C_2,D_2}\otimes_{\cO_A}\cO_{r_1}$ 
where $\cO_{r_1}:=c(r_1)=\bF_{-r_1}$ with the notation of  
\S~\ref{subsec:complex split objects}.  
Note $K(\cL_{C_2,D_2})=(D_2^{-1}Z_2+T_{22}D_2^{-1}Z_2)/(Z_2+T_{22}Z_2)$,  
$|K(\cL_{C_2,D_2})|=|D_2|^2$ and $A^t:=\Pic_{A/\bC}=A/K(\cL_{C_2,D_2})
\simeq V_2/D_2^{-1}Z_2+T_{22}D_2^{-1}Z_2$.  
If $D=mI_g$, then we set 
$\cL^0:=\cL_{C,I_g}$ and $\cM^0:=\cL_{C_2,I_{g_2}}$. 
Then $\cL_{C,D}=(\cL^0){\otimes m}$ and 
$\cL_{C_2,D_2}=(\cM^0){\otimes m}$. 
Moreover 
$(\sigma_{r_1}(\theta_a); a\in\Delta\cap Z_2)$ is a $\bC$-basis of 
$\Gamma(A,(\cM^0){\otimes m}\otimes_{\cO_A}\cO_{r_1})$.

\subsection{Complex analytic split objects}
\label{subsec:complex split objects}
Let $X,Y,Z_k,V,V_k$ and $\phi:Y\to X$ 
be the same as in \S~\ref{subsec:theta functions}. Replacing $\bC$ by 
$\Gamma(\cO_{\bD^{\times}})$, we set 
$T=(T_{i,j})\in M_{g,g}(\Gamma(\cO_{\bD^{\times}}))$ 
 as in \S~\ref{subsec:theta functions}.  Assume moreover 
\begin{enumerate}  
\item[(i)] $\omega:=T_{11}=B\log t+T_{11}^0(t),\ \omega^0=:T_{11}^0(t)$; 
\item[(ii)] $T_{11}^0$ and $T_{ij}$ are holomorphic over $\bD$ for 
$(i,j)\neq (1,1)$;
\item[(iii)] $B\in M_{g_1,g_1}(\bZ)$ is symmetric positive definite.
\end{enumerate} 

Then we define complex quotient manifolds over $\bD$:
\begin{equation}\label{eq:defn of G/tG/A/Gt/tGT/At}
\begin{aligned}
G&:=\bD\times V/Z\oplus TZ,\\
\tG&:=\bD\times V/Z\oplus T(0\oplus Z_2),\\
A&:=\bD\times V_2/Z_2\oplus T_{22}Z_2,\\
T&:=\bD\times V_1/Z_1,\\
G^t&:=\bD\times V/D^{-1}Z\oplus TD^{-1}Z,\\
\tG^t&:=\bD\times V/D^{-1}Z\oplus T(0\oplus D_2^{-1}Z_2),\\
A^t&:=\bD\times V_2/D_2^{-1}Z_2\oplus T_{22}D_2^{-1}Z_2,\\
T^t&:=\bD\times V_1/D_1^{-1}Z_1.
\end{aligned}
\end{equation}
Let $\bC^*:=\bC/\bZ$. Let   
$T_X:=\Hom(X,\bG_{m,\bZ})$ be a split torus over $\bZ$. 
Then by Eq.~(\ref{eq:defn of G/tG/A/Gt/tGT/At}), 
we have a natural isomorphism  
$T\simeq T_{X,\bD}:=\bD\times_{\bZ} T_X\simeq(\bC^*)^{g_1}_D$. 
It is clear that $A$ is an abelian $\bD$-scheme. 
Since any fiber $(G_t,\cL_t)$ $(t\neq 0)$ is a complex abelian variety, 
and $G_0\simeq(\bC^*)^{g_1}\times A_0$, $G$ is a semiabelian $\bD$-scheme. 
Meanwhile $\tG$ is a semiabelian $\bD$-scheme, which we call 
the (complex) Raynaud extension of $G$: 
\begin{equation}\label{eq:complex Raynaud extension}
0\to T\to\tG\to A\to 0.
\end{equation}

 Each $x\in X(T)$  
is a homomorphism $x:T\to\bG_m$. 
Then we obtain a holomorphic extension $\bF_x^{\times}$ 
of $A_t$ by $\bG_m$ from Eq.~(\ref{eq:complex Raynaud extension}) 
as the pushout of $\tG$ by $x$. Let $\bF_x$ 
be the holomorphic line bundle on $A$ associated with $\bF_x^{\times}$. 
Since $\bF_x^{\times}$ is a group $\bD$-scheme, 
$\bF_x\in\Pic^0_{A/\bD}(\bD)=A^t(\bD)$.  
So we define $c\in\Hom_{\bZ}(X,A^t(\bD))$ by setting 
$c(x):=\bF_x^{\otimes(-1)}$ 
with the sign change of $c$ and $c^t$ 
in \S~\ref{subsec:Fourier series general} taken into account. 
This $c$ is a complex analogue of 
the classifying morphism of $\tG$ 
in \S~\ref{subsec:Raynaud extensions split case}. 

Let $\cM:=\cL_{C_2,D_2}$ and let $\bM\in\Pic(A)$ 
be a holomorphic line bundle on $A$ 
with $\cO_{A}(\bM)=\cM$,  
$\lambda:A\to A^t$ the polarization morphism by $\bM$,  
$\bP$ the Poincar\`e bundle over $A$, 
$\tbL:=\pi^*\bM$ and $\tcL:=\cO_{\tG}(\tbL)$. 
Note that $\lambda:A\to A^t$ is the natural 
quotient morphism induced from Eq.~(\ref{eq:defn of G/tG/A/Gt/tGT/At}). 
Similarly we set $\lambda_T:T\to T^t$ and $\lambda_G:\tG\to\tG^t$ be the 
natural quotient morphisms induced 
from Eq.~(\ref{eq:defn of G/tG/A/Gt/tGT/At}). 
We have the Raynaud extension of $G^t$:
\begin{equation}\label{eq:Raynaud extension of Gt}
1\to T^t\to \tG^t\to A^t \to 0,
\end{equation}which defines the classifying morphism 
$c^t\in\Hom_{\bZ}(Y,A(\bD))$ by the pushouts via 
$X(T^t)=\Hom(T^t,\bG_{m,\bD})$. It is therefore clear that 
$\lambda_A\circ c^t=c\circ\phi$. 

Then we define 
\begin{gather*}\tau\in\Hom_{\bZ}(Y\times X,t^{\bZ}\left(\bP^{\otimes (-1)}(\bD)\right)^{\times}),\\
\psi\in\Hom(Y,t^{\bZ}\left(\bM^{\otimes -1}(\bD)\right)^{\times}),\ \   
\iota\in\Hom_{\bZ}(Y,\tG(\bD)),\\
\zeta:=(\tG,A,T,X,Y,c,c^t,\iota,\lambda,\phi,\tau,\tcL,\psi,\cM)
\end{gather*}by    
\begin{gather*} 
\tau(y,x):=\be(y^tT_{11}^tx),\ \  
\psi(y):=\be(W_1[y]/2),\\
\iota(y):=\bigoplus_{v\in X}\tau(y,v)\ \ (y\in Y, x\in X). 
\end{gather*}

\begin{defn}\label{defn:split obj in DD_ample}
Let $R:=\bC[[t]]$, $I=tR$ and  
$S=\Spec R$.  We say that $\zeta$ defined above is 
a {\it symmetric} split object over $\bD$  
if the pullback $\zeta_S$ of $\zeta$ to $S$ is 
a {\it symmetric} split object over $S$ in $\DD_{\ample}$. 
\end{defn}

\subsection{Relation with period matrices}
\label{subsec:relation with periods}
Let $\zeta$ be a symmetric split 
object over $\bD$. 
 Let $\be(x):=\exp(2\pi \sqrt{-1}x)$,     
$B:=v_t(\tau)$, $\tau^0:=t^{-B}\tau$ and 
$\psi^{0}:=t^{-v_t(\psi)}\psi$. 
 Let $N$, $\beta$ and $\mu$ be 
the same as in Definition~\ref{defn:miscellany}. 
By (the over-$\bD$-analogue of) Lemma~\ref{lemma:minimal Galois with extensions cte iotae psie taue}, we have extensions 
$c^{t,e}$ (resp. $\iota^e$, $\tau^e$, $\psi^e$) of  
$c^{t}$ (resp. $\iota$, $\tau$, $\psi$) 
from $Y$ to $X^{\vee}$ or from $Y\times X$ 
to $X^{\vee}\times X$.  
By Lemma~\ref{lemma:minimal Galois with 
extensions cte iotae psie taue}~(ii),
$v_t(\tau^e(u,x))=u(x)$ $(\forall u\in X^{\vee},\forall x\in X)$. 
The proof of 
Lemma~\ref{lemma:minimal Galois with extensions cte iotae psie taue} 
shows that the values of $c^{t,e}$ and 
$\iota^e$ (resp. $\tau^e$ and $\psi^e$) are 
{\it single-valued holomorphic} (resp. {\it single-valued 
meromorphic}) over $\bD$ without taking a finite cover of $\bD$.  
By $v_t(\tau^e(u,x))=u(x)$,   
there exists a $\Gamma(\cO_{\bD})$-valued 
bilinear form $\omega^{e,0}$ on $X^{\vee}\times X$ such that 
\begin{equation}\label{eq:extended periods}
\tau^e(u,x)=\be(\omega^e(u,x)),\ \omega^e(u,x)=\omega^{e,0}(u,x)+u(x)\log t.
\end{equation} 

Let 
$\omega:=\omega^e\circ(\beta\times\id_X)$ and 
$\omega^0:=\omega^{e,0}\circ(\beta\times\id_X)$. 
Then  
\begin{equation*}
\omega(y,x)=\omega^0(y,x)+B(y,x)\log t\  
(\forall y\in Y,\forall x\in X),
\end{equation*} where  
$B$ is {\it even} symmetric on $Y\times Y$ 
by \cite[3.5/3.7]{MN24} because $\zeta_S$ is symmetric. See 
Definition~\ref{defn:symm}.
Hence we have
\begin{gather*}
\psi^e(u)=\be(\omega^e(u,\mu(u))),\ \tau^e(u,x)=\be(\omega^e(u,x)),\\
\tau(y,x)=\tau^e(\beta(y),x)=\be(\omega(y,x)),\  
\psi(y)=\be(\omega(y,\phi(y))/2).
\end{gather*}

We define $\zeta^e_{4N^2l}$ (resp. $\xi^{\dagger}_{l}$) to be the 
$4N^2l$-th extension   
of $\zeta$ over $\bD$ as in Definition~\ref{defn:2nl th eFC} 
(resp. the $\NeFC$ $l$-th kit of $\zeta^e_{4N^2l}$).

\subsection{Twisted Mumford families over $\bD$}
\label{subsec:Mumford construction complex}
Let $\Fan(\xi^{\dagger}_{l})$ be a fan in $\tX^{\vee}$ with $e(\zeta)=1$
in \S~\ref{subsec:relation with torus embeddings},   
$U(\Fan(\xi^{\dagger}_{l}))$ the torus embedding over $\bD$ 
associated with  $\Fan(\xi^{\dagger}_{l})$ and $\Sigma^{\dagger}_{l}(0)$ 
the Voronoi polytope in Definition~\ref{defn:Btau Sigma_dagger_l0}.  
In what follows, {\it we assume that $\Sigma^{\dagger}_{l}(0)$ is integral. }
We define the contraction product of $U(\Fan(\xi^{\dagger}_{l}))$ 
 and $\tG$ by $T$ over $\bD$ 
via Eqs.~(\ref{eq:contraction product})/(\ref{eq:complex Raynaud extension}):
\begin{equation}\label{eq:contraction prod over D}
\tP_{l}:=U(\Fan(\xi^{\dagger}_{l}))\times^{T}_{\bD}\tG.
\end{equation}
The $\bD$-scheme $\tP_{l}$ is the normalization of an over-$\bD$-analogue of 
$\Proj \cR_{2,\Sigma^{\dagger}_{l}}$ in Eq.~(\ref{eq:R1Sigma,R2Sigma}). 
Then  $\cW:=(W_{l,\alpha,u};\alpha\in \Sigma^{\dagger}_{l},u\in X^{\vee})$ 
is an open covering of $\tP_{l}$, where $W_{l,\alpha,u}:=U(\tau_{l,\alpha,u})\times^{T}_{\bD}\tG.$  
Let $(w^x; x\in X)$ be the torus coordinates 
of $T\simeq T_{X,\bD}$ global over $D$.  
Let $\xi_{l,\alpha,v}:
=\epsilon^{\dagger}_{l}(v)\tau^e(v,\alpha)w^{\alpha+4Nl\mu(v)}$. 
Let $\tbL'_{l}$ be the holomorphic line bundle 
on $U(\Fan(\xi^{\dagger}_{l}))$ defined by a one cocyle 
$f_{l,\alpha,u;\beta,v}:=\xi_{l,\alpha,u}/\xi_{l,\beta,v}$ 
with respect to $\cW$. We define   
\begin{align}\label{eq:tbL}
\tbL_{l}&
:=\tbL'_{l}\boxtimes_{\cO_{\bD}}\pi^*\bM^{\otimes 4N^2l}.
\end{align}To be more precise, 
\begin{align*}\label{eq:tbL as sheaf}
\cO_{\tP_{l}}(\tbL_{l})&:=p_1^*\cO_{U(\Fan(\xi^{\dagger}_{l}))}(\tbL'_{l})
\otimes_{\cO_{\bD}}p_2^*\pi^*\cO_{A}(\bM^{\otimes 4N^2l})
\end{align*} where 
$p_i$ is the $i$-th projection of 
$U(\Fan(\xi^{\dagger}_{l}))\times^{T}_{\bD}\tG$.

We define an action of $X^{\vee}$ and $Y$  on $\tP_{l}$ and $\tbL_{l}$ 
as follows: Let $u\in X^{\vee}$ and let 
$\tc^{t,e}\in\Hom_{\bZ}(X^{\vee},\tA(D))$ be a lifting of $c^{t,e}$.
 Then we define a $\bD$-automorphism 
$\delta_u:\tP_{l}\to \tP_{l}$ and  $\tilde\delta_u:\tbL_{l}\to\tbL_{l}$ 
via Eqs.~(\ref{eq:contraction prod over D})-(\ref{eq:tbL}):
\begin{gather*}
\delta_u^*(t^aw^x)=t^{a}\tau^e(u,x)w^x,\ \delta_u^*(z)
=T^*_{\tA,\tc^{t,e}(u)}(z):=z+\tc^{t,e}(u),\\ 
\delta_u^*(f)=T^*_{\tA,\tc^{t,e}(u)}(f),\ 
\tilde\delta_u^*(\xi)
=\psi^e(u)^{-2Nl}w^{-4Nl\mu(u)}\xi,
\end{gather*}the last of which follows from 
$\tdelta_{l,u}\tilde\delta_u^*(\theta_{l})
=\epsilon^{\dagger}_{l}(u)w^{4Nl\mu(u)}
\theta_{l}$ by Eq.~(\ref{eq:tdelta_lu delta_lu*}). 
where $(a,x)\in\bZ\times X$ 
with $t^aw^x\in\cO_{\tP_{l}}$, $f\in \cO_A$, 
$(t,z)\in\tA\simeq\bD\times\bA^{g''}$ (as global coordinates)   
and $\xi$ is a fiber coordinate of $\tbL_{l}$. 
The action of $Y$  
on $\tP_{l}$ (via $\beta(Y)\subset X^{\vee}$) is 
 are properly discontinuous 
and fixed point free.  
Hence we have the quotient complex space 
$(P_{l},\bL_{l}):=(\tP_{l},\tbL_{l})/Y$ on which 
$X^{\vee}/\beta(Y)$ acts freely. 
Set $\cL_{l}:=\cO_{P_{l}}(\bL_{l})$. 

Let $\Fan^0:=\{(0),\sigma_{l,0,u}:=\bR_{\geq 0}(f_0+u)\ (u\in X^{\vee})\}$ 
be a fan in $\tX_{\bR}$ and $\Fan^{00}:=\{(0),\sigma_{l,0,\beta(y)}:
=\bR_{\geq 0}(f_0+\beta(y))\ (y\in Y)\}$. 
Let $U(\Fan^0)$ (resp. $U(\Fan^{00})$) be the torus 
embedding over $\bD$ associated with $\Fan^0$ (resp. $\Fan^{00}$). 
Let $\tU^0:=U(\Fan^0)\times_{\bD}\tG$ and 
$\tU^{00}:=U(\Fan^{00})\times_{\bD}\tG$. 
The action of $X^{\vee}$ (and hence of $Y$) 
on $\tP_{l}$ keeps $\tU^0$ and 
$\tU^{00}$ stable. Let 
$$\cG:=\tU^{0}/Y,\ \ G:=\tU^{00}/Y.$$  
The quotients $\cG$ and $G$ are  
open submanifolds of $P_{l}$. 
Thus we obtain
\begin{thm}\label{thm:cubical compactification complex case} 
Let $S:=\Spec\bC[[t]]$. Over $\bD$, the following are true: 
\begin{enumerate}
\item $(\cG,(\cL_{l})_{|\cG})$ and $(G,(\cL_{l})_{|G})$ 
are polarized semi-abelian $\bD$-schemes;  
\item $\cG_{\bD^{\times}}=G_{\bD^{\times}}$; $G_0$ is connected and 
$\cG_0/G_0\simeq X^{\vee}/\beta(Y)$;
\item $\cG_S$ is a N\'eron model of $G_{K}$ where 
$K$ is the fraction field of $\bC[[t]]$;
\item $(P_{l},\cL_{l})\times_{\bD}S$ is 
a cubical compactification of $\cG_S$ 
in \S~\ref{sec:twisted Mumford families}. 
\end{enumerate}
\end{thm}

We call $(P_{l},\cL_{l})$ a cubical compactification of $\cG$.  

\begin{rem}Let $t$ be a uniformizer of $\bD$.   
Then $\omega^e$ in Eq.~(\ref{eq:extended periods}) 
is recovered from $\cG$ and $G$  
as follows:   First  we can identify as follows:
\begin{align*}
Y&=\pi_1(\cG)/\pi_1(A)\simeq H_1(\cG,\bZ)/H_1(A,\bZ),\\
X^{\wedge}&=\pi_1(G)/\pi_1(A)\simeq H_1(G,\bZ)/H_1(A,\bZ).
\end{align*} Since $\cG$ is a smooth group $\bD$-scheme, 
the relative cotangent bundle $T^{\vee}_{\cG/\bD}$ 
of $\cG$ over $\bD$ is trivial. Hence $X$ is identified with 
the $\bZ$-module $V_{\log}$ 
consisting of holomorphic one forms 
$dw^x/w^x$ $(x\in X)$ on $\tG$. 
Any element of $V_{\log}$ 
descends to a holomorphic one form on $\cG$ and $G$.
Then we can define a multi-valued 
period $\omega^e$ of $V_{\log}$ on $X^{\wedge}$  
by integration over one cycles of a general fiber of $\pi$ and hence 
a single-valued pairing $\tau^e:X^{\wedge}\times X\to
\Gamma(\cO_{\bD})[1/t]$ such that $\tau^e(u,x)=\e(\omega^e(u,x))$. 
This recovers Eq.~(\ref{eq:extended periods}). 
See \cite[\S\S~2\,-3]{Nakamura77}. See \cite[9.6]{BLR90} 
for the case of Jacobian varieties.
\end{rem}

\appendix

\section{Algebraizability of $\hat{A}$}
\label{sec:appendix algebraizability of hatA}
In Appendix A, we shall prove   
\S~\ref{subsec:Raynaud extensions}~
Eq.~(\ref{eq:A is alg for nonsplit T0}), that is, 
$\hat{A}$ is algebraizable.  
\subsection{Galois descent}
\label{subsec:Galois descent}
We basically follow the notation and the notions from \cite[11.1-11.5]{MN24}. 
Let $T_0$ be the torus part of $G_0$ in \S~\ref{subsec:Raynaud extensions}. 
Then there exists an $S$-torus $T$ of $\tG$ such that 
$(T^{\wedge})_0=T_0$. If $T$ is not a split $S$-torus, then 
there exists a finite \'etale Galois cover 
$\varpi:S^*\to S$ of $S$ such that   
 $T_{S^*}$ is a split $S^*$-torus. Let $\Gamma:=\Aut(S^*/S)$.

Now let $K:=k(\eta)$, and let $L:=k(S^*)$ be the function field of $S^*$, 
and $A:=R_L=\Gamma(\cO_{S^*})$ the integral closure of $R$ in $L$. 
Then $\Gamma\simeq\Aut(L/K)=\Aut(A/R)$. 
For $\sigma\in\Gamma$, we define 
$s^*_{\sigma}\in\End_R(A)$ by $s^*_{\sigma}(a)=\sigma(a)$ $(\forall a\in A)$, 
and let $s_{\sigma}:S^*\to S^*$ be the $S$-morphism 
induced from $s_{\sigma}^*$. 
Set $t_{\sigma}:=s_{\sigma}^{-1}=s_{\sigma^{-1}}$. Then $s_{\tau}\circ s_{\sigma}=s_{\sigma\tau}$ and $t_{\sigma}\circ t_{\tau}=t_{\sigma\tau}$. 
Let $V$ be an $R_L$-module. Let $t_{\sigma}^*V:=V\otimes_{\sigma}A
:=V\otimes_{A_2;A_2\overset{t_{\sigma}^*}{\to}A_1}A_1$  
where $A_i:=A$ $(i=1,2)$. The $A\ (=1\otimes A_1)$-action 
on $t_{\sigma}^*V$ is given by  
$a(v\otimes 1):=(1\otimes a)(v\otimes 1)=v\otimes a=\sigma(a)v\otimes 1$\ $(\forall a\in A, \forall v\in V).$ Meanwhile we define an $A$-module 
$(V_{\sigma},\cdot_{\sigma})$ by $V_{\sigma}:=V$ 
and $a\cdot_{\sigma}v:=\sigma(a)v$ $(\forall a\in A_1,\forall v\in V)$. For brevity we denote $(V_{\sigma},\cdot_{\sigma})$ by $V_{\sigma}$. Then we define  an isomorphism of $A$-modules $
\kappa_{\sigma}:t_{\sigma}^*V\simeq V_{\sigma}\ (=V)$    
by $\kappa_{\sigma}(\sigma(a)v\otimes 1)=\kappa_{\sigma}(v\otimes a)
=a\cdot_{\sigma}v=\sigma(a)v$ $(\forall a\in A,\forall v\in V)$. 
\begin{lemma}\label{lemma:VGamma trivial}
Let $\varpi:T_{S^*}\to T$ be the natural morphism, 
$\Gamma:=\Gal(S^*/S)$, $C:=\Gamma(\cO_T)$, $L$ a projective 
$C$-module of rank one and $V:=\varpi^*L:=L\otimes_RA$.  Then 
$V_{\Gamma}:
=\bigotimes_{\sigma\in\Gamma}t_{\gamma}^*V$ 
descends to a trivial $C$-module.
\end{lemma}
\begin{proof}Let $D:=\Gamma(\cO_{T_{S^*}})$. 
Then $D=A[w^x;x\in X]$ for a lattice $X$. 
Since $T$ is affine and hence quasi-projective, 
there exist effective divisors $D_i$ of $T$ $(i=1,2)$ 
such that $\widetilde{L}\simeq \cO_{T}(D_1-D_2)$. Let $\cI_{D_i}$ 
be the ideal of $\cO_T$ defining $D_i$.  
Since $\Pic(T_{S^*})=\{1\}$ by \cite[6.1.2]{Nakamura16}, 
$\varpi^*\widetilde{L}\simeq\cO_{T_{S^*}}$ and 
there exists $f_i\in D$ 
such that $\varpi^*\cI_{D_i}=f_i\cO_{T_{S^*}}$. 
Let $I_i:=f_iD$ and $V=\varpi^*L$.  
Then as $D$-modules,  $V\simeq D=\Gamma(\cO_{T_{S^*}})$ and 
$V\simeq I_1^{\otimes (-1)}\otimes_DI_2$.  

Next we shall prove that 
$V_{\Gamma}$ descends to a $C$-module.  
Let $(\phi_{\sigma};\sigma\in\Gamma)$ be a descent datum of $V$ with respect to $\Gamma$. By definition $\phi_{\sigma}:V\to t_{\sigma}^*V$ is an $A$-isomorphism of $A$-modules such that $\phi_{\sigma\tau}=t_{\tau}^*\phi_{\sigma}\circ\phi_{\tau}$. Hence $\sigma:A\to t_{\sigma}^*A$ is lifted to  an $A$-isomorphism $\phi_{\sigma}:V\to t_{\sigma}^*V\simeq V_{\sigma}$. Hence $\sigma$ is lifted to an $\sigma$-semilinear $R$-automorphism $\psi_{\sigma}$ of $V=A[w^x;x\in X]$ such that $\phi_{\sigma}(w^x)=\psi_{\sigma}(w^x)\otimes_{\sigma} 1$ $(\forall x\in X)$. 
See \cite[12.2]{MN24}.

We define $V_{\Gamma}:=\bigotimes_{\sigma\in\Gamma}t_{\sigma}^*V$ 
and $F_i=\prod_{\gamma\in\Gamma}\psi_{\gamma}(f_i)\in D$. 
Since $\sigma:A\to t_{\sigma}^*A$
is lifted to $\psi_{\sigma}:D\to t_{\sigma}^*D$, we have 
$\psi_{\sigma}(F_i)=\prod_{\gamma\in\Gamma}\psi_{\sigma\gamma}(f_i)
=F_i\in D^{\Gamma}=C$.   
Let $\phi_{\Gamma,\sigma}:
=\bigotimes_{\gamma\in\Gamma}t_{\gamma}^*\phi_{\sigma}$. 
Since $\gamma\tau=\tau\gamma'$ for some $\gamma'\in\Gamma$,  
$\gamma'$ ranges over $\Gamma$ iff $\gamma$ ranges over $\Gamma$ 
when $\tau$ is fixed.
Since $t_{\tau\gamma'}^*=t_{\tau}^*t_{\gamma'}^*$, we obtain 
\begin{align*}
\phi_{\Gamma,\sigma\tau}
&=\bigotimes_{\gamma\in\Gamma}
t_{\gamma}^*(t_{\tau}^*\phi_{\sigma}\circ\phi_{\tau})
=\bigotimes_{\gamma\in\Gamma}t_{\tau\gamma}^*\phi_{\sigma}\circ 
t_{\gamma}^*\phi_{\tau}
=(\bigotimes_{\gamma\in\Gamma}t_{\tau\gamma}^*\phi_{\sigma})\circ
(\bigotimes_{\gamma\in\Gamma}t_{\gamma}^*\phi_{\tau})\\
&=(\bigotimes_{\gamma'\in\Gamma}t_{\gamma'\tau}^*\phi_{\sigma})\circ
(\bigotimes_{\gamma\in\Gamma}t_{\gamma}^*\phi_{\tau})
=t_{\tau}^*\phi_{\Gamma,\sigma}\circ \phi_{\Gamma,\tau}.
\end{align*} Hence $(V_{\Gamma},(\phi_{\sigma}:\sigma\in\Gamma))$ 
is a descent datum of $V_{\Gamma}$ with respect to $\Gamma$. 
By \cite[11.3]{MN24}, we have a descent $W_{\Gamma}$ 
of $V_{\Gamma}$. Since $F_i\in C$, $F_i$ 
determines an effective divisor $E_i$ of $T$ such that 
$\cO_T(-E_i)=F_i\cO_T\simeq\cO_T$. Hence  
$\widetilde{V}_{\Gamma}\simeq\varpi^*(\cO_{T}(E_1-E_2))$, so that  
$$\widetilde{W}_{\Gamma}=(\widetilde{V}_{\Gamma})^{\Gamma}
\simeq\left(\varpi^*(\cO_{T}(E_1-E_2))\right)^{\Gamma}
=\cO_T(E_1-E_2)\simeq\cO_T\simeq\widetilde{C}.$$
It follows $W_{\Gamma}\simeq C$. This completes the proof.
\end{proof}

\begin{cor}\label{cor:VGamma trivial}
Let $L$ be a projective 
$C^{\wedge}$-module of rank one and 
$V:=L\otimes_{R^{\wedge}}A^{\wedge}$.  Then 
$V_{\Gamma}:
=\bigotimes_{\sigma\in\Gamma}t_{\gamma}^*V$ 
descends to a trivial $C^{\wedge}$-module.
\end{cor}
\begin{proof}Similar to Lemma~\ref{lemma:VGamma trivial}.
\end{proof}

\subsection{Algebraizability of $\hat{A}$}
\label{subsec:alg of hatA}
Now we return to \S~\ref{subsec:Raynaud extensions}. 
Let $(G,\cL)$ be a semiabelian $S$-scheme with $T_0$ the torus part of $G_0$, 
and $0\to T\to\tG \to A\to 0$ the Raynaud extension of $G$ with 
$T$ an $S$-torus of $\tG$ such that the closed fiber of $T$ equal 
to the above $T_0$. Let $S_n=\Spec R/I^{n+1}$, 
$G_n:=G\times_SS_n\simeq\tG\times_SS_n=:\tG_n$, 
$T_n:=T\times_SS_n$, $A_n:=\tG_n/T_n$, 
$\pi_n:G_n\to A_n$ the natural projection and  
$\cL_n:=\cL\times_SS_n$.  Thus $G^{\wedge}\simeq\tG^{\wedge}$.
Since $A_0$ is an abelian $k(0)$-scheme, 
$A_n$ is a proper smooth group $S_n$-scheme $(\forall n\in\bN)$, 
so that it is an abelian $S_n$-scheme. 
Let $\hat{A}:=\projlim A_n$, $\pi^{\wedge}=\projlim \pi_n$ 
and $\cL^{\wedge}:=\projlim \cL_n$. Then 
$\pi^{\wedge}:\tG^{\wedge}\to \hat{A}$ is the natural 
morphism in Eq.~(\ref{eq:exact seq of Gwedge}) and 
$\cL^{\wedge}$ is an invertible sheaf on $\tG^{\wedge}$.   

\begin{lemma}\label{lemma:A is alg for nonsplit T0}
$\hat{A}$ is algebraizable.
\end{lemma}
\begin{proof}See \cite[II, p.~34, lines~2-3]{FC90}.
Let $B:=\Gamma(\cO_{S^*})$ for now.  Let  
$L=\Gamma(T^{\wedge},\cL^{\wedge}\otimes_{\cO_{\tG^{\wedge}}}
\cO_{T^{\wedge}})$, 
$V=L\otimes_{R^{\wedge}}B^{\wedge}$ 
and $V_{\Gamma}:=\prod_{\sigma\in\Gamma}t_{\sigma}^*V$. 
Let $\cN:=\cL^{\wedge}_{S^{*,\wedge}}$ and 
$\cN_{\Gamma}:=\prod_{\sigma\in\Gamma}t_{\sigma}^*\cN$. 
By Corollary~\ref{cor:VGamma trivial}, the $B^{\wedge}$-module 
$V_{\Gamma}$ descends to a trivial $R^{\wedge}$-module. Since 
$\cN_{\Gamma}=\widetilde{V}_{\Gamma}$, $\cN_{\Gamma}$ 
descends to a trivial $\cO_{T^{\wedge}}$-module. 
Let $i:T^{\wedge}\hookrightarrow\tG^{\wedge}$ 
(resp. $j:(T^{\wedge}_{S^{*,\wedge}})^{\Gamma}\hookrightarrow\tG^{\wedge}_{S^{*,\wedge}}$) be the closed immersion of $T^{\wedge}$ into $\tG^{\wedge}$ (resp. 
$(T^{\wedge}_{S^{*,\wedge}})^{\Gamma}$ into $\tG^{\wedge}_{S^{*,\wedge}}$) 
as $S^{\wedge}$-subschemes. 
Since $(T_{S^*})^{\Gamma}=T$, we have 
$(t_{\sigma}\circ j)^*\cN=j^*\cN\simeq i^*\cL^{\wedge}$. Hence we have 
\begin{align*}
j^*(\cN_{\Gamma})&=\prod_{\sigma\in\Gamma}j^*(t_{\sigma}^*\cN)
=\prod_{\sigma\in\Gamma}(t_{\sigma}\circ j)^*\cN
\simeq (i^*\cL^{\wedge})^{\otimes N_0}
\end{align*}where $N_0=|\Gamma|$. 
It follows $(i^*\cL^{\wedge})^{\otimes N_0}\simeq\cO_{T^{\wedge}}$. 
By \cite[I, 7.2.2]{MB85} 
there exists an ample invertible sheaf $\hat{\cM}$ on $\hat{A}$ 
such that $\cL^{\otimes N_0,\wedge}\simeq
(\pi^{\wedge})^*\hat{\cM}$.  Since $\cL$ is ample on $G$
in the sense of \cite[II, 4.5.2]{EGA}), 
so are $\cL^{\otimes N_0}_{n}$ on $G_n$ and $\hat{\cM}_n$ on $A_n$ 
with $\cL^{\otimes N_0}_n\simeq\pi_n^*\hat{\cM}_n$.   
It follows from \cite[III, 5.4.5]{EGA} 
that there exist an abelian $S$-scheme $A$ and 
 an ample invertible sheaf $\cM$ on $A$ such that 
$(A^{\wedge},\cM^{\wedge})\simeq(\hat{A},\hat\cM)$. This completes the proof. 
\end{proof}

\begin{example}\label{example:Weil res}
Let $k$ be a field of $\chara k\neq 2$, 
$k'$ a Galois extension of $k$ with $\Gamma:=\Gal(k'/k)\simeq\bZ/2\bZ$. Let $\Gamma=\{1,\sigma\}$, $k'=k(\epsilon)$ 
with $\sigma(\epsilon)=-\epsilon$ and $a:=-\epsilon^2\in k$. 
Let $H=\Spec A$, $H'=\Spec A'$ and $A':=k'[w, 1/w]$.  
Suppose that $H$ is not a split $k$-torus. 
Then 
\begin{equation}\label{defn of H}
A=(A')^{\Gamma}=k[x,y]/(x^2+ay^2-1),\ A'=A\otimes_kk',
\end{equation}where  $x:=(1/2)(w+1/w)$ and $y:=(1/2\epsilon)(w-1/w)$. 
The product of $H$ is given by $(x,y)\cdot(x',y')=(xx'-ayy',x'y+xy')$. 
In other words, $H'=\bG_{m,k'}$ and 
$H$ is the Weil restriction $\Res_{k'/k}H'$ of $H'$. 
\end{example}

\section{Extension of 
the degeneration data of $(G,\cL)$}
\label{sec:appendix deg data}

The purpose of Appendix~\ref{sec:appendix deg data} 
is to extend the degeneration data $\FC(G,\cL)$ of a 
semiabelian scheme $(G,\cL)$ 
from the lattices $(X,Y)$ to $(X,X)$.

\subsection{Fppf sheaves}
\label{subsec:fppf presheaf hW is sheaf}  
Let $S$ be a scheme, $\Scheme/S$ 
the category of schemes over $S$, and $\Scheme_S$ 
the category $\Scheme/S$ equipped with the fppf topology. 
For an $S$-scheme $W$, 
we define a presheaf $\underbarW$ by 
$\underbarW(U):=\Hom_{\Sch_S}(U,W)$ $(U\in\Sch_S)$. 
We denote $\Hom_{\Sch_S}(U,W)$ by $\Hom(U,W)$ or $W(U)$  
if no confusion is possible. 
This presheaf $\underbarW$ 
is an fpqc sheaf by descent theory 
\cite[VIII, 5.1~b), pp.~210-211]{SGA1}.  
See also \cite[Descent, 023Q]{Stacks}. 

To be more precise, 
$\underbarW$ satisfies the {\it sheaf condition} for the fpqc topology: 
for any fpqc covering $(U_i\to U; i\in I)$ of $U$, 
the following is exact:
\begin{equation*}
\Hom(U,W)\to \prod_{i\in I}\Hom(U_i,W)\rightrightarrows \prod_{i,j\in I}
\Hom(U_i\times_UU_j,W).
\end{equation*} 
Since any fppf covering is an fpqc covering,  
$\underbarW$ satisfies the sheaf condition for the fppf topology, {\it i.e.}, 
$\underbarW$ is an {\it fppf sheaf.}  We say that an fppf sheaf $\cF$ is 
{\it representable by a scheme $W$} if $\cF\simeq\underbarW$.

For a morphism $f:W\to Z$ of $S$-schemes, 
we define an {\it fppf sheaf homomorphism} $\underbarf:\underbarW\to \underbarZ$ by $\underbarf(U):=f(U):W(U)\to Z(U)$ $(U\in\Sch_S)$.   

\subsection{Extension of a split object}
Let $(G,\cL)$ be a {\it symmetric} semiabelian $S$-scheme 
with the torus part of $G_0$ split. 
Then we have a split object: \begin{equation}\label{eq:split obj zeta in DDample appendix}
\zeta:=\FC(G,\cL)=(\tG,A,T,X,Y,c,c^t,\iota,\lambda,\phi,\tau,\tcL,\psi,\cM)
\end{equation}in Eq.~(\ref{eq:split obj zeta in DDample}) 
associated with $(G,\cL)$.  
In \S\S~\ref{subsec:construction of c_t_ex}-\ref{subsec:construction of tau_ex,psi_ex,iota_ex}, we shall prove: 
\begin{lemma}\label{lemma:minimal Galois with ctex iotaex psiex tauex}
Let $(G,\cL)$ be a symmetric semiabelian $S$-scheme 
with $\zeta=\FC(G,\cL)$ and $T_0$ split. 
By taking a finite radical 
normal extension $K$ of $k(\eta)$ possibly ramified 
over the maximal ideal $I$ of $R$, we can find  
\begin{enumerate}
\item[(i)]$c^{t,\ex}\in\Hom_{\bZ}(X,A(R_K))$ with
$c^t=c^{t,\ex}\circ\phi$ and $c=\lambda\circ c^{t,\ex}$;
\item[(ii)] a trivialization of biextension 
$\tau^{\ex}:1_{X\times X}\simeq 
(c^{t,\ex}\times c)^*\cP_K^{\otimes (-1)}$, which is 
a symmetric blinear form $\tau^{\ex}\in\Hom_{\bZ}(X\times X,K^{\times})$ 
such that $\tau^{\ex}(\phi(y),x)=\tau(y,x)$\   
$(\forall x\in X,  \forall y\in Y)$;
\item[(iii)]$\iota^{\ex}(x)\in\Hom_{\bZ}(X,\tG(K))$ such that 
$\iota=\iota^{\ex}\circ\phi$ and $c^{t,\ex}=\pi\circ\iota^{\ex}$;
\item[(iv)] 
a cubical trivialization 
$\psi^{\ex}:1_X\simeq (c^{t,\ex})^*\cM_K^{\otimes (-1)}$, which is  
a $K^{\times}$-valued function compatible with 
$\tau^{\ex}$ in the sense that  
\begin{gather*}
\psi^{\ex}(x+z)=\psi^{\ex}(x)\psi^{\ex}(z)\tau^{\ex}(x,z),\ 
\psi^{\ex}(\phi(y))=\psi(y)\\ 
(\forall x,z\in X, \forall y\in Y);
\end{gather*} 
\end{enumerate}
\end{lemma}

Lemma~\ref{lemma:minimal Galois with ctex iotaex psiex tauex} is proved 
by Lemmas~\ref{lemma:cte}/\ref{lemma:tau_ex,psi_ex} 
and Eq.~(\ref{eq:defn of iota_ex}). 
By Lemma~\ref{lemma:minimal Galois with ctex iotaex psiex tauex}, 
we can construct a split object $\zeta^{\ex}$ extended from  
from $Y$ (resp. $Y\times X$) to $X$ (resp. $X\times X$). 
This fact was used in \cite[4.10]{Nakamura99} without proof.

\subsection{Construction of $c^{t,\ex}$}
\label{subsec:construction of c_t_ex}
Let $A_i=A$ and $B_i=B$ be abelian $S$-schemes $(i=1,2)$.  
Let $f_i:A_i\to B_i$ $(i=1,2)$, $g_1:A_1\to A_2$ and $g_2:B_1\to B_2$ 
be finite flat surjective $S$-morphisms. 
Let $H:=\ker(f_1)$ and  $K:=\ker(g_1)$. 
Then $A_2\simeq A_1/K$ and $B_1\simeq A_1/H$. 
Let $B_2:=A_1/(H+K)$ and let 
$f_2:A_2\to B_2$ (resp. $g_2:B_1\to B_2$) be the natural morphism. 

\begin{notation}
In what follows for a finite closed 
subgroup $S$-scheme $N$ of $A_1$, $U$ an $S$-scheme 
and $X,Y\in A_1(U)$, we denote $X-Y\in N(U)$ by $X\equiv Y\in (A_1/N)(U)$. 
We denote the class of $X$ in $(A_1/N)(U)$ by $X+N(U)$.
\end{notation}

\begin{lemma}\label{lemma:fiber prod}
$A_2\times_{B_2}B_1\simeq A_1/H\cap K$.
\end{lemma}
\begin{proof} 
Let $V=A_2\times_{B_2}B_1$ and $W=A_1/H\cap K$. 
We shall prove an isomorphism $V\simeq W$. 
Since there is a natural morphism $J:W\to V$, 
we shall construct its inverse $G:V\to W$. 
Let $T$ be any $S$-scheme and 
we shall define a map  
$G(T):V(T)\to W(T)$ functorial in $T$. 
Let $P:=(Q^0,R^0)\in V(T)$. Namely   
$Q^0\in A_2(T)$ and $R^0\in B_1(T)$ with  
$f_2(Q^0)=g_2(R^0)$. 
Hence there exist a scheme $U$ fppf over $T$, and  
$Q^1,R^1\in A_1(U)$ such that 
\begin{gather*}
Q^1\equiv Q^0_{U}\in A_2(U),\ R^1\equiv R^0_{U}\in B_1(U),\ 
Q^1\equiv R^1\in B_2(U)
\end{gather*} 

Therefore there exist $Q^2\in K(U)$ and 
$R^2\in H(U)$ such that 
$Q^1=R^1+R^2-Q^2\in A_1(U)$. Then we define
\begin{equation*}
P^1:=Q^1+Q^2=R^1+R^2\in A^1(U).
\end{equation*}
Suppose that there exist    
$E^1,F^1\in A_1(U)$ such that   
\begin{gather*}
E^1\equiv Q^0_{U}\in A_2(U),\  
F^1\equiv R^0_{U}\in B_1(U),\ 
E^1\equiv F^1\in B_2(U).
\end{gather*} 
Then there exist $E^2\in K(U)$ and 
$F^2\in H(U)$ such that $E^1_{U}=F^1_{U}+F^2-E^2\in A_1(U)$. 
Then we define $P^2:=E^1+E^2=F^1+F^2\in A^1(U).$  

Now we compare $P^1$ and $P^2$: 
\begin{align*}
P^1-P^2
&=(Q^1-E^1)+(Q^2-E^2)\in K(U),\\
P^1-P^2&=(R^1-F^1)+(R^2-F^2)\in H(U),
\end{align*}whence $P^1-P^2\in H(U)\cap K(U)$. 
Thus we have: 
\begin{claim}\label{claim:uniqueness of P1_in_W(U)}
Let $T$ be an $S$-scheme and 
$P=(Q^0,R^0)\in V(T)$. Then 
\begin{enumerate}
\item[(a)]there exist a $T$-scheme $U$ fppf over $T$ 
and $P^1\in A_1(U)$ 
such that $P^1\equiv Q^0_U\in A_2(U)$ and 
$P^1\equiv R^0_U\in B_1(U)$. Moreover 
the class $P^1+H(U)\cap K(U)\in W(U)$ is uniquely determined by 
$P\in V(T)$;
\item[(b)] suppose that $U'$ is a $T$-scheme fppf over $T$ 
with $P'\in A_1(U')$ such that $P'\equiv Q^0_{U'}\in A_2(U)$ and 
$P'\equiv R^0_{U'}\in B_1(U')$; then the class  
$P'+H(U')\cap K(U')\in W(U')$ is uniquely determined by $P\in V(T)$.
\end{enumerate}
\end{claim}
\begin{proof}(a) follows from the above argument. 
(b) is proved by applying the proof of (a) 
to $U'$ and $P'\in A_1(U')$.
\end{proof} 

Let $U$ be a $T$-scheme fppf over $T$ and $P^1\in A_1(U)$ 
in Claim~\ref{claim:uniqueness of P1_in_W(U)}~(a). 
Since the scheme $W$ is an fppf sheaf over $\Sch_S$, 
there is an exact sequence 
\begin{equation}\label{eq:W fppf}
0\to W(T)\to W(U)\underset{p_2^*}{\overset{p_1^*}{\rightrightarrows}} 
W(U\times_TU).
\end{equation}
Let $U^{(2)}:=U\times_TU$ and 
let  $p_i:U^{(2)}\to U$ be the $i$-th projection.  
Then we have $S^i:=p_i^*P^1\in A_1(U^{(2)})$. 
Then $S^i+K(U^{(2)})\in A_2(U^{(2)})$ and 
$S^i+H(U^{(2)})\in B_1(U^{(2)})$. Note that 
$U^{(2)}$ is a $T$-scheme fppf over $T$, and 
$S^i\equiv p_i^*Q^0_U\in A_2(U^{(2)})$ 
and $S^i\equiv p_i^*R^0_U\in B_1(U^{(2)})$. 
Hence we apply Claim~\ref{claim:uniqueness of P1_in_W(U)}~(b) to  
$(U',P')=(U^{(2)},S^i)$ $(i=1,2)$. It follows that  
$S^i+H(U^{(2)})\cap K(U^{(2)})$ is uniquely determined by $P\in V(T)$, 
and therfore, $S^1\equiv S^2\in W(U^{(2)})$. 
Hence $p_1^*P^1\equiv p_2^*P^1\in W(U^{(2)})$ 
in Eq.~(\ref{eq:W fppf}). Therefore there exists $P'\in W(T)$ such that 
$P'_U\equiv P^1\in W(U)$. 

Now we define a morphism of functors $G:V\to W$ by 
$$G(T)(P):=P'\in W(T).$$ 
Since $J(T):W(T)\to V(T)$ is given 
by defining $J(T)(P')=(P'+K(T),P'+H(T))\in V(T)$ for $P'\in W(T)$, 
it is easy to prove $J(T)\circ G(T)=\id_{V(T)}$ and 
$G(T)\circ J(T)=\id_{W(T)}$, so that $V(T)\simeq W(T)$ and 
$V\simeq W$ as $S$-schemes. This completes 
the proof of Lemma~\ref{lemma:fiber prod}.   
\end{proof}

\begin{cor}\label{cor:A1 to V surj}
The natural morphism $J:A_1\to A_2\times_{B_2}B_1$ is surjective. 
\end{cor}

Let $\Omega=\overline{k(\eta)}$ and $R_{\Omega}$ 
the integral closure of $R$ in $\Omega$.

\begin{lemma}\label{lemma:cte}
There exists a homomorphism 
$c^{t,e}\in\Hom_{\bZ}(X, A_1(R_{\Omega}))$ 
such that $\lambda\circ c^{t,\ex}=c\in\Hom_{\bZ}(X, A^t(R_{\Omega}))$ and 
$c^{t,\ex}\circ\phi=c^t\in\Hom_{\bZ}(Y, A^t(R_{\Omega}))$. 
\end{lemma}
\begin{proof}Let $X/\phi(Y)=\bigoplus_{i=1}^r(\bZ/e_i\bZ)$ with 
$e_i|e_{i+1}$ $(i\in [1,r-1])$. Then there exists 
a $\bZ$-basis $(x_i;i\in[1,r])$ of $X$ 
(resp. $(y_i;i\in[1,r])$ of $Y$) such that 
$e_ix_i=\phi(y_i)$ $(\forall i\in [1,r])$. 
For each $i$, we have a commutative diagram 
\begin{diagram}
A_1&\rTo^{\lambda}&A^t_1\\
\dTo^{e_i\id_A}&&\dTo_{e_i\id_{A^t}}\\
A_2&\rTo^{\lambda}&A^t_2
\end{diagram}where $A_i=A$ and $A^t_i=A^t$. 
 Since $\lambda\circ c^t(y_i)=c\phi(y_i)=e_ic(x_i)$, 
by Corollary~\ref{cor:A1 to V surj}, there exists 
$\alpha_i\in A_1(R_{\Omega})$ such that 
$J(\alpha_i)=(c^t(y_i),c(x_i))$, that is, 
$c^t(y_i)=e_i\alpha_i$ and $c(x_i)=\lambda(\alpha_i)$. 
Now we define $c^{t,\ex}(x_i)=\alpha_i$ $(\forall i\in[1,r])$ and 
extend $c^{t,\ex}$ to $X$ additively. 
Hence $\lambda\circ c^{t,\ex}=c$ and 
$c^{t,\ex}\phi=c^t$. 
\end{proof}

\subsection{Construction of $(\tau^{\ex}, \psi^{\ex}, \iota^{\ex})$}
\label{subsec:construction of tau_ex,psi_ex,iota_ex}
We can define $\iota^{\ex}$ by Lemma~\ref{lemma:minimal Galois with ctex iotaex psiex tauex}~(iii), equivalently, by a variant of Eq.~(\ref{eq:iotay}):
\begin{equation}\label{eq:defn of iota_ex}
\iota^{\ex}(x)^*:=\bigoplus_{z\in X}\tau^{\ex}(x,z)
\in (c^{t,\ex}(x)\times c)^*\cP_{\eta}^{\otimes (-1)}\simeq x\times X\times\Omega^{\times}.
\end{equation}

Therefore we suffice to construct $\tau^{\ex}$ and $\psi^{\ex}$.
In what follows, we use the same notation $e_i, x_i$ and $y_i$ 
as in the proof of 
Lemma~\ref{lemma:cte}. 
\begin{lemma}\label{lemma:cL symmetric, psi symmetric} 
If $\cL$ is symmetric, then  
$\psi(-y)=\psi(y)$ and 
$\psi(ey)=\psi(y)^{e^2}$ $(\forall e\in\bZ, \forall y\in Y)$.
\end{lemma}
\begin{proof}This is proved in the same manner as \cite[3.5]{MN24}.
Let $\theta\in\Gamma(G_{\eta},\cL_{\eta})$ and 
$\theta':=[-\id_G]^*\theta$. Then the Fourier expansion of 
$\theta$ is uniquely given by 
$\theta=\sum_{x\in X}\sigma_x(\theta)$, while 
$\theta'=\sum_{x\in X}\sigma_x(\theta')=\sum_{x\in X}\sigma_{-x}(\theta)$ with $\wt\sigma_{-x}=x$, whence $\sigma_x(\theta')=\sigma_{-x}(\theta)$. 
For any $y\in Y$, $\sigma_{x+\phi(y)}(\theta)=\psi(y)\tau(y,x)\sigma_x(\theta)$ and  
$\sigma_{x+\phi(y)}(\theta')=\psi(y)\tau(y,x)\sigma_x(\theta')$. 
Since $\sigma_{x+\phi(y)}(\theta)=\sigma_{-x-\phi(y)}(\theta')$, 
we obtain $\psi(-y)=\psi(y)$. The rest is easy.
\end{proof}

\begin{lemma}\label{lemma:tau_ex,psi_ex} 
By taking a finite radical normal extension $K$ of $k(\eta)$ possibly ramified over the maximal ideal $I$, we can find 
\begin{enumerate}
\item a $K^{\times}$-valued symmetric bilinear form 
$\tau^{\ex}\in\Hom_{\bZ}(X\times X,K^{\times})$;
\item a $K^{\times}$-valued function $\psi^{\ex}$ 
on $X$ such that 
\end{enumerate}
\begin{gather*}
\psi^{\ex}(x+z)=\psi^{\ex}(x)\psi^{\ex}(z)\tau^{\ex}(x,z),\\
\tau^{\ex}(\phi(y),x)=\tau(y,x),\ \psi^{\ex}(\phi(y))=\psi(y)\ 
(\forall x,z\in X, \forall y\in Y).
\end{gather*} 
\end{lemma}
\begin{proof}
To simplify the argument, we identify $Y$ with 
the subgroup $\phi(Y)$ of $X$. We first prove Lemma when $K=\Omega$.  
We choose $\bZ$-bases $(x_i;i\in[1,r])$ of $X$ 
and $(y_i;i\in[1,r])$ of $Y$ as before.
Let $s\geq 0$ and $X_s:=\sum_{i=1}^{s}\bZ x_i+Y$. 
Hence $X_0=Y$ and $X_r=X$. First we shall construct 
$\tau^{\ex}$ on $X\times X$ by the induction on $s$. 
We set $\tau_0=\tau$.  
Suppose that $\tau_{s-1}$ is defined and bilinear 
on $X_{s-1}\times X$ with $\tau_{s-1}=\tau$ on $Y\times X$. 
Then we define $\tau_{s}=\tau_{s-1}$ on $X_{s-1}\times X$ and 
$\tau_{s}(x_s,x_i)=\tau_s(x_i,x_s)$ if $i<s$. 
If $s\leq i\leq r$, we define $\tau_{s}(x_s,x_i)$ to be an $e_s$-th root of 
$\tau(y_s,x_i)$. And then we extend $\tau_{s}$ additively to $X_{s}\times X$. It is clear that $\tau_{s}$ is bilinear on 
$X_{s}\times X$ and $\tau_{s}=\tau_{s-1}$ on $X_{s-1}\times X$. Finally 
we set $\tau^{\ex}:=\tau_{r}$. \par
Next we define $\psi_s$ on $X_s$ by the induction on $s$. 
We set $\psi_0=\psi$. Suppose that $\psi_{s-1}$ is defined on $X_{s-1}$ 
so that $\psi_{s-1}(x+z)=\psi_{s-1}(x)\psi_{s-1}(z)\tau_{s-1}(x,z)$\  
$(\forall x\in X_{s-1},\forall z\in X)$. Then we define 
$\psi_{s}(x_s)$ to be a common solution of the equations 
$\psi_{s}(x_s)^2=\tau_{s}(x_s,x_s)$ and 
$\psi_{s}(x_s)^{e_s^2}=\psi_{s}(e_sx_s)=\psi_{s-1}(y_s)$. 
If $e_s$ is odd, then $\psi_{s}(x_s)$ is uniquely determined 
by these conditions. If $e_s$ is even, we can choose any of the 
square roots of $\tau_{s}(x_{s},x_{s})$ as $\psi_{s}(x_s)$.  
Then we define $\psi_{s}$ on $X_{s}$ by 
$$\psi_{s}(a_sx_s+x')=\psi_{s}(x_s)^{a_s^2}
\psi_s(x')\tau_{s}(a_sx_s,x')\ (a_s\in\bZ, x'\in X_s).  
$$

It is easy to check 
$\psi_s(ax_s+y+bx_s+z)=\psi_s(ax_s+y)\psi_s(bx_s+z)\tau_s(ax_s+y,bx_s+z)$ $(\forall a,b\in\bZ,\forall y,z\in X_{s-1})$.  
Finally we set $\psi^{\ex}:=\psi_{r}$.  The proof of Lemma is over by choosing a finite radical normal extension $K$ of $k(\eta)$ which is a subfield of $\Omega$ containing 
$\tau^{\ex}(x,y)$ and $\psi^{\ex}(z)$ $(\forall x,y,z\in X)$.
\end{proof}

This completes the proof 
of Lemma~\ref{lemma:minimal Galois with ctex iotaex psiex tauex}.

\subsection{A standard $K$-basis of $\Gamma(G_K,\cL_K)$}
\label{subsec:standard K basis}
Let $M$ be a positive integer, $\zeta_M$ a primitive 
$M$-th root of unity and $R$ an integral domain,  $K$ its fraction field
and $S=\Spec R$. Assume $R\supset \bZ[\zeta_M,1/M]$.  
Let $H$ be a finite abelian group such that 
$M=e_{\max}(H)$, the maximal order of elements in $H$. 
Let $H^{\vee}$ be the character group of $H$ 
and $K(H):=H\oplus H^{\vee}$. We define a {\it bilinear} form 
$e_H : K(H)\times K(H)\to \bG_{m}$ by 
\begin{equation}
e_H(z\oplus\alpha,w\oplus\beta)=\beta(z)\alpha(w)^{-1}
\end{equation} 
where $z,w\in H$, $\alpha,\beta\in H^{\vee}$. 
Let $V(H):=R[H]$ be the group algebra of $H$ over $R$.
Let $H_S$ be 
the constant finite abelian group 
$S$-scheme associated with $H$, $H^{\vee}_S:=D(H_S)$ 
the Cartier dual of $H_S$ over $S$ and 
$e_{H_S}$ the bilinear form on $K(H)_S$ induced from $e_H$. 
Let $\cG(H)$ be the central extension 
of $K(H)$ by $\bG_{m,S}$ associated with $e_H$, 
that is, the group $S$-scheme structure of $\cG(H)_S$ 
is defined by 
\begin{equation}\label{eq:symplectic pairing}
(a,z,\alpha)\cdot (b,w,\beta)
=(ab\beta(z),z+w,\alpha+\beta)
\end{equation}
where $a,b\in\bG_{m,S}$, $z,w\in H_S$ 
and $\alpha,\beta\in H^{\vee}_S$ are functorial points.	
Then $e_{H_S}$ is the commutator form of $\cG(H)_S$. 

We define $\rho_H\in\End_R(V(H))$ by 
\begin{equation}
\rho_H(a,x,\chi)(z)=a\chi(z)(x+z)
\end{equation}where $z\in H$ and $(a,x,\chi)\in\cG(H)$.
Hence $V(H)$ is a $\cG(H)_S$-module of weight one:
$\bG_{m,S}$ acts on $V(H)$ by scalar multiplication.
By \cite[V,\,2.4.2]{MB85}, $\rho_H$ 
is irreducible in the sense of Footnote~7.  
The following condition ($S_D$) is true over any integral 
domain $D$ over $R$: 
\begin{enumerate}
\item[($S_D$)]$(\rho_H)_D$ 
is irreducible and any $\cG(H)_D$-module $W$ 
of weight one is isomorphic to 
$V(H)_D^{\oplus n}$ for some $n\geq 0$,
\end{enumerate} 
where $Z_T:=Z_D:=Z\otimes_RD$ 
for $Z\in\{\rho(H), \cG(H), V(H)\}$ and 
$T:=\Spec D$.

\begin{defn}\label{defn cG(H) and V(H)}
Let $V$ be any irreducible $\cG(H)_K$-module. 
A $K$-basis  $(\theta_z;z\in H)$ of $V$ is called 
a {\it standard $K$-basis} of $V$ if 
$$\rho_H(a,x,\chi)(\theta_z)=a\chi(z)(\theta_{x+z})\ (\forall 
z\in H, \forall (a,x,\chi)\in\cG(H).$$  We call the isomorphism 
$\psi:V(H)\otimes_DK\to V_K$ sending 
$z\mapsto \theta_z$ a {\it standard $K$-isomorphism}, 
which is unique up to $K^{\times}$-multiples by Schur's lemma.
\end{defn}

\begin{cor}\label{cor:standard basis}
Suppose that $R$ is a CDVR. Let $(G,\cL)$ be a symmetric semi-abelian 
$S$-scheme, and $(A,\cM)$ the abelian part 
of the Raynaud extension $\tG$. 
If $(A,\cM)$ is a polarized abelian $S$-scheme with 
$\cM$ a separable polarization,
then by taking a finite radical normal extension $K$ 
of $k(\eta)$ if necessary, 
there exists a standard $K$-basis of 
$\Gamma(G_{K},\cL_{K})$. 
\end{cor}
\begin{proof}
Since $(A,\cM)$ is a polarized abelian $S$-scheme with 
$\cM$ a separable polarization, 
by \cite[Cor. of Th.~1, p.~294]{Mumford66} 
there exists a finite abelian group $H_2$ 
such that $\cG(A_{\Omega},\cM_{\Omega}))
\simeq\cG(H_2)$ and $\Gamma(A_{\Omega},\cM_{\Omega})\ 
(\simeq V(H)_{\Omega})$ has a standard $\Omega$-basis 
$(\vartheta_{z_2};z_2\in H_2)$.  
By Lemma~\ref{lemma:minimal Galois with ctex iotaex psiex tauex},  
we define:
$$\theta_{z_1,z_2}=\sum_{y\in Y}\psi^{\ex}(\tz_1+\phi(y))
T_{c^t(y)}^*(\vartheta_{z_2}), 
$$where $H_1=X/\phi(Y)$, $z_i\in H_i$\ $(i=1,2)$ and we choose 
$\tz_1\in X$ such that $\tz_1\mod Y=z_1$. Moreover 
$\psi^{\sharp}(\tz_1+\phi(y))$ is considered to have weight $\tz_1+\phi(y)$, in other words, in a down-to-earth manner (\S~\ref{subsec:down_to_earth tPl}), 
it is written as 
$\psi^{\sharp}(\tz_1+\phi(y))w^{\tz_1+\phi(y)}_U$ $(U\in\cF_0)$. 
Then ($\theta_{z_1,z_2};z_i\in H_i)$ 
is  a standard $K$-basis of 
$\Gamma(G_{K},\cL_{K})$.  
Indeed, let $H=H_1\oplus H_2$ and identify  
$\cG(G_{\Omega},\cL_{\Omega})\simeq\cG(H_1\oplus H_2)_{\Omega}$. 
By an easy computation, we obtain 
\begin{align*}
\rho&_{H_1\oplus H_2}(a,x_1\oplus x_2,\chi_1\oplus \chi_2)(\theta_{z_1,z_2})
=a\chi_1(z_1)\chi_2(z_2)\theta_{z_1+x_1,z_2+x_2}
\end{align*}where $\chi_i\in D(H_i)$\ $(i=1,2)$. Then by choosing any finite radical normal extension field $K$ of $k(\eta)$ which is a subfield of $\Omega$ 
 containing all $\tau^{\ex}(x,z)$ 
and $\psi^{\ex}(x)$ $(x,z\in X)$ and $\bZ[\zeta_M,1/M]$, 
we obtain Lemma.
\end{proof}

\subsection{$K$-rationality of $K_s^{\wedge}$-points}
\label{subsec:K rationaliy}
Let $D$ be a Dedekind domain of {\it characteristic zero,} 
$K$ its fraction field and $S=\Spec D$. 
Let $(G,\cL)$ be a semiabelian $S$-scheme   
 over $S$ with $\cL$ symmetric ample cubical and 
$\cG(G,\cL)$ the theta group of 
$(G,\cL)$. By Corollary~\ref{cor:standard basis}, 
 there exists a finite Galois extension $K'$ 
of $K$ and a finite abelian group $H$ 
such that $(\cG(G,\cL),\Gamma(G,\cL))\otimes_K K'
\simeq(\cG_H,V_H)\otimes_K K'$. 
We call the isomorphism $(\phi,\psi):(\cG_H,V_H)\otimes_K K'\simeq(\cG(G,\cL),\Gamma(G,\cL))\otimes_K K'$ a {\it standard isomorphism over $K'$}. Once we fix $\phi$, then $\psi$ is uniquely determined up to constant multiples. 
For simplicity we take $K'$ as $K$ 
(and then take the integral closure of $D$ in $K'$ for $D$).

Let $s\in \cB$, 
$D_s$, $D_s^{\wedge}$, $K_s$, $K_s^{\wedge}$ and   
$G^s:=G\times_SD_s$ and $\cL^s$ be 
the same as those in \S~\ref{subsec:plan}. 
Let $\Omega^{(s)}$ be an algebraic closure of $K_s^{\wedge}$ and 
$\phi^{s,\wedge}:=\phi\otimes_D\Omega^{(s)}$.
 
The following is now obvious: 
\begin{cor}
\label{cor:K rationality}
Let $(\phi,\psi)$ and $(\phi^{s,\wedge},\psi')$ be 
standard isomorphisms: 
\begin{gather*}
(\phi,\psi):(\cG(H),V(H))\simeq(\cG(G,\cL),\Gamma(G,\cL)),\\
(\phi^{s,\wedge},\psi'):(\cG(H),V(H))\otimes_D\Omega^{(s)}\simeq(\cG(G^{s,\wedge},\cL^{s,\wedge}),\Gamma(G^{s,\wedge},\cL^{s,\wedge}))\otimes_{D_s^{\wedge}}\Omega^{(s)}.
\end{gather*}
 Let $\theta_z:=\psi(z)$ and   
$\theta'_z:=\psi'(z)$.  
Then   
\begin{enumerate}
\item  $\theta_z/\theta_w=\theta'_z/\theta'_w$ 
on $G^{s,\wedge}_{\Omega^{(s)}}$\   
$(\forall z,w\in H)$;
\item if $\cL_K$ is very ample, then the following are equivalent:
\begin{enumerate}
\item[(2a)]$P\in G(K)$;
\item[(2b)]$P\in G(\Omega^{(s)})$ and 
$\theta'_z(P)/\theta'_w(P)\in K$\  
$(\forall z,w\in H)$ if $\theta'_w(P)\neq 0$, 
\end{enumerate}
\end{enumerate}where $G(L):=\Hom(\Spec L,G)$ for a field $L$.
\end{cor}
\begin{proof}Since $\Gamma(G,\cL)$ is $D$-free 
and $\rank_K\Gamma(G,\cL)=
\rank_\Omega^{(s)}\Gamma(G^{s,\wedge},\cL^{s,\wedge})
\otimes_{D_s^{\wedge}}\Omega^{(s)}$ by Theorem~\ref{thm:degeneration data2}~(\ref{item:Gamma Geta L general}), 
we have  $\Gamma(G,\cL)\otimes_{D}\Omega^{(s)}=\Gamma(G^{s,\wedge},\cL^{s,\wedge})\otimes_{D_s^{\wedge}}\Omega^{(s)}$.  
Since $(\theta_z;z\in H)$ (resp. 
$(\theta'_z;z\in H)$) is a standard $\Omega^{(s)}$-basis of 
$\Gamma(G,\cL)\otimes_K\Omega^{(s)}$ (resp. 
$\Gamma(G^{s,\wedge},\cL^{s,\wedge}))\otimes_{D^{s,\wedge}}\Omega^{(s)}$). 
By Schur's lemma, there exists $c\in(\Omega^{(s)})^{\times}$ such that 
$\theta_z=c\theta'_z$ $(\forall z\in H)$.  This proves (1) and (2).
\end{proof}

\end{document}